%% file: qcpaper4.tex
\documentclass[11 pt,reqno]{amsart}
\usepackage{amssymb}
\usepackage{amsmath}
\usepackage{amscd}
\usepackage[perpage,multiple]{footmisc}

\usepackage{diagrams}
\input{quasimacros}

\begin{document}


\input{qctitleabstract}

\setcounter{tocdepth}{1}
\tableofcontents


\input{introduction}

\part{Quasi--Coxeter algebras}


\input{DCPfinal}


\input{newassociahedron}


\input{quasiCoxeter}


\input{qcexamplesnew}


\input{dynkin}

\part{Quasi--Coxeter quasibialgebras}


\input{qcqtqb}


\input{dynkinbi}

\part{Quantum Weyl groups}


\input{transfer}


\input{rigidity}


\input{acknowledgements}


\input{quasibiblio}
\end{document}

%% file: quasimacros.tex
\newtheorem {theorem}{Theorem}[section]
\newtheorem {proposition}[theorem]{Proposition}
\newtheorem {corollary}[theorem]{Corollary}
\newtheorem {lemma}[theorem]{Lemma}
\newtheorem {conjecture}[theorem]{Conjecture}
\newtheorem {definition}[theorem]{Definition}

\numberwithin {equation}{section}
\renewcommand {\proof}{{\sc Proof.}\ }

\newcommand {\remark}
{\addtocounter{theorem}{1}{\bf Remark \thetheorem.}\ }
\newcommand {\rklabel}[1]
{\newcounter{#1.se}\setcounter{#1.se}{\value{section}}
 \newcounter{#1.rk}\setcounter{#1.rk}{\value{theorem}}}
\newcommand {\rkref}[1]{\arabic{#1.se}.\arabic{#1.rk}}

\newcommand {\halmos}{$\blacksquare$}

\newcommand {\IC}{\mathbb{C}}
\newcommand {\IN}{\mathbb{N}}                          
\newcommand {\IP}{\mathbb{P}}
\newcommand {\IR}{\mathbb{R}}
\newcommand {\IQ}{\mathbb{Q}}
\newcommand {\IZ}{\mathbb{Z}}

\newcommand {\bfA}{\mathbf A}
\newcommand {\bfI}{\mathbf I}

\newcommand {\A}{\mathcal A}
\newcommand {\B}{\mathcal B}
\newcommand {\C}{\mathcal C}
\newcommand {\D}{\mathcal D}

\newcommand {\F}{\mathcal F}
\newcommand {\G}{\mathcal G}
\renewcommand {\H}{\mathcal H}
\newcommand {\I}{\mathcal I}

\newcommand {\K}{\mathcal K}
\renewcommand {\L}{\mathcal L}
\newcommand {\N}{\mathcal N}

\renewcommand {\P}{\mathcal P}
\newcommand {\Q}{\mathcal Q}
\newcommand {\R}{\mathcal R}
\newcommand {\CS}{\mathcal S}

\newcommand {\T}{\mathcal T}
\newcommand {\U}{\mathcal U}
\newcommand {\V}{\mathcal V}
\newcommand {\X}{\mathcal X}
\newcommand {\Z}{\mathcal Z}

\renewcommand {\a}{\mathfrak a}

\newcommand {\g}{\mathfrak{g}}
\newcommand {\h}{\mathfrak h}
\renewcommand {\ll}{\mathfrak l}
\newcommand {\n}{\mathfrak n}
\newcommand {\p}{\mathfrak p}
\renewcommand {\r}{\mathfrak r}

\renewcommand {\SS}{\mathfrak S}

\newcommand {\sfA}{\mathsf{A}}
\newcommand {\sfB}{\mathsf{B}}
\newcommand {\sfC}{\mathsf{C}}
\newcommand {\sfD}{\mathsf{D}}
\newcommand {\sfE}{\mathsf{E}}
\newcommand {\sfG}{\mathsf{G}}

\newcommand {\ad}{\operatorname{ad}}
\newcommand {\Ad}{\operatorname{Ad}}

\newcommand {\Aut}{\operatorname{Aut}}
\newcommand {\End}{\operatorname{End}}
\newcommand {\Hom}{\operatorname{Hom}}
\newcommand {\Ker}{\operatorname{Ker}}
\newcommand {\Rep}{\operatorname{Rep}}

\newcommand {\id}{\operatorname{id}}
\newcommand {\rk}{\operatorname{rank}}

\newcommand {\gr}{\operatorname{gr}}

\newcommand {\gl}[1]{\mathfrak{gl}_{#1}}
\renewcommand {\sl}[1]{\mathfrak{sl}_{#1}}

\newcommand {\Ug}{U\g}

\newcommand {\Usl}[1]{U\sl{#1}}

\newcommand {\Uhg}{U_{\hbar}\g}
\newcommand {\Uhsl}[1]{U_{\hbar}\sl{#1}}

\newcommand {\kalpha}{\kappa_{\alpha}}

\newcommand {\calpha}{C_{\alpha}}
\newcommand {\nablac}{\nabla_{C}}
\newcommand {\nablak}{\nabla_{\kappa}}

\newcommand{\VA}{V_\A}

\newcommand {\reg}{_{\operatorname{reg}}}
\newcommand {\hreg}{\h\reg}

\newcommand {\Bg}{B_{\g}}

\newcommand {\ie}{{\it i.e., }}
\newcommand {\eg}{{\it e.g.}, }
\newcommand {\etal}{{\it et al.} }
\newcommand {\fd}{finite--dimensional }
\newcommand {\fg}{finitely--generated }
\newcommand {\lhs}{left--hand side }
\newcommand {\rhs}{right--hand side }

\newcommand {\wrt}{with respect to }

\newcommand {\ol}{\overline}
\newcommand {\wt}{\widetilde}
\newcommand {\wh}{\widehat}

\newcommand {\ns}{nested set }
\newcommand {\nss}{nested sets }
\newcommand {\mns}{maximal nested set }
\newcommand {\mnss}{maximal nested sets }
\newcommand {\fns}{fundamental nested set }
\newcommand {\fnss}{fundamental nested sets }
\newcommand {\fmns}{fundamental maximal nested set }
\newcommand {\fmnss}{fundamental maximal nested sets }

\newcommand {\supp}{\operatorname{supp}}
\newcommand {\zsupp}{\operatorname{\mathfrak{z}supp}}

\newcommand {\DCP}{De Concini--Procesi }
\renewcommand {\DJ}{Drinfeld--Jimbo }

\renewcommand {\DH}{Dynkin--Hochschild }

\newcommand {\FZ}{Fomin--Zelevinsky }

\newcommand {\HC}{Harish--Chandra }
\newcommand {\Ho}{Hochschild }
\newcommand {\IH}{Iwahori--Hecke }
\newcommand {\KM}{Kac--Moody }
\newcommand {\KR}{Kirillov--Reshetikhin }
\newcommand {\KZ}{Knizhnik--Zamolodchikov  }

\newcommand {\KKZ}{_{\scriptscriptstyle{\operatorname{KZ}}}}
\newcommand {\CKZ}{{\scriptscriptstyle{\operatorname{CKZ}}}}

\newcommand {\qW}{quantum Weyl group }

\newcommand {\fml}{[\negthinspace[\hbar]\negthinspace]}
\newcommand {\fmll}{[\negthinspace[h]\negthinspace]}
\newcommand {\ICh}{\IC\fml}
\newcommand {\rhoh}{{\rho_h}}

\newcommand {\bin}[3]{\begin{bmatrix}#1\\#2\end{bmatrix}_{#3}}
\newcommand {\half}[1]{\frac{#1}{2}}

\newcommand {\sqnorm}[1]{\frac{(#1,#1)}{2}}
\newcommand {\cor}[1]{\alpha_{#1}^{\vee}}	
\renewcommand {\root}[1]{\alpha_{#1}}
\newcommand {\cow}[1]{\lambda^{\vee}_{#1}}
\newcommand {\odots}[1]{#1\cdots #1}
\newcommand {\tens}[1]{^{\otimes #1}}

\newcommand {\aalpha}{\ul{\alpha}}
\newcommand {\bbeta}{\ul{\beta}}
\newcommand {\ggamma}{\ul{\gamma}}

\newcommand {\Ag}{\A_{\g}}
\newcommand {\old}{^{\scriptscriptstyle{\operatorname{old}}}}

\newcommand {\conn}[3]{{\complement}^{#1\setminus\{#2\}}_{#3}}
\newcommand {\con}[3]{{\complement}^{#1\setminus #2}_{#3}}
\newcommand {\aand}{\qquad\text{and}\qquad}
\newcommand {\oor}{\qquad\text{or}\qquad}
\newcommand {\dDD}{d_{DD}}
\renewcommand {\dh}{d_{H}}
\newcommand {\CDD}{CDD}

\newcommand {\ji}{_{j_{i}}}

\newcommand {\qc}{quasi--Coxeter }
\newcommand {\qca}{quasi--Coxeter algebra }
\newcommand {\qcas}{quasi--Coxeter algebras }
\newcommand {\ba}{bialgebra }
\newcommand {\bas}{bialgebras }
\newcommand {\qba}{quasibialgebra }
\newcommand {\qbas}{quasibialgebras }
\newcommand {\qtqba}{quasitriangular quasibialgebra }
\newcommand {\qtqbas}{quasitriangular quasibialgebras }

\newcommand {\qt}{quasitriangular }

\newcommand {\qcqtqba}{quasi--Coxeter quasitriangular quasibialgebra }
\newcommand {\qcqtqbas}{quasi--Coxeter quasitriangular quasibialgebras }

\newcommand {\Dg}{D_{\g}}

\newcommand {\dD}{d_{D}}
\newcommand {\gD}{\g_{D}}
\newcommand {\hD}{\h_{D}}
\newcommand {\lD}{\ll_{D}}

\newcommand {\AD}{\A_D}
\newcommand {\BD}{B_D}
\newcommand {\ND}{\N_D}

\newcommand {\ZD}{\Z_{D}}

\newcommand {\UgD}{U\gD}
\newcommand {\UhgD}{U_{\hbar}\gD}

\newcommand {\dH}{d_{H}}

\newcommand {\rg}{r_{\g}}
\newcommand {\rD}{r_{\gD}}

\newcommand {\Dh}{_{D,\hbar}}
\newcommand {\DKZ}{_{D}^{\KKZ}}
\newcommand {\ih}{_i^\hbar}
\newcommand {\jh}{_j^\hbar}
\newcommand {\ic}{_{i,\scriptscriptstyle{C}}}

\newcommand {\ul}[1]{\underline{#1}}

\newcommand {\Alt}{\operatorname{Alt}}
\newcommand {\Pent}{\operatorname{Pent}}

\newcommand {\eepsilon}{\varepsilon}
\newcommand {\veps}{\varepsilon}
\newcommand {\oveps}{\ol{\veps}}
\newcommand {\cell}{^{\scriptscriptstyle{\operatorname{cell}}}}
\newcommand {\sign}{\operatorname{sign}}
\newcommand {\bimod}{\operatorname{Bimod}}
\newcommand {\bimodD}{\bimod_D}

\newcommand {\nperp}{\perp\negthickspace\negthickspace\negthickspace\diagup}

\newcommand {\ruL}{\ruLine}
\newcommand {\rdL}{\rdLine}
\newarrow{Line}{-}{-}{-}{-}{-}
\newarrow{equal}{=}{=}{=}{=}{=}
\newarrow{Onecell}{o}{-}{-}{-}{o}

\newcommand {\op}{^{\scriptscriptstyle{\operatorname{op}}}}

\newcommand {\VV}{\mathbf{V}}

\newcommand {\Lpos}{\L^{*}}
\newcommand {\dprime}{^{\prime\prime}}
\newcommand {\BC}{_{B,C}}
\newcommand {\oPhi}{\ol{\Phi}}

\newcommand {\<}{\langle}
\renewcommand {\>}{\rangle}

\newcommand {\msp}{\medspace}

\newcommand {\Deltan}[1]{\Delta^{(#1)}}
\newcommand {\otimesn}[1]{^{\otimes #1}}
\newcommand {\otimesnp}[1]{^{\otimes (#1)}}
\newcommand {\Br}{\mathcal{B}r}

\newcommand {\Veck}{\mathbf{Vec_k}}
\newcommand {\Mod}[1]{\mathbf{Mod_{fd}}(#1)}
\newcommand {\bfF}{\mathbf{F}}

%% file: qctitleabstract.tex
\title
[Quasi--Coxeter algebras and Dynkin diagram cohomology]
{Quasi--Coxeter algebras,
Dynkin Diagram Cohomology
and quantum Weyl groups}
\author[V. Toledano Laredo]{Valerio Toledano Laredo}
\address{
Universit\'e Pierre et Marie Curie--Paris 6\\
Institut de Math\'ematiques de Jussieu, UMR 7586\\
Case 191\\
16 rue Clisson\\
F--75013 Paris}
\email{toledano@math.jussieu.fr}
\begin{abstract}
The author, and independently De Concini, conjectured that
the monodromy of the Casimir connection of a simple Lie
algebra $\g$ is described by the quantum Weyl group operators
of the quantum group $\Uhg$. The aim of this paper, and of its sequel \cite{TL4},
is to prove this conjecture. The proof relies upon the use of
{\it quasi--Coxeter algebras}, which are to generalised braid
groups what Drinfeld's \qtqbas are to the Artin braid groups
$B_n$. Using an appropriate deformation cohomology, we
reduce the conjecture to the existence of a quasi--Coxeter,
\qtqba structure on the enveloping algebra $\Ug$ which
interpolates between the \qca structure underlying the
Casimir connection and the \qtqba structure underlying
the \KZ equations. The existence of this structure will be
proved in \cite{TL4}.
\end{abstract}
\maketitle

%% file: introduction.tex

\section*{Introduction}

Let $\g$ be a complex, simple Lie algebra, $\h\subset\g$ a
Cartan subalgebra and $\Phi\subset\h^*$ the corresponding
root system. For each $\alpha\in\Phi$, let $\sl{2}^\alpha=\<e
_\alpha,f_\alpha,h_\alpha\>\subset\g$ be the corresponding
three--dimensional subalgebra and denote by
$$C_\alpha=\frac{(\alpha,\alpha)}{2}
\left(e_\alpha f_\alpha+f_\alpha e_\alpha+\frac{1}{2}h_\alpha^2\right)$$
its Casimir operator \wrt the restriction to $\sl{2}^\alpha$ of
a fixed non--degenerate, ad--invariant bilinear form $(\cdot,\cdot)$
on $\g$.
Let
$$\hreg=\h\setminus\bigcup_{\alpha\in\Phi}\Ker(\alpha)$$
be the set of regular elements in $\h$, $V$ a \fd $\g$--module,
and consider the following holomorphic connection on the
holomorphically trivial vector bundle $\VV$ over $\hreg$ with
fibre $V$
$$\nablac=
d-\frac{h}{2}\sum_{\alpha\in\Phi}\frac{d\alpha}{\alpha}\cdot C_\alpha
$$
where $d$ is the exterior derivative and $h$ is a complex number. 
The following result is due to the author and J. Millson \cite{MTL}
and was discovered independently by De Concini around 1995
(unpublished) and Felder \etal \cite{FMTV}.
\begin{theorem}
The connection $\nablac$ is flat for any $h\in\IC$.
\end{theorem}

Let $W\subset GL(\h)$ be the Weyl group of $\g$. It is well--known
that $W$ does not act on $V$ in general but that the triple exponentials
$$\wt{s}_i=
\exp(e_{\alpha_i})\exp(-f_{\alpha_i})\exp(e_{\alpha_i})\in GL(V)$$
corresponding to a choice $\alpha_1,\ldots,\alpha_n$ of simple roots
of $\g$ and of $\sl{2}$--triples $e_{\alpha_i},f_{\alpha_i},h_{\alpha_i}
\in\sl{2}^{\alpha_i}$ give rise to an action of an extension $\wt{W}$ of
$W$ by the sign group $\IZ^n$ \cite{Ti}. This action may be used to twist
$(\VV,\nablac)$ into a $W$--equivariant, flat vector bundle on $\hreg$
\cite[\S 2]{MTL}.
One therefore obtains an analytic, one--parameter family of monodromy
representations
$$\mu^h_V:B_W=\pi_1(\hreg/W)\longrightarrow GL(V)$$
of the generalised braid group $B_W$ which, for $h=0$, factors through
the action of $\wt{W}$ on $V$. Considering the deformation parameter
$h$ as formal and setting $\hbar=2\pi ih$, we regard this family as a
single representation
$$\mu_V:B_W\rightarrow GL(V\fml)$$
given by the formal Taylor series of $\mu_V^h$ at $h=0$.\\

Let now $\Uhg$ be the Drinfeld--Jimbo quantum group corresponding
to $\g$ and the inner product $(\cdot,\cdot)$. We regard $\Uhg$ as a
topological Hopf algebra over the ring of formal power series $\IC\fml$.
By a \fd representation of $\Uhg$ we shall mean a $\Uhg$--module
$\U$ which is free and \fg as $\IC\fml$--module.
The isomorphism class of such a representation
is uniquely determined by that of the $\g$--module $\U/\hbar\U$. Lusztig
\cite{Lu}, and independently \KR and Soibelman \cite{KR,So}, constructed
operators $S_1^\hbar,\ldots,S_n^\hbar$ labelled by the simple reflections
$s_i=s_{\alpha_i}$ in $W$, acting on any \fd representation of $\Uhg$.
These operators satisfy the braid relations
$$ S\ih S\jh\cdots=S\jh S\ih\cdots $$
for each $i\neq j$, where the number of terms on each side is equal to
the order of $s_is_j$ in $W$. As a consequence, each \fd representation
$\U$ of $\Uhg$ carries an action of $B_W$ called the {\it quantum Weyl
group} action. Its reduction mod $\hbar$ factors trough the action of the
Tits extension $\wt{W}$ on $\U/\hbar\U$.\\

Let $V$ be a \fd $\g$--module and let $\V$ be a quantum deformation
of $V$, that is a \fd $\Uhg$--module such that $\V/\hbar\V\cong V$ as
$\g$--modules. The following conjecture was formulated in \cite{TL1,TL2}
and independently by De Concini around 1995 (unpublished).

\begin{conjecture}\label{co:monodromy}
The monodromy action $\mu_V$ of $B_W$ on $V\fml$ is equivalent to
its quantum Weyl group action on $\V$.
\end{conjecture}

We note in passing the following interesting, and immediate consequence
of the above conjecture and of the fact that the \fd $\Uhg$--modules and
operators $S\ih$ are defined over $\IQ\fml$

\begin{corollary}
The monodromy representation $\mu_V$ is defined over $\IQ\fml$.
\end{corollary}

Conjecture \ref{co:monodromy} is proved in \cite{TL1} for all representations
of $\g=\sl{n}$ and in
\cite{TL2} for a number of pairs $(\g,V)$ including vector
and spin representations of classical Lie algebras and the adjoint representation
of all complex, simple Lie algebras. Its semi--classical analogue is proved for
all $\g$ in \cite{Bo1,Bo2}.\\

The aim of this paper, and of its sequel \cite{TL4}, is to give a proof of
this conjecture for any complex, simple Lie algebra $\g$. The strategy
we follow is very much inspired by Drinfeld's proof of the equivalence
of the monodromy of the \KZ (KZ) equations and the $R$--matrix
representations coming from $\Uhg$ \cite{Dr4}. It proceeds along
the following lines.
\begin{enumerate}
\item Define a suitable category of algebras carrying representations of
the generalised braid group $B_W$ on their finite--dimensional modules.
We call these algebras {\it quasi--Coxeter algebras}. They are the
analogues for $B_W$ of what Drinfeld's \qt \qbas are for the Artin
braid groups $B_n$.
\item Show that the monodromy representations $\mu_V$
arise from a \qca structure on the 
enveloping algebra $\Ug\fml$.
\item Show that the quantum Weyl group representations of $B_W$
arise from a \qca structure on the quantum group $\Uhg$.
\item Show that this \qca structure on $\Uhg$ is (cohomologically)
equivalent to one on $\Ug\fml$.
\item Construct a complex which controls the deformations of \qca
structures. We call the corresponding cohomology {\it Dynkin
diagram cohomology} $HD^*$. The infinitesimal of a \qca structure
on $A$ defines a canonical class in $HD^2(A)$.
\item Show, using the Dynkin complex of $\Ug$, that there is, up
to isomorphism at most one \qca structure on $\Ug\fml$ having
prescribed local monodromies.
\end{enumerate}

Steps (i)--(v) are carried out in Part I of this paper. Step (vi) however
hopelessly fails since it turns out that, for $\g$ of rank greater than
one, 
$HD^2(\Ug)$ is infinite--dimensional. To remedy this, one must
rigidify matters by taking into account the bialgebra structures on $\Uhg$
and $\Ug\fml$. This may perhaps be guessed from the fact that both
the quantum Weyl group operators $S\ih$ and the local monodromies
$S_i^C=\wt{s}_i\cdot\exp(\hbar C_{\alpha_i}/2)$ of the Casimir
connection $\nablac$ satisfy strikingly similar coproduct identities, 
namely
\begin{align*}
\Delta(S\ih)&=R_i^{21}\cdot S\ih\otimes S\ih\\
\Delta(S_i^C)&=
\exp(\hbar t_i)\cdot
S_i^C\otimes S_i^C
\end{align*}
where $R_i$ is the universal $R$--matrix of $\Uhsl{2}^{\alpha_i}\subset\Uhg$
and
$$t_i=(\Delta(C_{\alpha_i})-C_{\alpha_i}\otimes 1-1\otimes C_{\alpha_i})/2$$
Since the failure of $S\ih$ and $S_i^C$ to be group--like elements involves
R--matrices, a proper study of \qc \bas must in fact be concerned with \qc
\qt quasibialgebras, that is \bas which carry representations of both Artin's
braid groups and the group $B_W$ on the tensor products of their \fd modules.
We carry out the analogues of steps (i) and (v) for these algebras in Part II of this
paper and of steps (iii) and (iv) in Part III. 
We then show that, up 
to isomorphism, there exists a unique \qcqtqba structure on $\Ug\fml$ having
prescribed local monodromies as that coming from $\Uhg$ (step (vi)). The final
step needed to complete the proof of conjecture \ref{co:monodromy}, namely the
fact that the monodromy of the 
connection $\nablac$ and of the KZ connection
for $\g$ fit within a \qcqtqba structure on $\Ug\fml$ will be proved in \cite{TL4}.

\subsection*{Outline of the paper}

We turn now to an outline of the contents of the paper, referring the reader
to the introductory paragraphs of each section for more details. We begin in
section \ref{se:DCP} by reviewing the \DCP theory of asymptotic zones for
connections of KZ--type which provides a concise, combinatorial description
of their monodromy. This description forms, together with Drinfeld's theory
of \qt quasibialgebras, revisited through the author's duality between the
connection $\nablac$ for $\sl{n}$ and the KZ connection for $\sl{k}$ (see
\cite{TL1} and \S\ref{ss:qba examples}), the basis underpinning the definition
of a \qc algebra given in section \ref{se:qcoxeter}. Such algebras have a type
determined by a connected graph $D$ with labelled edges. For the examples
relevant to us, $D$ is the Dynkin diagram of the Lie algebra $\g$ but in most
of the paper we merely assume that $D$ is a connected graph and work in
this greater generality. Just as MacLane's coherence theorem for monoidal
categories is best proved using Stasheff's associahedra $K_n$ \cite{St1},
the most compact definition of a \qca has its relations labelled by the
two--dimensional faces of a regular CW--complex $\AD$, which is
defined and studied in section \ref{se:associahedron}. We call $\AD$
the {\it \DCP associahedron} since, when $D$ is the Dynkin diagram
of $\g$, $\AD$ is naturally realised inside their wonderful model of
$\hreg$ \cite{DCP2} and, when $\g=\sl{n+1}$, it coincides with $K_n
$\footnote{while
this paper was being written, Carr and Devadoss posted the preprint
\cite{CD} where the same $CW$--complex is introduced under the
name {\it graph--associahedron} and proved to be, as in our \S\ref
{ss:realisation}, a convex polytope.}. Section
\ref{se:qcexamples} describes several examples of \qc algebras. In
section \ref{se:dynkin}, we define the Dynkin complex of a \qca and
show that it controls its deformation theory. In section \ref{se:qcqtqb}
we define \qcqtqbas and show in section \ref{se:dynkinbi} that their
deformation theory is controlled by a suitable bicomplex, which we
call the {\it \DH} bicomplex. In section \ref{se:transfer}, we show that
Drinfeld's $R$--matrix and the quantum Weyl groups operators endow
the quantum group $\Uhg$ with the structure of a \qcqtqba and that
this structure may be cohomologically transferred to one on $\Ug
\fml$. Finally, section \ref{se:rigidity} is devoted to the proof of the rigidity
of \qcqtqba structures on the enveloping algebra $\Ug$. This result relies
on Drinfeld's rigidity of \qtqba deformations of $\Ug\fml$ \cite{Dr3} and
on the essential uniqueness of solutions of a certain relative twist
equation which is proved in \cite{TL3}.\\



%% file: DCPfinal.tex


\section{Asymptotic zones for hyperplane arrangements}\label{se:DCP}

This section is an exposition of the paper \cite{DCP2} and parts of \cite{DCP1}.
All results stated herein are due to De Concini and Procesi, with the minor exception
of the terminology of \S \ref{ss:DCP associators}--\S \ref{ss:elementary associators}
and \S \ref{ss:adapted family}--\S \ref{ss:forgetful} and of the results of \S \ref
{ss:support}--\S \ref{ss:forgetful}, which are however implicit in or immediate
consequences of \cite{DCP2}. The presentation in \S \ref{ss:chambers}--\S
\ref{ss:forgetful} is a little more general  than in \cite{DCP2} since we deal
with real arrangements endowed with a simplicial chamber rather than Coxeter
ones, but the proofs in \cite{DCP2} carry over verbatim to this context and
are therefore mostly omitted.\\

Let $\A$ be a hyperplane arrangement in a complex, \fd vector space $V$.
The aim of this section is to describe good solutions $\Psi$ of the {\it holonomy
equations} $\nabla\Psi=0$ on $\VA=V\setminus\A$. These are a generalisation
of the \KZ (KZ) equations, which are defined for the Coxeter arrangements of
type $\sfA_n$, to general hyperplane arrangements. The solutions in question are
universal, \ie take values in the {\it holonomy algebra} $A$ of the arrangement
(see \S \ref{ss:holonomy algebra}--\S \ref{ss:monodromy of holonomy} for
definitions) and have prescribed asymptotic behaviour, that is are such that
the monodromy of suitable commuting elements of $\pi_1(\VA)$ has a particularly
simple form. They generalise the ones constructed by Drinfeld for the KZ equations
\cite{Dr3} and by Cherednik for his generalisation of the KZ equations to Coxeter
arrangements \cite{Ch1,Ch2}. They are described in \S\ref{ss:residues}--\S \ref
{ss:fundamental solutions} in terms of local coordinates on a {\it wonderful model}
of $V$ that is, a smooth variety $Y_X$, proper over $V$ and in which the preimage
of $\A$ is a divisor with normal crossings. The construction and coordinatization of
$Y_X$ are reviewed in \S \ref{ss:wonderful model}--\S \ref{ss:charts}. In \S \ref
{ss:chambers}--\S \ref{ss:DCP associators}, we focus on real arrangements. This
allows one to consistently define multiplicative constants, specifically elements
of the holonomy algebra $A$, which relate different solutions. These are a
generalisation of Drinfeld's associator, which arises for the Coxeter arrangement
of type $\sfA_2$, and allow one to concisely describe the monodromy of $\nabla$.
The main properties of these \DCP associators, most notably their inductive
structure, are described in \S \ref{ss:elementary associators}--\S \ref{ss:forgetful}.
Finally, in \S \ref{ss:coxeter}, we specialise further to the case of Coxeter
arrangements and obtain that the associators satisfy, together with the local
monodromies of $\nabla$, the braid relations which give Brieskorn's presentation
of the corresponding generalised braid group.

\subsection{The holonomy algebra of a hyperplane arrangement}
\label{ss:holonomy algebra}

Let $V$ be a finite--dimensional, complex vector space and $\A$ a finite collection
of linear hyperplanes in $V$. Choose an equation $x\in H^\perp\setminus\{0\}$ for
each $H\in\A$ and let $X=\{x\}\subset V^*$ be the corresponding collection of linear
forms. Consider a connection on $\VA=V\setminus\A$ of the form
$$\nabla=d-\sum_{x\in X}\frac{dx}{x}\cdot t_x$$
Following Chen \cite{Chn1}, we do not regard the $t_x$ as acting on a vector space
but rather as formal variables, it being understood that any \fd representation
$\rho:F\rightarrow\End(U)$ of the free, associative algebra $F=\IC\<t_x\>_{x\in
X}$ generated by the $t_x$ gives rise to a holomorphic connection
\begin{equation}\label{eq:nabla rho}
\nabla_\rho=d-\sum_{x\in X}\frac{dx}{x}\cdot\rho(t_x)
\end{equation}
on the holomorphically trivial vector bundle over $\VA$ with fibre $U$.\\

By \cite{Ko2}, $\nabla$ is flat if, and only if the following relations hold for any
two--dimensional subspace $B\subset V^*$ spanned by elements of $X$ and
$x\in X\cap B$,
\begin{equation}\label{eq:rank 2}
[t_x,\sum_{y\in X\cap B}t_y]=0
\end{equation}
Let $I\subset F$ be the two--sided ideal generated by these relations. $I$
is homogeneous \wrt the grading on $F$ for which all $t_x$ have degree 1. 
The quotient $F/I$ is the universal enveloping algebra of the graded Lie
algebra with generators $t_x$, $x\in X$ defined by \eqref{eq:rank 2}. The
completion $A$ of $F/I$ with respect to its grading is a topological Hopf
algebra called the {\it holonomy algebra} of the arrangement $\A$. The
{\it holonomy equations} of $\A$ are the equations
\begin{equation}\label{eq:holonomy}
d\Psi=\sum_{x\in X}\frac{dx}{x}t_x\Psi
\end{equation}
for a locally defined function $\Psi:\VA\rightarrow A$. Such a $\Psi$ is
necessarily {\it holomorphic}, \ie such that each component $\Psi_m$,
$m\in\IN$ \wrt the $\IN$--grading on $A$ is a holomorphic function with
values in the \fd vector space $A_m$. Moreover, $\Psi$ takes values in
the group $A^\times$ of invertible elements of $A$ if, and only if its degree
zero term $\Psi_0:\VA\rightarrow A_0\cong\IC$ does not vanish. Since
\eqref{eq:holonomy} implies that $d\Psi_0=0$, this is the case if, and
only if $\Psi_0(v)\neq 0$ for some $v\in\VA$. Finally, since the $t_x$
are primitive elements of $A$, $\Psi$ takes values in the group like
elements of $A$ if, and only if $\Delta(\Psi(v))=\Psi(v)\otimes\Psi(v)$
for some $v\in\VA$. We denote by $N$ the group of group like elements
of $A$ with degree zero term equal to 1 and call a solution $\Psi$ of
the holonomy equations {\it unipotent} if it takes values in $N$.\\

Clearly, the connection \eqref{eq:nabla rho} determined by a \fd representation
$\rho:F\rightarrow\End(U)$ is flat if, and only if $\rho$ factors through $F/I$.
To avoid convergence issues, let $h$ be a formal variable and note that, due
to the homogeneity of the relations \eqref{eq:rank 2}, the representation
\begin{equation}\label{eq:rho_h}
\rho_h:F\rightarrow\End(U[h]),\qquad
\rho_h(t_x)=h\rho(t_x)
\end{equation}
factors through $F/I$ if, and only if, $\rho$ does. When that is the case, $\rho_h$
extends to a representation of $A$ on $U\fmll$ which we will denote by the same
symbol. Moreover, any invertible solution $\Psi$ of \eqref{eq:holonomy} gives rise
to a fundamental solution $\Psi_\rhoh=\rhoh(\Psi)$ of $\nabla_\rhoh\Psi_\rhoh=0$
with values in $\End(U)\fmll$.

\subsection{Monodromy of the holonomy equations}
\label{ss:monodromy of holonomy}

Throughout this paper, we adopt the convention that, for a topological space $X$,
the composition $\zeta\gamma$ of $\zeta,\gamma\in\pi_1(X,x_0)$ is given by $
\gamma$ {\it followed} by $\zeta$, so that the holonomy of a flat vector bundle
$(\V,\nabla)$ on $X$ at $x_0$ is a homomorphism $\pi_1(X,x_0)\rightarrow GL
(\V_{x_0})$.\\

Fix a basepoint $v_0\in\VA$ and let $(\wt{\VA},\wt{v_0})\stackrel{p}{\longrightarrow}
(\VA,v_0)$ be the pointed universal cover of $(\VA,v_0)$. $\wt{\VA}$ is endowed
with a canonical right action of the fundamental group $\pi_1(\VA,v_0)$ by deck
transformation. Any solution $\Psi$ of the holonomy equations \eqref{eq:holonomy}
defined in a neighborhood of $v_0$ lifts uniquely to, and will be regarded as, a global,
$A$--valued solution of $p^*\nabla\Psi=0$ on $\wt{\VA}$. If $\gamma\in\pi_1(\VA,v_0)$,
then $\gamma\bullet\Psi(\wt{v})=\Psi(\wt{v}\gamma)$ is another such solution. Thus,
if $\Psi$ is invertible, then
\begin{equation}\label{eq:mu_Psi}
\mu_\Psi(\gamma)=\Psi^{-1}\cdot\gamma\bullet\Psi
\end{equation}
is a locally constant, and therefore constant function on $\wt{\VA}$.

\begin{proposition}\label{pr:mu monodromy}\hfill
\begin{enumerate}
\item The map $\gamma\rightarrow\mu_\Psi(\gamma)$ is a homomorphism
$\pi_1(\VA,v_0)\rightarrow A$.
\item If $\Psi$ is unipotent, $\mu_\Psi$ takes values in $N$.
\item If $\Psi'=\Psi\cdot K$,  $K\in A^\times$ is another invertible solution,
then $$\mu_{\Psi'}=\Ad(K^{-1})\circ\mu_\Psi$$
\end{enumerate}
\end{proposition}

Note that if $\rho:F/I\rightarrow\End(U)$ is a \fd representation and $\rhoh:A
\rightarrow\End(U\fmll)$ is given by \eqref{eq:rho_h}, the composition $\rhoh
\circ\mu_\Psi:\pi_1(\VA,v_0)\rightarrow GL(U\fmll)$ is the monodromy
representation of the connection $\nabla_\rhoh$ expressed in the
fundamental solution $\Psi_\rhoh=\rhoh(\Psi)$.

\subsection{Chen--Kohno isomorphisms}\label{ss:Chen}
 
Set $\pi=\pi_1(\VA,v_0)$, let $\IC[\pi]$ be the group algebra of $\pi$ and
$J$ its augmentation ideal, that is the kernel of the counit $\veps:\IC[\pi]
\rightarrow\IC$ given by $\veps(\gamma)=1$ for any $\gamma\in\pi$.
Let $\wh{\IC[\pi]}$ be the prounipotent completion of $\IC[\pi]$, \ie the
inverse limit
$$\wh{\IC[\pi]}=
\lim_{\substack{\longleftarrow\\k\rightarrow\infty}}\IC[\pi]/J^k$$
Note that $\pi$ embeds in $\wh{\IC[\pi]}$ if, and only if $\IC[\pi]$ does and
that this is the case if, and only if $\pi$ is residually torsion free nilpotent
\cite[prop. 2.2.1]{Chn2}.\\

Let $\Psi:\wt{\VA}\rightarrow A$ be an invertible solution of the holonomy
equations \eqref{eq:holonomy}. Since their reduction by the ideal $A_+
\subset A$ consisting of elements of positive degree is the differential
equation $d\Psi=0$, the diagram
$$\begin{diagram}[height=2.5em,width=2.5em]
\IC[\pi]		&\rTo^{\mu_\Psi}	&A\\
\dTo^{\veps}	&				&\dTo\\
\IC			&\requal			&A/A_+
\end{diagram}$$
is commutative. Thus, $\mu_\Psi$ maps $J$ to $A_+$ and therefore factors
through an algebra homomorphism
$$\wh{\mu_\Psi}:\wh{\IC[\pi]}\longrightarrow A$$
which is a homomorphism of topological Hopf algebras if $\Psi$ is unipotent.

\begin{theorem}[Chen--Kohno]\label{th:Chen}
$\wh{\mu_\Psi}$ is an isomorphism.
\end{theorem}

A simple proof of theorem \ref{th:Chen} is given at the end of this section in
\S \ref{app:Chen}.

\subsection{The wonderful model $Y_X$}\label{ss:wonderful model}

Let $\L$ be the (atomic) lattice of subspaces of $V^*$ spanned by subsets of
$X$ and set $\Lpos=\L\setminus\{0\}$.

\begin{definition}
A decomposition of $U\in\Lpos$ is a family $U_1,\ldots,U_k$ of subspaces of $U$
lying in $\Lpos$ such that, 
for every subspace $W\subseteq U$ with $W\in\Lpos$ one has $W\cap U_i\in\L$
for every $i=1,\ldots,k$ and $$W=(W\cap U_1)\odots{\oplus}(W\cap U_k)$$
An element $U\in\Lpos$ is called reducible if it possesses such a decomposition
with $k\geq 2$ summands and irreducible otherwise.
\end{definition}

Any $U\in\Lpos$ possesses a unique decomposition into irreducible summands.
Let $\I\subset\Lpos$ be the collection of irreducible subspaces. For any $B\in\I$,
let
$$B^\perp=\bigcap_{x\in X\cap B}x^\perp\subset V$$
be the subspace orthogonal to $B$ and note that the rational map $\pi_B:V
\rightarrow V/B^\perp\rightarrow\IP(V/B^\perp)$ is regular outside $B^\perp$.

\begin{definition}
$Y_X$ is the closure of the image of $\VA$ under the embedding\footnote
{the notation $Y_X$ follows \cite{DCP2}. It should be manifest however that
$Y_X$ and all constructions to follow only depend upon the arrangement $
\A$ and not on the actual choice of the set $X$ of linear forms describing it.}
$$\VA\longrightarrow V\times\prod_{B\in\I}\IP(V/B^\perp)$$
\end{definition}

\begin{theorem}\hfill
\begin{enumerate}
\item $Y_X$ is an irreducible, smooth algebraic variety.
\item The natural projection $\pi:Y_X\rightarrow V$ is proper, surjective
and restricts to an isomorphism $Y_X\setminus\pi^{-1}(\A)\rightarrow\VA$.
\item $\pi^{-1}(\A)\subset Y_X$ is a divisor with normal crossings.
\end{enumerate}
\end{theorem}

\remark
The points of $Y_X$ may be parametrised by adapting Fulton and MacPherson's
notion of screens for configurations of points \cite[\S 1]{FM} as follows. Let $v\in V$
such that $x(v)=0$ for at least one $x\in X$, let $B$ be the maximal element of $
\Lpos$ such that $v\in B^\perp$ and let $B=B_1\odots{\oplus}B_k$ be its irreducible
decomposition. Define recursively a {\it sequence of screens} for $v$ to be given by
a non--zero vector $v_i$ in each $V/B_i^\perp$, defined up to multiplication by a scalar,
and, whenever $x(v_i)=0$ for at least an $x\in X\cap B_i$, a sequence of screens
for $v_i$ relative to the hyperplane arrangement on $V/B_i^\perp$ defined by $X
\cap B_i$. Then, a point $y\in Y_X$ is readily seen to determine a sequence of
screens for $v=\pi(y)$ and to be uniquely determined by it. Moreover, any sequence
of screens for a vector $v\in V$ arises in this way.

\subsection{Geometry and combinatorics of the divisor $\D$}\label{ss:divisor}

The following notion is crucial in describing the combinatorics of the
divisor $\D=\pi^{-1}(\A)$.

\begin{definition}
A family $\{U_i\}_{i\in I}$ of irreducible elements of $\Lpos$ is nested if,
for any subfamily $\{U_j\}_{j\in J}$ of pairwise non comparable elements,
the subspaces $U_j$ are in direct sum and are the irreducible summands
of $\bigoplus_{j\in J}U_j$.
\end{definition}

\begin{theorem}\label{th:divisor}\hfill
\begin{enumerate}
\item The irreducible components of $\D$ are smooth and labelled by
the irreducible elements of $\Lpos$, with $\D_B$ the unique component
of $\D$ such that $\pi(\D_B)=B^\perp$.
\item A family $\{\D_B\}_{B\in\CS}$ of such components has a non--zero
intersection if, and only if, $\CS$ is nested. In that case, the intersection
is transversal and irreducible.
\end{enumerate}
\end{theorem}

\remark
It is easy to see that the collection $\CS_y$ of subspaces appearing in the
sequence of screens attached to a point $y\in Y_X$ is nested. If $\CS\subset
\I$ is nested, the intersection $\bigcap_{B\in\CS}\D_B$ consists of those $y
\in Y_X$ such that $\CS_y$ contains $\CS$, its open locus being the set of
such $y\in Y$ for which $\CS_y=\CS$.

\subsection{Some properties of nested sets}\label{ss:lemma}

Let $\CS$ be a nested set, $B\in\CS$ and let $C_1,\ldots,C_m$ be the
maximal elements of $\CS$ properly contained in $B$. By nestedness,
the $C_i$ are in direct sum and
\begin{equation}\label{eq:i_S(B)}
i_\CS(B)=C_1\odots{\oplus}C_m
\end{equation}
is properly contained in $B$ since $B$ is irreducible. Set
$$n(B;\CS)=\dim(B/i_\CS(B))$$
so that $n(B;\CS)\geq 1$ and
$$n(\CS)=\sum_{B\in\CS}(n(B;\CS)-1)$$

Induction on $|\CS|$ readily shows the following

\begin{proposition}\label{pr:n(F)}
Let $B_1,\ldots,B_k$ be the maximal elements of $\CS$. Then, $$n(\CS)
=\sum_{i=1}^k\dim B_i-|\CS|$$
\end{proposition}

\begin{proposition}\label{pr:max nested}
Let $\CS$ be a maximal nested set. Then,
\begin{enumerate}
\item $n(B;\CS)=1$ for any $B\in\CS$.
\item The maximal elements of $\CS$ are the irreducible components
of the subspace $\X\subseteq V^*$ spanned by $X$.
\item $|\CS|=\dim\X$.
\end{enumerate}
\end{proposition}

\subsection{Coordinate charts on $Y_X$}\label{ss:charts}

\subsubsection{} 

To coordinatise $Y_X$, it is convenient to assume that the arrangement $\A$
is {\it essential}, that is satisfies
\begin{equation}\label{eq:essential}
\bigcap_{H\in\A}H=0
\end{equation}
so that $X$ spans $V^*$. This assumption, which we henceforth make, is not
truly restrictive since any linear isomorphism $V\cong V/X^\perp\oplus X^\perp$
induces an isomorphism $Y_X\cong Y_{\ol{X}}\times X^\perp$, where $\ol{X}
\subset(V/X^\perp)^*$ is the collection of linear forms induced by $X$. Note
that, by proposition \ref{pr:max nested}, \eqref{eq:essential} implies that
$$|\CS|=\dim V^*$$ for any \mns $\CS$. By theorem \ref{th:divisor}, the
corresponding intersection $\bigcap_{B\in\CS}\D_B$ therefore consists
of a single point. Such a point at infinity will be denoted $y_\CS$.

\subsubsection{} 

\begin{definition}
A basis $b$ of $V^*$ is adapted to a nested set $\CS$ if, for any $B\in\CS$,
$b\cap B$ is a basis of $B$\footnote{departing a little from \cite[\S 1.1]{DCP2},
but consistently with \cite[\S 1.3]{DCP1}, we do not assume that the elements
of an adapted basis lie in $X$, but denote them nonetheless by the letter $x$.}.
\end{definition}

If $\CS$ is a maximal nested set and $b$ an adapted basis, proposition \ref
{pr:max nested} (i) implies that, for any $B\in\CS$, there exists a unique $x
\in b$ such that
$$x\in B\setminus\sum_{\substack{C\in\CS,\\C\subsetneq B}}C$$
Such an element will be denoted $x_B$. Clearly, $x_B=x_C$ implies that $B
=C$ and any $x\in b$ is of this form since $|\CS|=\dim V^*$.

\subsubsection{} 

Let $U=\IC^\CS$ with coordinates $\{u_B\}_{B\in\CS}$ and let $\rho:U
\rightarrow V$ be the map defined in the coordinates $\{x_B\}_{B\in\CS}$
on $V$ by
$$x_B=\prod_{C\supseteq B}u_C$$
$\rho$ is birational, with inverse
\begin{equation}\label{eq:blowup}
u_B=\left\{\begin{array}{cl}
x_B&\text{if $B$ is maximal in $\CS$}\\
\frac{x_B}{x_{c(B)}}&\text{otherwise}
\end{array}\right.
\end{equation}
where $c(B)$ is the smallest element of $\CS$ properly containing $B$.
It restricts to an isomorphism between the open set of $U$ where all
the coordinates $u_B$ are different from zero and the open set of $V$
where all the coordinates $x_B$ are different from zero. Moreover,
$\rho$ maps the coordinate hyperplane defined by $u_B=0$ into the
subspace $B^\perp\subset V$.

\subsubsection{} 

For any subset $Z$ of $V^*$ containing a non--zero vector, the set of
elements of $\CS$ containing $Z$ is linearly ordered. Denote by $p_
\CS(x)$ its infimum if it is non--empty and set $p_\CS(Z)=V^*$
otherwise.

\begin{lemma}\label{le:p_S}\hfill
\begin{enumerate}
\item If $x\in X$, then $p_\CS(x)\in\CS$.
\item More generally, if $B\in\Lpos$ is irreducible, then $p_\CS(B)\in\CS$
and there is an $x\in B\cap X$ such that $p_\CS(B)=p_\CS(x)$.
\end{enumerate}
\end{lemma}
\proof
(i) $x$ lies in one of the irreducible components of the span of $X$. Since these
are contained in $\CS$ by proposition \ref{pr:max nested}, the set of elements
of $\CS$ containing $x$ is not empty.
(ii) Let $C_1,\ldots,C_m\in\CS$ be maximal among the $p_\CS(x)$, $x\in X\cap
B$. By nestedness, the $C_i$ are in direct sum. Since $B\subset C_1\odots{\oplus}
C_m$, the irreducibility of $B$ implies that $B\subset C_i$ for some $i$ \halmos

\subsubsection{} 

\begin{lemma}\label{le:P_x}
Let $x\in V^*\setminus\{0\}$ be such that $B=p_\CS(x)$ lies in $\CS$.
Then, $x=x_B\cdot P_x$ where $P_x$ is a polynomial in the variables
$u_C$, $C\subsetneq B$ such that $P_x(0)\neq 0$.
\end{lemma}
\proof
Since $x_C$, $C\subseteq B$ is a basis of $B\ni x$, we have
\begin{equation*}
\begin{split}
x
&=
\sum_{C\subseteq B} \alpha_C x_C
=
x_B(\alpha_B+\sum_{C\subsetneq B}\alpha_C\frac{x_C}{x_B})\\
&=
x_B(\alpha_B+\sum_{C\subsetneq B}\alpha_C\prod_{C\subseteq D\subsetneq B}u_D)
=
x_B P_x
\end{split}
\end{equation*}
where $P_x(0)=\alpha_B\neq 0$  or else $p_\CS(x)\subsetneq B$ \halmos\\

The following result explains the relevance of the polynomials $P_x$.

\begin{proposition}\label{pr:regularity}
Let $B\in\L^*$ be irreducible and $x\in B\cap X$ be such that $p_\CS(x)=p
_\CS(B)$ as in lemma \ref{le:p_S}. Then, the rational map
$$U\stackrel{\rho}{\longrightarrow}V\longrightarrow \IP(V/B^\perp)$$
restricts to a regular morphism on $U\setminus\{P_x=0\}$.
\end{proposition}
\proof Let $A=p_\CS(B)$. For any $y\in X\cap B$,
$$y=x_{p_\CS(y)}P_y=
x_A\prod_{p_\CS(y)\subseteq C\subsetneq A}u_C\cdot P_y=:x_A P_y^B$$
Complete $x$ to a basis $x=x_1,x_2,\ldots,x_n$ of $B$ whose elements lie
in $X$. Then, in the corresponding homogeneous coordinates on $\IP(V/B^
\perp)$, the composition above is given by
$$u\longrightarrow[x_A P_x(u),x_A P_{x_2}^B(u),\ldots,x_A P_{x_n}^B(u)]=
[P_x(u),P_{x_2}^B(u),\ldots,P_{x_n}^B(u)]$$
and is therefore regular on $\{P_x(u)\neq 0\}$ \halmos

\subsubsection{} 

\begin{definition}
Let $\U_\CS^b$ be the complement in $U$ of the zeros of the polynomials
$P_x$ as $x$ varies in $X$.
\end{definition}

By proposition \ref{pr:regularity}, the rational map $U\rightarrow V\times\prod
_{B\in\CS}\IP(V/B^\perp)$ restricts to a regular map $j_\CS^b$ on $\U_\CS
^b$. Since, for any $x\in X$
$$x=x_{p_\CS(x)}P_x=\prod_{B\supseteq p_\CS^b(x)}u_B\cdot P_x$$
$\rho$ maps the complement of the coordinate hyperplanes in $\U_\CS^b$
to $\VA$ so that $j_\CS^b(\U_\CS^b) \subset Y_X$.

\begin{proposition}
$j_\CS^b$ is an open embedding $\U_\CS^b\hookrightarrow Y_X$.
\end{proposition}

We shall henceforth identify $\U_\CS^b$ with its image in $Y_X$ under the
embedding $j_\CS^b$ and regard the functions $u_B$, $B\in\CS$ as local
coordinates on $Y_X$.

\subsubsection{} 

For every maximal nested set $\CS$ and $B\in\CS$, choose a basis $b_B$ of $B$
consisting of elements not lying in any $C\in\CS$, $C\subsetneq B$. Choosing an
element $x$ from each $b_B$, as $B$ varies in $\CS$, yields a basis of $V^*$
adapted to $\CS$. Varying the choice of such $x$ then gives rise to a finite set
$\B_\CS$ of such bases.

\begin{theorem}\label{th:covering}\hfill
\begin{enumerate}
\item $Y_X$ is covered by the open sets $\U_\CS^b$ as $\CS$ varies amongst
the maximal nested sets and $b$ varies in $\B_\CS$.
\item The intersection $\D_B\cap \U_\CS^b$ is non--zero if, and only if $B\in\CS$.
When that is the case, it is given by the equation $u_B=0$.
\end{enumerate}
\end{theorem}

Let $\CS$ be a \mns and $y_\CS=\bigcap_{B\in\CS}\D_B$ the corresponding point
at infinity in $Y_X$. By theorem \ref{th:covering}, $y_\CS$ lies in $\U_{\CS'}^b$ if,
and only if $\CS'=\CS$ and, when that is the case is the point with coordinates
$u_B=0$, $B\in\CS$.\\

\remark
Although the open set $\U_\CS^b$ depends upon the choice of the adapted
basis $b$, the union $\U_\CS=\bigcup_{b\in\B_\CS}\U_\CS^b$ has an intrinsic
description as the set of those points $y\in Y_X$ such that the subspaces involved
in the sequence of screens corresponding to $y$ are all contained in $\CS$\footnote
{one can show nonetheless that, if $b\subset X$, the open set $\U_\CS^b$ is independent of $b$ and equal to $\U_\CS$ \cite[\S 1.1]{DCP2}. In order to use
adapted families however (see \S \ref{ss:adapted family}) we shall need adapted
bases whose elements do not necessarily lie in $X$.}.

\subsection{Some properties of residues}
\label{ss:residues}

For any $B\in\Lpos$, set
$$t_B=\sum_{x\in X\cap B} t_x=\sum_{i=1}^m t_{B_i}$$
where $B_1,\ldots,B_m$ are the irreducible components of $B$.
\begin{lemma}\label{le:[t_B,t_C]}
Let $B,C$ be irreducible and nested, then $[t_B,t_C]=0$.
\end{lemma}
\proof Assume first that $B$ and $C$ are not comparable. By nestedness, $B\cap C
=0$ and $B,C$ are the irreducible summands of $B\oplus C$. It follows from this that
if $x\in B\cap X$ and $y\in C\cap X$, the subspace $\IC x\oplus\IC y\in\Lpos$ is reducible
and therefore cannot contain any other elements of $X$ other than $x$ and $y$. The
relations \eqref{eq:rank 2} therefore imply that $t_x$ and $t_y$ commute so that
$$[t_B,t_C]=\sum_{\substack{x\in B\cap X\\y\in C\cap X}}[t_x,t_y]=0$$
Assume now that $B\subset C$. Let $x\in B\cap X$ and define an equivalence relation
on $C\cap X\setminus x$ by $y_1\sim y_2$ if $y_1$ and $y_2$ span the same line
in $C/\IC x$. Let $\Xi$ be the set of equivalence classes, so that
$$t_C=t_x+\sum_{[y]\in\Xi}\sum_{y\in [y]}t_y$$
For a given equivalence class $[y]\in\Xi$, the span of $x$ and $y$, with $y\in[y]$, is
a two--dimensional subspace $C_{[y]}\subset C$ such that $C_{[y]}\cap X=\{x\}\cup
[y]$. Thus, $[t_x,\sum_{y\in[y]}t_y]=0$ by \eqref{eq:rank 2} and $[t_B,t_C]=0$
\halmos

\subsection{Fundamental solutions of the holonomy equations}
\label{ss:fundamental solutions}

For any \ns $\CS$ and $B\in\CS$, set
\begin{equation}\label{eq:R_B^S}
R_B^\CS=
\sum_{\substack{x\in X\cap B,\\p_\CS(x)=B}}t_x=
t_B-t_{i_\CS(B)}
\end{equation}
where $i_\CS(B)$ is defined by \eqref{eq:i_S(B)}.

\begin{theorem}\label{th:fundamental solutions}\hfill
\begin{enumerate}
\item The pull--back to $Y_X$ of $\nabla$ is a flat connection with
logarithmic singularities on the divisor $\D$.
\item Its residue of the irreducible component $\D_B$ of $\D$ is equal to $t_B$.
\item Let $\CS$ be a \mns and $b$ an adapted basis of $V^*$. Then, for any
simply--connected open set $\V\subset\U_\CS^b$ containing $y_\CS$, there
exists a unique holomorphic function $H_\CS^b:\V\rightarrow A$ such that
$H_\CS^b(y_\CS)=1$ and, for every determination of $\log(x_B)$, $B\in\CS$,
the multivalued function
$$\Psi_\CS^b=
H_\CS^b\cdot\prod_{B\in\CS} u_B^{t_B}=
H_\CS^b\cdot\prod_{B\in\CS} x_B^{R^\CS_B}$$
is a solution of the holonomy equation $\nabla\medspace\Psi_\CS^b=0$.
\item $\Psi_\CS^b$ is unipotent.
\end{enumerate}
\end{theorem}

\remark
Let $b'=\{x'_B\}_{B\in\CS}$ be another basis of $V^*$ adapted to
$\CS$. By lemma \ref{le:P_x}, $x'_B=x_B\cdot P_{B}$ where $P_B$ is a
polynomial in the variables $u_C$, $C\subsetneq B$, such that $P_B(0)
\neq 0$. Thus,
\begin{equation*}
\begin{split}
\Psi_\CS^{b'}
&=
H_\CS^{b'}\cdot\prod_{B\in\CS}(x'_B)^{R_B^\CS}\\
&=
H_\CS^{b'}\cdot\prod_{B\in\CS}\left(\frac{P_B}{P_B(0)}\right)^{R_B^\CS}
\cdot\prod_{B\in\CS}x_B^{R_B^\CS}
\cdot\prod_{B\in\CS}P_B(0)^{R_B^\CS}
\end{split}
\end{equation*}
Since the function $H_\CS^{b'}\cdot\prod_{B\in\CS}\left(\frac{P_B}{P_B(0)}
\right)^{R_B^\CS}$ is holomorphic in a neighborhood of $y_\CS$ and equal
to $1$ at $y_\CS$, it follows by the uniqueness statement of theorem \ref
{th:fundamental solutions} that 
\begin{equation}\label{eq:change of b}
\Psi_\CS^{b'}=\Psi_\CS^b\cdot\prod_{B\in\CS}P_B(0)^{R_B^\CS}
\end{equation}

\remark
The above solutions generalise those constructed by Drinfeld \cite{Dr3} for
the KZ equations, that is the holonomy equations for the Coxeter arrangement
of type $\sfA_n$, and those constructed by Cherednik for a general Coxeter
arrangement $\A_W$, see in particular \cite[\S 2]{Ch1} and \cite[\S 2]{Ch2}
where the blowup coordinates \eqref{eq:blowup} for $\A_W$ are introduced.

\subsection{Real arrangements and chambers}\label{ss:chambers}

\subsubsection{} 
\label{sss:chamber}

Assume henceforth that $\A$ is the complexification of a hyperplane arrangement
$\A_\IR$ in a real vector space $V_\IR\subset V$ and choose the linear forms $x
\in X$ so that they lie in $V_\IR^*\subset V^*$.  Fix a chamber $\C$ of $\A_\IR$,
that is a connected component of $V_\IR\setminus\A_\IR$. Up to replacing some
$x\in X$ by their opposites, we may assume that $x_{|_{\C}}>0$ for any $x\in X$.
Let $\Delta=\Delta(\C)\subset X$ be minimal for the property that any $x\in X$ is a
linear combination with non--negative coefficients of elements of $\Delta$. $\Delta$
is readily seen to consist of those $x\in X$ which are not linear combinations with
positive coefficients of two or more elements of $X$. It is therefore unique and
canonically determined by $\C$. An element $x$ of $V^*$ will be termed {\it real}
if $x\in V_\IR^*$ and {\it positive} if it is positive on $\C$. Note that $x\in V^*$ is
positive and real if, and only if it is a linear combination with non--negative
coefficients of the elements of $\Delta$. \\

We shall assume for simplicity that $\Delta$ is a basis of $V^*$, so that $\C$
is an open simplicial cone. The contents of \S \ref{ss:chambers}--\S \ref
{ss:adapted family} extends however, with suitable modifications, to the
case of a general polyhedral chamber.

\subsubsection{} 

Let $\I_\Delta\subset I$ be the set of irreducible elements $B$ which are
spanned by $\Delta_B=\Delta\cap B$.

\begin{definition}
A nested set $\CS$ is fundamental with respect to $\C$ if $\CS\subset\I_\Delta$.
\end{definition}

Let $\CS$ be a maximal nested set, $y_\CS\in Y_X$ the corresponding
point at infinity and identify the chamber $\C$ with its preimage in $Y_X$.
Our aim in this subsection is to prove the following

\begin{proposition}\label{pr:y_S in C}
$y_\CS$ lies in the closure of $\C$ in $Y_X$ if, and only if $\CS$ is
fundamental with respect to $\C$.
\end{proposition}

\subsubsection{} 

We shall need two preliminary results.

\begin{lemma}\label{le:y_S in C}\hfill
\begin{enumerate}
\item If $y_\CS$ lies in the closure of $\C$ then, for any positive, real basis
$b$ of $V^*$ adapted to $\CS$, the polynomials $P_x\in\IR[u_B]_{B\in\CS}$,
$x\in X$ defined by lemma \ref{le:P_x} satisfy $P_x(0)>0$.
\item Conversely, if $b$ is a positive, real basis of $V^*$ adapted to $\CS$
and $P_x(0)>0$ for any $x\in X$, then $y_\CS$ lies in the closure of $\C$.
\end{enumerate}
\end{lemma}
\proof (i) follows from the identity $x=x_{p_\CS(x)}\cdot P_x$. (ii) The sequence
of points $y_n\in\U_\CS^b$ with coordinates $u_B=1/n$, $B\in\CS$ converges
to $y_\CS$ and lies in $\C$ for $n$ large enough since for any $x\in X$,
$$x=x_{p_\CS(x)}\cdot P_x=\prod_{B\supseteq p_\CS(x)}u_B\cdot P_x$$
\halmos

\begin{lemma}\label{le:Px(0)}
If $b$ is a basis of $V^*$ adapted to $\CS$, then $P_x(0)>0$ for any $x\in
X$ if, and only if, $P_x(0)>0$ for any $x\in\Delta$.
\end{lemma}
\proof
Let $x\in X$ and write $x=\sum_{x'\in\Delta}k_{x'}x'$, where $k_{x'}\geq 0$. Let
$B_1,\ldots,B_m\in\CS$ be the maximal elements among the $p_\CS(x')$ with
$k_{x'}>0$ and rewrite the above as
$$x=\sum_{i=1}^m \sum_{x'\in\Delta_{B_i}}k_{x'}x'\in B_1\odots{\oplus}B_m$$
By nestedness of $\CS$, $x$ is contained in one of the $B_i$ whence $m=1$.
Thus, modulo $i_\CS(B_1)$,
$$x=
\sum_{\substack{x'\in\Delta,\\P_\CS(x')=B_1}}k_x' x'=
\sum_{\substack{x'\in\Delta,\\P_\CS(x')=B_1}}k_x' P_{x'}(0)x_{B_1}$$
whence $P_x(0)=\sum_{x'}k_x' P_{x'}(0)$ \halmos

\subsubsection{} 

{\sc Proof of proposition \ref{pr:y_S in C}.}
If $\CS$ is fundamental, then $\Delta$ is a positive, real basis of $V^*$ adapted
to $\CS$ and such that $P_x\equiv 1$ for any $x\in\Delta$. By lemmas \ref{le:Px(0)}
and \ref{le:y_S in C}, $y_\CS$ therefore lies in the closure of $\C$. Assume now
that $\CS$ is not fundamental, let $B\in\CS$ be minimal for the property that $B$
is not spanned by $\Delta_B$ and set $C=i_\CS(B)$. Let $x\in X\cap(B\setminus C)$
and write $x=\sum_{x'\in\Delta^x}k_{x'}x'$ where $\Delta^x\subseteq\Delta$ and $k_{x'}>0$ for any $x'\in\Delta^x$. Decompose
the sum as
\begin{equation}\label{eq:x1 x2}
x=
\sum_{x'\in\Delta^x\cap B}k_{x'}x'+
\sum_{x'\in\Delta^x\setminus B}k_{x'}x'=
x_1+x_2
\end{equation}
We claim that $x_2\neq 0$. Indeed, by minimality of $B$, $C$ is
spanned by $\Delta_C$ so that if $x_2=0$, then $B=\IC x\oplus C$ would
be spanned by $\Delta_B$, a contradiction. We claim next that there exists
an $x'\in\Delta^x$ such that $p_\CS(x')\supsetneq B$. Indeed, if $x'\in\Delta^x
\setminus B$, the nestedness of $\CS$ implies that either $B\subsetneq p_\CS
(x')$ or that the two are in direct sum. If the latter were the case for any $x'\in
\Delta^x\setminus B$, \eqref{eq:x1 x2} would yield $x_2=0$, a contradiction.
Let therefore $x'\in\Delta^x \setminus B$ be such that $A=p_\CS(x')\supsetneq
B$, and set $x_A=x'$, $x_B=x$ and complete this into a positive, real basis of
$V^*$ adapted to $\CS$.  Then, on $\U_\CS^b\cap\C$,
$$k_{x'}<\frac{x_B}{x_A}=\prod_{B\subseteq D\subsetneq A}u_D$$
Thus, the product of the coordinates $u_D$, $B\subseteq D\subsetneq A$ is
real and bounded below on $\C\cap\U_\CS^b$ so that $y_\CS$ does not lie
in the closure of $\C$ by lemma \ref{le:y_S in C} \halmos

\subsection{\DCP associators}\label{ss:DCP associators}

Let $\F$ be a \fmns and $b=\{x_B\}_{B\in\F}$ a positive, real basis of $V^*$
adapted to $\F$. Let $\V_\F^b\subset\U_\F^b$ be the complement of the
real codimension one semialgebraic subvarieties $\{x_B\leq 0\}$, $B\in\F$.
Note that the chamber $\C$ lies in $\V_\F^b$ since $x_B>0$ on $\C$ for any
$B\in\F$. We shall henceforth only consider the standard determination of the
function $\log(z)$ obtained by performing a cut along the negative real axis, so
that the functions $\log(x_B)$, $B\in\F$ are well--defined and single--valued
on $\V_\F^b$. Theorem \ref{th:fundamental solutions} then yields a
single--valued fundamental solution $\Psi_\F^b$ of the holonomy equations
\eqref{eq:holonomy} on the intersection of a neighborhood of $y_\F$ in $\U
_\F^b$ with $\V_\F^b$. Since $y_\F$ lies in the closure of $\C$ by proposition
\ref{pr:y_S in C}, $\Psi_\F^b$ may be continued to a single--valued solution
on $\C$ which we shall denote by the same symbol.\\

Let now $\F,\G$ be two \fmnss and $b,b'$ two positive, real bases of $V^*$
adapted to $\F$ and $\G$ respectively.

\begin{definition}
The \DCP associator $\Phi_{\G\F}^{b'b}$ is the element of $A$ defined by
$$\Phi_{\G\F}^{b'b}=
(\Psi_{\G}^{b'}(y))^{-1}\cdot\Psi_\F^b(y)$$
for any $y\in\C$.
\end{definition}

Note that $\Phi_{\G\F}^{b'b}$ is well--defined since the \rhs of the above
expression is a locally constant function on $\C$. The following properties
are immediate
\begin{itemize}
\item {\bf orientation} : for any pair $\F,\G$ and adapted bases $b,b'$
\begin{equation}\label{eq:DCP orientation}
\Phi_{\F\G}^{bb'}=(\Phi_{\G\F}^{b'b})^{-1}
\end{equation}
\item {\bf transitivity} : for any triple $\F,\G,\H$ and adapted
bases $b,b',b\dprime$
\begin{equation}\label{eq:DCP transitivity}
\Phi_{\H\F}^{b\dprime b}=
\Phi_{\H\G}^{b\dprime b'}\cdot\Phi_{\G\F}^{b'b}
\end{equation}
\end{itemize}


\subsection{Elementary associators}\label{ss:elementary associators}

The \DCP associators possess a number of other important properties
which will be given in \S \ref{ss:support}--\ref{ss:forgetful} and are easier
to formulate in terms of elementary associators.
\begin{definition}
A pair $(\G,\F)$ of \mnss (resp. an associator $\Phi_{\G\F}^{b'b}$)
is called elementary if $\G$ and $\F$ differ by one element.
\end{definition}

Elementary associators are sufficient to reconstruct general associators
in view of the transitivity relation \eqref{eq:DCP transitivity} and the
following result

\begin{proposition}\label{pr:connectivity}
For any pair $\F,\G$ of (fundamental) maximal nested sets, there exists
a sequence
$$\F=\H_1,\ldots,\H_m=\G$$
of (fundamental) \mnss such that $\H_i$ and $\H_{i+1}$ differ by an
element.
\end{proposition}

The transitivity of associators implies that the elementary ones satisfy the
following property :

\begin{itemize}
\item {\bf coherence (I)} : if $\H_1,\ldots,\H_\ell$ and $\K_1,\ldots,\K_m$
are two sequences of \fmnss such that $|\H_{i+1}\setminus\H_i|=1$ for any
$1\leq i\leq\ell-1$, $|\K_{j+1}\setminus\K_j|=1$ for any $1\leq j\leq m-1$,
$$\H_1=\K_1\qquad\text{and}\qquad\H_\ell=\K_m$$
then, 
$$
\Phi^{b_{\ell}b_{\ell-1}}_{\H_\ell\H_{\ell-1}}\cdots\Phi^{b_2b_1}_{\H_2\H_1}=
\Phi^{b'_mb'_{m-1}}_{\K_m\K_{m-1}}\cdots\Phi^{b'_2b'_1}_{\K_2\K_1}
$$
provided $b_1=b'_1$ and $b_\ell=b'_m$.
\end{itemize}

\subsection{2--dimensional reduction}\label{ss:rank 2 reduction}

Let $(\G,\F)$ be an elementary pair of fundamental maximal nested sets
and $b'=\{x'_B\}$, $b=\{x_B\}$ two positive, real bases of $V^*$ adapted
to $\G$ and $\F$ respectively. Assume that
\begin{equation}\label{eq:equal b}
x_B=x'_B\qquad\text{for any}\quad B\in\F\cap\G
\end{equation}
Then, as explained below, the associator $\Phi_{\G\F}^{b'b}$ coincides
with one obtained from a line arrangement in a two--dimensional vector
space determined by $\F$ and $\G$. This inductive structure of associators
is a consequence of the factorisation property of the solutions $\Psi_\F^b$
which was pointed out for the KZ equations by Tsuchiya--Kanie \cite{TK}
and extended to the holonomy equations for Coxeter arrangements by
Cherednik (Thm. 1 in \cite{Ch1,Ch2}).

\subsubsection{} 

Note first that, by propositions \ref{pr:n(F)} and \ref{pr:max nested}, there
is a unique $B\in\F\cap\G$ such that
$$n(B;\F\cap\G)=2\qquad\text{while}\qquad n(B';\F\cap\G)=1$$
for any other $B'\in\F\cap\G$. Set $C=i_{\F\cap\G}(B)$, so that $\dim(B/C)=2$
and let $X\BC$ be the quotient of $X\cap(B\setminus C)$ by the equivalence
relation $x\sim x'$ if $x$ is proportional to $x'$ mod $C$. Then, $X\BC$
defines an essential line arrangement $\A_{B,C}$ in the 2--dimensional
vector space $C^\perp/B^\perp$.

\subsubsection{} 

Let $C=C_1\odots{\oplus}C_m$ be the irreducible decomposition of $C$.

\begin{lemma}\label{le:C_x}
Let $x\in X\cap(B\setminus C)$ and consider the irreducible decomposition
$$\IC x\oplus C=C_x\oplus C_{x,1}\odots{\oplus}C_{x,p}$$
where $C_x$ is the summand containing $x$. Then,
\begin{enumerate}
\item The $C_{x,j}$ are exactly the irreducible summands $C_i$ of $C$ such
that $\IC x\oplus C_i$ is a decomposition in $\Lpos$.
\item $C_x=\IC x\bigoplus_i C_i$ where $C_i$ ranges over the irreducible
summands of $C$ such that $\IC x\oplus C_i$ is not a decomposition in $\Lpos$.
\item $C_x$ only depends upon the equivalence class of $x$ in $X\BC$.
\end{enumerate}
\end{lemma}
\proof (i) and (ii) follow easily by comparing the decompositions of $C$ and
$C\oplus\IC_x$. (iii) If $x'\sim x$, then
$$x'\in\IC x'\oplus C=\IC x\oplus C=C_x\oplus C_{x,1}\odots{\oplus}C_{x,p}$$
so that $x'\in C_x$ or, for some $j$, $x'\in C_{x,j}\subset C$ by (i). The latter
however is ruled out by the fact that $x'\notin C$ \halmos

\subsubsection{} 

For any equivalence class $\ol{x}\in X\BC$, set
$$t_{\ol{x}}=\sum_{x\in\ol{x}}t_x$$
and consider the connection on $C^\perp/B^\perp$ with logarithmic singularities
on $\A\BC$ defined by
\begin{equation}\label{eq:nablaBC}
\nabla\BC=d-\sum_{\ol{x}\in X\BC}\frac{d\ol{x}}{\ol{x}}\cdot t_{\ol{x}}
\end{equation}

For any $\ol{x}\in X\BC$, set $C_{\ol{x}}=C_x$ where $x\in\ol{x}$.

\begin{proposition}\label{pr:reduced connection}\hfill
\begin{enumerate}
\item For any $\ol{x}\in X\BC$ one has
$$t_{\ol{x}}=t_{C_{\ol{x}}}-\sum_{i:C_i\subset C_{\ol{x}}}t_{C_i}$$
\item $t_{\ol{x}}$ commutes with any $t_{x'}$, $x'\in X\cap C$.
\item The connection $\nabla\BC$ is flat.
\end{enumerate}
\end{proposition}
\proof
(i) We have
$$t_{C_{\ol{x}}}=
\sum_{y\in X\cap C_{\ol{x}}}t_y=
\sum_{y\in X\cap C_{\ol{x}}\setminus C}t_y+\sum_{y\in X\cap C_{\ol{x}}\cap C}t_y$$
The first summand is equal to $t_{\ol{x}}$. By lemma \ref{le:C_x}, $y\in X\cap C$
lies in $C_{\ol{x}}$ if, and only if the irreducible component $C_i$ of $C$ containing
$y$ is contained in $C_{\ol{x}}$. The second sum is therefore equal to $\sum_
{i:C\subset C_{\ol{x}}}t_{C_i}$ as claimed.
(ii) follows from (i) and lemma \ref{le:[t_B,t_C]}.
(iii) By \cite{Ko2}, the flatness of $\nabla\BC$ equivalent to the fact that each $t_{\ol
{x}}$ commutes with $$\sum_{\ol{x}\in X\BC}t_{\ol{x}}=t_B-t_C$$ which follows from
(ii) and lemma \ref{le:[t_B,t_C]} \halmos

\subsubsection{} 

We claim next that the lattice $\L\BC$ of subspaces of $B/C$ spanned by elements
of $X\BC$ contains two distinguished one--dimensional elements determined by $\F$
and $\G$ respectively. Indeed, let $B_1,B_2$ be the unique elements in $\F\setminus
\G$ and $\G \setminus\F$. Note that $B_1$ (resp. $B_2$) is one of the maximal elements of $\F$ (resp. $\G$) properly contained in $B$ since $n(B;\F)=1$ (resp.
$n(B;\G)=1$) while $n(B;\F\cap\G)=2$. Set, for $i=1,2$
$$\ol{B}_i=B_i/B_i\cap C\subseteq B/C$$
By nestedness of $\F$ (resp. $\G$), $B_1$ (resp. $B_2$) either contains, or is in
direct sum with, each $C_j$. In particular, the maximal elements of $\F$ (resp. $\G$)
properly contained in $B_1$ (resp. $B_2$) are exactly the $C_j$ contained in $B_1$
(resp. $B_2$) so that, by proposition \ref{pr:max nested},
$$\dim(\ol{B}_1)=1=\dim(\ol{B}_2)$$
Note that if $\A\BC$ contains at least three lines, so that $B/C=(C^\perp/B^\perp)^*$
is irreducible, then
$$\ol{\F}=\{\ol{B}_1,B/C\}\qquad\text{and}\qquad\ol{\G}=\{\ol{B}_2,B/C\}$$
are \mnss of irreducible elements of $\L\BC^*$.

\subsubsection{} 

Note next that $\A\BC$ is the complexification of a line arrangement in
$(C^\perp/B^\perp)_\IR$ endowed with a distinguished chamber $\C_
{B,C}$, namely the interior of the image of $\ol{\C}\cap C^\perp$ in $C
^\perp/B^\perp$. Moreover, if $x_i$ is the unique element in $\Delta_{B_i}
\setminus\Delta_C$, $i=1,2$, then any element of $X\BC$ is a linear
combination with non--negative coefficients of $\ol{x}_1,\ol{x}_2$ so
that the real lines $\ol{x}_1^\perp$ and $\ol{x}_2^\perp$ are two
contiguous walls of $\C\BC$. Since $\ol{B}_i=\IC\ol{x}_i$, $\ol{\F}$
and $\ol{\G}$ are \fmnss with respect to the chamber $\C\BC$ whenever
$|X\BC|\geq 3$.

\subsubsection{} 

Let now $b=\{x_A\}$ and $b'=\{x'_A\}$ be positive, real bases of $V^*$
adapted to $\F$ and $\G$ respectively and such that \eqref{eq:equal b}
holds. $b$ and $b'$ induce bases $\ol{b}=\{\ol{x}_B,\ol{x}_{B_1}\}$ and $\ol{b}'=\{\ol{x}_B,\ol{x'}_{B_2}\}$ of $C^\perp/B^\perp$ adapted to $\ol
{\F}$ and $\ol{\G}$ respectively which are positive and real with respect
to the chamber $\C\BC$.

\begin{theorem}\label{th:rank 2}\hfill
\begin{enumerate}
\item If the arrangement $\A\BC$ contains at least three lines, then
$$\Phi_{\G\F}^{b'b}=\Phi_{\ol{\G}\ol{\F}}^{\ol{b}'\ol{b}}$$
where the \rhs is the \DCP associator for the connection $\nabla\BC$
relative to the solutions $\Psi_{\ol{\G}}^{\ol{b}'}$ and
$\Psi_{\ol{\F}}^{\ol{b}}$.
\item Otherwise,
$$\Phi_{\G\F}^{b'b}=a_2^{t_{\ol{x}_{B_2}}}\cdot(a_1^{t_{\ol{x}_{B_1}}})^{-1}$$
where $a_1,a_2\in\IR_+^*$ are such that $x_B=a_1x_{B_1}+a_2x'_{B_2}$
mod $C$.
\end{enumerate}
\end{theorem}

\subsection{Adapted families}\label{ss:adapted family}

Recall that $\I_\Delta$ is the set of irreducible elements $B\in\L^*$ which are
spanned by $\Delta_B=\Delta\cap B$. If $\F\subset\I_\Delta$ is a \fns and $B
\in\F$, we set
$$\aalpha^B_\F=\Delta_B\setminus\Delta_{i_\F(B)}$$
so that $|\aalpha^B_\F|=n(B;\F)$. If $\F$ is maximal, we denote by $\alpha^B
_\F$ the unique element in $\aalpha^B_\F$. The following notion is useful to
obtain adapted bases which satisfy \eqref{eq:equal b}.

\begin{definition}
An adapted family is a collection $\beta=\{x_B\}_{B\in\I_\Delta}$ such that,
for any $B\in\I_\Delta$,
$$x_B\in B\setminus\bigcup_{C\in\I_\Delta,C\subsetneq B}C$$
\end{definition}

Clearly, if $\beta=\{x_B\}_{B\in\I_\Delta}$ is an adapted family then, for any
\fmns $\F$, $\beta_\F=\{x_B\}_{B\in\F}$ is a basis of $V^*$ adapted to $\F$
and all collections of adapted bases $b_\F=\{x_A^\F\}_{A\in\F}$ labelled by
\fmnss and satisfying \eqref{eq:equal b} are obtained in this way. If $\beta$
is a positive, real adapted family and $\F,\G$ are fundamental maximal
nested sets, set
$$\Phi_{\G\F}^\beta=\Phi_{\G\F}^{\beta_\G\beta_\F}$$
We shall henceforth only use the associators $\Phi_{\G\F}^\beta$ corresponding
to a fixed such $\beta$. These clearly satisfy the orientation property $\Phi
_{\F\G}^\beta=(\Phi_{\G\F}^\beta)^{-1}$ of \S \ref{ss:elementary associators}
and the following simpler version of the coherence property
\begin{itemize}
\item {\bf coherence (II)} : if $\H_1,\ldots,\H_\ell$ and $\K_1,\ldots,\K_m$
are two sequences of \fmnss such that $|\H_{i+1}\setminus\H_i|=1$ for any
$1\leq i\leq\ell-1$, $|\K_{j+1}\setminus\K_j|=1$ for any $1\leq j\leq m-1$,
$$\H_1=\K_1\qquad\text{and}\qquad\H_\ell=\K_m$$
then, 
\begin{equation}\label{eq:DCP elem transitivity 3}
\Phi^{\beta}_{\H_\ell\H_{\ell-1}}\cdots\Phi^{\beta}_{\H_2\H_1}=
\Phi^{\beta}_{\K_m\K_{m-1}}\cdots\Phi^{\beta}_{\K_2\K_1}
\end{equation}
\end{itemize}

\remark\rklabel{rk:gauge}
Let $\beta'=\{x'_B\}$ be another positive, real adapted family. For any $B\in\I_
\Delta$ and $\alpha\in\Delta_B$, let $c_{(B;\alpha)}\in\IR_+^*$ be such that
$x'_B=c_{(B;\alpha)}x_B$ modulo the span of $\Delta_B\setminus\alpha$.
Then, by \eqref{eq:change of b}, the associators $\Phi_{\G\F}^\beta$ and $\Phi_{\G\F}^{\beta'}$ are related by
\begin{equation}\label{eq:gauge}
\Phi_{\G\F}^{\beta'}=a_\G\cdot\Phi_{\G\F}^\beta\cdot a_\F^{-1}
\end{equation}
where, for any \fmns $\F$,
$$a_\F=\prod_{B\in\F}{c_{(B;\alpha^B_\F)}}^{-R_B^\F}$$

\subsection{Support properties of elementary associators}
\label{ss:support}

Let $(\F,\G)$ be an elementary pair of fundamental maximal nested sets.

\begin{definition}
The support of $(\F,\G)$ is the unique element $B=\supp(\F,\G)$ of $\F\cap
\G$ such that $n(B;\F\cap\G)=2$. The central support of $(\F,\G)$ is the
subspace $\zsupp(\F,\G)=i_{\F\cap\G}(\supp(\F,\G))$.
\end{definition}

For any $D\in\Lpos$, let
$$A_D\subseteq A$$
be the subalgebra topologically generated by the elements $t_x$, $x\in D\cap X$.
Let $\beta$ be a positive, real adapted family. It follows from theorem \ref{th:rank 2}
and (ii) of proposition \ref{pr:reduced connection} that the associator $\Phi
_{\G\F}^\beta$ satisfies the following property
\begin{itemize}
\item {\bf support} : $\Phi_{\G\F}^\beta$ lies in $A_{\supp(\F,\G)}$ and commutes
with $A_{\zsupp(\F,\G)}$.
\end{itemize}

\subsection{Forgetfulness properties of elementary associators}
\label{ss:forgetful}

\begin{definition}
Two elementary pairs $(\F,\G)$, $(\wt{\F},\wt{\G})$ of \fmnss are
equivalent if
\begin{gather*}
\supp(\F,\G)=\supp(\wt{\F},\wt{\G}),\\
\alpha^{\supp(\F,\G)}_\F=\alpha^{\supp(\wt{\F},\wt{\G})}_{\wt{\F}}
\qquad\text{and}\qquad
\alpha^{\supp(\F,\G)}_\G=\alpha^{\supp(\wt{\F},\wt{\G})}_{\wt{\G}}
\end{gather*}
\end{definition}

The following result guarantees that the equivalence of $(\F,\G)$
and $(\F',\G')$ implies the equality of the reduction data used in
\S \ref{ss:rank 2 reduction}.

\begin{proposition}\label{pr:combi=>practical}
Let $(\F,\G)$ be an elementary pair of fundamental, maximal nested
sets. Set $B=\supp(\F,\G)$ and let $B_1,B_2$ the unique elements in
$\F\setminus\G$ and $\G\setminus\F$ respectively. Then,
\begin{enumerate}
\item $\alpha^B_\F$ and $\alpha^B_\G$ are distinct and $\aalpha^B_{\F
\cap\G}=\{\alpha^B_\F,\alpha^B_\G\}$.
\item $\zsupp(\F,\G)$ is the span of $\Delta_{\supp(\F,\G)}\setminus\{\alpha^B_\F,\alpha^B_\G\}$.
\item $B_1$ is the irreducible component of $\<\Delta_B\setminus\alpha
^B_\F\>$ containing $\alpha^B_\G$. Moreover, $\alpha^{B_1}_\F=\alpha
^B_\G$.
\item $B_2$ is the irreducible component of $\<\Delta_B\setminus\alpha
^B_\G\>$ containing $\alpha^B_\F$. Moreover, $\alpha^{B_2}_\G=\alpha
^B_\F$.
\end{enumerate}
\end{proposition}
\proof
(i) Clearly, $\alpha^B_\F,\alpha^B_\G\in\aalpha^B_{\F\cap\G}$. Since the
maximal elements of $\F$ (resp. $\G$) properly contained in $B$ are the
irreducible components of the span of $\Delta_B\setminus\alpha^B_\F$
(resp. $\Delta_B\setminus\alpha^B_\G$), the equality $\alpha^B_\F=\alpha
^B_\G$ would imply that these components also lie in $\G$ (resp. $\F$)
and therefore that $B$ is saturated as an element of $\F\cap\G$, a contradiction. 
Since $n(B;\F\cap\G)=2$, this implies that $\aalpha^B_{\F\cap\G}=\{\alpha
^B_\F,\alpha^B_\G\}$.\\
(ii) follows from (i) and the fact that $\zsupp(\F,\G)$ is spanned by $\Delta_B\setminus\aalpha^B_{\F\cap\G}$.\\
(iii)--(iv) Since $\alpha^B_\G\neq\alpha^B_\F$, there exists an irreducible
component $B_1'\in\F$ of $\<\Delta_B\setminus\alpha^B_\F\>$ such that
$\alpha^B_\G\in B_1'$. However since $\alpha^B_\G\notin C$ for any $C
\in\G$, $C\subsetneq B$, $B_1'\in\F\setminus\G=\{B_1\}$. Since $B_1$ is
one of the proper maximal elements of $\F$ contained in $B$, we have
$\alpha^{B_1}_\F\in\aalpha^B_{\F\cap\G}=\{\alpha^B_\F,\alpha^B_\G\}$
whence $\alpha^{B_1}_\F=\alpha^B_\G$ since $\alpha^B_\F\notin B_1$.
Similarly, $B_2$ is the unique irreducible component of $\<\Delta_B\setminus
\alpha^B_\G\>$ containing $\alpha^B_\F$ and $\alpha^{B_2}_\G=\alpha
^B_\F$ \halmos\\

Let $\beta$ be a positive, real adapted family. Proposition \ref
{pr:combi=>practical} and theorem \ref{th:rank 2} imply that
the associators $\Phi_{\G\F}^\beta$ satisfy the following additional
property
\begin{itemize}
\item {\bf forgetfulness} : if $(\G,\F)$ and $(\wt{\G},\wt{\F})$
are equivalent,
$\Phi_{\G\F}^\beta=\Phi_{\wt{\G}\wt{\F}}^\beta$.
\end{itemize}

\subsection{Coxeter arrangements}\label{ss:coxeter}

\subsubsection{} 
\label{sss:root system}

Assume now that $\A_\IR$ is the arrangement of reflecting hyperplanes of a
finite (real) reflection group $W\subset GL(V_\IR)$. The set $X$ of defining
equations of $\A$ may then be chosen so that $\Phi=X\sqcup(-X)$ is invariant
under $W$. Thus, $\Phi$ is a (reduced) root system with respect to any $W
$--invariant Euclidean inner product $(\cdot,\cdot)$ on $V_\IR^*$, that is a
finite collection of non--zero vectors in $V_\IR^*$ satisfying, for any $\alpha
\in\Phi$
\begin{itemize}
\item[(R1)] $\Phi\cap\IR\alpha=\{\pm\alpha\}$
\item[(R2)] $s_\alpha\Phi=\Phi$
\end{itemize}
where $s_\alpha\in W^t$ is the orthogonal reflection determined by $\alpha$
and $W\cong W^t\subset O(V_\IR^*)$ is the group contragredient to $W$
\footnote{we follow here the terminology of \cite{Hu}. Thus, $\Phi$ need not
be crystallographic, \ie such that $2(\alpha,\beta)/(\alpha,\alpha)\in\IZ$ for
any $\alpha,\beta\in\Phi$.}. In accordance with \eqref{eq:essential} we assume
in addition that $\Phi$ spans $V_\IR^*$ so that no $v\in V\setminus\{0\}$
is fixed by $W$.

\subsubsection{} 
\label{sss:root terminology}

We shall need some mostly standard terminology. A root subsystem of $\Phi$
is a subset $\ol{\Phi}\subseteq\Phi$ satisfying (R1)--(R2) above and such that
the intersection of its linear span $\<\ol{\Phi}\>\subset V_\IR^*$ with $\Phi$ is
$\ol{\Phi}$\footnote{Thus, we do not regard the long roots in the root system
of type $\sfG_2$ as a root subsystem.}. A root subsystem $\oPhi$ is
reducible if it possesses a non--trivial
partition $\oPhi=\oPhi_1\sqcup\oPhi_2$ into mutually orthogonal subsets,
which are then necessarily root subsystems of $\oPhi$, and irreducible otherwise.
Two root subsystems $\oPhi_1,\oPhi_2\subseteq\Phi$ are said to be {\it completely orthogonal} if no element $\alpha\in\Phi$ is of the form 
$\alpha=a_1\alpha_1+a_2\alpha_2$ with $\alpha_i\in\oPhi_i$ and $a_i\in
\IR^*$\footnote{When $\Phi$ is a crystallographic root system and $\oPhi_i
=\{\pm\beta_i\}$, this notion is more stringent than the strong orthogonality
of $\beta_1$ and $\beta_2$ \ie the requirement that $\beta_1\pm\beta_2
\notin\Phi$. For example, if $\Phi$ is the root system of type $\sfG_2$ and
$\alpha_1$, $\alpha_2$ are the short and long simple roots respectively,
then $\beta_1=2\alpha_1+\alpha_2$ and $\beta_2=\alpha_2$ are strongly
orthogonal but $\{\pm\beta_1\}$ and $\{\pm\beta_2\}$ are not completely
orthogonal.}. This implies in particular that $\oPhi_1\perp\oPhi_2$.\\

Let now $\C$ be a chamber of $\A_\IR$, $\Phi=\Phi_+\sqcup\Phi_-$ the
corresponding partition into positive and negative roots and $\Delta=\{
\alpha_i\}_{i\in\bfI}\subset\Phi_+$ the basis of $V^*$ consisting of
indecomposable elements of $\Phi_+$. We shall say that a root
subsystem $\oPhi\subseteq\Phi$ is {\it fundamental} if $\<\oPhi\>$
is spanned by $\Delta\cap\oPhi$. It is easy to see that two fundamental
root subsystems $\oPhi_1,\oPhi_2\subseteq\Phi$ are completely
orthogonal if, and only if they are orthogonal.

\subsubsection{} 
\label{sss:subspace dictionary}

In accordance with \S \ref{sss:chamber}, we assume that $X=\Phi_+$. Let
$\L$ be the lattice of subspaces of $V^*$ spanned by the elements of $X$,
as in \S  \ref{ss:wonderful model}, and $\R$ the lattice of root subsystems
of $\Phi$. The following result provides a dictionary between the terminology
of \S \ref{ss:wonderful model}--\S \ref{ss:divisor} and that of \S \ref
{sss:root terminology}.

\begin{proposition}\label{pr:dictionary}\hfill
\begin{enumerate}
\item The map $\ol{\Phi}\rightarrow\<\ol{\Phi}\>$ is a bijection between
$\R$ and $\L$, with inverse given by $B\rightarrow \Phi_B=B\cap\Phi$.
\item A root subsystem $\ol{\Phi}\subset\Phi$ admits a partition $\ol
{\Phi}=\ol{\Phi}_1\sqcup\ol{\Phi}_2$ into two orthogonal subsets if,
and only if, $\<\ol{\Phi}\>=\<\ol{\Phi}_1\>\oplus\<\ol{\Phi}_2\>$ is a
decomposition in $\L$. In particular, $\ol{\Phi}$ is irreducible if, and
only if $\<\ol{\Phi}\>$ is an irreducible element of $\Lpos$.
\item A collection $\CS$ of irreducible elements of $\Lpos$ is nested
if, and only if the root subsystems $\{\Phi_B\}_{B\in\CS}$ are pairwise
completely orthogonal when non--comparable. In particular, if each
$\Phi_B$ is fundamental, then $\CS$ is nested if, and only if, the $\Phi
_B$ are pairwise orthogonal when non--comparable.
\end{enumerate}
\end{proposition}

\subsubsection{} 
\label{sss:graph dictionary}

Recall that the Coxeter graph $D$ of $\Phi$ is the graph with vertex
set $\Delta$ and an edge between $\alpha_i$ and $\alpha_j$ if, and
only if $\alpha_i\nperp\alpha_j$. It follows from proposition \ref{pr:dictionary}
that the map $B\rightarrow B\cap\Delta$ induces a bijection between
\fnss of irreducible elements of $\Lpos$ and sets of connected subgraphs
of $D$ which are pairwise {\it compatible} that is such that $D'\subseteq
D\dprime$, $D\dprime\subseteq D'$ or $D'\perp D\dprime$, where the
last statement means that no vertex of $D'$ is connected to a vertex
of $D\dprime$ by an edge of $D$.\\

\remark\rklabel{rk:Kn}
Such collections admit the following well--known alternative description
when $D$ is the Dynkin diagram of type $\sfA_{n-1}$. Identify for this
purpose $D$ with the interval $[1,n-1]$ and its connected subdiagrams
with the subintervals $[i,j]$, with $1\leq i\leq j\leq n-1$. This induces a
bijection between the sets of pairwise compatible connected subdiagrams
of $D$ and consistent bracketings on the non--associative monomial
$x_1\cdots x_n$ by attaching to $B=[i,j]$ the bracket
$x_1\cdots x_{i-1}(x_i\cdots x_{j+1})x_{j+2}\cdots x_n$.

\subsubsection{} 

Fix a basepoint $v_0\in\VA$, let $[v_0]$ be its image in $\VA/W$ and let
$$P_W=\pi_1(\VA,v_0)
\qquad\text{and}\qquad
B_W=\pi_1(\VA/W,[v_0])$$
be the generalised (topological) pure and full braid groups of type $W$
respectively\footnote{to distinguish $\pi_1(\VA/W,[v_0])$ from the isomorphic
abstract group introduced in \S \ref{sss:monodromy for Brieskorn}, we denote
them by $B_W$ and $B_D$ and refer to them as the topological and algebraic
braid groups of $W$ respectively.}. Since $W$ acts freely on $\VA$ \cite{Ste},
the quotient map $\VA\xrightarrow{\pi}\VA/W$ is a covering and gives rise to
an exact sequence
\begin{equation}\label{eq:fibration ses}
1\longrightarrow P_W \longrightarrow B_W
\longrightarrow W\longrightarrow 1
\end{equation}
where the rightmost arrow is obtained by associating to $\gamma\in B_W$ the
unique $w\in W$ such that $w^{-1}v_0=\wt{\gamma}(1)$, with $\wt{\gamma}$
the unique lift of $\gamma$ to a path in $\VA$ such that $\wt{\gamma}(0)=v_0$.

\subsubsection{} 

Let $(\wt{\VA},\wt{v_0})\stackrel{p}{\longrightarrow}(\VA,v_0)$ be the universal
covering space of $\VA$. Then, $(\wt{\VA},\wt{v_0})$ is also the universal covering
space of $(\VA/W,[v_0])$ via $\pi\circ p$ and we get a canonical right action of
$B_W$ on $\wt{\VA}$ by deck transformations extending that of $P_W$. The
group $W$, and therefore $B_W$, act on $A$ by
\begin{equation}\label{eq:W on A}
w(t_\alpha)=t_{|w\alpha|}
\end{equation}
where $|w\alpha|=\pm w\alpha$ depending on whether $w\alpha\in\Phi_{\pm}$.
If $b\in B_W$ and $\Psi:\wt{\VA}\rightarrow A$ is a solution of the holonomy
equations $p^*\nabla\Psi=0$, one readily checks that
$b\bullet\Psi(\wt{v})=b(\Psi(\wt{v}b))$
is another solution. Thus, if $\Psi$ is invertible, then $\mu_\Psi=\Psi^{-1}\cdot b\bullet\Psi$ is a constant element of $A$. When $b\in P_W$, $\mu_\Psi(b)$
coincides with the element defined by \eqref{eq:mu_Psi} since $P_W$ acts
trivially on $A$. Set
$$\nu_\Psi(b)=\mu_\Psi(b)\cdot b\in A\rtimes B_W$$

\begin{proposition}\label{pr:nu monodromy}\hfill
\begin{enumerate}
\item The map $b\rightarrow \nu_\Psi(b)$ is a homomorphism $B_W\rightarrow
A\rtimes B_W$.
\item If $\Psi$ is unipotent, $\mu_\Psi$ takes values in $N\rtimes B_W$.
\item If $\Psi'$ is another invertible solution, then $\nu_{\Psi'}=\Ad(K^{-1})\circ
\nu_\Psi$, where $K=\Psi^{-1}\cdot\Psi'\in A^\times$. 
\end{enumerate}
\end{proposition}

\subsubsection{} 
\label{sss:relevance}

The relevance of the map $\nu_\Psi$ is the following. To any \fd representation
$\wt{\rho}:F/I\rtimes B_W\rightarrow\End(U)$, we may associate a flat holomorphic
vector bundle $(\U_{\wt{\rho}},\nabla_{\wt{\rho}})$ with fibre $U$ over $\VA/W$
as follows. Let $\rho$ be the restriction of $\wt\rho$ to $F/I$,
$$\nabla_\rho=
d-\sum_{\alpha\in\Phi_{+}}\frac{d\alpha}{\alpha}\cdot\rho(t_\alpha)$$
the corresponding flat connection on $\VA\times U$ defined by \eqref
{eq:nabla rho} and $p^*\nabla_\rho$ its pull--back to $\wt{\VA}\times U$.
Then, $(\U_{\wt{\rho}},\nabla_{\wt{\rho}})$ is the quotient
$$(U_{\wt{\rho}},\nabla_{\wt{\rho}})=(\wt{\VA}\times U,p^*\nabla_\rho)/B_W$$
where $B_W$ acts on $U$ via $\wt{\rho}$.
Note that if $\wt{\rho}$ factors through $F/I\rtimes W$, then $(\U_{\wt{\rho}},
\nabla_{\wt{\rho}})$ is simply the quotient of $(\VA\times U,\nabla_\rho)$
by $W$. As in \S \ref{ss:holonomy algebra}, we associate to $\wt{\rho}$
a representation $\wt{\rho}_h:A\rtimes B_W
\rightarrow\End(U\fmll)$ by setting, for $b\in B_W$ and $\alpha\in\Phi_+$,
\begin{equation}\label{eq:extended rho_h}
\wt\rho_h(b)=\wt\rho(b)
\qquad\text{and}\qquad
\wt\rho_h(t_\alpha)=h\wt\rho(t_\alpha)
\end{equation}
Then, for any invertible solution $\Psi:\wt{\VA}\rightarrow A$ of $p^*\nabla
\Psi=0$, $\Psi_{\wt{\rho}_h}=\wt\rho_h(\Psi)$ is a fundamental solution of
$\nabla_{\wt{\rho}_h}\medspace\Psi_{\wt{\rho}_h}=0$ and, for any $b\in
B_W$, $\wt\rho_h(\nu_\Psi(b))\in GL(U\fmll)$ is the monodromy of $b$
expressed in that solution.

\subsubsection{Extended Chen--Kohno isomorphisms}\label{sss:extended Chen}

Let $I_W\subset\IC[B_W]$ be the kernel of the epimorphism $\IC[B_W]\rightarrow
\IC[W]$. $I_W$ is readily seen to be the ideal generated by the augmentation ideal
$J$ of $\IC[P_W]$ inside $\IC[B_W]$. It follows from this that $I_W^m\cap\IC[P_W]
=J^m$ for any $m\geq 0$ and therefore that 
\begin{equation}\label{eq:exact completions}
0\longrightarrow
\wh{\IC[P_W]}\longrightarrow
\wh{\IC[B_W]}\longrightarrow
\IC[W]\longrightarrow 0
\end{equation}
is exact, where
$$\wh{\IC[P_W]}=\lim_{\substack{\longleftarrow\\m\rightarrow\infty}}\IC[P_W]/J^m
\qquad\text{and}\qquad
\wh{\IC[B_W]}=\lim_{\substack{\longleftarrow\\m\rightarrow\infty}}\IC[B_W]/I_W^m$$
are the prounipotent completion of $\IC[P_W]$ and the completion of $\IC[B_W]$
relative to the homomorphism $\IC[B_W]\rightarrow\IC[W]$ respectively \cite{Ha}.\\

Let $\Psi:\wt{\VA}\rightarrow A$ be an invertible solution of the holonomy equations
\eqref{eq:holonomy} and let $\Theta_\Psi$ be the composition of the monodromy map
$\nu_\Psi$ with the projection to the quotient $A\rtimes W$ of $A\rtimes B_W$. One readily checks that $\Theta_\Psi$ maps $I_W$ into the ideal $A_+\rtimes W$ of positive
elements \wrt the grading given by $\deg(t_\alpha)=1$ for $\alpha\in\Phi_+$ and $\deg(w)=0$ for $w\in W$. $\Theta_\Psi$ therefore factors through a map $\wh
{\Theta}_\Psi:\IC[B_W]\rightarrow A\rtimes W$ which fits into the commutative
diagram

$$\begin{diagram}
0&\rTo&\wh{\IC[P_W]}		&\rTo&\wh{\IC[B_W]}		 	&\rTo&\IC[W]	  &\rTo&0\\
  &	  &\dTo^{\wh{\mu_\Psi}}	&	&\dTo^{\wh{\Theta}_\Psi}	&	&\dTo^{\id} &      &  \\	
0&\rTo&A					&\rTo&A\rtimes W		 	&\rTo&\IC[W]	  &\rTo&0
\end{diagram}$$

where $\wh{\mu_\Psi}$ is the monodromy map of \S \ref{ss:Chen}. The exactness
of the rows and theorem \ref{th:Chen} readily yield the following

\begin{theorem}\label{th:extended Chen}
The monodromy map $\wh{\Theta}_\Psi$ is an isomorphism. In particular,
the exact sequence \eqref{eq:exact completions} is (non--canonically) split.
\end{theorem}

\subsubsection{} 
\label{sss:monodromy for Brieskorn}

Recall that $W$ possesses a presentation on the reflections $s_i=s_{\alpha_i}$,
$i\in\bfI$ corresponding to the walls of the chamber $\C$ with relations  $s_i^2=
1$ and, for any $i\neq j\in\bfI$,
$$\underbrace{s_is_j\cdots}_{m_{ij}}=\underbrace{s_js_i\cdots}_{m_{ij}}$$
where the number $m_{ij}$ of factors on each side is the order of $s_is_j$
in $W$. Assume henceforth that the base point $v_0$ lies in $\C$. Then, by
Brieskorn's theorem \cite{Br}, $B_W$ is canonically isomorphic to the algebraic
braid group of $W$, that is the group $B_D$ presented on generators $S_i$,
$i\in\bfI$ with relations
\begin{equation}\label{eq:Brieskorn}
\underbrace{S_{i}S_{j}\cdots}_{m_{ij}}=
\underbrace{S_{j}S_{i}\cdots}_{m_{ij}}
\end{equation}
for any $i\neq j\in\bfI$. The image of $S_i$ in $B_W$ is a generator of
monodromy around the image of the hyperplane $\Ker(\alpha_i)$ in $
\VA/W$. The isomorphism is compatible with the diagrams
$$\begin{diagram}[height=1.7em,width=2em]
B_W&		&\rTo&		&B_D\\
       &\rdTo 	&	&\ldTo_q	&	\\
       &		&W	&		&
\end{diagram}
\qquad\text{and}\qquad
\begin{diagram}[height=1.7em,width=2em]
\pi_1(\VA/W,[v_0])	&		&\rTo^\imath&		&\pi_1(\VA/W,[v_0'])\\
			   	&\rdTo 	&	 	   &\ldTo	&	\\
			   	&		&B_D	   &		&
\end{diagram}$$
where $q$ maps $S_i$ to $s_i$, $v_0'\in\C$ is another basepoint and $\imath$
is the canonical identification induced by the contractibility of $\C$.

\subsubsection{}

Let $\F$ be a \fmnss and $b$ a positive, real basis of $V^*$ adapted to $\F$.
Note that the solution $\Psi_\F^b$ lifts uniquely to $\wt{\VA}$ since it is defined
on $\C\ni v_0$.

\begin{theorem}\label{th:monodromy factor}
If $\F$ contains $\IC\alpha_i$, then
$$\nu_{\Psi_\F^b}(S_i)=\exp(\pi\sqrt{-1}\cdot t_{\alpha_i})\cdot S_i$$
\end{theorem}

\subsubsection{} 
\label{sss:DCP braid relations}

For any $i\in\bfI$, set
$$S_i^\nabla=\exp(\pi\sqrt{-1}\cdot t_{\alpha_i})\cdot S_i\in A\rtimes B_W$$
Theorem \ref{th:monodromy factor} and proposition \ref{pr:nu monodromy}
allow one to concisely describe the monodromy of the holonomy equations
in terms of the elements $S_i^\nabla$ and the associators $\Phi^\beta_{\G
\F}$ corresponding to a fixed positive, real adapted family $\beta$\footnote
{If $\Phi$ is a crystallographic root system, there are several preferred ways
to obtain such families $\beta=\{x_B\}$ such that $\beta\subset\Phi_+$. One
may for example set $x_B=\theta_B$ where the latter is the highest root of
the irreducible subroot system $B\cap\Phi$ relative to the basis $\Delta_B=
\Delta\cap B$, or $x_B=\sum\alpha_i$ where the sum ranges over those $i
\in\bfI$ such that $\alpha_i\in\Delta_B$. These two choices coincide for root
systems of type $\sfA_n$. The latter is the choice adopted by Drinfeld \cite{Dr3}
for root systems of type $\sfA_n$ and by \DCP for general root systems.}. Indeed,
let $\F$ be a fixed fundamental maximal nested set. If $\IC\alpha_i\in\F$, then
$$\nu_{\Psi_\F^\beta}(S_i)=S_i^\nabla$$
Otherwise, if $\G$ is a \fmns such that $\IC\alpha_i\in\G$, then
$$\nu_{\Psi_\F^\beta}(S_i)=
\nu_{\Psi_\G^\beta\cdot\Phi_{\G\F}^\beta}(S_i)=
\Phi_{\F\G}^\beta\cdot S_i^\nabla\cdot\Phi_{\G\F}^\beta$$

It follows in particular from \eqref{eq:Brieskorn} that the elementary associators
$\Phi_{\G\F}^\beta$ satisfy the following additional property

\begin{itemize}
\item {\bf braid relations} : if $(\F_i,\F_j)$ is a pair of \fmnss such that $\IC\alpha
_i\in\F_i$ and $\IC\alpha_j\in\F_j$, then
\begin{equation}\label{eq:DCP braid relations}
\underbrace{
S_i^\nabla\cdot\Ad(\Phi_{\F_i\F_j}^\beta)(S_j^\nabla)
\cdots}_{m_{ij}}=
\underbrace{
\Ad(\Phi_{\F_i\F_j}^\beta)(S_j^\nabla)\cdot S_i^\nabla
\cdots}_{m_{ij}}
\end{equation}
\end{itemize}

\subsection{Appendix : proof of theorem \ref{th:Chen}}\label{app:Chen}

We prove below that the map
$$\wh{\mu_\Psi}:\wh{\IC[\pi]}\longrightarrow A$$
defined in \S \ref{ss:Chen} is an isomorphism. The surjectivity of $\wh{\mu_\Psi}$
is due to Chen \cite[thm 3.4.1]{Chn1} and its injectivity to Kohno \cite{Ko1,Ko3}.
We merely repeat here Chen's simple proof and give an alternative approach
to Kohno's based on the observation, used by Bar--Natan in the case of the
Coxeter arrangements of type $\sfA_n$ \cite[prop. 3.6]{BN}, that the defining
relations \eqref{eq:rank 2} of the holonomy algebra $A$ may be obtained by
linearising suitable commutation relations in $\pi$.\\

Note first that $\wh{\IC[\pi]}$ and $A$ are endowed with decreasing $\IN$--filtrations.
It therefore suffices to show that $\gr(\wh{\mu_\Psi}):\gr(\wh{\IC[\pi]})\rightarrow
\gr(A)=A$ is an isomorphism. Let $x\in X$ and $\gamma_x\in\pi$ a generator
of monodromy around the hyperplane $x^\perp$. Picard iteration readily shows
that, mod $A_+^2$
\begin{equation}\label{eq:Picard}
\mu_\Psi(\gamma_x-1)=
\int_{\gamma_x}\sum_{x'\in X}\frac{dx'}{x'}\cdot t_{x'}=
2\pi i t_x
\end{equation}
so that $\gr(\wh{\mu_\Psi})(\gamma_x-1)=2\pi i t_x$. In particular, $\gr(\wh{\mu_
\Psi})$ is surjective. To construct an inverse to $\gr(\wh{\mu_\Psi})$ note first that,
for any $\gamma,\zeta\in\pi$, the identity
$$(\gamma\zeta\gamma^{-1}-1)-(\zeta-1)=
\left((\gamma-1)(\zeta-1)-(\zeta-1)(\gamma-1)\right)\gamma^{-1}$$
shows that the image of $\zeta-1$ in $J/J^2$ only depends upon the
conjugacy class of $\zeta$ in $\pi$. In particular, if $x\in X$, the class of $\gamma
_x-1$ in $J/J^2$ does not depend upon the choice of the generator of monodromy $\gamma_x$ around $x^\perp$ since any two such choices are conjugate in $\pi$. Define now $\nu:A_1\rightarrow J/J^2$ by
$$t_x\longrightarrow(2\pi i)^{-1}\cdot(\gamma_x-1)+J^2$$
for any $x\in X$. We claim that $\nu$ extends to a homomorphism $A\rightarrow
\gr(\wh{\IC[\pi]})$. It suffices to show that the elements $\delta_x=\gamma_x-1$
satisfy the relations \eqref{eq:rank 2} modulo $J^3$.  Let $B\subset V^*$ be a
2--dimensional subspace spanned by elements of $X$. We shall need the
following result whose proof is given below

\begin{lemma}\label{le:central element}
Let $x_1,\ldots,x_m$ be an enumeration of $B\cap X$. Then, there exist
generators of monodromy $\gamma_i\in\pi$ around each $x_i^\perp$ such
that the product $\gamma_1\cdots\gamma_m$ commutes with each
$\gamma_i$.
\end{lemma}

Let now $\gamma_i$, $i=1,\ldots,m$, be as in lemma \ref{le:central element} and
set $\delta_i=\gamma_i-1\in J$. Replacing each $\gamma_j$ by $\delta_j+1$
in $\gamma_i\cdot\gamma_1\cdots\gamma_m=\gamma_1\cdots\gamma_m
\cdot\gamma_i$ yields $[\delta_i,\sum_j\delta_j]=0\mod J^3$ as required. Thus,
$\nu$ extends to a homomorphism $A\rightarrow\gr(\wh{\IC[\pi]})$ satisfying $\mu\circ\nu=\id$. Moreover, $\nu\circ\mu=\id$ since this holds on any generator
of monodromy $\gamma_x$, $x\in X$ and these generate $\pi$ \cite[prop A.2]
{BMR} \halmos\\

{\sc Proof of lemma \ref{le:central element}.}\footnote{I owe the proof of this
lemma to D. Bessis and J. Millson.}
Let $\A_B=\bigcup_{x\in X\cap B}x^\perp\subset V$ be the arrangement
determined by $B$ and $D\subset V$ an open ball centered at $v\in B^\perp$
such that $D$ does not intersect any hyperplane $y^\perp$, $y\in X\setminus
B$. Since the composition
$$D\setminus\A_B\hookrightarrow V\setminus\A\hookrightarrow V\setminus\A_B$$
is a homotopy equivalence, $\pi_1(D\setminus\A_B)$ embeds in $\pi$. The
elements $\gamma_i$ will be chosen in $\pi_1(D\setminus\A_B)\cong\pi_1
(V\setminus\A_B)$. Let $V'\subset V$ be a complementary subspace to $
B^\perp$. The corresponding projection $V\rightarrow V'$ induces a homotopy
equivalence $V\setminus\A_B\sim V'\setminus\A_B$ and therefore an
isomorphism of $\pi_1(V\setminus\A_B)$ with the fundamental group
of the complement in $V'$ of the lines $L_x=V'\cap x^\perp$,
$x\in X\cap B$. Consider now the Hopf fibration
$$\IC^*\longrightarrow
V'\setminus\bigcup_i L_{x_i}\longrightarrow
\IP^1\setminus\{z_1,\ldots,z_m\}$$
where $z_i=[L_{x_i}]\in\IP(V')\cong\IP^1$. Since $m\geq 1$, the fibration is trivial
and the image of $\pi_1(\IC^*)$ in $\pi_1(V'\setminus\bigcup_i L_{x_i})$ is central.
Thus, if $\ol{\gamma}_1,\ldots,\ol{\gamma}_m\in\pi_1(\IP^1\setminus\{z_1,\ldots,
z_m\})$ are small loops around $z_1,\ldots,z_m$ such that $\ol{\gamma}_1\cdots
\ol{\gamma}_m=1$ and each $\ol{\gamma}_i$ is lifted to a generator of monodromy $\gamma_i$ around $L_{x_i}$, the product $\gamma_1\cdots\gamma_m$ is central
in $\pi_1(V'\setminus\bigcup_i L_{x_i})$\halmos\\

\remark
By \eqref{eq:Picard}, the map $\gr(\mu_\Psi)$ does not depend upon
the choice of $\Psi$. Thus, the monodromy of $\nabla$ yields a canonical
isomorphism
identification $\gr(\wh{\IC[\pi]})\cong A$.


%% file: newassociahedron.tex

\section{The \DCP associahedron $\AD$}\label{se:associahedron}

Let $D$ be a connected graph. The aim of this section is to construct
a regular cell complex $\AD$ whose face poset is that of nested sets
of connected subgraphs of $D$, ordered by reverse inclusion. When
$D$ is the Dynkin diagram of type $\sfA_{n-1}$, $\AD$ is isomorphic
to Stasheff's associahedron $K_n$ \cite{St1}. More generally, if $D$
is the Coxeter graph of an irreducible, finite Coxeter group $W$, $\AD$
is isomorphic to the cell complex constructed by \DCP inside their
wonderful model of the reflection arrangement of $W$ \cite[\S 3.2]{DCP2}
(see also \cite{Ga}). For this reason, we call $\AD$ the {\it \DCP associahedron}
corresponding to $D$. When $D$ is the affine Dynkin diagram of type
$\sfA_{n-1}$, $\AD$ is isomorphic to Bott and Taubes' cyclohedron $W_n$
\cite {BT}.\footnote{The reader should  be cautioned that the \DCP associahedra
corresponding to Dynkin diagrams of finite type differ from the generalised
associahedra defined by Fomin and Zelevinsky \cite{FZ,CFZ}, since
the former do not depend upon the multiplicities of the edges of the
diagram. For example, the \DCP associahedra of type $\sfA_{n-1},
\sfB_{n-1},\sfC_{n-1}$ are all isomorphic to the associahedron $K_n$
while the \FZ associahedra of types $\sfB_n,\sfC_n$ are homeomorphic
to the cyclohedron $W_n$.}\\

After defining the face poset of $\AD$ in \S \ref{ss:nested}, we realise
$\AD$ as a convex polytope in \S \ref{ss:realisation}, thereby settling
its existence and proving its contractibility. The simple connectedness
of $\AD$ will be used in section \ref{se:qcoxeter} to prove an analogue
for \qcas of Mac Lane's coherence theorem for monoidal categories.
We then prove in \S \ref{ss:faces of AD} that the faces of $\AD$ are
isomorphic to products of associahedra corresponding to subquotients
of $D$, a well--known fact for the associahedron $K_n$. In particular,
each facet of $\AD$ is a product of two associahedra, one corresponding
to a proper, connected subgraph $B$ of $D$ and the other to the quotient
graph $D/B$. Finally, in \S \ref{ss:edges}--\S \ref{ss:2 faces} we describe
the edges and 2--faces of $\AD$ explicitly. Interestingly perhaps, the latter
turn out to be squares, pentagons or hexagons.

\subsection{The poset $\ND$ of nested sets on $D$}
\label{ss:nested}

By a {\it diagram} we shall mean a non--empty undirected graph $D$ with no
multiple edges or loops. We denote the set of vertices of $D$ by $V(D)$ and
set $|D|=|V(D)|$. A {\it subdiagram} $B\subset D$ is a full subgraph of $D$, that
is a graph consisting of a subset $V(B)$ of vertices of $D$, together with all
edges of $D$ joining any two elements of $V(B)$. We will often abusively identify
such a $B$ with its set of vertices and write $\alpha\in B$ to mean $\alpha\in V(B)$.
The {\it union} $B_1\cup B_2$ of two subdiagrams $B_1,B_2\subset D$ of $
D$ is the subdiagram having $V(B_1)\cup V(B_2)$ as set of vertices.\\

Two subdiagrams $B_1,B_2\subseteq D$ are {\it orthogonal} if no two vertices
$\alpha_1\in B_1,\alpha_2\in B_2$ are joined by an edge in $D$. $B_1$ and
$B_2$ are {\it compatible} if either one contains the other or they are orthogonal.
Assume henceforth that $D$ is connected.

\begin{definition}\label{de:nested}
A {\it nested set} on $D$ is a collection $\H$ of pairwise compatible, connected
subdiagrams of $D$ which contains $D$.\footnote{This is the opposite of Stasheff's
convention in which the faces of the associahedron $K_n$ are labelled by
consistent bracketings of a monomial $x_1\cdots x_n$ which do {\it not} contain
the big bracket $(x_1\cdots x_n)$, but is better suited to our needs.}\footnote
{Such collections should perhaps be called {\it fundamental} nested sets on $D$
since, when $D$ is the graph of a finite Coxeter group $W$ they correspond, via
the dictionary of \S \ref{sss:subspace dictionary}--\S \ref{sss:graph dictionary}, to
\fnss of subspaces spanned by the roots of $W$. Since general nested sets of
subspaces do not seem to have an analogue for diagrams however, we prefer
to omit the adjective fundamental when speaking about diagrams.}
\end{definition}

We denote by $\ND$ the partially ordered set of nested sets on $D$, ordered
by reverse inclusion. $\ND$ has a unique maximal element $\hat{1}=\{D\}$.
Its minimal elements are the maximal nested sets. When $D$ is the Dynkin
diagram of type $\sfA_{n-1}$, $\ND$ is the face poset of the associahedron
$K_n$ by remark \rkref{rk:Kn}.


\subsection{Convex realisation of $\AD$}
\label{ss:realisation}

Recall that a $CW$--complex $X$ is {\it regular} if all its attaching maps
are homeomorphisms \cite[\S IX.6]{Ms}. In this case, an induction on
the skeleton of $X$ shows that its cellular isomorphism type is uniquely
determined by its face poset. This justifies the following
\begin{definition}
The \DCP associahedron $\AD$ is the regular $CW$--complex whose
poset of (non--empty) faces is $\ND$.
\end{definition}

We shall prove the existence of $\AD$ by realising it as a convex polytope of
dimension $|D|-1$. Our construction follows the Shnider--Sternberg--Stasheff
realisation of the associahedron $K_n$ as a truncation of the $(n-2)$--simplex
\cite{SS}, as presented in \cite[Appendix B]{St2}, and coincides with Stasheff
and Markl's convex realisation of the cyclohedron $W_n$ \cite{St2,Ma} when
$D$ is the affine Dynkin diagram of type $\sfA_{n-1}$.\\

Let $c$ be a function on the set of connected subdiagrams of $D$ with values
in $\IR_+^*$ such that
\begin{equation}\label{eq:c}
c(B_1\cup B_2)>c(B_1)+c(B_2)
\end{equation}
whenever $B_1$ and $B_2$ are not compatible. An example of such a $c$
is given by $c(B)=3^{|B|}$. Let $\{t_\alpha\}_{\alpha\in D}$ be the canonical
coordinates on $\IR^{|D|}$ and consider, for any connected $B\subseteq D$,
the linear hyperplane
$$\L_B^c=\{t\in\IR^{|D|}|\sum_{\alpha\in B}t_\alpha=c(B)\}\subset\IR^{|D|}$$
Consider next the convex polytope
$$P_D^c=\{t\in\IR^{|D|}|\msp \sum_{\alpha\in D}t_\alpha=c(D),\msp\sum_
{\alpha\in B}t_{\alpha}\geq c(B)\quad\text{for any connected $B\subsetneq D$}\}$$

\begin{theorem}\label{th:realisation}\hfill
\begin{enumerate}
\item The polytope $P_D^c$ has nonempty interior in the hyperplane $\L_D^c$.
\item For any connected subdiagrams $B_1,\ldots,B_m\subsetneq D$, the
intersection
$$P_{D,B_1,\ldots,B_m}^c=P_D^c\cap\bigcap_{i=1}^{m}\L_{B_i}^c$$
is nonempty if, and only if $B_1,\ldots,B_m$ are pairwise compatible.
\item If $B_1,\ldots,B_m$ are pairwise compatible and distinct, $P^c_
{D,B_1,\ldots,B_m}$ is a face of $P_D^c$ of dimension $|D|-1-m$.
\item All nonempty faces of $P_D^c$ are obtained in this way.
\end{enumerate}
\end{theorem}
\proof The proof given in \cite[appendix B]{St2} in the case when $D$ is the
Dynkin diagram of type $\sfA_{n-1}$ carries over easily to the general case
\halmos

\begin{corollary}\label{co:realisation}
The map
$$\H\longrightarrow P_{\H}^c=P_{D,B_1,\ldots,B_m}^c$$
where $\H=\{D,B_1,\ldots,B_m\}$ is a nested set on $D$ is an isomorphism
between $\ND$ and the poset of nonempty faces of $P_D^c$.
\end{corollary}

Thus, for any function $c$ satisfying \eqref{eq:c}, the polytope $P^c_D$
gives a convex realisation of the associahedron $\AD$. In particular

\begin{corollary}\label{co:contractibility}
The \DCP associahedron $\AD$ is contractible.
\end{corollary}



\remark\rklabel{rk:1-skeleton}
By theorem \ref{th:realisation}, the \mnss on $D$, which label the vertices
of $\AD$, are of cardinality $|D|$ and any $\H\in\ND$ of cardinality $|D|-1$
is contained in exactly two maximal nested sets. Thus, the 1--skeleton of
$\AD$ may equivalently be described as having a $0$--cell for each \mns
$\F$ on $D$ and a $1$--cell between $\F$ and $\G$ if, and only if $\F$
and $\G$ differ by an element. In particular, the connectedness of $\AD$
gives another proof of proposition \ref{pr:connectivity} for Coxeter
arrangements.\\

\remark
When $D$ is the graph of an irreducible, finite Coxeter system $(W,S)$, the
associahedron $\AD$ may be obtained more geometrically as follows  \cite
[\S 3.2]{DCP2}. Let $\A\subset V$ be the complexified reflection arrangement
of $W$ and $Y_X$ the wonderful model of $\VA=V\setminus\A$ described in
section \ref{se:DCP}. The irreducible component $\D_{V^*}\subset Y_X$ of
the exceptional divisor corresponding to $V^*$ is a smooth projective variety
and $Y_X$ is the total space of a line bundle over $\D_{V^*}$ in such a way
that the corresponding action of $\IC^*$ agrees, on $Y_X\setminus\D\cong
\VA$ with its natural action on $V$. Let $\A_\IR\subset V_\IR$ be the real
reflection arrangement of $W$, $\C\subset V_\IR\setminus\A_\IR$ the
chamber corresponding to $S$ and $\ol{\C}$ the closure of $\C$ in $Y_X$.
Then, the intersection $\ol{\C}\cap\D_{V^*}$ possesses a regular cellular
structure with corresponding face poset given by $\ND$.\\

We shall often identify a \ns $\H\in\ND$ with the corresponding face of $\AD$
and speak of the dimension $\dim(\H)=|D|-|\H|$ of $\H$ to mean the dimension
of that face.

\subsection{The rank function of $\ND$}\label{ss:rank}

For any nested set $\H\in\ND$ and $B\in\H$, set
\begin{equation}\label{eq:i_H(B)}
i_\H(B)=B_1\odots{\cup}B_m
\end{equation}
where the $B_i$ are the maximal elements of $\H$ properly contained in $B$.

\begin{definition}
Set
$$\aalpha^{B}_{\H}=B\setminus i_\H(B)
\aand
n(B;\H)=|\aalpha^{B}_{\H}|$$
\end{definition}
Note that $n(B;\H)\geq 1$. Indeed, if $m=1$, then $n(B;\H)=|B\setminus B_1|
\geq 1$. Otherwise, $B_1,\ldots,B_m$ are necessarily pairwise orthogonal,
their union is disconnected and cannot be equal to $B$. Note in passing that
$B_1,\ldots,B_m$ are the connected components of $B\setminus\aalpha^B_\H$.
In particular, the latter lie in $\H$. Set now
$$n(\H)=\sum_{B\in\H}(n(B;\H)-1)$$
The following is an analogue of propositions \ref{pr:n(F)} and \ref{pr:max nested}

\begin{proposition}\label{pr:n(H)}\hfill
\begin{enumerate}
\item For any \ns $\H\in\ND$, $$n(\H)=|D|-|\H|=\dim(\H)$$
\item If $\H$ is a maximal nested set, then $n(B;\H)=1$ for any $B\in\H$.
\item Any maximal nested set is of cardinality $|D|$.
\end{enumerate}
\end{proposition}
\proof (i) If $|\H|=1$, then $\H=\{D\}$ and $n(\H)=|D|-1$ as required.
Assume now that $|\H|\geq 2$ and let $D_1,\ldots,D_m\subsetneq D$
be the proper, maximal elements in $\H$ so that
$$\H=\{D\}\sqcup\H_1\odots{\sqcup}\H_m$$
where $\H_i$ is a \ns on $D_i$. By induction,
$n(\H_i)=|D_i|-|\H_i|$ whence
\begin{equation*}
\begin{split}
n(\H)
&=
(n(D;\H)-1)+\sum_{i}n(\H_i)\\
&= 
(n(D;\H)-1)+\sum_{i}(|D_i|-|\H_i|)\\
&=
|D|-|\H|
\end{split}
\end{equation*}
(ii) Let $B\in\H$ and $\alpha\in\aalpha^{B}_{\H}$. Since the connected
components of $B\setminus\alpha$ are compatible with the elements
of $\H$, they lie in $\H$ by maximality whence $|\aalpha^{B}_{\H}|=1$.
(iii) follows from (ii) and (i) \halmos\\

If $\F$ is a maximal nested set and $B\in\F$, we denote the unique
element of $\aalpha^B_\F$ by $\alpha^B _\F$.

\subsection{Quotienting by a subdiagram}\label{ss:quotienting}

We define in this subsection the quotient $D/B$ of $D$ by a proper
subdiagram $B$ and relate nested sets on $D/B$ with those on $D$.
Let $B_1,\ldots,B_m$ be the connected components of $B$.

\begin{definition}\label{de:D/D}
The set of vertices of the diagram $D/B$ is $V(D)\setminus V(B)$.
Two vertices $\alpha\neq\beta$ of $D/B$ are linked by an edge if,
and only if the following holds in $D$
$$\alpha\nperp\beta
\oor
\alpha,\beta\nperp B_i
\quad\text{for some}\quad
i=1,\ldots,m$$
\end{definition}

For any connected subdiagram $C\subseteq D$ not contained in $B$,
we denote by $\ol{C}\subseteq D/B$ the connected subdiagram with
vertex set $V(C)\setminus V(B)$. We shall need the following

\begin{lemma}\label{le:quotient compatibility}
Let $C_1,C_2\nsubseteq B$ be two connected subdiagrams of $D$
which are compatible. Then
\begin{enumerate}
\item $\ol{C}_1,\ol{C}_2$ are compatible unless $C_1\perp C_2$ and
$C_1,C_2\nperp B_i$ for some $i$.
\item If $C_1$ is compatible with every $B_i$, then $\ol{C}_1$ and $
\ol{C}_2$ are compatible.
\end{enumerate}
In particular, if $\F$ is a \ns on $D$ containing each $B_i$, then $\ol
{\F}=\{\ol{C}\}$, where $C$ runs over the elements of $\F$ such that
$C\nsubseteq B$, is a \ns on $D/B$.
\end{lemma}
\proof
(i) Clearly, if $C_1\subset C_2$ or $C_2\subset C_1$ then $\ol{C}_1
\subset\ol{C}_2$ or $\ol{C}_2\subset\ol{C}_1$ respectively. If, on the
other hand, $C_1\perp C_2$ then $\ol{C}_1\perp\ol{C}_2$ if, and only
if, for any connected component $B_i$ of $B$ at least one of $C_1,
C_2$ is perpendicular to $B_i$. (ii) We may assume by (i) that $C_1
\perp C_2$. Since $C_1$ and $B_i$ are compatible for any $i$ and
$C_1\nsubseteq B_i$, either $B_i\subset C_1$, in which case $B_i
\perp C_2$, or $B_i\perp C_1$. $\ol{C}_1$ and $\ol{C}_2$ are therefore
compatible by (i) \halmos\\

Let now $A$ be a connected subdiagram of $D/B$ and denote by
$\wt{A}\subseteq D$ the connected sudbdiagram with vertex set
\begin{equation}\label{eq:A tilde}
V(\wt{A})=V(A)\bigcup_{i:B_i\nperp V(A)} V(B_i)
\end{equation}
Clearly, $A_1\subseteq A_2$ or $A_1\perp A_2$ imply $\wt{A_1}
\subseteq\wt{A_2}$ and $\wt{A_1}\perp\wt{A_2}$ respectively, so
the lifting map $A\rightarrow\wt{A}$ preserves compatibility.\\

For any connected subdiagrams $A\subseteq D/B$ and $C\subseteq
D$, we have
\begin{equation}\label{eq:bar tilde}
\ol{\wt{A}}=A
\aand\wt{\ol{C}}=C\bigcup_{i:B_i\nperp C} B_i
\end{equation}
In particular, $\wt{\ol{C}}=C$ if, and only if, $C$ is compatible with $B_1,
\ldots,B_m$ and not contained in $B$. The applications $C\rightarrow
\ol{C}$ and $A\rightarrow\wt{A}$ therefore yield a bijection between the
connected subdiagrams of $D$ which are either orthogonal to or strictly
contain each $B_i$ and the connected subdiagrams of $D/B$. By lemma
\ref{le:quotient compatibility}, this bijection preserves compatibility and
therefore induces an embedding $\N_{D/B}\hookrightarrow\ND$. This
yields an embedding
\begin{equation}\label{eq:embedding}
\N_{B_1}\odots{\times}\N_{B_m}\times\N_{D/B}\hookrightarrow\ND
\end{equation}
with image the poset of nested sets on $D$ containing each $B_i$.

\subsection{Unsaturated elements and the faces of $\AD$}
\label{ss:faces of AD}

We show below that the faces of $\AD$ are products of associahedra
corresponding to subquotient diagrams of $D$. Let $\H$ be a \ns on
$D$.
\begin{definition}
An element $B\in\H$ is called unsaturated if $n(B;\H)\geq 2$.
\end{definition}

Let $\AD^\H$ be the face of the associahedron $\AD$ corresponding
to $\H$. The face poset of $\AD^\H$ is the poset $\ND^\H$ of \nss on
$D$ containing $\H$.

\begin{proposition}\label{pr:ND^H}
As posets,
$$\ND^\H\cong
\prod_{C\in\H} \N_{C/ i_\H(C)}\cong
\prod_{j=1}^p \N_{D_j/i_\H(D_j)}$$
where $D_1,\ldots,D_p$ are the unsaturated elements of $\H$. In particular,
$$\AD^\H\cong\prod_{j=1}^p \A_{D_j/i_\H(D_j)}$$
as $CW$--complexes.
\end{proposition}
\proof Let $B_1,\ldots,B_m$ be the proper maximal elements of $\H$,
so that $i_\H(D)=B_1\odots{\cup}B_m$ and let $\H_i$ be the \ns on
$B_i$ induced by $\H$. The embedding \eqref{eq:embedding} yields
an isomorphism
$$\ND^\H\cong
\N_{D/i_\H(D)}\times\N_{B_1}^{\H_1}\odots{\times}\N_{B_m}^{\H_m}$$
The first isomorphism now follows from an easy induction, the second
from the fact that $\N_{C/i_\H(C)}$ consists of a single element if $C$
is saturated. The corresponding description of $\AD^\H$ follows from
the fact that a regular $CW$--complex is determined by its face poset
\halmos\\

\remark The isomorphism $\prod_{j=1}^p \N_{D_j/i_\H(D_j)}
\longrightarrow\ND^\H$ is explicitly given by
\begin{equation}\label{eq:explicit iso}
\{\K_j\}_{j=1}^p\longrightarrow
\H\thickspace\bigcup_{j=1}^p\thickspace(\wt{\K_j}\setminus\{D_j\})
\end{equation}
where, for a \ns $\K_j$ on $D_j/i_\H(D_j)$, $\wt{\K_j}$ is the \ns on $D_j
$ obtained by lifting the elements of $\K_j$ to connected subdiagrams of
$D_j$.\\




By theorem \ref{th:realisation}, the facets of $\AD$ are labelled by the \nss on
$D$ of the form $\H=\{D,B\}$ where $B$ is a proper, connected subdiagrams
of $D$.
\begin{corollary}\label{co:facets of AD}
The facet of $\AD$ corresponding to $B$ is isomorphic, as cell complex, to the
product $\A_{B}\times\A_{D/B}$.
\end{corollary}
When $D$ is the finite or the affine Dynkin diagram of type $\sfA_{n-1}$, we
recover from corollary \ref{co:facets of AD} the familiar fact that each facet
of the associahedron $K_n$ or of the cyclohedron $W_n$ is the product $K_r
\times K_s$ or $K_r\times W_s$, with $r+s=n+1$ of two smaller associahedra
or an associahedron and a cyclohedron respectively.\\

\remark The set of Dynkin diagrams of is not closed under quotienting. For
example, if $D$ is the Dynkin diagram of type $\sfD_n$ and $\alpha$ is the
trivalent node of $D$, then $D/\alpha$ is the affine Dynkin diagram of type
$\sfA_2$ if $n=4$ and a tadpole if $n\geq 5$. Thus, if $D$ is a Dynkin
diagram with a trivalent node other than $\sfD_4$, the faces of $\AD$
are products of associahedra some of which correspond to non--Dynkin
diagrams.

\subsection{An alternative description of the lifting map}
\label{ss:explicit lift}

We shall need an alternative description of the map \eqref{eq:explicit iso}.
Let $\H$ be a \ns on $D$ with $|\H|<|D|$ and let $D_1,\ldots,D_p$ be the
unsaturated elements of $\H$. For any $1\leq j\leq p$, let
$\aalpha_j=\aalpha^{D_j}_{\H}$
and, for any subset $\emptyset\neq\bbeta_j\subseteq\aalpha_j$, set
\begin{equation}\label{eq:Dbeta}
\bbeta_j^c=\aalpha_j\setminus\bbeta_j
\qquad\text{and}\qquad
D_{\bbeta_j}=\con{D_j}{\bbeta_j^c}{\bbeta_j}
\end{equation}
where the latter denotes the connected component of $D_j\setminus\bbeta_j^c$
containing $\bbeta_j$ if one such exists and the empty set otherwise.

\begin{lemma}\label{le:explicit lift}
Let $B_j\subset D_j/i_\H(D_j)$ be the subdiagram with vertex set $\bbeta_j$.
\begin{enumerate}
\item $D_{\bbeta_j}$ is non--empty if, and only if $B_j$ is connected.
\item When that is the case,
$$\wt{B_j}=D_{\bbeta_j}\aand\ol{D_{\bbeta_j}}=B_j$$
\end{enumerate}
where $\ol{\cdot}$ and $\wt{\cdot}$ are the quotient and lifting maps for the
quotient $D_j/i_\H(D_j)$.
\end{lemma}
\proof
If $D_{\bbeta_j}$ is non--empty, $\ol{D_{\bbeta_j}}$ is a subdiagram of $D_j
/i_\H(D_j)$ with vertex set $\bbeta_j$. Thus, $\ol{D_{\bbeta_j}}=B_j$ and the
latter is connected since $D_{\bbeta_j}$ is. Conversely, if $B_j$ is connected,
$\wt{B_j}$ is a connected subdiagram of $D_j\setminus\bbeta_j^c$ containing $\bbeta_j$. Thus, $D_{\bbeta_j}$ is non--empty and, by \eqref{eq:bar tilde}, $D_{\bbeta_j}=\wt{\ol{D_{\bbeta_j}}}=\wt{B_j}$ \halmos

\subsection{Edges of the associahedron $\AD$}
\label{ss:edges}

Let $\H\in\ND$ be a nested set of dimension 1. By proposition \ref{pr:n(H)},
$\H$ has a unique unsaturated element $B$ and $\aalpha^B_\H$ consists
of two vertices $\alpha_1,\alpha_2$. Thus, $\ol{B}=B/i_\H(B)$ is the connected
diagram with vertices $\ol{\alpha}_1,\ol{\alpha}_2$ and $\A_{\ol{B}}$ is the
interval
\begin{equation}
\begin{diagram}[width=3em]\label{eq:quotient edge}
\{\ol{B},\ol{\alpha}_1\}&&\rOnecell^{\{\ol{B}\}}&&\{\ol{B},\ol{\alpha}_2\}
\end{diagram}
\end{equation}
Setting $B_1=\con{B}{\alpha_2}{\alpha_1}$ and $B_2=\con{B}{\alpha_1}
{\alpha_2}$, we see that, by proposition \ref{pr:ND^H} and lemma \ref
{le:explicit lift}, the edge of $\AD$ corresponding to $\H$ is of the form
\begin{diagram}[width=3em]
\H\cup\{B_1\}&&\rOnecell^{\H}&&\H\cup\{B_2\}
\end{diagram}

\subsection{Two--faces of the associahedron $\AD$}
\label{ss:2 faces}

We work out below the two--faces of $\AD$ and show that they are squares,
pentagons or hexagons. Let $\H\in\ND$ be a nested set of dimension 2. By
proposition \ref{pr:n(H)}, $\H$ either has two unsaturated elements $B_1,B_2$,
and $|\aalpha^{B_1}_{\H}|=2=|\aalpha^{B_2}_{\H}|$, or a unique unsaturated
element $B$ and $|\aalpha^B_\H|=3$. We treat these two cases separately.

\subsubsection{Square 2--faces}\label{sss:square}

Assume first that $\H$ has two unsaturated elements $B_1,B_2$. By proposition
\ref{pr:ND^H}, $\AD^\H\cong\A_{B_1/i_\H(B_1)}\times\A_{B_2/i_\H(B_2)}$ is the
product of two intervals of the form \eqref{eq:quotient edge}. Thus, setting for
$i,j\in\{1,2\}$
$$\aalpha^{B_i}_\H=\{\alpha_i^1,\alpha_i^2\}
\qquad\text{and}\qquad
B_i^j=\con{B_i}{\alpha_i^{3-j}}{\alpha_i^j}\subset B_i$$
we see that $\AD^\H$ is the square
\begin{diagram}
\H^{1,1}&\hLine	&\H^{1,2}\\
\vLine   &\H		&\vLine   \\
\H^{2,1}&\hLine	&\H^{2,2}
\end{diagram}
where $\H^{j,k}=\H\cup\{B_1^j,B_2^k\}$.

\subsubsection{Pentagonal and hexagonal 2--faces}\label{sss:pentagon}

Assume now that $\H$ has a unique unsaturated element $B$ and set $\aalpha
^B_\H=\{\alpha_1,\alpha_2,\alpha_3\}$. Then $\ol{B}=B/i_\H(B)$ is a connected
diagram with vertices $\ol{\alpha}_1,\ol{\alpha}_2,\ol{\alpha}_3$ and is therefore,
up to a relabelling of the $\alpha_i$, one of the following diagrams
$$\begin{diagram}[inline,abut,height=.8em,width=1em]
\bullet	    &\rLine	&\bullet		&\rLine&\bullet\\
\ol{\alpha}_1&		&\ol{\alpha}_2	&	   &\ol{\alpha}_3
\end{diagram}
\qquad\qquad\qquad
\begin{diagram}[inline,abut,height=.7em,width=.5em]
			&		& 	&	& 	&\ol{\alpha}_2	&	&	&	&		&\\
			&		& 	&	& 	&\bullet		&	&	&	&		&\\
			&		& 	&	&\ruL&			&\rdL&	&	&		&\\
			&		& 	&\ruL& 	&			&	&\rdL&	&		&\\
			&		&\ruL&	& 	& 			&	&	&\rdL&		&\\
\ol{\alpha}_1	&\bullet	&	&	& 	&\rLine		&	&	&	&\bullet	&\thinspace\ol{\alpha}_3
\end{diagram}$$

For any $1\leq i\leq 3$ and $1\leq j\neq k\leq 3$, set
$$B_i=\con{B}{(\aalpha^B_\H\setminus\alpha_i)}{\alpha_i} 
\qquad\text{and}\qquad
B_{jk}=\con{B}{(\aalpha^B_\H\setminus\{\alpha_j,\alpha_k\})}{\{\alpha_j,\alpha_k\}}$$

If $\ol{B}$ is the affine Dynkin diagram of type $\sfA_2$, proposition \ref{pr:ND^H}
and lemma \ref{le:explicit lift} imply that the vertices of $\A_D^\H$ are of the form
$\H\cup\{B_i,B_{ij}\}$ with $1\leq i\neq j\leq 3$, so that $\A_D^\H$ is the hexagon
\begin{diagram}[vtrianglewidth=2.5em]
&\H\cup\{B_1,B_{13}\}&&\rLine&&\H\cup\{B_3,B_{13}\}&\\
\H\cup\{B_1\ruLine(1,1),B_{12}\}&&&\H&&&\H\cup\{\luLine(1,1)B_3,B_{23}\}\\
&\rdLine(1,1)\H\cup\{B_2,B_{12}\}&&\rLine&&\H\cup\{B_2,B_{23}\}\ruLine(1,1)&
\end{diagram}
If, on the other hand $\ol{B}$ is the Dynkin diagram of type $\sfA_3$, the
vertices of $\AD^\H$ are $\H\cup\{B_i,B_{ij}\}$, with $1\leq i\neq j\leq 3$ and
$(i,j)\neq(1,3)$, and $\H\cup\{B_1,B_3\}$ so that $\AD^\H$ is the pentagon
\begin{diagram}[vtrianglewidth=2.5em]
&&\H\cup\{B_1,B_3\}&&\\
\H\cup\{B_1,B_{12}\}&\ruLine(2,1)&\H&\luLine(2,1)&\H\cup\{B_3,B_{23}\}\\
&\H\cup\rdLine(1,1)\{B_2,B_{12}\}&\rLine&\H\cup\{B_2,\ruLine(1,1)B_{23}\}&
\end{diagram}

\remark If $B$ is linear, $B_{ij}$ is empty for a unique pair $(i,j)$ and $\H$ is
therefore a pentagon. Thus, the associahedron $K_n$, and more generally
the associahedra corresponding to linear Dynkin diagrams, only have squares
and pentagons as 2--faces. The associahedra of type $\sfD_n,\sfE_6,\sfE_7,
\sfE_8$ and those corresponding to the affine Dynkin diagrams of type $\sfA_n,
\sfB_n,\sfD_n,\sfE_6,\sfE_7,\sfE_8$ on the other hand all have some hexagonal
2--faces.


%% file: quasiCoxeter.tex


\section{$D$--algebras and quasi--Coxeter algebras}\label{se:qcoxeter}

The aim of this section is to define the category of \qc algebras. We begin
in \S \ref{ss:D algebra}--\S \ref{ss:D morphisms} by describing the underlying
notion of $D$--algebras. We then give three equivalent definitions of \qc
algebras. The first two, in \S \ref{ss:qca 1} and \S \ref{ss:qca 2} respectively,
are more closely modelled on the \DCP theory of asymptotic zones reviewed
in section \ref{se:DCP} as well as Drinfeld's theory of \qbas \cite{Dr3}.
The first definition is better suited to the study of examples which are
considered in section \ref{se:qcexamples}, while the second is more
convenient to show that \qcas define representations of braid groups,
as explained in \S \ref{ss:qca braid}. The third definition, given in \S
\ref{ss:qca 3}, is the most compact one and will be used in section
\ref{se:dynkin} to study the deformation theory of \qc algebras. The
equivalence of this definition with the first two is the analogue for
\qcas of Mac Lane's coherence theorem for monoidal categories. It
relies on the simple connectedness of the \DCP associahedron $\AD$
introduced in section \ref{se:associahedron}. In \S \ref{ss:qca twisting},
we define the twisting of a \qca and show that it yields equivalent braid
group representations.

\subsection{D--algebras}\label{ss:D algebra}

Let $k$ be a fixed commutative ring with unit. By an algebra we
shall henceforth mean a unital, associative $k$--algebra. All algebra
homomorphisms will be tacitly assumed to be unital.
Let $D$ be a connected diagram. 

\begin{definition}
A $D$--algebra 
is an algebra $A$ endowed with subalgebras
$A_{D'}$ labelled by the non--empty connected subdiagrams $D'$ of $D$ such
that the following holds
\begin{itemize}
\item $A_{D'}\subseteq A_{D\dprime}$ whenever $D'\subseteq D\dprime$.
\item $A_{D'}$ and $A_{D\dprime}$ commute whenever $D'$ and $D\dprime$
are orthogonal.
\end{itemize}
\end{definition}

If $A$ is a $D$--algebra and $\alpha_i$ is a vertex of $D$ we denote $A_{\alpha
_i}$ by $A_i$. If $D_1,D_2\subseteq D$ are subdiagrams with $D_1$ connected,
we denote by $A_{D_1}^{D_2}$ the centraliser in $A_{D_1}$ of the subalgebras
$A_{D_2'}$ where $D_2'$ runs over the connected components of $D_2$.

\subsection{Examples}

Most, but not all examples of $D$--algebras arise in the following way.
An algebra $A$ is endowed with subalgebras $A_i$ labelled by the
vertices $\alpha_i$ of $D$ with $[A_i,A_j]=0$ whenever $\alpha_i$ and
$\alpha_j$ are orthogonal. In this case, letting $A_{D'}\subseteq A$ be
the subalgebra generated by the $A_i$ corresponding to the vertices
of $D'$ endows $A$ with a $D$--algebra structure.

\subsubsection{} 
\label{sss:Coxeter D}

Let $W$ be an irreducible Coxeter group with system of generators $S=
\{s_i\}_{i\in\bfI}$ and let $D$ be the Coxeter graph of $W$. For any $i\in
\bfI$, let $\IZ_2\cong W_i\subset W$ be the subgroup generated by $s_i$.
Then, $(k[W],k[W_i])$ is a $D$--algebra. Similarly, let $q_i\in k$ be invertible
elements such that $q_i=q_j$ whenever $s_i$ and $s_j$ are conjugate
in $W$ and let $\H(W)$ be the Iwahori--Hecke algebra of $W$, that is the
algebra with generators $\{S_i\}_{i\in\bfI}$ and relations
\begin{gather*}
(S_i-q_i)(S_i+q_i^{-1})=0\\
\underbrace{S_iS_j\cdots }_{m_{ij}}=
\underbrace{S_jS_i\cdots }_{m_{ij}}
\end{gather*}
where $m_{ij}$ is the order of $s_is_j$ in $W$. Then, $(\H(W),\H(W_i))$
is a $D$--algebra.

\subsubsection{} 
 
Let $\bfA=(a_{ij})_{i,j\in\bfI}$ be an irreducible, generalised Cartan matrix, $\g
=\g(\bfA)$ the corresponding Kac--Moody algebra and $\g'=[\g,\g]$ its derived
subalgebra with generators $e_i,f_i,h_i$, $i\in\bfI$ \cite{Ka}. Let $D=D(\bfA)$
be the Dynkin diagram of $\g$, that is the connected graph having $\bfI$ as its
vertex set and an edge between $i$ and $j$ if $a_{ij}\neq 0$. For any $i\in\bfI$,
let $\sl{2}^i\subseteq\g'$ be the three--dimensional subalgebra spanned by $
e_i,f_i,h_i$. Then $(U\g',U\sl{2}^i)$ is a $D$--algebra over $k=\IC$. Similarly,
if $\bfA$ is symmetrisable and $\Uhg'$ is the corresponding quantum
enveloping algebra (see \cite{Lu}, or \S \ref{ss:qWeyl}), then $(\Uhg',
\Uhsl{2}^i)$ is a $D$--algebra over the ring $\ICh$ of formal power series
in $\hbar$.

\subsection{Strict morphisms of $D$--algebras}

In example \ref{sss:Coxeter D}, set $k=\IC\fml$ and $q_i=\exp(\hbar k_i)$,
where $k_i\in\IC^*$ and $k_i=k_j$ if $s_i$ and $s_j$ are conjugate in $W$.
It is well--known in this case that if $W$ is finite and $\g=\g'$ is a complex,
simple Lie algebra, $\H(W)$ and $\Uhg$ are isomorphic to $\IC[W]\fml$ and
$\Ug\fml$ respectively. We will need to use such isomorphisms to compare
the corresponding structures of $D$--algebras. The following result shows
that the na\"{\i}ve notion of isomorphism between $D$--algebras is too
restrictive for this purpose however.

\begin{proposition}\label{pr:no graded}
Assume that $|\bfI|\geq 2$. Then,
\begin{enumerate}
\item There exists no algebra isomorphism $\Psi:\H(W)\rightarrow\IC[W]\fml$
such that $\Psi(\H[W_i])=\IC[W_i]\fml$ for any $i\in\bfI$. 
\item There exists no algebra isomorphism $\Psi:\Uhg\rightarrow\Ug
\fml$ equal to the identity mod $\hbar$ such that $\Psi(\Uhsl{2}^i)=
\Usl{2}^{\alpha_i}\fml$ for any $i\in\bfI$. 
\end{enumerate}
\end{proposition}
\proof
We may assume that $|\bfI|=2$.
(i) Let $\Psi:\H(W)\rightarrow\IC[W]\fml$ be an isomorphism such that $\Psi
(\H(W_i))=\IC[W_i]\fml$ for any $i\in\bfI$. Then, $\Psi(S_i)=x_i s_i+y_i$ for
some $x_i,y_i\in\ICh$ with $x_i\neq 0$. Since $S_i^2=(q_i-q_i^{-1})S_i+1$,
we get $y_i=(q_i-q_i^{-1})/2$. Equating the coefficients of $s_1s_2\cdots$
($m_{12}-2$ factors) in
$$\Psi(\underbrace{S_1S_2\cdots}_{m_{12}})=
(x_1s_1+y_1)(x_2s_2+y_2)\cdots$$
and in
$$\Psi(\underbrace{S_2S_1\cdots}_{m_{12}})=
(x_2s_2+y_2)(x_1s_1+y_1)\cdots$$
yields $2y_iy_j=y_iy_j$, a contradiction.\\

(ii) Recall that any algebra isomorphism $\Psi:\Uhg\rightarrow\Ug\fml$
equal to the identity mod $\hbar$ canonically identifies the centre of $\Uhg$
with $Z(\Ug)\fml$ \cite[p. 331]{Dr2}. Let $C_i=e_if_i+f_ie_i+1/2\thinspace
h_i^2\in U\sl{2}^i$ be the Casimir operator of $\sl{2}^i$ and $C_i^\hbar$ the
corresponding element of $Z(\Uhsl{2}^i)$. If $\Psi(\Uhsl{2}^i)=\Usl{2}^i\fml$
for any $i$, then $\Psi(C_i^\hbar)=C_i$. We will prove that such a $\Psi$
does not exist by showing that $C_1C_2=\Psi(C_1^\hbar C_2^\hbar)$
and $C_1^\hbar C_2^\hbar$ have different eigenvalues on the adjoint
representation $V$ of $\g$ and its quantum deformation $\V$ respectively.
Lusztig has given an explicit presentation of $\V$ \cite[\S 2.1]{Lu2}. Its zero
weight space $\V[0]$ is spanned by $t_1,t_2$ and
$$E_i\medspace t_j=-[a_{ji}]_j X_i$$
with $E_i\medspace X_i=0$.\footnote{for a definition of the $q$--numbers
$[n]_i$, see \S \ref{ss:qWeyl}.} Thus, for $i\neq j$, $t_i$ and $t_j-[a_{ji}]_j
/[2]_i\medspace t_i$ lie in the zero weight spaces of the quantum
deformations of the simple $\sl{2}^i$--modules $V^i_2$ and $V^i_0$ of
highest weights $2$ and $0$ respectively. Since $C_i$ acts as multiplication
by $m(m+2)/2$ on $V^i_m$, we get
$$C_i^\hbar t_i=4 t_i
\qquad\text{and}\qquad
C_i^\hbar(t_j-[a_{ji}]_j/[2]_i t_i)=0$$
so that, on $\V[0]$,
$$C_1^\hbar=\begin{pmatrix}4&4[a_{21}]_2/[2]_1\\0&0\end{pmatrix}
\qquad\text{and}\qquad
C_2^\hbar=\begin{pmatrix}0&0\\4[a_{12}]_1/[2]_2&4\end{pmatrix}$$
It follows that
$$C_1^\hbar C_2^\hbar =\begin{pmatrix}
16[a_{12}]_1[a_{21}]_2/[2]_1[2]_2&16[a_{21}]_2/[2]_1\\0&0\end{pmatrix}$$
has eigenvalues $0$ and $16[a_{12}]_2[a_{21}]_1/[2]_1[2]_2\in\ICh\setminus
\IC$ on $\V[0]$. Since the isomorphism class of a \fd $\Uhg$--module $\U$
is uniquely determined by that of the $\g$--module $\U/\hbar\U$, $\V$ is
isomorphic as $\Uhg$--module to $V\fml$, where $\Uhg$ acts on the latter
via $\Psi:\Uhg\rightarrow\Ug\fml$. The eigenvalues of $C_1^\hbar C_2^\hbar$
on $\V$ cannot therefore depend upon $\hbar$ \halmos

\subsection{Morphisms of $D$--algebras}\label{ss:D morphisms}

The following gives the correct notion of morphism of $D$--algebras.

\begin{definition}
A morphism of $D$--algebras $A,A'$ is a collection of algebra
homomorphisms $\Psi_\F:A\rightarrow A'$ labelled by the maximal
nested sets $\F$ on $D$ such that for any $\F$ and $D'\in\F$, $\Psi
_\F(A_{D'})\subseteq A'_{D'}$.
\end{definition}

\remark We will prove in theorems \ref{th:Ch iso} and \ref{th:transfer}
that $\H(W)$ and $\Uhg$ are isomorphic, as $D$--algebras to $\IC[W]
\fml$ and $\Ug\fml$ respectively.

\subsection{Completion with respect to \fd representations}

Let $\Veck$ be the category of finitely--generated, free $k$--modules and $
\Mod{A}$ that of \fd $A$--modules, that is the $A$--modules whose underlying
$k$--module lies in $\Veck$. Consider the forgetful functor
$$\bfF:\Mod{A}\rightarrow\Veck$$
By definition, the completion of $A$ with respect to its \fd representations is the
algebra $\wh{A}$ of endomorphisms of $\bfF$. Thus, an element of $\wh{A}$ is
a collection $\Theta=\{\Theta_V\}$, with $\Theta_V\in\End_{k}(V)$ for any $V\in
\Mod{A}$, such that for any $U,V\in\Mod{A}$ and $f\in\Hom_A(U,V)$
$$\Theta_V\circ f=f\circ\Theta_U$$

There is a natural homomorphism $A\rightarrow\wh{A}$ mapping $a\in A$
to the element $\Theta(a)$ which acts on a \fd representation $\rho:A\rightarrow
\End_k(V)$ as $\rho(a)$. The following is a straightforward consequence of
the above definitions
\begin{proposition}\label{pr:completion}\hfill
\begin{enumerate}
\item $A\rightarrow\wh{A}$ is a functor.
\item For any algebra $A$, the natural map $\wh{A}\rightarrow\wh{{\wh{A}}}$
is an isomorphism.
\item Let $\hbar$ be a formal variable and $A\fml$ the formal power series
in $\hbar$ with coefficients in $A$, regarded as an algebra over $k\fml$.
Then, the natural homomorphism $\wh{A}\fml\rightarrow\wh{A\fml}$ induced
by mapping $U\in\Mod{A}$ to $U\fml\in\Mod{A\fml}$ is an isomorphism
if \fd $A$--modules do not admit non--trivial deformations.
\item If $A$ is a bialgebra, $\wh{A}$ is a bialgebra with coproduct and counit
$$\Delta(\Theta)_{U\otimes V}=\Theta_U\otimes\Theta_V
\qquad\text{and}\qquad
\veps(\Theta)=\Theta_k$$
where $A$ acts on $U\otimes V$ via the coproduct $\Delta:A\rightarrow A^
{\otimes 2}$ and on $k$ via the counit $\veps:A\rightarrow k$.
\end{enumerate}
\end{proposition}

\subsection{Elementary pairs of \mnss}\label{ss:elementary pairs on D}

The terminology below corresponds, via the dictionary of \S \ref
{sss:graph dictionary}, to that of \S \ref{ss:elementary associators}
and \S \ref{ss:support}--\S \ref{ss:forgetful}.

\begin{definition}
An ordered pair $(\G,\F)$ of \mnss on $D$ is called elementary if $\G$
and $\F$ differ by an element. A sequence $\H_1,\ldots,\H_m$ of \mnss
on $D$ is called elementary if $|\H_{i+1}\setminus\H_i|=1$ for any
$i=1,\ldots,m-1$.
\end{definition}

By remark \rkref{rk:1-skeleton}, elementary pairs correspond to
oriented edges of the associahedron $\AD$ and elementary
sequences to edge--paths in $\AD$.

\begin{definition}\label{de:support on D}
The support $\supp(\G,\F)$ of an elementary pair of \mnss on $D$
is the unique unsaturated element of $\F\cap\G$. The central support
$\zsupp(\G,\F)$ of $(\G,\F)$ is the union of the maximal elements
of $\F\cap\G$ properly contained in $\supp(\G,\F)$. Thus,
$$\zsupp(\G,\F)=
\supp(\G,\F)\setminus\aalpha^{\supp(\G,\F)}_{\G\cap\F}$$
\end{definition}

\begin{definition}\label{de:DCP equiv}
Two elementary pairs $(\F,\G)$, $(\F',\G')$ of \mnss on $D$
are equivalent if
\begin{gather*}
\supp(\F,\G)=\supp(\F',\G'),\\
\alpha^{\supp(\F,\G)}_\F=\alpha^{\supp(\F',\G')}_{\F'}
\qquad\text{and}\qquad
\alpha^{\supp(\F,\G)}_\G=\alpha^{\supp(\F',\G')}_{\G'}
\end{gather*}
\end{definition}

\remark As in proposition \ref{pr:combi=>practical}, one readily shows
that for an elementary pair $(\F,\G)$,
$$\aalpha^{\supp(\F,\G)}_{\F\cap\G}=
\{\alpha^{\supp(\F,\G)}_\F,\alpha^{\supp(\F,\G)}_\G\}$$
In particular, two equivalent elementary pairs have the same central
support.

\subsection{Labelled diagrams and Artin braid groups}

\begin{definition}
A labelling of the diagram $D$ is the assignement of an integer $m_{ij}
\in\{2,3,\ldots,\infty\}$ to any pair $\alpha_i,\alpha_j$ of distinct vertices
of $D$ such that
$$m_{ij}=m_{ji}\qquad\text{and}\qquad m_{ij}=2$$
if, and only if $\alpha_i$ and $\alpha_j$ are orthogonal.
\end{definition}

Let $D$ be a labelled diagram.

\begin{definition}[Brieskorn--Saito \cite{BS}]
The Artin group $B_D$ is the group generated by elements $S_i$ labelled
by the vertices $\alpha_i$ of $D$ with relations
$$\underbrace{S_iS_j\cdots}_{m_{ij}}=
\underbrace{S_jS_i\cdots}_{m_{ij}}$$
for any $\alpha_i\neq\alpha_j$ such that $m_{ij}<\infty$.
\end{definition}

We shall also refer to $\BD$ as the braid group corresponding to $D$.

\subsection{Quasi--Coxeter algebras}\label{ss:qca 1}

Let $D$ be a labelled diagram.

\begin{definition}\label{de:qcoxeter}
A quasi--Coxeter algebra 
of type $D$ is a $D$--algebra $A$ endowed with the following additional
data
\begin{itemize}
\item {\bf Local monodromies : }%
for each $\alpha_i\in D$, an invertible element $$S^A_i\in\wh{A_i}$$
where $\wh{A}_i$ is the completion of $A_i$ with respect to its
finite--dimensional representations.
\item {\bf Elementary Associators : }%
for each elementary pair $(\G,\F)$ of \mnss on $D$, an invertible
element $$\Phi_{\G\F}\in A$$
\end{itemize}
satisfying the following axioms
\begin{itemize}
\item {\bf Orientation : }%
for any elementary pair $(\G,\F)$ of \mnss on $D$,
$$\Phi_{\F\G}=\Phi_{\G\F}^{-1}$$
\item {\bf Coherence : }%
for any pair of elementary sequences $\H_1,\ldots,\H_m$ and
$\K_1,\ldots,\K_\ell$ of \mnss on $D$ such that $\H_1=\K_1$
and $\H_m=\K_\ell$,
$$\Phi_{\H_m\H_{m-1}}\cdots\Phi_{\H_2\H_1}=
\Phi_{\K_\ell\K_{\ell-1}}\cdots\Phi_{\K_2\K_1}$$
\item {\bf Support : }%
for any elementary pair $(\G,\F)$ of \mnss on $D$,
$$\Phi_{\G\F}\in A_{\supp(\G,\F)}^{\zsupp(\G,\F)}$$
\item {\bf Forgetfulness : }%
for any equivalent elementary pairs $(\G,\F)$, $(\G',\F')$ of
\mnss on $D$,
$$\Phi_{\G\F}=\Phi_{\G'\F'}$$
\item {\bf Braid relations : }%
for any pair $\alpha_i,\alpha_j$ of distinct vertices of $D$ such that
$m_{ij}<\infty$, and elementary pair $(\G,\F)$ of \mnss on $D$ such
that $\alpha_i\in\F$ and $\alpha_j\in\G$,
$$\Ad(\Phi_{\G\F})(S^A_i)\cdot S^A_j\cdots=
S^A_j\cdot\Ad(\Phi_{\G\F})(S^A_i)\cdots$$
where the number of factors on each side is equal to $m_{ij}$.
\end{itemize}
\end{definition}

\remark The braid relations for $m_{ij}=2$ follow from the other axioms.
Indeed, let $\alpha_i\neq\alpha_j$ be two vertices of $D$ and $(\F_i,\F_j)$
an elementary pair of \mnss on $D$ such that $\alpha_i\in\F_i$ and $\alpha
_j\in\F_j$. If $\alpha_i\perp\alpha_j$, then either $\alpha_i\in\F_j$ or $\alpha
_j\in\F_i$ since, by maximality of $\F_i$ and $\F_j$, $\alpha_i\in\F_i\setminus
\F_j$ and $\alpha_j\in\F_j\setminus\F_i$ imply the incompatibility of $\alpha_i$
and $\alpha_j$. Assuming therefore that $\alpha_i\in\F_i\cap\F_j$ so that either
$\alpha_i\in\zsupp(\F_i,\F_j)$ or $\alpha_i\perp\supp(\F_i,\F_j)$. In either case,
the support axiom implies that $\Ad(\Phi_{\F_j\F_i})S^A_i=S^A_i$. The corresponding braid relation therefore reduces to $S^A_iS^A_j=S^A_jS^A_i$ and holds because
$S^A_i\in\wh{A_i},S^A_j\in\wh{A_j}$ and $[A_i,A_j]=0$.

\begin{definition}
A morphism of quasi--Coxeter algebras $A,A'$ of type $D$ is a morphism
$\{\Psi_\F\}$ of the underlying $D$--algebras such that
\begin{itemize}
\item for any $\alpha_i\in D$ and \mns $\F$ on $D$ with $\{\alpha_i\}\in\F$, $$\Psi_\F(S_i^A)=S_i^{A'}$$
\item for any elementary pair $(\G,\F)$ of \mnss on $D$,
$$\Psi_\G\circ\Ad(\Phi^A_{\G\F})=\Ad(\Phi^{A'}_{\G\F})\circ\Psi_\F$$
\end{itemize}
\end{definition}

subsection{The symmetric difference $\F\Delta\G$}

Let $\F,\G$ be \mnss on $D$. We characterise below the unsaturated
elements $B_1,\ldots,B_m$ of $\F\cap\G$ and the subsets of vertices
$\alpha^{B_i}_{\F\cap\G}$ in terms of the symmetric difference
$$\F\Delta\G=\F\setminus\G\cup\G\setminus\F$$
These results will be used in \S \ref{ss:D support}--\S
\ref{ss:D good sequences}. For any (maximal) \ns $\H$ on $D$ and
$B\in\H$, let
$$\H_B=\{B'\in\H|\thinspace B'\subseteq B\}$$
be the (maximal) \ns on $B$ induced by $\H$.

\begin{lemma}\label{le:unsat}
Let $C$ be an unsaturated element of $\F\cap\G$. Then,
\begin{equation}\label{eq:unsat}
C=\bigcup_{B\in\F_C\Delta\G_C}B
\end{equation}
\end{lemma}
\proof The \rhs is clearly contained in the left--hand side.
Let $\alpha_1=\alpha^C_\F$, $\alpha_2=\alpha^C_\G$ and note
that $\alpha_1\neq\alpha_2$ since $C$ is unsaturated. Set $C_1=\con
{C}{\alpha_1}{\alpha_2}\in\F\setminus\G$ and $C_2=\con{C}{\alpha_2}
{\alpha_1}\in\G\setminus\F$. Since $C_1$ and $C_2$ are not compatible,
$C_1\cup C_2$ is connected and properly contains $C_1,C_2$. Let
$\ol{C}$ be the connected component of the \rhs of \eqref{eq:unsat}
containing $C_1,C_2$. We claim that $\ol{C}$ is compatible with any
$B\in\F$. This is clear if $B\perp C$, if $B\supseteq C$ or if $B\subset
C$ and $C\notin\G$. If, on the other hand, $B\in\G$ then $B$ is compatible
with any element in $\F\Delta\G$ and therefore with $\ol{C}$. By maximality
of $\F$, $\ol{C}\in\F$ whence $\ol{C}=C$ since $C_1\subsetneq\ol{C}
\subset C$ \halmos

\begin{proposition}\label{pr:max unsat}
The connected components of
$$\bigcup_{B\in\F\Delta\G}B$$
are the maximal unsaturated elements of $\F\cap\G$.
\end{proposition}
\proof We claim that any $B\in\F\Delta\G$ is contained in an unsaturated
element of $\F\cap\G$. Let $\ol{B}$ be the minimal element in $\F\cap\G$
containing $B$. If $B\in\F\setminus\G$ (resp. $\G\setminus\F$), $B$ is
contained in a connected components $\ol{B}'$ of $\ol{B}\setminus\alpha
^{\ol{B}}_\F$ (resp. $\ol{B}\setminus\alpha^{\ol{B}}_\G$). If $\ol{B}$ were
saturated as an element of $\F\cap\G$, then $\alpha^{\ol{B}}_\F=\alpha^
{\ol{B}}_{\F\cap\G}=\alpha^{\ol{B}}_\G$ and $\ol{B}'\in\F\cap\G$, in
contradiction with the minimality of $\ol{B}$. It follows that
$$\bigcup_{B\in\F\Delta\G}B=
\bigcup_{\substack{C\in\F\cap\G,\\\text{$C$ unsaturated}}}
\bigcup_{\substack{B\in\F\Delta\G,\\ B\subseteq C}}B=
\bigcup_{\substack{C\in\F\cap\G,\\\text{$C$ unsaturated}}}
C$$
where the second equality holds by lemma \ref{le:unsat} \halmos

\begin{proposition}\label{pr:prop max}
Let $C$ be an unsaturated element of $\F\cap\G$. Then, the connected
components of
\begin{equation}\label{eq:prop max}
\bigcup_{\substack{B\subsetneq C,\\
\text{$B$ compatible with any $B'\in\F\Delta\G$}}} B
\end{equation}
are the maximal elements of $\F\cap\G$ properly contained in $C$.
\end{proposition}
\proof Any $B\in\F\cap\G$ with $B\subsetneq C$ is clearly contained in
\eqref{eq:prop max}. Let now $B\subsetneq C$ be connected and compatible
with any $B'\in\F\Delta\G$. Assume that $B$ is not contained in one of the
maximal elements $C_1,\ldots,C_m$ of $\F\cap\G$ properly contained
in $C$, set $\ul{C}=C_1\odots{\sqcup}C_m$ and let $\ol{B}$ be the image
of $B$ in $\ol{C}=C/\ul{C}$. Let $\ol{\F}_C$ and $\ol
{\G}_C$ be the \mnss on $\ol{C}$ induced by $\F_C,\G_C$ respectively.
By (ii) of lemma \ref{le:quotient compatibility}, $\ol{B}$ is compatible with
$$\ol{\F}_C\setminus\ol{\G}_C=\ol{\F}_C\setminus\{\ol{C}\}
\aand\ol{\G}_C\setminus\ol{\F}_C=\ol{\G}_C\setminus\{\ol{C}\}$$ and
therefore with $\ol{\F}_C,\ol{\G}_C$. By maximality of $\ol{\F}_C$ and
$\ol{\G}_C$, $\ol{B}$ lies in $\ol{\F}_C\cap\ol{\G}_C=\{\ol{C}\}$.
We claim that this is a contradiction. It suffices for this to prove the
existence of a $B'\in\F_C\Delta\G_C$ such that $B\subseteq B'$, for
then $\ol{C}=\ol{B}\subset\ol{B}'\subset\ol{C}$, whence, by \eqref{eq:bar tilde},
$B'=\wt{\ol{B}'}=C\in\F\cap\G$. Assume for this purpose that $B\nsubseteq B'$
for any $B'\in\F_C\Delta\G_C$ and set
$$C_\perp=\bigcup_{\substack{B'\in\F_C\Delta\G_C,\\ B'\perp B}}B'
\qquad\text{and}\qquad
C_{\subset}=\bigcup_{\substack{B'\in\F_C\Delta\G_C,\\B'\subset B}}B'$$
By lemma \ref{le:unsat}, $C=C_\perp\cup C_\subset$. Since $C$ is
connected and $C_\perp$ and $C_\subset$ are orthogonal, one has
$C_\perp=\emptyset$ or $C_\subset=\emptyset$ and therefore $C
\subset B$ or $C\perp B$ respectively, both of which contradict $B
\subsetneq C$ \halmos\\

The following is a direct consequence of propositions \ref{pr:max unsat}
and \ref{pr:prop max}.

\begin{corollary}\label{co:unsat Delta}
The unsaturated elements $B_1,\ldots,B_m$ of $\F\cap\G$ and subsets
of vertices $\aalpha^{B_i}_{\F\cap\G}\subset B_i$ only depend upon the
symmetric difference $\F\Delta\G$.
\end{corollary}

\subsection{Support and central support of a pair of maximal nested sets}
\label{ss:D support}

We extend below the notions of support and central support to a general
pair $(\F,\G)$ of \mnss on $D$.

\begin{definition}\label{de:gen supp}
The support of $(\F,\G)$ is the union
$$\supp(\F,\G)=\bigcup_{B\in\F\Delta\G}B$$
\end{definition}

By proposition \ref{pr:max unsat}, $\supp(\F,\G)$ is the union of the
maximal unsaturated elements of $\F\cap\G$. In particular, definition
\ref{de:gen supp} is consistent with definition \ref{de:support on D}
when $\F$ and $\G$ differ by an element.\\

For any collection $\C$ of connected subdiagrams of $D$, set now
$$\kappa(\C)=
\{B\subseteq D|\medspace B\perp C\thickspace\text{or}\thickspace
B\subseteq C\thickspace\text{for any $C\in\C$}\}$$
One readily checks that if $B_1,B_2\in\kappa(\C)$ are incompatible,
then $B_1\cup B_2\in\kappa(\C)$. In particular, the maximal elements
$B_1,\ldots,B_m$ of $\kappa(\C)$ are pairwise orthogonal and
$$\bigcup_{B\in\kappa(\C)}B=B_1\odots{\sqcup}B_m$$
Note that if $B\in\kappa(\F\Delta\G)$, then either $B\subset\supp(\F,\G)$
or $B\perp\supp(\F,\G)$.

\begin{definition}\label{de:gen zsupp}
The central support $\zsupp(\F,\G)$ is the union $\bigcup B$, where $B$
ranges over the elements of $\kappa(\F\Delta\G)$ contained in $\supp(\F,
\G)$.
\end{definition}

The following result shows that definition \ref{de:gen zsupp} is consistent
with definition \ref{de:support on D}.

\begin{proposition}\label{pr:elem support}
Assume that $\G$ and $\F$ differ by an element and let $B=\supp(\F,\G)$
be the unique unsaturated element of $\F\cap\G$. Then,
$$\bigcup_{\substack{C\in\kappa(\F\Delta\G),\\C\subseteq B}}
C=B\setminus\aalpha^B_{\F\cap\G}$$
\end{proposition}
\proof Set $\alpha_1=\alpha^B_\F$ and $\alpha_2=\alpha^B_\G$, so that
$\aalpha^B_{\F\cap\G}=\{\alpha_1,\alpha_2\}$, $B_1=\con{B}{\alpha_1}
{\alpha_2}\in\F\setminus\G$ and $B_2=\con{B}{\alpha_2}{\alpha_1}\in\G
\setminus\F$. Since any $C\in\kappa(\{B_1,B_2\})$ does not contain $
\alpha_1$ and $\alpha_2$, the \lhs is contained in the right--hand side.
The opposite inclusion is easy to check  \halmos\\

\remark\rklabel{rk:disjoint}
Note that $\kappa(\F\Delta\G)$ and $\F\Delta\G$ are disjoint. Indeed,
any $B\in\kappa(\F\Delta\G)\cap\F\Delta\G$ is compatible with $\G=
\G\setminus\F\cup(\F\cap\G)$ and $\F=\F\setminus\G\cup(\F\cap
\G)$ and therefore lies in $\F\cap\G$ by maximality of $\F$ and $\G$,
a contradiction. Thus, if $B\in\kappa(\F\Delta\G)$ and $C\in\F\Delta\G$,
then
$$B\perp C\oor B\subsetneq C$$

\subsection{Equivalence of pairs of \mnss}
\label{ss:D equivalence}

We shall need to extend the notion of equivalence to general pairs of \mnss
on $D$. We begin by giving an alternative characterisation of the equivalence
of two elementary such pairs.

\begin{proposition}\label{pr:new equivalence}
Two elementary pairs $(\F,\G)$ and $(\F',\G')$ of \mnss on $D$ are equivalent
if, and only if
$$\F\setminus\G=\F'\setminus\G'
\aand
\G\setminus\F=\G'\setminus\F'$$
\end{proposition}

Proposition \ref{pr:new equivalence} is an immediate corollary of the following.

\begin{proposition}\label{pr:combi<=>practical}
Let $(\F,\G)$ be an elementary pair of maximal nested sets on $D$. Let
$B=\supp(\F,\G)$ be the unique unsaturated element of $\F\cap\G$ and
$B_1,B_2$ the unique elements in $\F\setminus\G$ and $\G\setminus
\F$ respectively. Then,
\begin{enumerate}
\item $\alpha^B_\F$ and $\alpha^B_\G$ are distinct and
$\aalpha^B_{\F\cap\G}=\{\alpha^B_\F,\alpha^B_\G\}$.
\item $B_1=\con{B}{\alpha^B_\F}{\alpha^B_\G}$ and
$B_2=\con{B}{\alpha^B_\G}{\alpha^B_\F}$.
\item $\alpha^{B_1}_\F=\alpha^B_\G$ and
$\alpha^{B_2}_\G=\alpha^B_\F$.
\item $B_1,B_2$ are not compatible and $B_1\cup B_2=B$.
\item $\alpha^B_\F$, $\alpha^B_\G$ are uniquely determined by (ii).
\end{enumerate}
\end{proposition}
\proof
(i)--(iii) are proved exactly as in proposition \ref{pr:combi=>practical}.
(iv) the incompatibility of $B_1,B_2$ is a direct consequence of (ii).
The fact that $B=B_1\cup B_2$ follows by lemma \ref{le:unsat}.
(v) Let $\alpha_1$ a vertex of $B$ such that $B_1$ is a connected
component of $B\setminus\alpha_1$. Then, $\alpha_1\nperp B_1$
so that, if $\alpha_1\neq\alpha^B_\F$, $\alpha_1$ and $B_1$ lie in
the same connected component of $B\setminus\alpha^{B}_\F$. By
(ii), this implies that $\alpha_1\in B_1$, a contradiction. Similarly, $
\alpha^B_\G$ is the unique vertex of $B$ such that $B_2$ is a
connected component of $B\setminus\alpha^B_\G$ \halmos\\

Proposition \ref{pr:new equivalence} ensures that the following is consistent
with definition \ref{de:DCP equiv}.

\begin{definition}\label{de:equiv pairs}
Two ordered pairs $(\F,\G)$, $(\F',\G')$ of \mnss on $D$ are equivalent if
$$\F\setminus\G=\F'\setminus\G'
\aand
\G\setminus\F=\G'\setminus\F'$$
\end{definition}

Note that the equivalence of $(\F,\G)$ and $(\F',\G')$ implies that
$$\supp(\F,\G)=\supp(\F',\G')
\aand
\zsupp(\F,\G)=\zsupp(\F',\G')$$

\subsection{Existence of good elementary sequences}
\label{ss:D good sequences}

\begin{proposition}\label{pr:elem seq}\hfill
\begin{enumerate}
\item For any pair $(\F,\G)$ of \mnss on $D$, there exists an elementary
sequence $$\F=\H_1,\H_2,\ldots,\H_m=\G$$ such that, for any $i=1,\ldots,
m-1$,
\begin{gather*}
\H_i\cap\H_{i+1}\supseteq\F\cap\G,\\
\supp(\H_i,\H_{i+1})\subseteq\supp(\F,\G)
\end{gather*}
and, for any component $B$ of $\zsupp(\F,\G)$, either
$$B\perp\supp(\H_i,\H_{i+1})
\quad\text{or}\quad
B\subseteq\zsupp(\H_i,\H_{i+1})$$
\item If $(\G,\F)$ and $(\G',\F')$ are equivalent pairs of \mnss
on $D$, the corresponding elementary sequences
$$\F=\H_1,\H_2,\ldots,\H_m=\G
\quad\text{and}\quad
\F'=\H'_1,\H'_2,\ldots,\H'_\ell=\G'$$
may be chosen such that $\ell=m$ and such that,  for any $i=1\ldots m-1$,
$(\H_i,\H_{i+1})$ is equivalent to $(\H_i',\H_{i+1}')$.
\end{enumerate}
\end{proposition}
\proof
(i) By the connectedness of the face of the associahedron $\A_D$ corresponding
to $\K=\F\cap\G$, there exists an elementary sequence $\F=\H_1,\ldots,\H_m$
$=\G$
such that $\K\subset\H_i$ for any $i$. Let $B_1,\ldots,B_p$ be the unsaturated
elements of $\K$ and set  $\aalpha_j=\aalpha^{B_j}_\K$. By proposition
\ref{pr:ND^H} and lemma \ref{le:explicit lift}, each $\H_i$ is the union of
$\K$ and of a compatible family of diagrams of the form
$D_{\bbeta_j}=\con{B_j}{(\aalpha_j\setminus\bbeta_j)}{\bbeta_j}$ for some $
1\leq j\leq p$ and $\emptyset\neq\bbeta_j\subsetneq\aalpha_j$. For any $i=1
,\ldots,p-1$, set
$$D_{\bbeta\ji}=\H_i\setminus\H_{i+1}
\qquad\text{and}\qquad
D_{\ggamma_{k_i}}=\H_{i+1}\setminus\H_i$$
Since $D_{\bbeta\ji}$ and $D_{\ggamma_{k_i}}$, one has $j_i=k_i$ whence,
$$\supp(\H_i,\H_{i+1})=
D_{\bbeta\ji}\cup D_{\ggamma_{k_i}}\subseteq
B\ji\subseteq\supp(\F,\G)$$

Let now $B$ be a component of $\zsupp(\F,\G)$. We shall need the following
\begin{lemma}
For any $1\leq j\leq p$, one has $B\perp B\ji$ or $B\subseteq B\ji\setminus
\aalpha\ji$.
\end{lemma}
\proof By lemma \ref{le:unsat}, $B\ji=\bigcup_{B'\in\F_{B\ji}\Delta\G_{B\ji}}
B'$. Since $B$ is compatible with any such $B'$ and does not contain it by
remark \rkref{rk:disjoint}, either $B\perp B\ji$ or $B\subsetneq B\ji$. In
the latter case, $B\subseteq\bigcup B'$ where $B'$ now ranges over
the proper connected subdiagrams of $B\ji$ which are compatible with
any element of $\F_{B\ji}\Delta\G_{B\ji}$, whence $B\subseteq B\ji
\setminus\aalpha\ji$ by proposition \ref{pr:prop max} \halmos\\

If $B\perp B\ji$, then $B\perp\supp(\H_i,\H_{i+1})$ as required. If, on the
other hand $B\subseteq B\ji\setminus\aalpha\ji$, then $B$ is compatible
with $D_{\bbeta\ji}$ and $D_{\ggamma\ji}$ and contains neither since $
\aalpha\ji\cap D_{\bbeta\ji},\aalpha\ji\cap D_{\ggamma\ji}\neq\emptyset$.
Thus, either $B\perp D_{\bbeta\ji},D_{\ggamma\ji}$, in which case $B\perp
D_{\bbeta\ji}\cup D_{\ggamma\ji}=\supp(\H_i,\H_{i+1})$, or $B\subseteq
\bigcup B'$, where the union ranges over the connected subdiagrams of
$\supp(\H_i,\H_{i+1})$ compatible with, but not containing either of $D_
{\bbeta\ji},D_{\ggamma\ji}$ and therefore $B\subseteq\zsupp(\H_i,\H_
{i+1})$.\\

(ii) Let $\H_i=\K\cup\{D_{\bbeta_j}\}_{j\in J_i}$ be the elementary sequence
obtained in (i). By corollary \ref{co:unsat Delta}, $\K'=\F'\cap\G'$ and $\K$
have the same unsaturated elements $B_1,\ldots,B_p$ and $\aalpha^{B_i}
_{\K'}=\aalpha^{B_i}_\K$ for any $i=1,\ldots,p$. It follows from this, proposition
\ref{pr:ND^H} and lemma \ref{le:explicit lift} that $\H'_i=\K'\cup\{D_{\bbeta_j}
\}_{j\in J_i}$ is a \mns on $D$ \halmos

\subsection{General associators}\label{ss:qca 2}

Let $A$ be a quasi--Coxeter algebra of type $D$. By the connectedness
of the associahedron $\AD$, there exists, for any pair $\G,\F$ of \mnss
on $D$, an elementary sequence $\H_1,\ldots,\H_m$ such that $\H_1=\F$
and $\H_m=\G$. Set
$$\Phi_{\G\F}=\Phi_{\H_m\H_{m-1}}\cdots\Phi_{\H_2\H_1}$$
The coherence axiom implies that this definition is independent of the
choice of the elementary sequence and that $\Phi_{\G\F}$ is the elementary
associator corresponding to $(\G,\F)$ if $\F$ and $\G$ differ by an element.
The following result summarises the main properties of the general associators
$\Phi_{\G\F}$ and gives an equivalent characterisation of \qcas in terms
of them.

\begin{theorem}\label{th:gen ass}
The associators $\Phi_{\G\F}$ satisfy the following properties
\begin{itemize}
\item {\bf Orientation : }%
for any pair $(\G,\F)$ of \mnss on $D$,
$$\Phi_{\F\G}=\Phi_{\G\F}^{-1}$$
\item {\bf Transitivity : }%
for any triple $\H,\G,\F$ of \mnss on $D$,
$$\Phi_{\H\F}=\Phi_{\H\G}\cdot\Phi_{\G\F}$$
\item {\bf Forgetfulness : }%
for any equivalent pairs $(\G,\F)$ and $(\G',\F')$ of \mnss on $D$,
$$\Phi_{\G\F}=\Phi_{\G'\F'}$$
\item {\bf Support : }%
for any pair $(\G,\F)$ of \mnss on $D$,
$$\Phi_{\G\F}\in A_{\supp(\F,\G)}^{\zsupp(\F,\G)}$$
\item {\bf Braid relations : }%
for any pair $\alpha_i,\alpha_j$ of distinct vertices of $D$ such
that $m_{ij}<\infty$, and pair $(\G,\F)$ of \mnss on $D$ such
that $\alpha_i\in\F$ and $\alpha_j\in\G$,
$$\Ad(\Phi_{\G\F})(S^A_i)\cdot S^A_j\cdots=
S^A_j\cdot\Ad(\Phi_{\G\F})(S^A_i)\cdots$$
where the number of factors on each side is equal to $m_{ij}$.
\end{itemize}
Conversely, if $A$ is a $D$--algebra endowed with invertible elements
$S^A_i\in\wh{A}_i$ for any $\alpha_i\in D$ and $\Phi_{\G\F}\in A$ for
any pair $(\G,\F)$ of \mnss on $D$ which satisfy the above properties,
then the $S^A_i$ and associators $\Phi_{\G\F}$ corresponding to
elementary pairs give $A$ the structure of a quasi--Coxeter algebra
of type $D$.
\end{theorem}
\proof Orientation and transitivity follow at once from the orientation
and coherence axioms satisfied by the elementary associators. If
$(\G,\F)$ and $(\G',\F')$ are two equivalent pairs of \mnss on $D$
and
$$\F=\H_1,\ldots,\H_m=\G\qquad\text{and}
\qquad\F'=\H_1',\ldots,\H_m'=\G'$$
are two elementary sequences such that $(\H_i,\H_{i+1})$ and
$(\H_i',\H_{i+1}')$ are equivalent for any $i=1,\ldots,m-1$ as in
proposition \ref{pr:elem seq} (ii), then, by the forgetfulness axiom
satisfied by elementary associators
$$\Phi_{\G\F}=
\Phi_{\H_m\H_{m-1}}\cdots\Phi_{\H_2\H_1}=
\Phi_{\H'_m\H'_{m-1}}\cdots\Phi_{\H'_2\H'_1}=
\Phi_{\G'\F'}$$
Similarly, if $\F=\H_1,\ldots,\H_m=\G$ is an elementary sequence of
\mnss on $D$ as in proposition \ref{pr:elem seq} (i), then, by the support
properties of elementary associators,
$$\Phi_{\H_{i+1}\H_{i}}\in
A_{\supp(\H_{i+1},\H_i)}^{\zsupp(\H_{i+1},\H_i)}\subseteq
A_{\supp(\G,\F)}^{\zsupp(\G,\F)}$$
so that
$$\Phi_{\G\F}=
\Phi_{\H_m\H_{m_1}}\cdots\Phi_{\H_2\H_1}\in
A_{\supp(\G,\F)}^{\zsupp(\G,\F)}$$
Let now $\alpha_i\neq\alpha_j$ be such that $m_{ij}<\infty$, $\G,\F$
such that $\alpha_i\in\F$, $\alpha_j\in\G$ and choose an elementary
pair $(\F',\G')$ such that $\alpha_i\in\F'$ and $\alpha_j\in\G'$. Since
$\alpha_i\in\F\cap\F'$ and $\alpha_i$ does not contain any elements
of $\F_i\Delta\F'_i$, $\alpha_i$ is either contained in $\zsupp(\F,\F')$
or perpendicular to $\supp(\F,\F')$. In either case, 
$$\Ad(\Phi_{\F'\F})\medspace(S^A_i)=S^A_i
\quad\text{and, similarly,}\quad
\Ad(\Phi_{\G'\G})\medspace(S^A_j)=S^A_j$$
Thus
\begin{equation*}
\begin{split}
\Ad(\Phi_{\G\F})(S^A_i)\cdot S^A_j\cdots
&=
\Ad(\Phi_{\G\G'}\cdot\Phi_{\G'\F'}\cdot\Phi_{\F'\F})(S^A_i)\cdot S^A_j\cdots\\
&=
\Ad(\Phi_{\G\G'})\left(\Ad(\Phi_{\G'\F'})(S^A_i)\cdot S^A_j\cdots\right)
\end{split}
\end{equation*}
and
$$S^A_j\cdot\Ad(\Phi_{\G\F})(S^A_i)\cdots=
\Ad(\Phi_{\G\G'})\left(S^A_j\cdot\Ad(\Phi_{\G'\F'})(S^A_i)\cdots\right)$$
so the two are equal because of the braid relations satisfied by
the elementary associators. The converse implication is clear
\halmos\\

\remark Unlike the case of elementary associators, the braid relations
involving general associators do not follow from the other axioms when
$m_{ij}=2$. For example, if $D$ is the Dynkin diagram of type $\sfA_3$,
which we identify with the interval $[1,3]$, $\F=\{[1,1],[1,2],[1,3]\}$ and
$\G=\{[3,3],[2,3],[1,3]\}$, then $\supp(\F,\G)=[1,3]$ and $\zsupp(\F,\G)
=\emptyset$. Thus, $\Phi_{\G\F}$ does not centralise $A_1$ or $A_3$
and neither of the relations
$$\Ad(\Phi_{\G\F})(S_1^A)\cdot S_3^A=S_3^A\cdot \Ad(\Phi_{\G\F})(S^A_1)
\thickspace\thickspace\text{or}\thickspace
\Ad(\Phi_{\F\G})(S_3^A)\cdot S_1^A=S_1^A\cdot \Ad(\Phi_{\F\G})(S_1^A)$$
can be deduced from the fact that $S_1^A\cdot S_3^A=S_3^A\cdot S_1^A$.\\

We record for later use the following consequence of proposition \ref
{pr:elem seq} and of the definition of general associators

\begin{proposition}\label{pr:stability}
If $A$ is a \qc algebra of type $D$ then, for any pair $(\G,\F)$ of \mnss
on $D$ and $B\in\G\cap\F$,
$$\Ad(\Phi_{\G\F})(A_B)\subseteq A_B$$
\end{proposition}
\proof Assume first that $\G$ and $\F$ differ by an element and let $C=
\supp(\F,\G)$ be the unique unsaturated element of $\G\cap\F$, so that
$\Phi_{\G\F}\in A_C^{C\setminus\aalpha^C_{\G\cap\F}}$. If $B\perp C$
or $B\subseteq C$, in which case $B\subseteq C\setminus\aalpha^C_{
\G\cap\F}$, then $[\Phi_{\G\F},A_B]=0$ and the result follows. If, on the
other hand $B\supseteq C$, then $\Phi_{\G\F}\in A_C\subseteq A_B$
and the result follows again. If $(\G,\F)$ is a general pair of \mnss on
$D$, proposition \ref{pr:elem seq} implies the existence of an elementary
sequence $\F=\H_0,\cdots,\H_m=\G$ such that, for any $i=1,\dots,m-1$,
$B\in\H_i\cap\H_{i+1}$. The result now follows from our previous analysis
\halmos

\subsection{Braid group representations}
\label{ss:qca braid}

By mimicking the monodromy computations of \S \ref{se:DCP} (in particular
\S \ref{sss:DCP braid relations}), we show below that a quasi--Coxeter algebra
$A$ of type $D$ defines representations of the braid group $B_D$ on any
\fd $A$--module, with isomorphic quasi--Coxeter algebras defining equivalent
representations of $\BD$. It is worth keeping in mind that, just as the action
of Artin's braid group $B_n$ on the $n$--fold tensor product $V^{\otimes n}$
of an object in a braided tensor category depends upon the choice of a
complete bracketing on $V^{\otimes n}$, the procedure described below
yields not one, but a {\it family} of canonically equivalent representations
$$\pi_\F:B_D\longrightarrow\wh{A}$$
labelled by the \mnss $\F$ on $D$.\\

Let $\F$ be a \mns on $D$. For any $\alpha_i\in D$, choose a \mns
$\G_i$ such that $\alpha_i\in\G_i$ and set
$$\pi_\F(S_i)=\Phi_{\F\G_i}\cdot S^A_i\cdot\Phi_{\G_i\F}$$
\begin{theorem}\label{th:monodromy maps}\hfill
\begin{enumerate}
\item The above assignement is independent of the choice of $\G_i$
and extends to a homomorphism $\pi_\F:\BD\rightarrow\wh{A}$.
\item If $\alpha_i\in\F$, then
$$\pi_\F(S_i)=S^A_i$$
\item  For any $D'\in\F$,
$$\pi_\F(B_{D'})\subset \wh{A_{D'}}$$
\item If $\G$ is another \mns on $D$ then, for any $b\in\BD$,
$$\pi_\G(b)=\Phi_{\G\F}\cdot \pi_\F(b)\cdot \Phi_{\F\G}$$
so that $\pi_\F$ and $\pi_\G$ are canonically equivalent.
\item If $\G$ is another \mns and $D'\in\G\cap\F$ is such
that the induced \mnss $\F_{D'},\G_{D'}$ on $D'$ coincide,
the restrictions of $\pi_\F,\pi_\G$ to $B_{D'}$ are equal.
\item If $\{\Psi_\F\}_\F:A\rightarrow A'$ is a morphism of
quasi--Coxeter algebras, then for any \mns $\F$ and $b
\in\BD$,
$$\Psi_\F(\pi^A_\F(b))=\pi^{A'}_\F(b)$$
In particular, isomorphic quasi--Coxeter algebras yield
equivalent representations of $\BD$.
\end{enumerate}
\end{theorem}
\proof
(i) If $\G_i'$ is such that $\alpha_i\in\G_i'$, then either $\alpha_i
\perp\supp(\G_i,\G_i')$ or $\alpha_i\subseteq\zsupp(\G_i,\G_i')$
so that $\Phi_{\G_i'\G_i}$ centralises $A_i$ by theorem \ref
{th:gen ass}. Thus,
$$\Phi_{\F\G_i}\cdot S^A_i\cdot\Phi_{\G_i\F}=
\Phi_{\F\G_i'}\cdot\Phi_{\G_i'\G_i}\cdot S^A_i
\cdot\Phi_{\G_i\G_i'}\cdot\Phi_{\G_i'\F}=
\Phi_{\F\G_i'}\cdot S^A_i\cdot\Phi_{\G_i'\F}$$
Let now $\alpha_i\neq\alpha_j$ be such that $m_{ij}<\infty$. Then,
\begin{equation*}
\begin{split}
\pi_\F(S_i)\pi_\F(S_j)\cdots
&=
\Phi_{\F\G_i}\cdot S^A_i\cdot\Phi_{\G_i\F}\cdot
\Phi_{\F\G_j}\cdot S^A_j\cdot\Phi_{\G_j\F}\cdots\\
&=
\Phi_{\F\G_i}\cdot\left(
S^A_i\cdot\Phi_{\G_i\G_j}\cdot S^A_j\cdot\Phi_{\G_j\G_i}\cdots
\right)\cdot\Phi_{\G_i\F}
\end{split}
\end{equation*}
and
\begin{equation*}
\begin{split}
\pi_\F(S_j)\pi_\F(S_i)\cdots
&=
\Phi_{\F\G_j}\cdot S^A_j\cdot\Phi_{\G_j\F}\cdot
\Phi_{\F\G_i}\cdot S^A_i\cdot\Phi_{\G_i\F}\cdots\\
&=
\Phi_{\F\G_i}\cdot\left(
\Phi_{\G_i\G_j}\cdot S^A_j\cdot\Phi_{\G_j\G_i}\cdot S^A_i\cdots
\right)\cdot\Phi_{\G_i\F}
\end{split}
\end{equation*}
so that the two coincide by theorem \ref{th:gen ass}.\\
(ii) follows by choosing $\G_i=\F$.\\
(iii) For any $\alpha_i\in D'$, let $\G_{D'}$ be a \mns on $D'$ such that
$\alpha_i\in\G_{D'}$ and set
$$\G_i=(\F\setminus\F_{D'})\cup\G_{D'}$$
Then, $\supp(\F,\G_i)\subseteq D'$ so that, by theorem \ref{th:gen ass}
$$\pi_\F(S_i)=\Phi_{\F\G_i}\cdot S^A_i\cdot\Phi_{\G_i\F}\in\wh{A_{D'}}$$
(iv) we have
$$\pi_\G(S_i)=
\Phi_{\G\G_i}\cdot S^A_i\cdot\Phi_{\G_i\G}=
\Phi_{\G\F}\cdot\Phi_{\F\G_i}\cdot S^A_i\cdot\Phi_{\G_i\F}\cdot\Phi_{\F\G}=
\Phi_{\G\F}\cdot\pi_\F(S_i)\cdot\Phi_{\F\G}$$
(v) By assumption, either $D'\perp\supp(\F,\G)$ or $D'\subseteq\zsupp(\F,\G)$.
It follows that $\Phi_{\G\F}$ centralises $A_{D'}$ whence, by (iii)
$$\pi_\G(S_i)=\Phi_{\G\F}\cdot\pi_\F(S_i)\cdot\Phi_{\F\G}=\pi_\F(S_i)$$
(vi) By definition,
\begin{equation*}
\begin{split}
\pi^{A'}_\F(S_i)
&=
\Phi^{A'}_{\F\G_i}\cdot S_i^{A'}\cdot\Phi^{A'}_{\G_i\F}\\
&=
\Phi^{A'}_{\F\G_i}\cdot\Psi_{\G_i}(S_i^A)\cdot\Phi^{A'}_{\G_i\F}\\
&=
\Psi_\F(\Phi^A_{\F\G_i}\cdot S_i^A\cdot\Phi^A_{\G_i\F})\\
&=
\Psi_\F(\pi^A_\F(S_i))
\end{split}
\end{equation*}
\halmos\\

\remark\rklabel{rk:univeral}
The group algebra $A'=k[B_D]$ of $B_D$ may be regarded as a \qca of type
$D$ by setting
$$A'_{D'}=k[B_{D'}],\qquad
S_i^{A'}=S_i\qquad\text{and}\qquad
\Phi^{A'}_{\G\F}=1$$
Theorem \ref{th:monodromy maps} may then be rephrased as saying that the
collection $\{\pi_\F\}$ is a morphism of \qcas $k[B_D]\rightarrow\wh{A}$ which
is functorial in $A$. 

\subsection{Generalised pentagon relations}

The coherence relations satisfied by the elementary associators of a 
\qca are convenient for most applications but somewhat redundant.
In this subsection, we use the simple connectedness of the associahedron
$\AD$ to reduce them to a smaller number of identities labelled by
the pentagonal and hexagonal faces of $\AD$.\\

Let $A$ be a $D$--algebra endowed with invertible elements
$\Phi_{\G\F}$ labelled by elementary pairs of \mnss on $D$.
For any $2$--face $\H$ of $\AD$, orientation $\oveps$ of $\H$
and \mns $\F_0$ on the boundary of $\H$, let $\F_0,\F_1,\ldots,
\F_{k-1},\F_k=\F_0$ be the vertices of $\H$ listed in their order
of appearance along $\partial\H$ when the latter is endowed
with the orientation $\oveps$. Set
$$\mu(\H;\F_0,\oveps)=
\Phi_{\F_0\F_{k-1}}\cdots\Phi_{\F_1\F_0}\in A$$
The following is immediate
\begin{lemma}\label{le:choices}\hfill
\begin{enumerate}
\item For any $i=0,\ldots,k-1$,
$$\mu(\H;\F_i,\oveps)=
\Ad(\Phi_{\F_i\F_{i-1}}\cdots\Phi_{\F_1\F_0})\medspace
\mu(\H;\F_0,\oveps)$$
\item If $\Phi_{\F\G}=\Phi_{\G\F}^{-1}$ for any elementary pair
$(\G,\F)$ of \mnss on $D$, then
$$\mu(\H;\F_0,-\oveps)=
\mu(\H;\F_0,\oveps)^{-1}$$
where $-\oveps$ is the opposite orientation to $\oveps$.
\end{enumerate}
\end{lemma}

By lemma \ref{le:choices},  the identity $\mu(\H;\F_0,\oveps)=1$, regarded
as an identity in the variables $\Phi_{\G\F}$ does not depend upon the
choice of $\oveps$ and $\F_0$ provided the $\Phi_{\G\F}$ satisfy the
orientation axiom of definition \ref{de:qcoxeter}. We shall henceforth
denote this identity by $\mu(\H)=1$.

\begin{proposition}\label{pr:zero on squares}
Assume that the elements $\Phi_{\G\F}$ satisfy the orientation,
forgetfulness and support axioms of definition \ref{de:qcoxeter}.
Then, for any square 2--face $\H$ of the associahedron $\AD$, 
$\mu(\H)=1$.
\end{proposition}
\proof Let $D_1,D_2$ be the unsaturated elements of $\H$ and
set, for $i,j,k\in\{1,2\}$,
$$\aalpha_i=\aalpha^{D_i}_\H=\{\alpha_i^1,\alpha_i^2\},\quad
D_{i,j}=\con{D_i}{\alpha_i^{3-j}}{\alpha_i^j}
\quad\text{and}\quad
\H_{j,k}=\H\cup\{D_{1,j},D_{2,k}\}$$
By \S \ref{sss:square}, $\H$ is given by
\begin{diagram}
\H_{1,1}&\hLine&\H_{1,2}\\
\vLine   &   \H   &\vLine   \\
\H_{2,1}&\hLine&\H_{2,2}
\end{diagram}
Set $\F_0=\H_{1,1}$ and let $\oveps$ be the clockwise
orientation of $\H$ so that
$$\mu(\H;\F_0,\oveps)=
\Phi_{\H_{1,1}\H_{2,1}} \Phi_{\H_{2,1}\H_{2,2}}
\Phi_{\H_{2,2}\H_{1,2}} \Phi_{\H_{1,2}\H_{1,1}}$$
For $j=1,2$, the unsaturated element of $\H_{j,1}\cap\H_{j,2}=\H\cup\{D_{1,j}\}$
is $D_2$ with $$\aalpha^{D_2} _{\H_{j,1}}=\alpha_2
\qquad\text{and}\qquad
\alpha^{D_2}_{\H_{j,2}}=\alpha_1$$
It follows from the forgetfulness and support axioms that
$$\Phi_{\H_{1,2}\H_{1,1}}=
\Phi_{\H_{2,1}\H_{2,2}}^{-1}\in
A_{D_2}^{D_2\setminus\aalpha_2}$$
and similarly that
$$\Phi_{\H_{2,2}\H_{1,2}}=
\Phi_{\H_{1,1}\H_{2,1}}^{-1}\in
A_{D_1}^{D_1\setminus\aalpha_1}$$
Since $D_1,D_2\in\H$ are compatible, $[A_{D_2}^{D_2\setminus\aalpha_2},
A_{D_1}^{D_1\setminus\aalpha_1}]=0$ and
\begin{multline*}
\Phi_{\H_{1,1}\H_{2,1}} \Phi_{\H_{2,1}\H_{2,2}}
\Phi_{\H_{2,2}\H_{1,2}}\Phi_{\H_{1,2}\H_{1,1}}\\
=
\Phi_{\H_{1,1}\H_{2,1}}\Phi_{\H_{2,2}\H_{1,2}}
\Phi_{\H_{2,1}\H_{2,2}}\Phi_{\H_{1,2}\H_{1,1}}=1
\end{multline*}
as claimed \halmos\\

The following result is the analogue for quasi--Coxeter algebras of
Mac Lane's coherence theorem for monoidal categories.

\begin{theorem}\label{th:coherence}
Let $A$ be a $D$--algebra and $\{\Phi_{\G\F}\}$ a collection of invertible
elements of $A$ labelled by elementary pairs of \mnss on $D$. Assume
that $\Phi_{\G\F}$ satisfy the orientation, forgetfulness and support
axioms of definition \ref{de:qcoxeter}. Then, the coherence axiom of
definition \ref{de:qcoxeter} is equivalent to the identities
$$\mu(\H)=1$$
for any pentagonal or hexagonal 2--face $\H$ of the associahedron
$\AD$.
\end{theorem}
\proof By the simple--connectedness of $\AD$, the coherence axiom
is equivalent to the identities $\mu(\H)=1$ for any two--face $\H$. The
conclusion now follows since, by \S \ref{ss:2 faces}, the two--faces of
$\AD$ are either squares, pentagons or hexagons and, by proposition
\ref{pr:zero on squares}, $\mu(\H)=1$ for any square 2--face $\H$
\halmos\\

\remark The identities $\mu(\H)=1$ corresponding to the pentagonal
and hexagonal 2--faces of $\AD$ are analogous to the pentagon identity
satisfied by the associator of a quasi--bialgebra. We shall refer to them
as {\it generalised pentagon relations}. These will be spelled out in
\S \ref{ss:qca 3}.

\subsection{Diagrammatic notation for elementary pairs}

\begin{proposition}\label{pr:diagrammatic}
The map
$$\imath:(\G,\F)\longrightarrow
(\supp(\G,\F);\alpha^{\supp(\G,\F)}_\G,\alpha^{\supp(\G,\F)}_\F)$$
induces a bijection between equivalence classes of elementary pairs of
\mnss on $D$ and triples $(B;\alpha,\beta)$ consisting of a connected
subdiagram $B\subseteq D$ and an ordered pair $(\alpha,\beta)$ of
distinct vertices of $B$.
\end{proposition}
\proof $\imath$ is injective by definition of equivalence. To show that it is
surjective, let $B\subseteq D$ be connected and let $\alpha_1\neq\alpha
_2$ be two vertices of $B$. Set $B_1=\con{B}{\alpha_1}{\alpha_2}$ and
$B_2=\con{B}{\alpha_2}{\alpha_1}$. Let $B^1,\ldots,B^k$ be the connected
components of $B\setminus\{\alpha_1,\alpha_2\}$ and choose a \mns $
\H^j$ on each $B^j$. Since $B^j$ is either contained in, or orthogonal to,
each of $B_1,B_2$,
$$\ol{\F}=\H^1\odots{\sqcup}\H^k\sqcup\{B_1,B\}
\qquad\text{and}\qquad
\ol{\G}=\H^1\odots{\sqcup}\H^k\sqcup\{B_2,B\}$$
are an elementary pair of \mnss on $B$ such that $\supp(\ol{\F},\ol{\G})=
B$, $\alpha^{B}_{\ol{\F}}=\alpha_1$ and $\alpha^{B}_{\ol{\G}}=\alpha_2$.
Choose next an increasing sequence $B=D_1\odots{\subset}D_m=D$
of connected subdiagrams such that $|D_{j+1}\setminus D_j|=1$
for any $j=1,\ldots,m-1$ and set
$$\F=\ol{\F}\sqcup\{D_2,\ldots,D_m\}
\qquad\text{and}\qquad
\G=\ol{\G}\sqcup\{D_2,\ldots,D_m\}$$
Then, $\imath(\G,\F)=(B;\alpha_2,\alpha_1)$ \halmos

\subsection{Diagrammatic notation for elementary associators}\label{ss:qca 3}

Let $A$ be a quasi--Coxeter algebra of type $D$. For any connected
subdiagram $B\subseteq D$ and ordered pair $(\alpha_i,\alpha_j)$
of distinct vertices of $B$, there exists by proposition \ref{pr:diagrammatic}
an elementary pair $(\G,\F)$ of \mnss on $D$ such that $B=\supp(\F,
\G)$, $\alpha^{B}_\F=\alpha_i$ and $\alpha^{B}_\G=\alpha_j$. Set
$$\Phi_{(B;\alpha_j,\alpha_i)}=\Phi_{\G\F}$$
The forgetfulness axioms implies that this definition is independent of
the choice of $\G,\F$. 

\begin{theorem}\label{th:comb ass}
The associators $\Phi_{(B;\alpha_j,\alpha_i)}$ satisfy the following
properties
\begin{itemize}
\item {\bf Orientation : }%
$$\Phi_{(B;\alpha_i,\alpha_j)}=\Phi_{(B;\alpha_j,\alpha_i)}^{-1}$$
\item {\bf Generalised pentagon relations : }%
For any connected $B\subseteq D$ and triple $(\alpha_i,\alpha_j,
\alpha_k)$ of distinct vertices of $B$, set
$$B_i=\conn{B}{\alpha_j,\alpha_k}{\alpha_i}\aand
B_{jk}=\con{B}{\alpha_i}{\{\alpha_j,\alpha_k\}}$$
Then, if $B_{jk}'=\emptyset$,
$$\phantom{Theorem}
\Phi_{(B		;\alpha_k,\alpha_i)}\cdot
\Phi_{(B		;\alpha_i,\alpha_j)}\cdot
\Phi_{(B_{ik}	;\alpha_i,\alpha_k)}\cdot
\Phi_{(B		;\alpha_j,\alpha_k)}\cdot
\Phi_{(B_{ij} 	;\alpha_j,\alpha_i)}=1$$
whereas, if $B_{ij},B_{jk},B_{ik}\neq\emptyset$,
\begin{multline*}
\phantom{Theore||}
\Phi_{(B		;\alpha_k,\alpha_i)}\cdot
\Phi_{(B_{jk}	;\alpha_k,\alpha_j)}\cdot
\Phi_{(B		;\alpha_i,\alpha_j)}\cdot\\
\cdot
\Phi_{(B_{ik}	;\alpha_i,\alpha_k)}\cdot
\Phi_{(B		;\alpha_j,\alpha_k)}\cdot
\Phi_{(B_{ij}	;\alpha_j,\alpha_i)}=1
\end{multline*}
\item {\bf Support : }%
$$\Phi_{(B;\alpha_j,\alpha_i)}\in A_{B}^{B\setminus\{\alpha_i,\alpha_j\}}$$
\item {\bf Braid relations : }%
if $m_{ij}<\infty$ and $B$ is the connected diagram with vertices
$\{\alpha_i,\alpha_j\}$, then
$$\Ad(\Phi_{(B;\alpha_i,\alpha_j)})(S^A_i)\cdot S^A_j\cdots=
S^A_j\cdot\Ad(\Phi_{(B;\alpha_i,\alpha_j)})(S^A_i)\cdots$$
where the number of factors on each side is equal to $m_{ij}$.
\end{itemize}
Conversely, if $A$ is a $D$--algebra endowed with invertible elements
$S_i^A\in\wh{A}_i$ for any $\alpha_i\in D$ and $\Phi_{(B;\alpha_j,\alpha_i)}
\in A$ for any connected subdiagram $B\subseteq D$ and ordered pair
of distinct vertices $(\alpha_i,\alpha_j)\in B$ which satisfy the above
properties, then the $S_i^A$ and the associators
\begin{equation}\label{eq:comb to pairs}
\Phi_{\G\F}=
\Phi_{(\supp(\G,\F);\alpha^{\supp(\G,\F)}_\G,\alpha^{\supp(\G,\F)}_\F)}
\end{equation}
corresponding to elementary pairs of \mnss on $D$ give $A$ the
structure of a quasi--Coxeter algebra of type $D$.
\end{theorem}
\proof Orientation and support are equivalent to the orientation and
support axioms satisfied by the elementary associators $\Phi_{\G\F}$.
Let $B\subseteq D$ be connected and $(\alpha_i,\alpha_j,\alpha_k)$ a
triple of distinct elements of $B$. Let $\H$ be a 2--face of $\AD$ such
that $\H$ has $B$ as its unique unsaturated element and $\aalpha^{B}
_\H=\{\alpha_i,\alpha_j,\alpha_k\}$. A simple exercise, using \S \ref
{sss:pentagon} shows that the identity $\mu(\H)=1$ corresponding
to $\H$ is the first or the second of the two generalised pentagonal
relations above depending on whether $B_{ij},B_{jk},B_{ki}$ are
all non--empty or one of them, which up to a relabelling we may assume
to be $B_{jk}$, is empty. Let now $\alpha_i\neq\alpha_j$ be such that
$m_{ij}<\infty$ and let $B$ be the connected diagram with vertices
$\alpha_i$ and $\alpha_j$. Let $(\F,\G)$ be an elementary pair of
\mnss on $D$ such that $\alpha_i\in\F$ and $\alpha_j\in\G$. Since
$\alpha_i$ and $\alpha_j$ are not compatible, they are the unique
elements in $\F\setminus\G$ and $\G\setminus\F$ respectively. Thus,
$\supp(\F,\G)=B$ and clearly $\alpha^{B}_\F=\alpha_j$ and $\alpha^
{B}_\G=\alpha_i$. It follows that $\Phi_{\G\F}=\Phi_{(B;\alpha_i,\alpha_j)}$
and the braid relations for $\Phi_{\G\F}$ coincide with those for $\Phi_
{(B;\alpha_i,\alpha_j)}$. The converse follows from the fact that the
associators $\Phi_{\G\F}$ defined by \eqref{eq:comb to pairs} clearly
satisfy the forgetfulness axiom and the fact that, by theorem \ref
{th:coherence}, the coherence axiom for the associators $\Phi_{\G\F}$
is equivalent to the relations $\mu(\H)=1$ which, as pointed out,
coincide with the generalised pentagonal relations above \halmos

\subsection{Twisting of quasi--Coxeter algebras}\label{ss:qca twisting}

Let $A$ be a quasi--Coxeter algebra of type $D$.

\begin{definition}\label{de:a-twist}
A twist $a=\{a_{(B;\alpha)}\}$ of $A$ is a collection of invertible elements of
$A$ labelled by pairs $(B;\alpha)$ consisting of a connected subdiagram
$B\subseteq D$ and a vertex $\alpha$ of $B$ such that
$$a_{(B;\alpha)}\in A_B^{B\setminus\alpha}$$
\end{definition}

Let $a$ be a twist. For any connected $B\subseteq D$ and \mns $\F$ on
$B$, set
\begin{equation}\label{eq:a_F}
a_{\F}=\prod_{B'\in\F}a_{(B';\alpha^{B'}_\F)}
\end{equation}
Note that the product does not depend upon the order of the factors.
Indeed, if $B'\neq B\dprime\in\F$, then either $B'\perp B\dprime$ in which
case $A_{B'}\ni a_{(B';\alpha^{B'}_\F)}$ and $A_{B\dprime}\ni a_{(B\dprime;
\alpha^{B\dprime}_\F)}$ commute, or, up to a permutation, $B'\subsetneq
B\dprime$ so that $B'\subseteq B\dprime\setminus\alpha^{B\dprime}_\F$
and $a_{(B\dprime;\alpha^{B\dprime}_\F)}$ commutes with $A_{B'}\ni a_
{(B';\alpha^{B'}_\F)}$.\\

We shall need the following

\begin{lemma}\label{le:twist of x}
Let $a$ be a twist, $\F$ a \mns on $D$ and $B\in\F$. Then,
\begin{enumerate}
\item for any $x\in A_B$ : $\Ad(a_\F)\medspace x=\Ad(a_{\F_B})
\medspace x$.\\[-.6em]
\item For any $x\in A^B$ : $\displaystyle{\Ad(a_\F)\medspace x=
\Ad(\prod_{\substack{B'\in\F,\\ B'\nsubseteq B}}
a_{(B';\alpha^{B'}_\F)})\medspace x}$.
\end{enumerate}
\end{lemma}
\proof
(i) If $B'\in\F$ and $B'\nsubseteq B$, then either $B'\perp B$,
so that $a_{(B';\alpha^{B'}_\F)}\in A_{B'}$ commutes with $A
_B$, or $B'\supsetneq B$, so that $B'\setminus\alpha^{B'}_\F
\supseteq B$ and again $a_{(B';\alpha^{B'}_\F)}\in A^{B'
\setminus\alpha^{B'}_\F}$ commutes with $A_B$. Thus, for
any $x\in A_B$,
$$\Ad(a_\F)\medspace x=
\prod_{B'\subseteq B}\Ad(a_{(B';\alpha^{B'}_\F)})\medspace x=
\Ad(a_{\F_B})\medspace x$$
(ii) If $B'\subseteq B$, then $a_{(B';\alpha^{B'}_\F)}\in A_{B'}
\subseteq A_B$ commutes with $x$. The result follows \halmos\\

For any $\alpha_i\in D$ and elementary pair $(\G,\F)$ of \mnss
on $D$, set
\begin{align}
S_i^{A^a}		&= a_{(\alpha_i;\alpha_i)}\cdot
				S^A_i\cdot a_{(\alpha_i;\alpha_i)}^{-1}
\label{eq:a twist of S}\\
\Phi^a_{\G\F}	&= a_\G\cdot\Phi_{\G\F}\cdot a_\F^{-1}
\label{eq:a twist of Phi}
\end{align}
Let $B=\supp(\F,\G)$ and set $\alpha_1=\alpha^{B}_{\F}$, $\alpha
_2=\alpha^{B}_\G$, $B_1=\con{B}{\alpha_1}{\alpha_2}$ and $B_2
=\con{B}{\alpha_2}{\alpha_1}$. 

\begin{lemma}\label{le:Phi^a}
$$\Phi_{\G\F}^a=
a_{(B;\alpha_2)}\cdot a_{(B_2;\alpha_1)}\cdot
\Phi_{\G\F}\cdot
a_{(B_1;\alpha_2)}^{-1}\cdot a_{(B;\alpha_1)}^{-1}$$
\end{lemma}
\proof Let $B'\in\G\cap\F=\F\setminus\{B_1\}=\G\setminus\{B_2\}$ be
distinct from $B$ so that $\alpha^{B'}_{\F}=\alpha^{B'}_{\F\cap\G}=
\alpha^{B'}_{\G}$. Lemma \ref{le:twist of x} readily implies that $a_
{(B';\alpha^{B'}_\G)}=\alpha_{(B';\alpha^{B'}_\F)}$ commutes with
$\Phi_{\G\F}\in A_B^{B\setminus\{\alpha_1,\alpha_2\}}$. The identity
above now follows from this and the fact that, by proposition \ref
{pr:combi<=>practical}, $\alpha^{B_1}_\F=\alpha_2$ and $\alpha^
{B_2}_\G=\alpha_1$ \halmos 

\begin{proposition}
$A^a=(A,A_{D'},S_i^{A^a},\Phi^a_{\G\F})$ is a quasi--Coxeter algebra
of type $D$ called the twist of $A$ by $a$.
\end{proposition}
\proof The elements $\Phi_{\G\F}^a$ clearly satisfy the orientation and
coherence axioms of definition \ref{de:qcoxeter}. Lemma \ref{le:Phi^a}
implies that they also satisfy the support and forgetfulness axiom.
Let now $\alpha_i\neq\alpha_j\in D$ be such that $m_{ij}<\infty$ and
$(\F,\G)$ an elementary pair of \mnss on $D$ such that $\alpha_i\in\F$
and $\alpha_j\in\G$. By lemma \ref{le:twist of x}, $\Ad(a_\F)(S^A_i)=
\Ad(a_{(\alpha_i;\alpha_i)})(S^A_i)=S_i^{A^a}$ and similarly $\Ad(a_\G)
(S^A_j)=S_j^{A^a}$. Thus,
$$\Ad(\Phi^a_{\G\F})(S^{A^a}_i)\cdot S^{A^a}_j\cdots=
\Ad(a_\G)\left(\Ad(\Phi_{\G\F})(S_i^A)\cdot S_j^A\cdots\right)$$
and 
$$S_j^{A^a}\cdot \Ad(\Phi^a_{\G\F})(S_i^{A^a})\cdots=
\Ad(a_\G)\left(S_j^A\cdot \Ad(\Phi_{\G\F})(S_i^A)\cdots\right)$$
so that $\Phi^a_{\G\F}$ satisfies the braid relations with respect to
$S_i^{A^a}$, $S_j^{A^a}$ \halmos\\

The following result shows that twisting does not change the isomorphism
class of a \qc algebra.

\begin{proposition}
If $a$ is a twist of $A$, then $\Psi_\F=\Ad(a_\F)$ is an isomorphism of
$A$ onto the \qca $A^a$. In particular, $A$ and $A^a$ define equivalent
representations of the Artin group $B_D$.
\end{proposition}
\proof If $\alpha_i\in\F$, then, by lemma \ref{le:twist of x}, $\Ad(a_\F)(S_i^A)
=S_i^{A^a}$. Moreover, for any elementary pair $(\G,\F)$ of \mnss on $D$, 
$$\Ad(a_\G)\circ\Ad(\Phi_{\G\F})=\Ad(\Phi^a_{\G\F})\circ\Ad(a_\F)$$
\halmos\\

\remark By lemma \ref{le:Phi^a}, the twist by $a=\{a_{(B;\alpha)}\}_{\alpha\in
B\subseteq D}$ of the associators of a quasi--Coxeter algebra reads, in
diagrammatic notation
\begin{equation}\label{eq:Phi^a_D}
\Phi^a_{(B;\alpha_j,\alpha_i)}=
a_{(B;\alpha_j)}\cdot a_{(\con{B}{\alpha_j}{\alpha_i};\alpha_i)}\cdot
\Phi_{(B;\alpha_j,\alpha_i)}\cdot
a_{(\con{B}{\alpha_i}{\alpha_j};\alpha_j)}^{-1}\cdot a_{(B;\alpha_i)}^{-1}
\end{equation}


%% file: qcexamplesnew.tex

\newcommand {\iC}{\imath_C}
\newcommand {\ik}{\imath_\kappa}

\section{Examples of \qcas}\label{se:qcexamples}

This section is devoted to the study of several examples of quasi--Coxeter
algebras. In \S \ref{ss:quantum examples} we consider 'quantum' examples,
the defining feature of which is that their associators are all trivial. We begin
with the universal one given by the braid group of a Coxeter group, then
consider the corresponding Hecke algebra and finally Lusztig's quantum
Weyl group operators for a symmetrisable \KM algebra. In \S \ref
{ss:diff examples}, we consider examples which underlie the monodromy
representations of several flat connections, specifically the holonomy
equations \eqref{eq:holonomy} of a Coxeter arrangement, Cherednik's
KZ connection and the Casimir connection described in the Introduction.
The study of these examples relies heavily on the \DCP theory of asymptotic
zones described in section \ref{se:DCP}. Finally, in \S \ref{ss:qba examples}
we show how to obtain \qcas of type $\sfA_n$ as commutants of \qbas and
of quasitriangular quasibialgebras.\\

Throughout this section, $W$ denotes an irreducible Coxeter group with
system of generators $S=\{s_i\}_{i\in\bfI}$. We denote by $D$ the Coxeter
graph of $(W,S)$ and label the pair $i\neq j\in D$ with the order $m_{ij}$
of $s_is_j$ in $W$.

\subsection{Quantum examples}\label{ss:quantum examples}

\subsubsection{Universal example}\label{ss:q universal}

For any connected subgraph $D'\subseteq D$ with vertex set $\bfI'\subseteq
\bfI$, let $W_{D'}\subseteq W$ be the parabolic subgroup generated by $s_i$,
$i\in\bfI'$ and $B_{D'}$ the (algebraic) braid group of $W_{D'}$, that is the
group with generators $S_i$, $i\in\bfI'$ and relations \eqref{eq:Brieskorn}
for any $i\neq j\in\bfI'$. Then, as noted in remark \rkref{rk:univeral}, the
assignement
$$A_{D'}=k[B_{D'}],\qquad
S_i^A=S_i\qquad\text{and}\qquad
\Phi_{\G\F}^A=1$$
endows $A=k[B_D]$ with the structure of a \qca of type $D$ and if $A'$ is a
\qca of type $D$, the maps $\pi_\F:B_D\rightarrow\wh{A'}$ given by theorem
\ref{th:monodromy maps} define a morphism of \qcas $k[B_D]\rightarrow\wh{A'}$.

\subsubsection{Hecke algebras}\label{ss:Hecke}

Let $q_i\in k$ be invertible elements such that $q_i=q_j$ whenever $s_i$
is conjugate to $s_j$ in $W$. For any $D'\subseteq D$ with vertex set $\bfI'
\subseteq\bfI$, let $\H(W_{D'})$ be the \IH algebra of $W_{D'}$,
that is the quotient of $k[B_{D'}]$ by the quadratic relations
\begin{equation}\label{eq:Hecke}
(S_i-q_i)(S_i+q_i^{-1})=0
\end{equation}
Then, the assignement
$$A_{D'}=\H(W_{D'}),\qquad
S_i^A=S_i\qquad\text{and}\qquad
\Phi_{\G\F}^A=1$$
endows $A=\H(W)$ with a structure of \qca of type $D$. The corresponding
maps $\pi_\F:B_D\rightarrow A$ are all equal to the quotient map $k[B_D]
\rightarrow\H(W)$.

\subsubsection{Quantum Weyl groups}\label{ss:qWeyl}

Assume now that $W$ is the Weyl group of a \KM algebra $\g=\g(\bfA)$
with generalised Cartan matrix $\bfA=(a_{ij})_{i,j\in\bfI}$ \cite{Ka}. Let
$(\h,\Delta,\Delta^\vee)$ be the unique realisation of $\bfA$. Thus, $\h$
is a complex vector space of dimension $2|\bfI|-\rk(\bfA)$,
$$\Delta=\{\alpha_i\}_{i\in\bfI}\subset\h^*
\qquad\text{and}\qquad
\Delta^\vee=\{\alpha^\vee_i\}_{i\in\bfI}\subset\h$$
are linearly independent sets of cardinality $|\bfI|$ and, for any $i,j\in\bfI$,
$a_{ij}=\<\alpha_i^\vee,\alpha_j\>$. Then, $W$ is the subgroup of $GL(\h)$
generated by the reflections $s_i$ acting on $t\in\h$ by
$$s_i(t)=t-\alpha_i(t)\cdot \alpha_i^\vee$$
and $m_{ij}=2,3,4,6,\infty$ according to whether $a_{ij}a_{ji}=0,1,2,3,\geq 4$.\\

Assume further that $\bfA$ is symmetrisable, that is that there exists relatively
prime integers $d_i\geq 1$, $i\in\bfI$ such that $d_ia_{ij}=d_ja_{ji}$ for any $i,j\in
\bfI$. Then, there exists a non--degenerate bilinear form $(\cdot,\cdot)$ on $\h^
*$, unique up to a scalar, such that $a_{ij}=2\frac{(\alpha_i,\alpha_j)}{(\alpha_i,
\alpha_i)}$. Let $\g'=[\g,\g]$ be the derived subalgebra of $\g$, $k=\ICh$ the
ring of formal power series in the variable $\hbar$ and $\Uhg'$ the \DJ quantum
group corresponding to $\bfA$ and $(\cdot,\cdot)$, that is the algebra over
$\IC\fml$ topologically generated by elements $E_{i},F_{i},H_{i}$, $i\in\bfI$,
subject to the relations\footnote{we follow here the conventions of \cite{Lu}.}
\begin{gather*}
[H_{i},H_{j}]=0
\\[1.1 ex]
[H_{i},E_{j}]= a_{ij}E_{j}
\qquad
[H_{i},F_{j}]=-a_{ij}F_{j}
\\[1.1 ex]
[E_{i},F_{j}]=
\delta_{ij}\frac
{q_i^{H_{i}}-q_i^{-H_{i}}}
{q_i-q_i^{-1}}
\intertext{where}
q_i=q^{(\root{i},\root{i})/2}
\qquad\text{with}\qquad
q=e^{\hbar}
\intertext{and the $q$--Serre relations}
\sum_{k=0}^{1-a_{ij}}(-1)^{k}
\bin{1-a_{ij}}{k}{i}
E_{i}^{k}E_{j}E_{i}^{1-a_{ij}-k}=0
\\[1.1 ex]
\sum_{k=0}^{1-a_{ij}}(-1)^{k}
\bin{1-a_{ij}}{k}{i}
F_{i}^{k}F_{j}F_{i}^{1-a_{ij}-k}=0
\intertext{where for any $k\leq n$,}
[n]_{i}=\frac{q_i^{n}-q_i^{-n}}{q_i-q_i^{-1}}\\
[n]_{i}!=[n]_{i}[n-1]_{i}\cdots[1]_{i}
\qquad\text{and}\qquad
\bin{n}{k}{i}=\frac{[n]_{i}!}{[k]_{i}![n-k]_{i}!}
\end{gather*}

For any connected $D'\subseteq D$ with vertex set $\bfI'\subseteq\bfI$,
let $\Uhg'_{D'}\subseteq\Uhg'$ be the subalgebra topologically generated
by the elements $E_i,F_i$ and $H_i$ with $i\in\bfI'$. Then, $\Uhg'_{D'}$ is
the \DJ quantum group corresponding to the Cartan matrix $(a_{ij})_{i,j
\in \bfI'}$ and the unique non--degenerate bilinear form on its realisation
which coincides with the restriction of $(\cdot,\cdot)$ on the span of the
$\alpha_i$, $i\in \bfI'$. If $\bfI'=\{i\}$, we denote $\Uhg'_{D'}$ by $\Uhsl{2}
^i$.\\

For any $i$, let $\ol{S}\ih$ be the operator acting on an integrable
$\Uhg$--module $\V$ as\footnote{the element $\ol{S}\ih$ is, in the
notation of \cite[\S 5.2.1]{Lu}, the operator $T_{i,+1}''$.}
\begin{equation}\label{eq:barSih}
\ol{S}\ih\medspace v=
\sum_{\substack{a,b,c\in\IZ : \\a-b+c=-\lambda(\cor{i})}}
(-1)^{b}q_i^{b-ac}
E_i^{(a)}F_i^{(b)}E_i^{(c)}
v
\end{equation}
where
$$E_{i}^{(a)}=\frac{E_i^{a}}{[a]_i!}
\qquad
F_{i}^{(a)}=\frac{F_i^{a}}{[a]_i!}$$
and $v\in V$ if of weight $\lambda\in\h^*$. Set
\begin{equation}\label{eq:sih}
S\ih=
\ol{S}\ih\cdot q_i^{H_i^2/4}=
q_i^{H_i^2/4}\cdot\ol{S}\ih
\in
\wh{\Uhsl{2}^i}
\end{equation}
We shall refer to $S_i$ as the {\it quantum Weyl group operator}
corresponding to $i$. The following result is due to Lusztig, \KR
and Soibelman \cite{Lu,KR,So}.

\begin{proposition}\label{pr:qW ids1}
If the order $m_{ij}$ of $s_is_j$ in $W$ is finite, then
$$\underbrace{S\ih S\jh\cdots}_{m_{ij}}=
\underbrace{S\jh S\ih\cdots}_{m_{ij}}$$
\end{proposition}
\proof this is an immediate consequence of the braid relations
satisfied by the operators $\ol{S}\ih$ \cite [\S 39.4]{Lu} \halmos\\

It follows from proposition \ref{pr:qW ids1} that the assignement
$$A_{D'}=\Uhg_{D'}',\qquad
S_i^A=S\ih\qquad\text{and}\qquad
\Phi^A_{\G\F}=1$$
endows $A=\Uhg'$ with the structure of a \qca of type $D$. The
corresponding representations $B_D\rightarrow GL(\V)$, with
$\V$ a \fd $\Uhg'$--module are called {\it quantum Weyl group
representations}.

\subsection{Differential examples}\label{ss:diff examples}

Assume now that $W$ is finite. Let $V_\IR$ be its reflection representation,
$\A_\IR\subset V_\IR$ the corresponding arrangement of reflecting
hyperplanes of $W$ and $V,\A$ their complexifications. Retain the notation
of section \ref{se:DCP}, particularly \S \ref{ss:coxeter}. Thus, $\Phi\subset
V_\IR^*$ is a root system for $W$, $\C\subset V_\IR\setminus\A_\IR$ a chamber,
which we choose to be the one bound by the reflecting hyperplanes of the
generators $s_i$ of $W$, $\Phi=\Phi_+\sqcup\Phi_-$ the corresponding
partition into positive and negative roots, $\Delta=\{\alpha_i\}_{i\in\bfI}\subset
\Phi_+$ the basis of $V_\IR^*$ consisting of indecomposable elements of
$\Phi_+$, and we choose $X=\Phi_+$ as the set of defining equations for
$\A$. Fix $v_0\in\C$ and identify $B_W=\pi_1(\VA/W,v_0)$ with $B_D$ via
the presentation \eqref{eq:Brieskorn}.\\

Recall that, by proposition \ref{pr:dictionary}, there is a bijection between
\nss of connected subdiagrams of $D$ and \fnss of irreducible subspaces
of $V^*$ which contain $V^*$. This correspondence maps $D'\subseteq D$
with vertex set $\bfI'\subseteq\bfI$ to the subspace $\<D'\>$ spanned by
the roots $\alpha_i$, $i\in \bfI'$ and $B\subseteq V^*$ to the connected
subdiagram with vertex set $B\cap\Delta$.

\subsubsection{Universal example}\label{ss:diff universal}

For any subdiagram $D'\subseteq D$, let $F_{D'}\subseteq F$ be the subalgebra
generated by the elements $t_\alpha$, $\alpha\in\<D'\>$, $I_{D'}=I\cap
F_{D'}$ the ideal defined by the relations \eqref{eq:rank 2} for $B\subseteq
\<D'\>$ and $A_{D'}$ the completion of $F_{D'}/I_{D'}$ with respect to its $\IN
$--grading. The braid group $B_{D'}$ acts on $F/I$ by \eqref{eq:W on A} via the
homomorphism $B_{D'}\rightarrow W_{D'}$ and leaves $F_{D'}/I_{D'}$ invariant.\\

Let $h$ be a formal variable and $\imath_h:A\rtimes B_D\rightarrow F/I\fmll
\rtimes B_D$ the embedding given on $b\in B_D$ and $t_\alpha\in A$ by
\begin{equation}\label{eq:i_h}
\imath_h(b)=b
\qquad\text{and}\qquad
\imath_h(t_\alpha)=ht_\alpha
\end{equation}
Note that if $\wt{\rho}:F/I\rtimes B_D\rightarrow\End(U)$ is a \fd representation,
the action $\wt{\rho}_h$ of $A\rtimes B_D$ on $U\fmll$ defined by \eqref
{eq:extended rho_h} is given by $\rho\circ\imath_h$.

\begin{proposition}\label{pr:universal diff} Set $k=\IC\fmll$.
\begin{enumerate}
\item The assignement
$$A^\nabla_{D'}=F_{D'}/I_{D'}\fmll\rtimes B_{D'}$$
endows $A^\nabla=F/I\fmll\rtimes B_D$ with the structure of a $D$--algebra
\mbox{over $k$}.
\item Let $\beta\subset V^*$ be a positive, real adapted family. Then, the
elements
$$S_i^\nabla=\exp(\pi\sqrt{-1}\cdot h t_{\alpha_i})\cdot S_i
\quad\text{and}\quad
\Phi^\nabla_{\G\F}=\imath_h(\Phi^\beta_{\G\F})$$
where the latter are the \DCP associators corresponding to $\beta$,
endow $A^\nabla$ with the structure of a \qca of type $D$.
\item If $\wt{\rho}:F/I\rtimes B_D\longrightarrow End(U)$ is a finite--dimensional
representation, the action $\wt{\rho}\circ\pi_\F$ of $B_D$ on $U\fmll$ induced
by the \qca structure on $A^\nabla$ coincides with the monodromy of the flat
vector bundle $(\wt{\VA}\times U,p^*\nabla_{\wt{\rho}_h})/B_D$ over $\VA/W$,
where
$$\nabla_{\wt{\rho}_h}=
d-h\sum_{\alpha\in\Phi_+}\frac{d\alpha}{\alpha}\wt{\rho}(t_\alpha)$$
expressed in the fundamental solution $\wt{\rho}_h(\Psi_\F^{\beta_\F})$.
\item The \qca structures on $A^\nabla$ given by two positive, real adapted
families $\beta,\beta'$ differ by a canonical twist.
\end{enumerate}
\end{proposition}
\proof
(i) Clearly, $A^\nabla_{D'}\subset A^\nabla_{D\dprime}$ whenever $D'\subset
D\dprime$. If, on the other hand $D'\perp D\dprime$, then $\{\<D'\>,\<D\dprime
\>\}$ is nested by proposition \ref{pr:dictionary} so that $[A_{D'},A_{D\dprime}]
=0$ by lemma \ref{le:[t_B,t_C]}. It follows that $[A^\nabla_{D'},A^\nabla_{
D\dprime}]=0$ since, whenever $D'\perp D\dprime$, $B_{D'}$ fixes $A_{D
\dprime}$ and commutes with $B_{D\dprime}$.

(ii) The orientation, transitivity, forgetfulness and braid relations axioms have
been proved in \S \ref{ss:elementary associators}, \S \ref{ss:adapted family},
\S \ref {ss:forgetful} and \S \ref{sss:DCP braid relations} respectively. To check
the support axiom, let $(\G,\F)$ be an elementary pair of \mnss on $D$, $D'
=\supp(\G,\F)$ and $D\dprime=\zsupp(\F,\G)$. By \S \ref{ss:support}, the
associator $\Phi_{\G\F}^\beta$ lies in $A_{D'}$ and commutes with $A_{D
\dprime}$. We therefore need only prove that it is invariant under $W_{D
\dprime}$. It is sufficient to show that this is the case for the coefficients
$t_{\ol{x}}$ of the connection $\nabla_{\<D'\>,\<D\dprime\>}$ defined by
\eqref{eq:nablaBC}. However, if $x\in X\cap(\<D'\>\setminus\<D\dprime\>)$
and $\alpha_i\in D\dprime$, then
$$s_i x=x-2(x,\alpha_i)/(\alpha_i,\alpha_i)\alpha_i=x\mod\<D\dprime\>$$
so that
$$s_i(t_{\ol{x}})=s_i(\sum_{x\in\ol{x}}t_x)=\sum_{x\in\ol{x}}t_{s_i(x)}=t_{\ol{x}}$$

(iii) Follows from the discussion in \S \ref{sss:relevance}.

(iv) By remark \rkref{rk:gauge}, $\Phi_{\G\F}^{\beta'}=a_\G\cdot\Phi_{\G\F}
^\beta\cdot a_\F^{-1}$, where
\begin{equation}\label{eq:explicit twisting}
a_\F=\prod_{D'\in\F}{c_{(D';\alpha^{D'}_\F)}}^{-R_{\<D'\>}^\F}
\end{equation}
for some $c_{(D';\alpha)}\in\IR_{+}^*$ and
$$R_{\<D'\>}^\F=
t_{\<D'\>}-t_{i_\F\<D'\>}=
t_{\<D'\>}-t_{\<D'\setminus\alpha^{D'}_\F\>}\in A_{D'}$$
is defined by \eqref{eq:R_B^S}. The claim now follows since, by lemma
\ref{le:[t_B,t_C]}, $R_{\<D'\>}^\F$, commutes with $A_{D'\setminus\alpha
^{D'}_\F}$ and $t_{\alpha_i}$ is fixed by $s_i$, so that $\Ad(c_{(\alpha_i;
\alpha_i)}^{-t_{\alpha_i}})S_i^\nabla=S_i^\nabla$ \halmos\\

\remark Since the elements $S_i^\nabla$ lie in $A^\nabla$, the maps $\pi_\F
:B_D\rightarrow\wh{A^\nabla}$ given by theorem \ref{th:monodromy maps}
factor through $A^\nabla$. Moreover, by theorem \ref{th:extended Chen},
they induce isomorphism
$$\wh{\pi_\F}:\wh{\IC[B_D]\fmll}\longrightarrow F/I\fmll\rtimes W$$
where $\wh{\IC[B_D]\fmll}$ is the completion of $\IC[B_D]\fmll$ with respect
to the kernel of the epimorphism $\IC[B_D]\fmll\rightarrow\IC[W]\fmll$.

\subsubsection{Cherednik's rational KZ connection}\label{s s:CKZ}

We recast below the monodromy of Cherednik's KZ connection for $W$ \cite
{Ch1,Ch2} in the language of \qc algebras. For any $\alpha\in\Phi$, let $s_\alpha
\in W$ be the corresponding orthogonal reflection and consider the connection
on $\VA$ with values in $\IC[W]$ given by\footnote{The letter $C$ is used to
distinguish $\nabla_\CKZ$ from the \KZ connection $\nabla\KKZ$ which arises
in Conformal Field Theory. The latter depends upon the choice of a complex,
reductive Lie algebra $\r$ and takes values in the $n$--fold tensor product $
V^{\otimes n}$ of a \fd $\r$--module $V$. It coincides with $\nabla_\CKZ$ for
the Coxeter group $W=\SS_n$ when $\r=\gl{n}$ and $V=\IC^n$ is the vector
representation, via the identication $V^{\otimes n}[0]\cong\IC[\SS_n]$ (see,
\eg \cite[\S 4]{TL2}). In \cite{TL2}, the name Coxeter--KZ connection was
used for $\nabla_\CKZ$. The name Cherednik--KZ connection seems far
more appropriate however.}
$$\nabla_\CKZ=
d-\sum_{\alpha\in\Phi_+}k_\alpha\frac{d\alpha}{\alpha}s_\alpha$$
where the $k_\alpha$ are complex numbers such that $k_{\alpha}=k_{\alpha'}$
whenever $s_\alpha$ and $s_{\alpha'}$ are conjugate in $W$.

\begin{theorem}[\cite{Ch1,Ch2}]\label{th:Cherednik}
The connection $\nabla_\CKZ$ is flat.
\end{theorem}

By theorem \ref{th:Cherednik}, the $W$--equivariant homomorphism $F
\rightarrow\IC[W]$ given by $t_\alpha\rightarrow k_\alpha s_\alpha$ factors
through $F/I$. Composing with the embedding \eqref{eq:i_h} yields a
homomorphism $\imath_\CKZ:A\rtimes B_D\rightarrow \IC\fmll[W]$ restricting
to the canonical projection $B_D\rightarrow W$. The following result is
an immediate consequence of proposition \ref{pr:universal diff}.

\begin{proposition}\label{pr:Cherednik diff}\hfill
\begin{enumerate}
\item The assignement $A^\CKZ_{D'}=\IC\fmll[W_{D'}]$ endows $A^\CKZ=
\IC\fmll[W]$ with the structure of a $D$--algebra over $k=\IC\fmll$.
\item Let $\beta\subset V^*$ be a positive, real adapted family. Then,
the elements
$$S_i^\CKZ=\exp(\pi\sqrt{-1}\cdot h k_{\alpha_i}s_i)\cdot s_i
\quad\text{and}\quad
\Phi^\CKZ_{\G\F}=\imath_\CKZ(\Phi^\beta_{\G\F})$$
give $A^\CKZ$ the structure of a quasi--Coxeter algebra of type $D$.
\item If $\rho:W\rightarrow GL(U)$ is a finite--dimensional representation,
the action $\rho\circ\pi_\F$ of $B_D$ on $U\fmll$ induced by the \qca
structure on $A^\CKZ$ coincides with the monodromy representation of
the flat connection
$$\nabla_\CKZ^\rho=
d-h\sum_{\alpha\in\Phi_+}k_\alpha\frac{d\alpha}{\alpha}\rho(s_\alpha)$$
on the holomorphic vector bundle $\U_{\rho}=\VA\times_{W}U\fmll$ over
$\VA/W$ expressed in the fundamental solution $\rho\circ\imath_\CKZ
(\Psi_\F^{\beta_\F})$.
\item The \qca structures on $A^\CKZ$ given by two positive, real adapted
families $\beta,\beta'$ differ by a twist.
\end{enumerate}
\end{proposition}

Let $\H(W)$ be the \IH algebra of $W$ over $k=\IC\fml$ defined by 
\eqref{eq:Hecke} with $q_i=\exp(\pi\sqrt{-1}\cdot h k_{\alpha_i})$.
The following is a reformulation of a well--known result of Cherednik
\cite{Ch1,Ch2}.

\begin{theorem}\label{th:Ch iso}
For any \mns $\F$ on $D$, the map $\pi_\F:B_D\rightarrow A^\CKZ$
induced by the \qca structure on $A^\CKZ$ factors through $\H(W)$
and induces an isomorphism $\pi_\F:\H(W)\longrightarrow\IC\fmll[W]$
which restricts to an isomorphism
\begin{equation}\label{eq:parabolic}
\H(W_{D'})\longrightarrow\IC\fmll[W_{D'}]
\quad
\text{for any $D'\in\F$}
\end{equation}
\end{theorem}
\proof 
That $\pi_\F$ factors through $\H(W)$ follows at once from theorem \ref
{th:monodromy factor} and proposition \ref{pr:nu monodromy}. $\mu_\F$
is surjective by Nakayama's lemma, since its reduction mod $h$ is the
canonical projection $B_D\rightarrow W$ and it is injective since the
dimension of $\H(W)$ is bounded above by that of $\IC[W]$ \halmos\\

\remark The only significant difference between theorem \ref{th:Ch iso}
and theorem 2 in \cite{Ch1,Ch2} is that we consider a {\it collection} of
isomorphisms $\mu_\F:\H(W)\rightarrow\IC\fmll[W]$. These are labelled
by \mnss on the Coxeter graph $D$ of $W$, respect the parabolic structure
on $\H(W)$ and $\IC\fmll[W]$ in the sense of \eqref{eq:parabolic} and
are related by
$$\mu_\G=\Ad(\Phi^\CKZ_{\G\F})\circ\mu_\F$$

\remark Cherednik has given an explicit computation of the monodromy
of $\nabla_\CKZ$ when $W$ is of type $\sfA,\sfB,\sfC,\sfD$ in terms of
the classical hypergeometric function \cite{Ch3}.

\subsubsection{The Casimir connection}\label{ss:Casimir}

Assume now that $W$ is the Weyl group of a complex, simple Lie algebra
$\g$ with Cartan subalgebra $\h\subset\g$, identify $\VA$ with the set
$$\hreg=\h\setminus\bigcup_{\alpha\in\Phi}\Ker(\alpha)$$
of regular elements in $\h$ and choose $\Phi\subset\h^*$ to be the root
system of $\g$. For any $\alpha\in\Phi$, let $\sl{2}^\alpha=
\<e_\alpha,f_\alpha,h_\alpha\>\subset\g$ be the corresponding
three--dimensional subalgebra and denote by
$$C_\alpha=\frac{(\alpha,\alpha)}{2}
\left(e_\alpha f_\alpha+f_\alpha e_\alpha+\frac{1}{2}h_\alpha^2\right)$$
its Casimir operator with respect to the restriction to $\sl{2}^\alpha$ of
a a fixed non--degenerate, invariant bilinear form $(\cdot,\cdot)$ on $\g$.
Note that $C_\alpha$ is independent of the choice of the root vectors
$e_\alpha,f_\alpha$ and satisfies $C_{-\alpha}=C_{\alpha}$. Let $V$
a \fd $\g$--module and consider the following holomorphic connection
on the holomorphically trivial vector bundle over $\hreg$ with fibre $V$
\begin{equation}\label{eq:nablac}
\nablac=
d-\sum_{\alpha\in\Phi_+}\frac{d\alpha}{\alpha}\cdot C_\alpha
\end{equation}
The following result is due to the author and J. Millson \cite{MTL}
and was discovered independently by De Concini around 1995
(unpublished) and by Felder \etal \cite{FMTV}.
\begin{theorem}
The connection $\nablac$ is flat.
\end{theorem}

It will be convenient to use the following, closely related auxiliary connection
$$\nablak=d-\sum_{\alpha\in\Phi_+}\frac{d\alpha}{\alpha}\cdot\kalpha
\qquad\text{where}\qquad
\kalpha=\frac{(\alpha,\alpha)}{2}(e_\alpha f_\alpha+f_\alpha e_\alpha)$$
is the truncated Casimir operator of $\sl{2}^\alpha$. $\nablak$ is also flat
\cite{MTL} and therefore determines a homomorphism $F/I\rightarrow \Ug$
given by $t_\alpha\rightarrow\kalpha$. We extend it to a morphism $F/I\rtimes
B_D\rightarrow\wh{\Ug}$ by mapping the generator $S_i\in B_D$ to the
triple exponential
\begin{equation}\label{eq:triple exp}
\wt{s}_i=\exp(e_i)\exp(-f_i)\exp(e_i)
\end{equation}
where $e_i=e_{\alpha_i}$ and $f_i=f_{\alpha_i}$ are a fixed choice of simple
root vectors.\\

By \cite{Ti}, the assignement $S_i\rightarrow\wt{s}_i$ extends
to a homomorphism $\sigma:B_D\rightarrow N(H)\subset G$ where $G$ is
the connected and simply--connected complex Lie group corresponding to
$\g$, $H\subset G$ the maximal torus with Lie algebra $\h$ and $N(H)$ its
normaliser. The image of $\sigma$ is an extension $\wt{W}$ of $W=N(H)/H$
by the sign group $\IZ_2^{\dim(\h)}$ which we shall call the {\it Tits extension}.
We shall however regard the elements $\wt{s}_i$ as lying in the completion
$\wh{\Ug}$ of $\Ug$ \wrt to its \fd representations rather than in $N(H)$.\\

Composing with the embedding \eqref{eq:i_h} we therefore obtain a
homomorphism
$\ik:A\rtimes B_D\rightarrow\wh{\Ug}\fmll=\wh{\Ug\fmll}$ given by
$$\ik(t_\alpha)=h\kalpha\aand\ik(S_i)=\wt{s}_i$$
Similarly, we have a homomorphism $\iC:A\rtimes B_D\rightarrow\wh{\Ug}\fmll$
given by $\iC(t_\alpha)=h\calpha$ and $\iC(S_i)=\wt{s}_i$.  For any
subdiagram $D'\subseteq D$ with vertex set $\bfI'\subseteq\bfI$, let $\g_{D'}\subseteq\g_D$
be the subalgebra generated by $e_i,f_i,h_i$, $i\in\bfI'$.

\begin{theorem}\label{th:qca on Ug}\hfill
\begin{enumerate}
\item The assignement $D'\rightarrow U\g_{D'}\fmll$ defines a $D$--algebra
structure on $\Ug\fmll$.
\item Let $\beta$ be a positive, real adapted family. Then, the elements
$$S_i^C=\wt{s}_i\cdot\exp(\pi\sqrt{-1}\cdot h C_{\alpha_i})
\qquad\text{and}\qquad
\Phi^C_{\G\F}=\ik(\Phi^\beta_{\G\F})$$
define a \qca structure of type $D$ on $\Ug\fmll$.
\item If $\rho:\g\rightarrow\End(V)$ is a \fd $\g$--module, the corresponding
action $\pi_\F:B_D\rightarrow GL(V\fmll)$ given by the \qca structure
is equivalent to the monodromy representation of the flat vector bundle
$$(\wt{\hreg}\times V,p^*\nablac^h)/B_D
\quad\text{where}\quad
\nablac^h=
d-h\sum_{\alpha\in\Phi_+}\frac{d\alpha}{\alpha}\cdot C_\alpha$$
expressed in the fundamental solution $\rho\circ\iC(\Psi_\F^\beta)$.
\item The \qca structures corresponding to two adapted families
$\beta,\beta'$ differ by a canonical twist.
\end{enumerate}
\end{theorem}
\proof
(i) is clear. (ii) For any $i\in\bfI$, set
$$S_i^\kappa=
\wt{s}_i\cdot\exp(\pi\sqrt{-1}\cdot h\kappa_{\alpha_i})=
\ik(S_i\cdot\exp(\pi\sqrt{-1}\cdot t_{\alpha_i}))$$
We first prove that the $S_i^\kappa$ and associators $\Phi^C_{\G\F}$
endow $\Ug\fmll$ with the structure of a \qca of type $D$. This follows by
transport of structure from proposition \ref{pr:universal diff} except for the
statement that, for an elementary pair $(\G,\F)$ of \mnss on $D$ with $D'
=\supp(\G,\F)$ and $D\dprime=\zsupp(\G,\F)$ the associator $\Phi^C_{\G\F}$
commutes with $\Ug_{D\dprime}$. To see this, it is sufficient to prove
that the coefficients $\ik(t_{\ol{x}})$ of the connection $\nabla_{\<D'\>,
\<D\dprime\>}$ defined by \eqref{eq:nablaBC} are invariant under the
adjoint action of $\g_{D\dprime}$\footnote{this is not true of the coefficients
$\iC(t_{\ol{x}})$ which is why we prefer to work with the connection
$\nablak$. $\nablac$ on the other hand has better local monodromies
since their squares $(S_i^C)^2=(-1)^{h_i}\exp(2\pi h\sqrt{-1}
C_{\alpha_i})$ are central in the corresponding $\wh{\Uhsl{2}^i}\fmll$.}.
Let $\alpha\in\Phi_+\cap(\<D'\>\setminus \<D\dprime\>)$
and $\<D\dprime\>_\alpha$ the irreducible component of $\<D\dprime\>
\oplus\IC\alpha$ containing $\alpha$ as in lemma \ref{le:C_x}. Set
$\<D\dprime\>^\alpha=\bigoplus_i\<D\dprime\>_i$ where the latter
are the irreducible components of $\<D\dprime\>$ contained in
$\<D\dprime\>_\alpha$. By proposition \ref{pr:reduced connection},
$$t_{\ol{\alpha}}=t_{\<D\dprime\>_\alpha}-t_{\<D\dprime\>^\alpha}$$
For any subdiagram $B\subset D$, with corresponding subalgebra
$\g_B\subset\g$, one has
$$\ik(t_{\<B\>})=
h\sum_{\alpha\in\Phi_+\cap\<B\>}\kappa_\alpha=
h\left(C_{\g_B}-\sum_j t_j^2\right)$$
where $C_{\g_B}$ is the Casimir operator of $\g_B$ relative to the restriction of
$(\cdot,\cdot)$ to it and $t_j$ is an orthonormal basis of the Cartan subalgebra
$\g_B\cap\h=\<B\>$ of $\g_B$. Thus,
$$\ik(t_{\ol{\alpha}})=
h\left(C_{\g_{\<D\dprime\>_\alpha}}-C_{\g_{\<D\dprime\>^\alpha}}
-t^2\right)$$
where $t$ is a vector in $\<D\dprime\>_\alpha$ which is orthogonal to
$\<D\dprime\>^\alpha$ and of norm one so that $\ik(t_{\ol{\alpha}})$
is invariant under $\g_{D\dprime}$ as claimed. To prove that the
$S_i^C$ and $\Phi^C_{\G\F}$ endow $\Ug\fmll$ with a \qca structure,
it suffices to show that the braid relations hold. This follows easily
from the fact that $S_i^C=S_i^\kappa\cdot\exp(\pi\sqrt{-1}\cdot
h_{\alpha_i}^2/2)$ and that the elements $\Phi^\kappa_{\G\F}$
are of weight $0$.\\

(iii) Consider the multivalued function $\Theta:\hreg\longrightarrow\Ug\fmll$
given by
$$\Theta=
\prod_{\alpha\in\Phi_+}\alpha^{h\frac{\<\alpha,\alpha\>}{4}h_\alpha^2}$$
One readily checks that $\Psi$ satisfies $\nablak\Psi=0$ if, and only if, $\nablac(\Psi\Theta)=0$. The claim follows easily from this and from the
discussion in \S \ref{sss:relevance}.\\

(iv) This again follows by transport of structure from proposition \ref{pr:universal
diff}, except for the statement that the explicit twisting element $\ik(c_{(D';\alpha
^{D'}_\F)}^{-R^\F_{<D'>}})$ of \eqref{eq:explicit twisting} commutes with $U\g_{D'
\setminus\alpha^{D'}_\F}$. This however follows from the already noted fact that
$$\ik(R^\F_{<D'>})=h(C_{\g_{D'}}-C_{\g_{D'\setminus\alpha^{D'}_\F}}-t^2)$$
where $t$ is a unitary vector in $\<D'\>$ which is orthogonal to
$\<D'\setminus\alpha^{D'}_\F\>$ \halmos

\subsection{From \qbas to \qcas}\label{ss:qba examples}

We show below that the commutant in the $n$--fold tensor product $A^{\otimes n}$
of a \qba $A$ with a coassociative coproduct, or of a \qtqba $A$ with a cocommutative
and coassociative coproduct, has the structure of a \qca of type $\sfA_{n-1}$. In the
latter case, the corresponding \qc representations of Artin's braid group $B_n$ coincide
with the $R$--matrix representations obtained from $A$. This construction abstracts
the author's duality between the \KZ connection for $\gl{k}$ and the Casimir connection
$\nablac$ for $\gl{n}$ \cite[\S 3]{TL1}.

\subsubsection{}\label{sss:qba recap}

Recall \cite[\S 1]{Dr3} that a \qba $(A,\Delta,\veps,\Phi)$ is an algebra
$A$ endowed with algebra homomorphisms $\Delta:A\rightarrow A^{\otimes
2}$ and $\veps:A\rightarrow k$ called the coproduct and counit, and an
invertible element $\Phi\in A^{\otimes 3}$ called the associator which
satisfy, for any $a\in A$
\begin{align}
\id\otimes\Delta(\Delta(a))
&=
\Phi\cdot\Delta\otimes\id(\Delta(a))\cdot\Phi^{-1}
\label{eq:coass}
\\
\id^{\otimes 2}\otimes\Delta(\Phi)\cdot
\Delta\otimes\id^{\otimes 2}(\Phi)
&=
1\otimes\Phi\cdot
\id\otimes\Delta\otimes\id(\Phi)\cdot
\Phi\otimes 1\\
\veps\otimes\id\circ\Delta
&=\id\\
\id\otimes\veps\circ\Delta
&=\id\\
\id\otimes\veps\otimes\id(\Phi)
&=1
\end{align}

Recall also that a twist of a \qba $A$ is an invertible element $F\in A^
{\otimes 2}$ satisfying $$\veps\otimes\id(F)=1=\id\otimes\veps(F)$$
Given such an $F$, the twisting of $A$ by $F$ is the \qba $(A,\Delta_F,
\veps,\Phi_F)$ where the coproduct $\Delta_F$ and associator $\Phi
_F$ are given by
\begin{align*}
\Delta_F(a)&=F\cdot\Delta(a)\cdot F^{-1}\\
\Phi_F
&=
1\otimes F\cdot\id\otimes\Delta(F)\cdot
\Phi\cdot
\Delta\otimes\id(F^{-1})\cdot F^{-1}\otimes 1
\end{align*}

A na\"{\i}ve, or strict morphism $\Psi:A\rightarrow A'$ of \qbas is an
algebra homomorphism satisfying
$$\veps=\veps'\circ\Psi,\qquad
\Psi^{\otimes 2}\circ\Delta=\Delta'\circ\Psi
\qquad\text{and}\qquad
\Psi^{\otimes 3}(\Phi)=\Phi'$$
A morphism $A\rightarrow A'$ of \qbas is a pair $(\Psi,F')$ where
$F'$ is a twist of $A'$ and $\Psi$ is a na\"{\i}ve morphism of $A$
to the twisting of $A'$ by $F'$.

\subsubsection{}\label{sss:qba abstract nonsense}

Let $A$ be a quasibialgebra. For any $n\geq 1$, let $\Br_n$ be
the set of complete bracketings on the non--associative monomial
$x_1\cdots x_n$. If $n\geq 2$, we require that such bracketings
contain the parentheses $(x_1\cdots x_n)$. Define, for any $\B\in
\Br_n$, a homomorphism $\Delta_\B :A\rightarrow A^{\otimes n}$
in the following way. If $n=1$, set $\Delta_{x_1}=\id$. Otherwise,
let
\begin{equation}\label{eq:B i}
1\leq i_1<i_2-1<\cdots<i_k-(k-1)\leq n-k
\end{equation}
be the indices $i$ such that $\B$ contains the bracket $(x_ix_{i+1})$.
Let $\ol{\B}$ be the bracketing on $\ol{x}_1\cdots\ol{x}_{n-k}$ obtained
from $\B$ by performing the substitutions
\begin{equation}\label{eq:B substitution}
x_\ell\longrightarrow\left\{\begin{array}{cl}
\ol{x}_\ell&\text{if $1\leq\ell<i_1$}\\
\ol{x}_{\ell-(j-1)}&\text{if $i_j+1<\ell<i_{j+1}$}\\
\ol{x}_{\ell-(k-1)}&\text{if $i_k+1<\ell\leq n$}
\end{array}\right.
\quad\text{and}\quad
(x_{i_j}x_{i_j+1})\longrightarrow \ol{x}_{i_j}
\end{equation}
and set
\begin{equation*}
\Delta_\B=
\id\otimesnp{i_1-1}\otimes\Delta\otimes
\id\otimesnp{i_2-i_1-2}
\odots{\otimes}
\id\otimesnp{i_k-i_{k-1}-2}\otimes\Delta\otimes\id\otimesnp{n-1-i_k}
\circ\Delta_{\ol{\B}}
\end{equation*}

Let now $\B,\B'\in\Br_n$ be two bracketings differing by the change of
one pair of parentheses. To such a data one associates an element $
\Phi_{\B'\B}$ of $A^{\otimes n}$ such that, for any $a\in A$,
$$\Phi_{\B'\B}\cdot\Delta_\B(a)=\Delta_{\B'}(a)\cdot\Phi_{\B'\B}$$
in the following way. Up to a permutation, $\B$ and $\B'$ differ by the
replacement of a monomial $((\B_1\B_2)\B_3)$ by $(\B_1(\B_2\B_3)$
where $\B_1,\B_2,\B_3$ are parenthesised monomials in the variables
$x_i\cdots x_{i+n_1-1}$, $x_{i+n_1}\cdots x_{i+n_1+n_2-1}$ and
$x_{i+n_1+n_2}\cdots x_{i+n_1+n_2+n_3-1}$ respectively. Set
\begin{equation}\label{eq:Phi_B}
\Phi_{\B'\B}=
1\otimesnp{i-1}\otimes
\Delta_{\B_1}\otimes\Delta_{\B_2}\otimes\Delta_{\B_3}(\Phi)\otimes
1\otimesnp{n-(i+n_1+n_2+n_3)+1}
\end{equation}

For any bracketing $\B\in\Br_n$, and twist $F$ of $A$, one also defines
an element $F_\B\in A^{\otimes n}$ such that
\begin{equation}\label{eq:Delta_B F_B}
(\Delta_F)_{\B}=\Ad(F_\B)\circ\Delta_\B
\end{equation}
as follows. If $n=1$, $F_{x_1}=1$. Otherwise, let $1\leq i_1\odots{<}
i_k-(k-1)\leq n-k$ and $\ol{\B}$ be as in \eqref{eq:B i} and \eqref
{eq:B substitution} respectively and set
\begin{equation}\label{eq:F_B}
\begin{split}
F_\B
&=
1\otimesnp{i_1-1}\otimes F\otimes 1\otimesnp{i_2-i_1-2}
\odots{\otimes}
1\otimesnp{i_k-i_{k-1}-2}\otimes F\otimes 1\otimesnp{n-1-i_k}\\
&\phantom{=}\cdot
\id\otimesnp{i_1-1}\otimes\Delta\otimes\id\otimesnp{i_2-i_1-2}
\odots{\otimes}
\id\otimesnp{i_k-i_{k-1}-2}\otimes\Delta\otimes\id\otimesnp{n-1-i_k}
(F_{\ol{\B}})
\end{split}
\end{equation}
Then, one also has
\begin{equation}\label{eq:F_B Phi_B}
(\Phi_F)_{\B'\B}=F_{\B'}\cdot\Phi_{\B'\B}\cdot F_\B^{-1}
\end{equation}

\subsubsection{}\label{sss:functor qba qca}

Assume henceforth that the coproduct $\Delta$ of $A$ is coassociative,
\ie satisfies $\Delta\otimes\id\circ\Delta=\id\otimes\Delta\circ\Delta$, so
that $(A,\Delta,\veps)$ is a bialgebra. Then, $\Delta_\B=\Delta_{\B'}$
for any $\B,\B'\in\Br_n$. Denote their common value by $\Delta^{(n)}$
and set
$$(A\otimesn{n})^A=
\{\alpha\in A\otimesn{n}|\medspace [\alpha,\Delta^{(n)}(a)]=0
\thickspace\thickspace\text{for any $a\in A$}\}$$

Fix $n\in\IN$, with $n\geq 2$ and let $D$ be the Dynkin diagram of
type $\sfA_{n-1}$. We wish to define a \qca structure of type $D$ on
$(A\otimesn{n})^A$. Identify for this purpose $D$ with the interval
$[1,n-1]$, its connected subdiagrams with subintervals $[i,j]\subseteq
[1,n-1]$ with integral endpoints and \mnss on $D$ with elements
in $\Br_n$ by attaching to $D'=[i,j]$ the bracket $x_1\cdots x_{i-1}
(x_i\cdots x_{j+1})x_{j+2}\cdots x_n$ as in remark \rkref{rk:Kn}.

\begin{theorem}\label{th:qba to qca}
Let $A$ be a \qba with coassociative coproduct and set $\T_n(A)=
(A\otimesn{n})^A$. Then,
\begin{enumerate}
\item The assignement
$$[i,j]\longrightarrow \T_n(A)_{[i,j]}=
1\otimesnp{i-1}\otimes (A\otimesnp{j-i+2})^A\otimes 1\otimesnp{n-1-j}$$
endows $\T_n(A)$ with the structure of a $D$--algebra.
\item If the edges of $D$ are each given an infinite multiplicity, the elements
$$S_i^{\T_n(A)}=1\otimesn{n},\thickspace i=1,\ldots,n-1
\qquad\text{and}\qquad
\Phi_{\B'\B}$$
give $\T_n(A)$ the structure of a \qca of type $D$.
\item If $\Psi:A\rightarrow A'$ is a na\"{\i}ve epimorphism of quasibialgebras,
then $\Psi_\B=\Psi\otimesn{n}$, $\B\in\Br_n$, is a morphism of \qcas
$\T_n(A)\rightarrow\T_n(A')$.
\item If $F$ is a twist of $A$ such that $\Delta_F$ is coassociative\footnote
{this is equivalent to the requirement that $\id\otimes\Delta(F^{-1})\cdot 1
\otimes F^{-1}\cdot F\otimes 1\cdot\Delta\otimes\id(F)\in(A^{\otimes 3})^A$.},
then $\Psi_\B=\Ad(F_\B)$, $\B\in\Br_n$, is an isomorphism of \qcas $\T_n(A)\rightarrow\T_n(A_F)$.
\item If, in addition, $F$ is an invariant twist, then $\{\Ad(F_\B)\}_{\B\in\Br_n}$
is the isomorphism induced by a canonical twist of $\T_n(A)$.
\end{enumerate}
\end{theorem}

The following is an immediate consequence of theorem \ref{th:qba to qca}
\begin{corollary}\label{co:qba to qca}
Let $\Q$ be the category of \qbas with coassociative coproducts and
morphisms $(\Psi,F')$ with $\Psi$ surjective. Then, the assignements
$A\rightarrow\T_n(A)$ and $(\Psi,F')\rightarrow\{\Ad(F'_\B)\circ\Psi
\otimesn{n}\}_{\B\in\Br_n}$ define a functor from $\Q$ to the category
of \qcas of type $\sfA_{n-1}$.
\end{corollary}

\remark Note that $\T_n(A)_{[i,j]}$ is not generated by the subalgebras
$\T_n(A)_{k}$ corresponding to the vertices $k$ of $D'=[i,j]$.

\subsubsection{}

{\sc Proof of theorem \ref{th:qba to qca}}. (i) The identity $$\Delta^{(n)}=
\id\otimesnp{i-1}\otimes\Delta^{(j-i+2)}\otimes\id\otimesnp{n-1-j}\circ\Delta
^{(n-1-(j-i))}$$ shows that $\T_n(A)_{[i,j]}$ is a subalgebra of $\T_n(A)$
and that $\T_n(A)_{[i,j]}\subseteq\T_n(A)_{[i',j']}$ if $[i,j]\subseteq[i',j']$.
It is moreover clear that $[\T_n(A)_{[i,j]},\T_n(A)_{[i',j']}]=0$ if either $j<
i'-1$ or $j'<i-1$.\\
(ii) By definition, the associators $\Phi_{\B'\B}$ satisfy $\Phi_{\B\B'}=\Phi
_{\B'\B}^{-1}$. When $\Delta$ is coassociative, \eqref{eq:Phi_B} reduces
to
\begin{equation}\label{eq:new Phi_B}
\Phi_{\B'\B}=
1\otimesnp{i-1}\otimes
\Delta^{(n_1)}\otimes\Delta^{(n_2)}\otimes\Delta^{(n_3)}(\Phi)\otimes
1\otimesnp{n-(i+n_1+n_2+n_3)+1}
\end{equation}
This element lies in $1\otimesnp{i-1}\otimes(A\otimesnp{n_1+n_2+n_3})^A
\otimes 1\otimesnp{n-(i+n_1+n_2+n_3)+1}$, and therefore in $\T_n(A)$,
because of the identity
$$\Delta^{(n_1+n_2+n+3)}=
\Delta^{(n_1)}\otimes\Delta^{(n_2)}\otimes\Delta^{(n_3)}\circ\Delta^{(3)}$$
and of the fact that, when $\Delta$ is coassociative, $\Phi$ lies in $(A\otimesn
{3})^A$. \eqref{eq:new Phi_B} also shows that $\Phi_{\B'\B}$ satisfies the
forgetfulness axiom as well as the support axiom since it commutes with
\begin{multline*}
\T_n(A)_
{[i,i+n_1+n_2+n_3-1]\setminus\alpha^{[i,i+n_1+n_2+n_3-1]}_{\B\cap\B'}}\\
=
1\otimesnp{i-1}\otimes(A\otimesn{n_1})^A\otimes(A\otimesn{n_2})^A\otimes (A\otimesn{n_3})^A\otimes1\otimesnp{n-(i+n_1+n_2+n_3)+1}
\end{multline*}
The coherence axiom follows from MacLane's coherence theorem and
the braid relations are void in this case since $m_{ij}=\infty$ for any $i,j$.\\
(iii) The surjectivity of $\Psi$ guarantees that $\Psi^{\otimes n}$ maps
$(A^{\otimes n})^{A}$ to $(A'^{\otimes n})^{A'}$, the rest of the claim
is clear.\\
(iv) When $\Delta$ is coassociative, the relation \eqref{eq:Delta_B F_B}
reduces to $\Delta_F^{(n)}=\Ad(F_\B)\circ\Delta^{(n)}$ and shows that
$\Ad(F_\B)$ maps $\T_n(A)$ to $\T_n(A_F)$. Set now, for $1\leq i\leq
k\leq j\leq n-1$,
\begin{equation}\label{eq:F ijk}
\begin{split}
F_{([i,j];k)}
&=
1\otimesnp{i-1}\otimes
\Deltan{k-i+1}\otimes\Deltan{j-k+1}(F)
\otimes 1\otimesnp{n-1-j}\\
&\in
1\otimesnp{i-1}\otimes
A\otimesn{(j-i+2)}\otimes 1\otimesnp{n-1-j}
\end{split}
\end{equation}
An induction based on \eqref{eq:F_B} shows that
\begin{equation}\label{eq:factored F_B}
F_\B=
\stackrel{\longrightarrow}{\prod_{[i,j]\in\B}}F_{([i,j];\alpha^{[i,j]}_\B)}
\end{equation}
where the product is taken with $F_{([i,j];\alpha^{[i,j]}_\B)}$ written to the
left of $F_{([i',j'];\alpha^{[i',j']}_\B)}$ whenever $[i,j]\subset[i',j']$\footnote
{this does not specify the order of the factors uniquely, but any two orders
satisfying this requirement are easily seen to yield the same result.}. Since
$F_{([i,j];k)}$ commutes with $\T_n(A)_{[i,k-1]\cup[k+1,j]}$ we get, for any
$[i,j]\in\B$
$$\Ad(F_\B)\left(\T_n(A)_{[i,j]}\right)=
\Ad(F_{\B_{[i,j]}})\left(\T_n(A)_{[i,j]}\right)=
\T_n(A_F)_{[i,j]}$$
where $\B_{[i,j]}$ is the bracketing on $x_i\cdots x_{j+1}$ induced by $\B$.
Thus $\{\Ad(F_\B)\}$ is a $D$--algebra morphism $\T_n(A)\rightarrow\T_n
(A_F)$ and therefore a morphism of \qcas by \eqref{eq:F_B Phi_B}.\\
(iv) Assume that $F\in (A^{\otimes 2})^A$, so that $\Delta_F=\Delta$. The
identity $\Delta\otimesnp{j-i+2}=\Delta\otimesnp{k-i+1}\otimes\Delta
\otimesnp{j-k+1}\circ\Delta$ shows that the element $F_{([i,j];k)}$ defined
by \eqref{eq:F ijk} is invariant under $A$. In other words,
$$a_{([i,j];k)}=F_{([i,j];k)}\in\T_n(A)^{[i,j]\setminus k}_{[i,j]}$$
is a twist of $\T_n(A)$ and $\T_n(A_F)$ is obtained from $\T_n(A)$ 
by twisting by $\{a_{([i,j];k)}\}$ \halmos

\subsubsection{}\label{sss:qtqba recap}

We give next a similar construction for \qt \qbas with a
coassociative and cocommutative coproduct. Recall first \cite[\S 3]{Dr3}
that a \qba $(A,\Delta,\veps,\Phi)$ is \qt if it is endowed with an invertible
element $R\in A^{\otimes 2}$ satisfying, for any $a\in A$,
\begin{align}
\Delta\op(a)&=R\cdot\Delta(a)\cdot R^{-1}
\label{eq:Delta R}\\
\Delta\otimes\id(R)&=
\Phi_{312}\cdot R_{13}\cdot\Phi_{132}^{-1}\cdot R_{23}\cdot\Phi_{123}
\label{eq:R12}\\
\id\otimes\Delta(R)&=
\Phi_{231}^{-1}\cdot R_{13}\cdot\Phi_{213}\cdot R_{12}\cdot\Phi_{123}^{-1}
\label{eq:R23}
\end{align}

A twist $F$ of a \qtqba $A$ is a twist of the underyling quasibialgebra.
The twisting of $A$ by $F$ is the \qtqba $(A,\Delta_F,\veps,\Phi_F,R_F)$
where 
$$R_F=F_{21}\cdot R\cdot F^{-1}$$
A morphism $(\Psi,F'):A\rightarrow A'$ of \qt \qbas is a morphism of the
underlying \qbas such that $\Phi^{\otimes 2}(R)=R'_{F'}$.

\subsubsection{}\label{ss:functor qtqba qc}
 
Let $A$ be a \qba with cocommutative and coassociative coproduct. Then
$\sigma\circ\Delta^{(n)}=\Delta^{(n)}$ for any $\sigma\in\SS_n$, so that
$(A\otimesn{n})^A$ is invariant under $\SS_n$. Set $\wt{\T}_n(A)=\T_n(A)\rtimes\SS_n=(A\otimesn{n})^{A}\rtimes\SS_n$ and,
for any $1\leq i\leq j\leq n-1$, let $\SS_{[i,j]}\subseteq\SS_n$ be the subgroup
generated by the transpositions $(k\medspace k+1)$, $k=i,\ldots,j$.

\begin{theorem}\label{th:qtqba to qca}\hfill
\begin{enumerate}
\item The assignement
$$\wt{\T}_n(A)_{[i,j]}=
1\otimesnp{i-1}\otimes (A\otimesnp{j-i+2})^A\otimes 1\otimesnp{n-1-j}
\rtimes\SS_{[i,j]}$$
endows $\wt{\T}_n(A)$ with the structure of a $D$--algebra.
\item If the edges of $D$ are each given the multiplicity $3$, the elements
$$S_i^{\wt{\T}_n(A)}
=(i\medspace i+1)\cdot 1\otimesnp{i-1}\otimes R\otimes1\otimesnp{n-i-1}$$
and the associators $\Phi_{\B'\B}$ give $\wt{\T}_n(A)$ the structure of a \qca
of type $D$.
\item If $(\Psi,F'):A\rightarrow A'$ is a surjective morphism of \qt quasibialgebras,
$\Psi_\B=\Ad(F'_\B)\circ\Psi\otimesn{n}$, $\B\in\Br_n$, is a \qca morphism $\wt{\T}_n(A)\rightarrow\wt{\T}_n(A')$.
\item If $F$ is an invariant twist of $A$, then $\wt{T}_n(A)$ and $\wt{T}_n
(A_F)$ differ by a twist.
\end{enumerate}
\end{theorem}
\proof (i) is clear. (ii) Note first that the cocommutativity of $\Delta$ and \eqref
{eq:Delta R} imply that $R\in(A\otimesn{2})^A$ so that $S^{\wt{\T}_n(A)}_i\in\wt
{\T}_n(A)_i$ for any $i=1,\ldots,n-1$. By theorem \ref{th:qba to qca}, the associators
$\Phi_{\B'\B}$ satisfy the orientation, forgetfulness and coherence axioms.
The support axioms follows from the fact that, by \eqref{eq:new Phi_B}, 
$\Phi_{\B'\B}$ commutes with
$\SS_{n_1}\times\SS_{n_2}\times\SS_{n_3}$
since $\Delta$ is cocommutive. The braid relations follow from the hexagon
relations \eqref {eq:R12}--\eqref{eq:R23} in the usual way.\\
(iii) The claim is obvious if $F'=1$. Assume now that $\Psi=\id_A$ and
that $F$ is a twist of $A$ such that $\Delta_F$ is cocommutative and coassociative\footnote{the cocommutativity of $\Delta_F$ is equivalent to the
requirement that $F_{21}^{-1}F\in(A^{\otimes 2})^A$.}. By \eqref{eq:Delta_B F_B},
$\Ad(F_\B)$ maps $\T_n(A)$ to $\T_n(A_F)$. Moreover, if $\sigma\in\SS_n$
and $a\in A$, then
\begin{equation*}
\begin{split}
F_\B\sigma F_\B^{-1}\sigma^{-1}\Delta_F^{(n)}(a)
&=
F_\B\sigma F_\B^{-1}\Delta_F^{(n)}(a)\sigma^{-1}
=
F_\B\sigma \Delta^{(n)}(a)F_\B^{-1}\sigma^{-1}\\
&=
F_\B \Delta^{(n)}(a)\sigma F_\B^{-1}\sigma^{-1}
=
\Delta_F^{(n)}(a)F_\B \sigma F_\B^{-1}\sigma^{-1}
\end{split}
\end{equation*}
where the first and third equalities follow from the cocommutativity of $\Delta_F$
and $\Delta$ respectively and the second and fourth from \eqref{eq:Delta_B F_B}. It
follows from this that $\Ad(F_\B)$ maps $\wt{\T}_n(A)$ to $\wt{\T}_n(A_F)$.
If $[i,j]\in\B$, the factorisation \eqref{eq:factored F_B} together with the
cocommutativity of $\Delta$ readily imply that
\begin{equation}\label{eq:local F_B}
\Ad(F_\B)\left(\wt{\T}_n(A)_{[i,j]}\right)=
\Ad(F_{\B_{[i,j]}})\left(\wt{\T}_n(A)_{[i,j]}\right)
\end{equation}
Thus, $\Ad(F_\B)$ defines a $D$--algebra morphism $\wt{\T}_n(A)\rightarrow
\wt{\T}_n(A_F)$ and therefore a morphism of \qcas by \eqref{eq:F_B Phi_B}
since, by \eqref{eq:local F_B}, for any $i=1,\ldots,n-1$ such that $[i,i]\in\B$,
\begin{equation*}
\begin{split}
\Ad(F_\B)S_i^{\wt{\T}_n(A)}
&=
\Ad(F_{([i,i];i)})S_i^{\wt{\T}_n(A)}\\
&=
(i\medspace i+1) 1\otimesn{i-1}\otimes F_{21}RF^{-1}\otimes 1\otimesn{n-1-i}\\
&=
S_i^{\wt{\T}_n(A_F)}
\end{split}
\end{equation*}
(iv) follows as in the proof of theorem \ref{th:qba to qca} \halmos

\begin{corollary}
Let $\wt{\Q}$ be the category of \qt \qbas with cocommutative and coassociative
coproduct and morphisms $(\Psi,F')$ where $\Psi$ is surjective. Then, for any
$n\geq 2$,
\begin{enumerate}
\item The assignement $A\rightarrow\wt{\T}_n(A)$ is a functor from $\wt{\Q}$
to the category of \qcas of type $\sfA_{n-1}$.
\item The $R$--matrix representation of Artin's braid group $B_n$
corresponding to an $A$--module $V$ and a bracketing $\B\in\Br_n$
coincides with the \qca representation $\pi_\B$ of $B_n$ on the $\wt
{\T}_n(A)$--module $V^{\otimes n}$.
\end{enumerate}
\end{corollary}


%% file: dynkin.tex

\section{The Dynkin complex and deformations of \qc algebras}\label{se:dynkin}

Let $D$ be a connected diagram and $A$ a $D$--algebra. We define in \S
\ref{ss:Dynkin complex} the {\it Dynkin complex} of $A$ and study some of its
elementary properties in \S\ref{ss:elementary}--\S\ref{ss:Dynkin Hochschild}.
Its main property, which we establish in \S \ref{ss:Dynkin controls qc}, is that
it controls the deformation theory of \qca structures on $A$. This is
obtained by showing in \S\ref{ss:Dynkin cellular} that, in degrees greater
or equal to 2, the Dynkin complex of $A$ embeds into the cellular cochain
complex of the \DCP associahedron $\AD$. In turn, this embedding is
obtained from an explicit presentation of the cellular chain complex of
$\AD$ in terms of the poset $\ND$ of nested sets on $D$ which is
described in \S \ref{ss:Hasse complex} and \S \ref{ss:cellular complex}. 

\subsection{The Dynkin complex of $A$}\label{ss:Dynkin complex}

We begin by defining the category of coefficients of the Dynkin complex of $A$.

\begin{definition}\label{de:D-bimod}
A $D$--bimodule over $A$ is an $A$--bimodule $M$, with left and right actions
denoted by $am$ and $ma$ respectively, endowed with a family of subspaces
$M_{D_1}$ indexed by the connected subdiagrams $D_1\subseteq D$ of $D$
such that the following properties hold
\begin{itemize}
\item for any $D_1\subseteq D$,
$$A_{D_1}\msp M_{D_1}\subseteq M_{D_1}
\qquad\text{and}\qquad
M_{D_1}\msp A_{D_1}\subseteq M_{D_1}$$
\item For any pair $D_2\subseteq D_1\subseteq D$,
$$M_{D_2}\subseteq M_{D_1}$$
\item For any pair of orthogonal subdiagrams $D_1,
D_2$ of $D$, $a_{D_1}\in A_{D_1}$ and $m_{D_2}\in
M_{D_2}$,
$$a_{D_1}m_{D_2}=m_{D_2}a_{D_1}$$
\end{itemize}
A morphism of $D$--bimodules $M,N$ over $A$ is an
$A$--bimodule map $T:M\rightarrow N$ such that $T(M_D)
\subseteq N_D$ for all $D\subseteq D$.
\end{definition}

Clearly, $A$ is a $D$--bimodule over itself. We denote by $\bimodD(A)$ the
abelian subcategory of $\bimod(A)$ consisting of $D$--bimodules over $A$.
If $M\in\bimodD(A)$, and $D_1\subseteq D$ is a subdiagram, we set
$$M^{D_1}=\{m\in M|\medspace a m=ma\quad
        \text{for any $a\in A_{D_1^i}$}\}$$
where $D_1^i$ are the connected components of $D_1$. In particular, if $D_1,D_2\subseteq D$ are orthogonal, and $D_1$ is connected, then $M_{D_1}\subseteq M^{D_2}$.\\

Let $M\in\bimodD(A)$. For any integer $0\leq p\leq n=|D|$, set
$$C^{p}(A;M)=
\bigoplus_{\substack{\aalpha\subseteq D_1\subseteq D,\\|\aalpha|=p}}
M_{D_1}^{D_1\setminus\aalpha}$$
where the sum ranges over all connected subdiagrams $D_1$ of $D$ and
ordered subsets $\aalpha=\{\alpha_1,\ldots,\alpha_p\}$ of cardinality $p$ of
$D_1$ and $M_{D_1}^{D_1\setminus\aalpha}=(M_{D_1})^{D_1\setminus
\aalpha}$. We denote the component of $m\in C^{p}(A;M)$ along $M_{D_1}^
{D_1\setminus\aalpha}$ by $m_{(D_1;\aalpha)}$.

\begin{definition} The group of Dynkin $p$--cochains on $A$ with coefficients
in $M$ is the subspace $CD^{p}(A;M)\subset C^{p}(A;M)$ of elements $m$
such that $$m_{(D_1;\sigma\aalpha)}=(-1)^{\sigma}m_{(D_1;\aalpha)}$$
where, for any $\sigma\in\SS_{p}$, $\sigma\{\alpha_1,\ldots,\alpha_p\}=
\{\alpha_{\sigma(1)},\ldots,\alpha_{\sigma(p)}\}$.
\end{definition}

Note that
$$
CD^{0}(A;A)=\bigoplus_{D_1\subseteq D}Z(A_{D_1})
\aand
CD^{n}(A;M)\cong M_D$$
For $1\leq p\leq n-1$, define a map $\dD^{p}:C^{p}(A;M)\rightarrow
C^{p+1}(A;M)$ by
$$\dD m_{(D_1;\aalpha)}=
\sum_{i=1}^{p+1}(-1)^{i-1}\left(
m_{(D_1;\aalpha\setminus\root{i})}-
m_{(\con{D_1}{\root{i}}{\aalpha\setminus\root{i}};
\aalpha\setminus\root{i})}
\right)$$
where $\aalpha=\{\root{1},\ldots,\root{p+1}\}$,
$\con{D_1}{\root{i}}{\aalpha\setminus\root{i_j}}$ is
the connected component of $D_1\setminus\root{i}$
containing $\aalpha\setminus\root{i}$ if such a component
exists and the empty set otherwise, and we set $m_{
(\emptyset;-)}=0$. For $p=0$, define $\dD^0:
C^{0}(A;M)\rightarrow C^{1}(A;M)$ by
$$\dD^0 m_{(D_1;\alpha_i)}=m_{D_1}-m_{D_1\setminus\root{i}}$$
where $m_{D_1\setminus\alpha_i}$ is the sum of $m_{D_2}$
with $D_2$ ranging over the connected components of
$D_1\setminus\root{i}$. Finally, set $\dD^n=0$. It is easy to see that
the Dynkin differential $\dD$ is well--defined and that it leaves $CD^
{*}(A;M)$ invariant.

\begin{theorem}\label{th:Dynkin complex}
$(CD^{*}(A;M),\dD)$ is a complex. The cohomology groups $$HD^{p}
(A;M)=\Ker(\dD^{p})/\Im(\dD^{p-1})$$
for $p=0,\ldots,|D|$ are called the Dynkin diagram cohomology groups
of $A$ with coefficients in $M$.
\end{theorem}
\proof
Assume first that $m$ is a zero--cochain. Let $B\subseteq D$
be a connected diagram, $\root{1},\root{2}$ two distinct
elements of $D$ and set
$$B_1=\con{B}{\root{2}}{\root{1}}
\qquad\text{and}\qquad
B_2=\con{B}{\root{1}}{\root{2}}$$
Then,
\begin{equation*}
\begin{split}
\dD^2\msp m_{(B;\root{1},\root{2})}
&=
 \dD m_{(B;\root{2})}
-\dD m_{(B_2;\root{2})}
-\dD m_{(B;\root{1})}
+\dD m_{(B_1;\root{1})}\\
&=
 m_B-m_{B\setminus\root{2}}
-m_{B_2}+m_{B_2\setminus\root{2}}
-m_B+m_{B\setminus\root{1}}
+m_{B_1}-m_{B_1\setminus\root{1}}\\
&=
 (m_{B\setminus\root{1}}-m_{B_2}+m_{B_2\setminus\root{2}})
-(m_{B\setminus\root{2}}-m_{B_1}+m_{B_1\setminus\root{1}})\\
&=
m_{B\setminus\{\root{1},\root{2}\}}-
m_{B\setminus\{\root{1},\root{2}\}}\\
&=
0
\end{split}
\end{equation*}

To treat the general case, we shall need the following

\begin{lemma}\label{le:cc}
Let $\aalpha=\{\root{1},\ldots,\root{k}\}\subseteq B$
be a subset of cardinality $k\geq 3$. For any $1\leq
i\neq j\leq k$, set
$$B_{ij}=
\con{\con{B}{\root{i}}{\aalpha\setminus\root{i}}}
{\root{j}}{\aalpha\setminus\{\root{i},\root{j}\}}$$
Then,
\begin{enumerate}
\item if $\con{B}{\root{j}}{\aalpha\setminus\root{j}}
=\emptyset$, then,
$B_{ij}=
\con{B}{\root{j}}{\aalpha\setminus\{\root{i},\root{j}\}}$.
\item If $\con{B}{\root{i}}{\aalpha\setminus\root{i}}
\neq\emptyset$, then,
$B_{ij}=
\con{B}{\{\root{i},\root{j}\}}
{\aalpha\setminus\{\root{i},\root{j}\}}$.
\end{enumerate}
\end{lemma}
\proof
(i) The \lhs is contained in the \rhs since, when it
is not empty, it is connected, contained in $B\setminus
\root{j}$ and contains $\aalpha\setminus\{\root{i},
\root{j}\}\neq\emptyset$. Similarly, the \rhs is
contained in $\con{B}{\root{i}}{\aalpha\setminus
\root{i}}$ since, when non--empty, it is connected,
does not contain $\root{i}$ by assumption, and
contains $\aalpha\setminus\{\root{i},\root{j}\}$. 
Since the \rhs does not contain $\root{j}$, it is
therefore contained in the left--hand side. (ii) is
proved in a similar way \halmos\\

Write $\dD=d_1-d_2$ where, for $\aalpha=\{\root{1},
\ldots,\root{k}\}$,
$${d_1 m}_{(B;\aalpha)}=
\sum_{i=1}^{k}(-1)^{i-1}m_{(B;\aalpha\setminus\root{i})}
\qquad
{d_2 m}_{(B;\aalpha)}=
\sum_{i=1}^{k}(-1)^{i-1}
m_{(\con{B}{\root{i}}{\aalpha\setminus\root{i}};
\aalpha\setminus\root{i})}
$$
Note that
$$d_1^2\msp m_{(B;\aalpha)}=
\sum_{1\leq i\neq j\leq k}(-1)^{i+j}\sign(i-j)\msp
m_{(B;\aalpha\setminus\{\root{i},\root{j}\})}=
-d_1^2\msp m_{(B;\aalpha)}
$$
whence $d_1^2=0$. We next have
\begin{equation*}
\begin{split}
d_1d_2\msp m_{(B;\aalpha)}
&=
\phantom{\leq}
\sum_{i=1}^{k}(-1)^{i-1}
{d_2 a}_{(B;\aalpha\setminus\root{i})}\\
&=
\sum_{1\leq i\neq j\leq k}(-1)^{i+j}\sign(i-j)\msp
m_{(\con{B}{\root{j}}{\aalpha\setminus\{\root{i},\root{j}\}};
\aalpha\setminus\{\root{i},\root{j}\})}
\end{split}
\end{equation*}
and
\begin{equation*}
\begin{split}
d_2d_1\msp m_{(B;\aalpha)}
&=
\phantom{\leq}
\sum_{j=1}^{k}(-1)^{j-1}
{d_1 a}_{(\con{B}{\root{j}}{\aalpha\setminus\root{j}};
\aalpha\setminus\root{j})}\\
&=
\sum_{1\leq i\neq j\leq k}(-1)^{i+j}\sign(j-i)\msp
m_{(\con{B}{\root{j}}{\aalpha\setminus\root{j}};
\aalpha\setminus\{\root{i},\root{j}\})}
\end{split}
\end{equation*}
so that
\begin{multline*}
(d_1d_2+d_2d_1)\msp m_{(B;\aalpha)}\\
=
\sum_{1\leq i\neq j\leq k}(-1)^{i+j}\sign(i-j)\msp
\delta_{\con{B}{\root{j}}{\aalpha\setminus\root{j}}=\emptyset}\msp\cdot
m_{(\con{B}{\root{j}}{\aalpha\setminus\{\root{i},\root{j}\}};
\aalpha\setminus\{\root{i},\root{j}\})}
\end{multline*}
Finally,
\begin{equation*}
\begin{split}
d_2^2\msp m_{(B;\aalpha)}
&=
\phantom{\leq}
\sum_{i=1}^{k}(-1)^{i-1}
{d_2 a}_{(\con{B}{\root{i}}{\aalpha\setminus\root{i}};
\aalpha\setminus\root{i})}\\
&=
\sum_{1\leq i\neq j\leq k}(-1)^{i+j}\sign(i-j)\msp
m_{\left(\con{\con{B}{\root{i}}{\aalpha\setminus\root{i}}}
{\root{j}}{\aalpha\setminus\{\root{i},\root{j}\}};
\aalpha\setminus\{\root{i},\root{j}\}\right)}\\
&=
\sum_{1\leq i\neq j\leq k}(-1)^{i+j}\sign(i-j)\msp
\delta_{\con{B}{\root{j}}{\aalpha\setminus\root{j}}=\emptyset}\msp\cdot
m_{\left(\con{\con{B}{\root{i}}{\aalpha\setminus\root{i}}}
{\root{j}}{\aalpha\setminus\{\root{i},\root{j}\}};
\aalpha\setminus\{\root{i},\root{j}\}\right)}\\
&+
\sum_{1\leq i\neq j\leq k}(-1)^{i+j}\sign(i-j)\msp
\delta_{\con{B}{\root{j}}{\aalpha\setminus\root{j}}\neq\emptyset}\msp\cdot
m_{\left(\con{\con{B}{\root{i}}{\aalpha\setminus\root{i}}}
{\root{j}}{\aalpha\setminus\{\root{i},\root{j}\}};
\aalpha\setminus\{\root{i},\root{j}\}\right)}
\end{split}
\end{equation*}
We claim that the second summand is zero. Indeed,  it is equal to
\begin{multline*}
\sum_{1\leq i\neq j\leq k}(-1)^{i+j}\sign(i-j)\msp
\delta_{\con{B}{\root{j}}{\aalpha\setminus\root{j}}\neq\emptyset}\msp\cdot
\delta_{\con{B}{\root{i}}{\aalpha\setminus\root{i}}\neq\emptyset}\msp\cdot
m_{\left(\con{\con{B}{\root{i}}{\aalpha\setminus\root{i}}}
{\root{j}}{\aalpha\setminus\{\root{i},\root{j}\}};
\aalpha\setminus\{\root{i},\root{j}\}\right)}\\
=
\sum_{1\leq i\neq j\leq k}(-1)^{i+j}\sign(i-j)\msp
\delta_{\con{B}{\root{j}}{\aalpha\setminus\root{j}}\neq\emptyset}\msp\cdot
\delta_{\con{B}{\root{i}}{\aalpha\setminus\root{i}}\neq\emptyset}\msp\cdot
m_{(\conn{B}{\root{i},\root{j}}
{\aalpha\setminus\{\root{i},\root{j}\}};
\aalpha\setminus\{\root{i},\root{j}\})}
\end{multline*}
where we used (ii) of lemma \ref{le:cc}, and it therefore vanishes since
the summand is antisymmetric in the interchange $i\leftrightarrow j$.
Thus, by (i) of lemma \ref{le:cc},
\begin{equation*}
\begin{split}
d_2^2\msp m_{(B;\aalpha)}
&=
\sum_{1\leq i\neq j\leq k}(-1)^{i+j}\sign(i-j)\msp
\delta_{\con{B}{\root{j}}{\aalpha\setminus\root{j}}=\emptyset}\msp
m_{(\con{B}{\root{j}}{\aalpha\setminus\{\root{i},\root{j}\}};
\aalpha\setminus\{\root{i},\root{j}\})}\\
&=
(d_1d_2+d_2d_1)m_{(B;\aalpha)}
\end{split}
\end{equation*}
so that $\dD^2=d_1^2-(d_1d_2+d_2d_1)+d_2^2=0$ \halmos

\subsection{Elementary Properties of the Dynkin complex}\label{ss:elementary}

\subsubsection{Functoriality with respect to restriction to subdiagrams}

Let $M\in\bimodD(A)$, $B\subseteq D$ a connected subdiagram
and consider $M_B$ as an $B$--bimodule over $A_B$. Then, the
map $$p_{B,D}:CD^*(A;M)\longrightarrow CD^*(A_B;M_B)$$
given by $p_{B,D}m_{(B;\aalpha)}=m_{(B;\aalpha)}$ for any $\aalpha
\subseteq B\subseteq B$ is a chain map satisfying $$p_{C,B}\cdot
p_{B,D}=p_{C,D}$$ for any $C\subseteq B\subseteq D$. More
generally, the Dynkin complex is functorial \wrt na\"{\i}ve
morphisms of $D$--algebras but not \wrt general morphisms.

\subsubsection{Low--dimensional cohomology groups}

\begin{proposition}\label{pr:low dimensional}
Let $M$ be a $D$--bimodule over $A$. Then
\begin{enumerate}
\item $H^0(A;M)=0$.
\item If, for each connected $B\subseteq D$, the algebra $A_B$ is
generated by the subalgebras $A_\alpha$, $\alpha\in B$, the map
$$\bigoplus_{\alpha\in D}p_{\alpha,D}:H^1(A;M)\longrightarrow
\bigoplus_{\alpha\in D}H^1(A_\alpha;M_\alpha)$$
is injective.
\end{enumerate}
\end{proposition}
\proof (i) We need to prove that the differential $\dD^0$ is injective. Let $m=\{m_B\}
_{B\subseteq D}$ be a zero--cocycle and assume by induction that $m_B=0$ whenever
$|B|\leq k$. Let $B\subseteq D$ be a connected subdiagram of cardinality $k+1$ and
let $\alpha\in B$. Then $$0=dm_{(B;\alpha)}=m_B-m_{B\setminus\alpha}$$
implies that $m_B=0$ whence $m=0$.\\
(ii) Note that $H^1(A_\alpha;M_\alpha)=M_\alpha/M_\alpha^\alpha$. Let $m=\{m_
{(B;\alpha)}\}_{\alpha\in B\subseteq D}$ be a one--cocycle such that $m_{(\alpha;
\alpha)}\in M_\alpha^\alpha$ for any $\alpha\in D$. Replacing $m$ by $m-\dD n$,
where $n\in CD^0(A;M)$ is given by
$$n_B=\left\{\begin{array}{cl}
m_{(\alpha;\alpha)}&\text{if $B=\alpha$}\\
0&\text{if $|B|\geq 2$}
\end{array}\right.$$
we may assume that $m_{(\alpha;\alpha)}=0$ for any $\alpha\in D$. Assume
therefore that, up to the addition of a 1--coboundary, $m_{(B;\alpha)}=0$ for any
$\alpha\in B\subseteq D$ such that $|B|\leq k$. For any $B$ of cardinality $k+1$
and $\alpha\neq\beta\in B$, we have
$$0=\dD m_{(B;\alpha,\beta)}=
m_{(B;\beta)}-m_{(\con{B}{\alpha}{\beta};\beta)}-
m_{(B;\alpha)}+m_{(\con{B}{\beta}{\alpha};\alpha)}$$
whence, given that $|\con{B}{\alpha}{\beta}|\leq k$,
$$m_{(B;\alpha)}=m_{(B;\beta)}$$
This implies that $m_{(B;\alpha)}\in M_B^{B\setminus\alpha}\cap M_B
^{B\setminus\beta}$ which, by assumption is equal to $M_B^B$, and
that $m_{(B;\alpha)}$ is independent of $\alpha\in B$. Thus, replacing
$m$ by $m-\dD n$ where
$$n_B=\left\{\begin{array}{cl}
m_{(B;\alpha)}&\text{if $|B|=k+1$ and $\alpha\in B$}\\
0			  &\text{if $|B|\neq k+1$}
\end{array}\right.$$
we find that $m_{(B;\alpha)}=0$ whenever $|B|\leq k+1$ \halmos

\begin{corollary}
If $(W,S)$ is a Coxeter system, then $HD^i(k[W];k[W])=0$ for $i=0,1$.
\end{corollary}
\proof For any simple reflection $s_i\in S$, $H^1(A_{s_i};A_{s_i})=A_{s_i}/
A_{s_i}^{s_i}=0$ since $A_{s_i}\cong k[\IZ_2]$ is commutative. The result
now follows from proposition \ref{pr:low dimensional} \halmos

\subsection{Dynkin cohomology and Hochschild cohomology}
\label{ss:Dynkin Hochschild}

Let $n\geq 2$, $D$ the Dynkin diagram of type $\sfA_{n-1}$, $A$ a
\ba and $\T_n(A)=(A\otimesn{n})^A$ the $D$--algebra constructed
in \S \ref{sss:functor qba qca}. We relate below the cobar complex
of $A$ with the Dynkin complex of $\T_n(A)$.\\

Recall first (see \eg \cite[\S XVIII.5]{Ks}) that if $(A,\Delta,\veps)$ is a
coalgebra endowed with an element $1\in A$ such that $\Delta(1)=1
\otimes 1$, the cobar complex $(C^k(A),\dH)$ of $A$ is defined by
setting $C^k(A)=A\otimesn{k}$ for $k\geq 0$, $\dh=0$ in degree 0
and, for $a\in A\otimesn{k}$, $k\geq 1$,
$$\dh a=
1\otimes a
+\sum_{i=1}^k(-1)^i
\id^{\otimes (i-1)}\otimes\Delta\otimes\id^{\otimes (k-i)}(a)
+(-1)^{k+1}1\otimes a$$
If $A$ is a bialgebra with unit $1$, the subspaces $(A\otimesn{k})^A$
are readily seen to form a subcomplex of $C^k(A)$. For $k=1,\ldots,n$,
define a map $\phi_k:(A\otimesn{k})^A\rightarrow\CDD^{k-1}(\T_n(A);
\T_n(A))$ by
$$\phi_1 a_{([i,j])}=
\sum_{\ell=0}^{j-i+1}
1\otimesnp{i-1+\ell}\otimes a\otimes 1\otimesnp{n-i-\ell}-
1\otimesnp{i-1}\otimes\Delta^{(j-i+2)}(a)\otimes 1\otimesnp{n-1-j}$$
and, for $k\geq 2$,
\begin{multline*}
\phi_k a_{([i,j];\ell_1,\ldots,\ell_{k-1})}=(-1)^k\cdot
1\otimesnp{i-1}\otimes\Delta^{(\ell_1-i+1)}\otimes\Delta^{(\ell_2-\ell_1)}
\otimes\cdots\\
\phantom{==================}\cdots\otimes
\Delta^{(\ell_{k-1}-\ell_{k-2})}\otimes\Delta^{(j-\ell_k+1)}
(a)\otimes 1\otimesnp{n-i-j}
\end{multline*}
where $1\leq i\leq\ell_1\odots{<}\ell_{k-1}\leq j\leq n-1$.

\begin{proposition}
For any $k=1,\ldots,n-2$,
$$\dDD\circ\phi_k=\dH\circ\phi_{k+1}$$
so that $\phi$ is a degree $-1$ chain map from the invariant, truncated
cobar complex $((A\otimesn{k})^A,\dH)_{k=1}^{n}$ of $A$ to the Dynkin 
complex of $\T_n(A)$.
\end{proposition}
\proof This follows by a straightforward computation \halmos

\subsection{The chain complex $C_*(\ND)$}\label{ss:Hasse complex}

We define in this subsection the chain complex $C_*(\ND)$ of {\it oriented}
nested sets on $D$. We then show in \S \ref{ss:cellular complex} that $C_*
(\ND)$ is isomorphic to the cellular chain complex of the \DCP associahedron
$\AD$.

\begin{definition}
Let $\H\in\ND$ be a nested set with $|\H|<D$. An orientation $\eepsilon$
of $\H$ is a choice of
\begin{enumerate}
\item an enumeration $D_1,\ldots,D_m$ of the unsaturated elements of $\H$,
\item a total order on each $\aalpha^{D_i}_{\H}$.
\end{enumerate}
\end{definition}

By convention, a \mns has a unique orientation.

\begin{definition}\label{de:C(ND)}
Let $C_k(\ND)$ be the free $\IZ$--module generated by symbols $\H_\eepsilon$
where $\H$ is a nested set of dimension $k=0,\ldots,|D|$ and $\eepsilon$ is an
orientation of $\H$, modulo the following relations if $k\geq 1$
\begin{equation}\label{eq:or 1}
\H_{\eepsilon'}=
(-1)^{(|\aalpha^{D_i}_{\H}|-1)\cdot(|\aalpha^{D_{i+1}}_{\H}|-1)}\cdot\H_{\eepsilon}
\end{equation}
if $\eepsilon'$ is obtained from $\eepsilon$ by permuting the unsaturated elements
$D_i$, $D_{i+1}$ of $\H$ while leaving the total order on each $\aalpha^{D_j}_{\H}$
unchanged and
\begin{equation}\label{eq:or 2}
\H_{\eepsilon'}=(-1)^{\sigma_i}\cdot\H_{\eepsilon}
\end{equation}
if $\eepsilon'$ is obtained from $\eepsilon$ by changing the order on $\aalpha
^{D_i}_{\H}$ by a permutation $\sigma_i$.
\end{definition}

Let $\H_\eepsilon$ be an oriented nested set of positive dimension. Let $D_1,\ldots,D_m$
be the unsaturated elements of $\H$ and set $\aalpha_i=\aalpha^{D_i}_{\H}$. Let $
\G$ a boundary facet of $\H$, that is $\G\supset\H$ and $|\G\setminus\H|=1$. By
proposition \ref{pr:ND^H} and lemma \ref{le:explicit lift}, $\G=\H\cup D_{\bbeta_i}$
for a unique $i\in[1,m]$ and $\emptyset\neq\bbeta_i\subsetneq\aalpha_i$. The unsaturated elements of $\G$ are $D_j$, with $j\neq i$, and $D_{\bbeta_i}$,
$D_i$, provided $|\bbeta_i|\geq 2$ and $|\bbeta_i|\leq|\aalpha_i|-2$ respectively.
The corresponding subsets of vertices are
$$\aalpha^{D_j}_\G=\aalpha_j,\qquad
\aalpha^{D_{\bbeta_i}}_\G=\bbeta_i\qquad\text{and}\qquad
\aalpha^{D_i}_\G=\aalpha_i\setminus\bbeta_i$$

\begin{definition}
The orientation $\ol{\veps}$ of $\G=\H\cup D_{\bbeta_i}$ induced by
$\eepsilon$ is obtained by enumerating the unsaturated elements of
$\G$ as $$D_1,\ldots,D_{i-1},D_{\bbeta_i},D_i,D_{i+1},\ldots,D_{m},$$
ordering $\aalpha_j$, $j\neq i$ as prescribed by $\varepsilon$ and endowing $\bbeta_i$
and $\aalpha_i\setminus\bbeta_i$ with the restriction of the total order on $\aalpha_i$
prescribed by $\eepsilon$.
\end{definition}

Define the {\it shuffle number} $s(\bbeta_i;\aalpha_i)$ of $\bbeta_i$ with respect
to $\aalpha_i$ to be the number of elementary transpositions required to move
all elements of $\bbeta_i$ to the left of the first element of $\aalpha_i$. In other
words, if
$$\aalpha_i=\{\alpha_i^1,\ldots,\alpha_i^{n_i}\}
\qquad\text{and}\qquad
\bbeta_i=\{\alpha_i^{j_1},\ldots,\alpha_i^{j_p}\}$$
for some $1\leq j_1\odots{<}j_p\leq n_i$, then
$$s(\bbeta_i;\aalpha_i)=(j_1-1)+(j_2-2)+\cdots+(j_p-p)$$

\begin{proposition}
Let $\partial_k:C_k(\ND)\longrightarrow C_{k-1}(\ND)$ be the operator given by
$$\partial_k\H_\varepsilon=
\sum_{\substack{1\leq i\leq m,\\\emptyset\neq\bbeta_i\subsetneq\aalpha_i}}
(-1)^{(|\aalpha_1|-1)+\cdots+(|\aalpha_{i-1}|-1)}\cdot
(-1)^{|\bbeta_i|-1}\cdot
(-1)^{s(\bbeta_i;\aalpha_i)}\cdot
(H\cup D_{\bbeta_i})_{\ol{\eepsilon}}$$
if $k=1,\ldots,|D|$ and $\partial_0=0$. Then, $\partial_k$ is well--defined and
$\partial_{k-1}\circ\partial_k=0$.
\end{proposition}
\proof By proposition \ref{pr:ND^H} and lemma \ref{le:explicit lift}, the summands
which arise in writing $\partial^2\H_\eepsilon$ are of the form $\H\cup\{D_{\bbeta_i},
D_{\ggamma_j}\}$ with $D_{\bbeta_i}$ and $D_{\ggamma_j}$ compatible.
We must therefore prove that the sign contributions corresponding to the
two sides of the diamond
\begin{diagram}[height=2em]
					&		&\H		&		&						\\	
					&\ldTo	&		&\rdTo	&						\\
\H\cup\{D_{\bbeta_i}\}&		&		&		&\H\cup\{D_{\ggamma_j}\}\\
					&\rdTo	&		&\ldTo	&						\\
					&		&\H\cup\{D_{\bbeta_i},D_{\ggamma_j}\}	&&
\end{diagram}
are opposite to each other. We consider 
the various cases corresponding to the relative position of $\bbeta_i$ and $\ggamma_j$.

\subsubsection{$i\neq j$}

In this case, the orientations induced on $\H\cup\{D_{\bbeta_i},D_{\ggamma_j}\}$ by
the two sides of the diamond are the same. We may assume, up to a permutation of
$i$ and $j$ that $i<j$. The sign contribution of the left side is then
\begin{equation*}
\begin{split}
&
(-1)^{(|\aalpha_1|-1)\odots{+}(|\aalpha_{i-1}|-1)}
\cdot(-1)^{|\bbeta_i|-1}\cdot(-1)^{s(\bbeta_i;\aalpha_i)}\\
\cdot&
(-1)^{(|\aalpha_1|-1)\odots{+}(|\aalpha_{i-1}|-1)+
(|\bbeta_i|-1)+(|\aalpha_i\setminus\bbeta_i|-1)+
(|\aalpha_{i+1}|-1)\odots{+}(|\aalpha_{j-1}|-1)}\\
\cdot&
(-1)^{|\ggamma_j|-1}\cdot(-1)^{s(\ggamma_j;\aalpha_j)}\\
=&-
(-1)^{(|\aalpha_i|-1)\odots{+}(|\aalpha_{j-1}|-1)}\cdot
(-1)^{|\bbeta_i|+|\ggamma_j|}\cdot(-1)^{s(\bbeta_i;\aalpha_i)+s(\ggamma_j;\aalpha_j)}
\end{split}
\end{equation*}
while that of the right side is
\begin{equation*}
\begin{split}
&
(-1)^{(|\aalpha_1|-1)\odots{+}(|\aalpha_{j-1}|-1)}
\cdot(-1)^{|\ggamma_j|-1}\cdot(-1)^{s(\ggamma_j;\aalpha_j)}\\
\cdot&
(-1)^{(|\aalpha_1|-1)\odots{+}(|\aalpha_{i-1}|-1)}\cdot
(-1)^{|\bbeta_i|-1}\cdot(-1)^{s(\bbeta_i;\aalpha_i)}\\
=&
(-1)^{(|\aalpha_i|-1)\odots{+}(|\aalpha_{j-1}|-1)}\cdot
(-1)^{|\bbeta_i|+|\ggamma_j|}\cdot(-1)^{s(\bbeta_i;\aalpha_i)+s(\ggamma_j;\aalpha_j)}
\end{split}
\end{equation*}
as required.

\subsubsection{$i=j$ and $\bbeta_i\subset\ggamma_i$ or $\ggamma_i\subset\bbeta_i$}

Up to a permutation, we may assume that $\bbeta_i\subset\ggamma_i$. In this case
again, the orientations induced on $\H\cup\{D_{\bbeta_i},D_{\ggamma_j}\}$ by the two
sides of the diamond are the same. The sign contribution from the left side is
\begin{equation*}
\begin{split}
&
(-1)^{(|\aalpha_1|-1)\odots{+}(|\aalpha_{i-1}|-1)}\cdot(-1)^{|\bbeta_i|-1}\cdot
(-1)^{s(\bbeta_i;\aalpha_i)}\\
\cdot&
(-1)^{(|\aalpha_1|-1)\odots{+}(|\aalpha_{i-1}|-1)+(|\bbeta_i|-1)}
\cdot
(-1)^{|\ggamma_i\setminus\bbeta_i|-1}\cdot
(-1)^{s(\ggamma_i\setminus\bbeta_i;\aalpha_i\setminus\bbeta_i)}\\
=&
-(-1)^{|\bbeta_i|+|\ggamma_i|}\cdot
(-1)^{s(\bbeta_i;\aalpha_i)+s(\ggamma_i\setminus\bbeta_i;\aalpha_i\setminus\bbeta_i)}
\end{split}
\end{equation*}
while that of the right one is
\begin{equation*}
\begin{split}
&
(-1)^{(|\aalpha_1|-1)\odots{+}(|\aalpha_{i-1}|-1)}\cdot
(-1)^{|\ggamma_i|-1}\cdot(-1)^{s(\ggamma_i;\aalpha_i)}\\
\cdot&
(-1)^{(|\aalpha_1|-1)\odots{+}(|\aalpha_{i-1}|-1)}\cdot
(-1)^{|\bbeta_i|-1}\cdot
(-1)^{s(\bbeta_i;\ggamma_i)}\\
=&
(-1)^{|\bbeta_i|+|\ggamma_i|}\cdot
(-1)^{s(\bbeta_i;\ggamma_i)+s(\ggamma_i;\aalpha_i)}
\end{split}
\end{equation*}
These are opposite to each other in view of the following

\begin{lemma}
$(-1)^{s(\bbeta_i;\aalpha_i)+
s(\ggamma_i\setminus\bbeta_i;\aalpha_i\setminus\bbeta_i)}=
(-1)^{s(\bbeta_i;\ggamma_i)+s(\ggamma_i;\aalpha_i)}$
\end{lemma}
\proof The \lhs is the parity of the number of elementary transpositions
required to shuffle $\bbeta_i$ to the left of $\aalpha_i$ and then $\ggamma_i
\setminus\bbeta_i$ to the left of $\aalpha_i\setminus\bbeta_i$, thus arriving
at the ordered configuration $\bbeta_i,\ggamma_i\setminus\bbeta_i,\aalpha_i
\setminus\ggamma_i$. The \rhs on the other hand is the parity of the number
of transpositions needed to shuffle $\ggamma_i$ to the left of $\aalpha_i$
and then $\bbeta_i$ to the left of $\ggamma_i$ resulting in the very same
configuration \halmos

\subsubsection{$i=j$ and $\bbeta_i\cap\ggamma_i=\emptyset$}
	
In this last case, the sign contribution of the left side of the diamond is			
\begin{equation*}
\begin{split}
&
(-1)^{(|\aalpha_1|-1)\odots{+}(|\aalpha_{i-1}|-1)}\cdot
(-1)^{|\bbeta_i|-1}\cdot(-1)^{s(\bbeta_i;\aalpha_i)}\\
\cdot&
(-1)^{(|\aalpha_1|-1)\odots{+}(|\aalpha_{i-1}|-1)+(|\bbeta_i|-1)}
\cdot
(-1)^{|\ggamma_i|-1}\cdot
(-1)^{s(\ggamma_i;\aalpha_i\setminus\bbeta_i)}\\
=&
(-1)^{|\ggamma_i|-1}\cdot
(-1)^{s(\bbeta_i;\aalpha_i)+s(\ggamma_i;\aalpha_i\setminus\bbeta_i)}
\end{split}
\end{equation*}
so that, by symmetry, that of the right side is
$$(-1)^{|\bbeta_i|-1}\cdot
(-1)^{s(\ggamma_i;\aalpha_i)+s(\bbeta_i;\aalpha_i\setminus\ggamma_i)}$$
In this case however, the orientation induced on $\H\cup\{D_{\bbeta_i},
D_{\ggamma_i}\}$ by each side may differ. Indeed, the left side leads to
enumerating the unsaturated elements of $\H\cup\{D_{\bbeta_i},
D_{\ggamma_i}\}$ as
$$D_1,\ldots,D_{i-1},D_{\bbeta_i},D_{\ggamma_i},D_i,\ldots,D_m$$
while the right branch leads to enumerating them as
$$D_1,\ldots,D_{i-1},D_{\ggamma_i},D_{\bbeta_i},D_i,\ldots,D_m$$
In view of the relation \eqref{eq:or 1} we must therefore prove that
\begin{multline*}
(-1)^{|\ggamma_i|}\cdot
(-1)^{s(\bbeta_i;\aalpha_i)+s(\ggamma_i;\aalpha_i\setminus\bbeta_i)}\\
=
-
(-1)^{(|\bbeta_i|-1)(|\ggamma_i|-1)}\cdot
(-1)^{|\bbeta_i|}\cdot
(-1)^{s(\ggamma_i;\aalpha_i)+s(\bbeta_i;\aalpha_i\setminus\ggamma_i)}
\end{multline*}
which is settled by the following

\begin{lemma}
$(-1)^{s(\bbeta_i;\aalpha_i)-s(\bbeta_i;\aalpha_i-\ggamma_i)}=
(-1)^{|\bbeta_i||\ggamma_i|}
(-1)^{s(\ggamma_i;\aalpha_i)-s(\ggamma_i;\aalpha_i\setminus\bbeta_i)}$
\end{lemma}
\proof The \lhs is the parity of the set $N_{\bbeta_i,\ggamma_i}$ of pairs
$(\beta,\gamma)\in\bbeta_i\times\ggamma_i$ which are permuted when
$\bbeta_i$ is shuffled to the left of $\aalpha_i$. Similarly, $(-1)^{s(\ggamma
_i;\aalpha_i)-s(\ggamma_i;\aalpha_i\setminus\bbeta_i)}$ is the parity of
the set $N_{\ggamma_i,\bbeta_i}$ of pairs $(\beta,\gamma)\in\bbeta_i
\times\ggamma_i$ which are permuted when $\ggamma_i$ is shuffled
to the left of $\aalpha_i$. Since $N_{\bbeta_i,\ggamma_i}\cup N_{\ggamma
_i,\bbeta_i}=\bbeta_i\times\ggamma_i$, the product of these parities is
equal to $(-1)^{|\bbeta_i||\ggamma_i|}$ \halmos

\subsection{The cellular chain complex of $\AD$}\label{ss:cellular complex}

We construct below an isomorphism between the complex $C_*(\ND)$ and
the cellular chain complex of the associahedron $\AD$ by using its realisation
as a convex polytope $P_D^c$ given in section \ref{ss:realisation}. We begin
by explaining how an orientation of a nested set $\H\in\ND$ determines one
of the corresponding face $P_{\H}^c$ of $P_D^c$.\\

For any $t\in P_D^c\subset\IR^{|D|}$ with coordinates $\{t_{\alpha}\}_{\alpha
\in D}$ and subset $B\subset D$, set $$t_B=\sum_{\gamma\in B}t_\gamma$$
If $t\in P_{\H}^c$, then, for any $A\in\H$
\begin{equation}\label{eq:constraint}
t_{\aalpha^A_{\H}}=
t_A-t_{A\setminus\aalpha^A_{\H}}=
c(A)-c(A\setminus\aalpha^A_{\H})
\end{equation}
where we extend the function $c$ to non--connected subdiagrams $B\subset
D$ with connected components $B_1,\ldots,B_m$ by setting $c(B)=c(B_1)+
\cdots+c(B_m)$. It follows that if $D_1,\ldots,D_m$ are the unsaturated elements
of $\H$, a redundant system of coordinates on $P_{\H}^c$ is given by the
components $t_{\gamma}$ with $\gamma$ ranging over $\aalpha^{D_1}
_{\H}\odots{\cup}\aalpha^{D_m}_{\H}$. These coordinates are only subject
to the constraints that equation \eqref{eq:constraint} should hold whenever
$A=D_i$ for some $i$.\\

Assume that $\H$ is of positive dimension. Let $\eepsilon$ be an orientation of $\H$
and let
$$D_1,\ldots,D_m
\qquad\text{and}\qquad
\aalpha_i=\aalpha^{D_i}_{\H}=\{\alpha_i^1,\ldots,\alpha_i^{n_i}\}\subset D_i$$
be the corresponding enumeration of the unsaturated elements of $\H$ and
ordered subsets of vertices respectively.

\begin{definition}
The orientation of the face $P_{\H}^c$ induced by $\eepsilon$ is the one
determined by the volume element
$$\omega_{\H_\eepsilon}=
dt_{\alpha_1^1}\odots{\wedge}dt_{\alpha_1^{n_1-1}}
\odots{\wedge}
dt_{\alpha_m^1}\odots{\wedge}dt_{\alpha_m^{n_m-1}}$$
\end{definition}

Note that the assignement $\H_{\eepsilon}\longrightarrow\omega_{\H_\eepsilon}$
is consistent with relations \eqref{eq:or 1}--\eqref{eq:or 2}. This is clear if one
permutes $D_i$ and $D_{i+1}$ or $\alpha_i^{j}$ and $\alpha_i^{j+1}$ within
$\aalpha_i$, so long as $1\leq j\leq n_i-2$. If $j=n_i-1$, the new contribution
of $D_i$ to the volume form on $P_{\H}^c$ is
\begin{equation*}
\begin{split}
dt_{\alpha_i^1}\odots{\wedge}dt_{\alpha_i^{n_i-2}}&\wedge dt_{\alpha_i^{n_i}}\\
&=
dt_{\alpha_i^1}\odots{\wedge}dt_{\alpha_i^{n_i-2}}\wedge
d(c(D_i)-c(D_i\setminus\aalpha_i)-\sum_{\ell=1}^{n_i-1}t_{\alpha_i^{\ell}})\\
&=
-dt_{\alpha_i^1}\odots{\wedge}dt_{\alpha_i^{n_i-2}}\wedge dt_{\alpha_i^{n_i-1}}
\end{split}
\end{equation*}
as required. We next work out the orientation given by the volume form $\omega
_{\H_\eepsilon}$ more explicitly. Assume first that $\H$ is of dimension 1 so that $m=1
$ and $\aalpha_1=\{\alpha_1^1,\alpha_1^2\}$. By \S\ref{ss:edges}, the corresponding edge $P_\H^c$ has boundary points labelled by $\H\cup D_{\alpha_1^1}$ and $\H\cup D_{\alpha_1^2}$.

\begin{lemma}\label{le:oriented edge}
The orientation of $P_\H^c$ induced by $\eepsilon$ is given by
$$\begin{diagram}
\H\cup D_{\alpha_1^1}&&\rTo^{\H}&&\H\cup D_{\alpha_1^2}
\end{diagram}$$
\end{lemma}
\proof
In this case $\omega_{\H_\eepsilon}=dt_{\alpha_1^1}$. By \eqref{eq:constraint},
$$t_{\alpha_1^1}(P_{\H\cup D_{\alpha_1^1}})=
c(D_{\alpha_1^1})-c(D_{\alpha_1^1}\setminus{\alpha_1^1})
\quad\text{and}\quad
t_{\alpha_1^1}(P_{\H\cup D_{\alpha_1^2}})=
c(D_1)-c(D_1\setminus{\alpha_1^1})$$
Since the connected components of $D_1\setminus{\alpha_1^1}$ not containing
${\alpha_1^2}$ are the connected components of $D_{\alpha_1^1}\setminus
\alpha_1^1$ which are orthogonal to $\alpha_1^2$, we have
$$c(D_1)-(c(D_1\setminus{\alpha_1^1})-
c(D_{\alpha_1^1}\setminus{\alpha_1^1}))-c(D_{\alpha_1^1})\geq
c(D_1)-c(D_{\alpha_1^2})-c(D_{\alpha_1^1})>0$$
where the inequality follows from \eqref{eq:c} and the fact $D_{\alpha_1^1}$
and $D_{\alpha_1^2}$ are not compatible and such that $D_{\alpha_1^1}\cup
D_{\alpha_1^2}=D_1$ \halmos\\

Assume now that $\H$ is of dimension greater or equal to two and let $\G=\H
\cup\bbeta_i$ be a boundary facet of $\H$. Note that the function $t_{\bbeta_i}
=\sum_{\gamma\in\bbeta_i}t_\gamma$ is identically equal to $c(D_{\bbeta_i})-c(D_{\bbeta_i}\setminus\bbeta_i)$ on $P_{\G}^c\subset\partial P_{\H}^c$ and
strictly greater than that value on the interior of $P_{\H}^c$ since $D_{\bbeta_i}
\notin\H$ while the connected components of $D_{\bbeta_i}\setminus\bbeta_i$
lie in $\H$. The orientation of $P_{\G}^c$ as a boundary component of
$P_{\H}^c$ induced by the volume form $\omega_{\H_\eepsilon}$ is therefore
given by the form $\partial_{\bbeta_i}\omega_{\H_\eepsilon}$ such that
$$\omega_{\H_\eepsilon}=
 -dt_{\bbeta_i}\wedge\partial_{\bbeta_i}\omega_{\H_\eepsilon}$$
We wish to relate $\partial_{\bbeta_i}\omega_{\H_\eepsilon}$ to the volume
form $\omega_{\G_{\ol{\eepsilon}}}$ where $\ol{\eepsilon}$ is the orientation
of $\G$ induced by $\eepsilon$.

\begin{lemma}\label{le:induced}
$$\partial_{\bbeta_i}\omega_{\H_\eepsilon}
=-
(-1)^{(|\aalpha_1|-1)+\cdots+(|\aalpha_{i-1}|-1)}\cdot
(-1)^{|\bbeta_i|-1}\cdot
(-1)^{s(\bbeta_i;\aalpha_i)}\cdot
\omega_{\G_{\ol{\eepsilon}}}$$
\end{lemma}
\proof  For any ordered set $\ggamma=\{\alpha_1,\ldots,\alpha_k\}\subset D$,
set $$\omega_{\ggamma}=dt_{\alpha_1}\odots{\wedge}dt_{\alpha_{k-1}}$$
so that
\begin{align*}
\omega_{\H_\eepsilon}
&=
\omega_{\aalpha_1}\odots{\wedge}\omega_{\aalpha_m}\\
\omega_{\G_{\ol{\eepsilon}}}
&=
\omega_{\aalpha_1}\odots{\wedge}\omega_{\aalpha_{i-1}}\wedge
\omega_{\bbeta_i}\wedge\omega_{\aalpha_i\setminus\bbeta_i}\wedge
\omega_{\aalpha_{i+1}}\odots{\wedge}\omega_{\aalpha_m}
\end{align*}

If $\bbeta_i=\{\alpha_i^{j_1},\ldots,\alpha_i^{j_p}\}\subset\{\alpha_i^1,
\ldots,\alpha_i^{n_i}\}=\aalpha_i$, then assuming first $j_p<n_i$,
\begin{equation*}
\begin{split}
\omega_{\aalpha_i}
&=
(-1)^{s(\bbeta_i;\aalpha_i)}\cdot
dt_{\alpha_i^{j_1}}\odots{\wedge}dt_{\alpha_i^{j_p}}\wedge
\omega_{\aalpha_i\setminus\bbeta_i}\\
&=
(-1)^{s(\bbeta_i;\aalpha_i)}\cdot
dt_{\alpha_i^{j_1}}\odots{\wedge}dt_{\alpha_i^{j_p-1}}\wedge
d(t_{\alpha_i^{j_1}}+\cdots+t_{\alpha_i^{j_p}})\wedge
\omega_{\aalpha_i\setminus\bbeta_i}\\
&=
(-1)^{s(\bbeta_i;\aalpha_i)}\cdot
(-1)^{|\bbeta_i|-1}\cdot
dt_{\bbeta_i}\wedge
\omega_{\bbeta_i}\wedge\omega_{\aalpha_i\setminus\bbeta_i}
\end{split}
\end{equation*}
If, on the other hand, $j_p=n_i$, then, denoting the maximal element of
$\aalpha_i\setminus\bbeta_i$ by $\alpha_\ell$, we get
\begin{equation*}
\begin{split}
\omega_{\aalpha_i}
&=
(-1)^{s(\bbeta_i;\aalpha_i)-(|\aalpha_i|-|\bbeta_i|)}\cdot
\omega_{\bbeta_i}\wedge
\omega_{\aalpha_i\setminus\bbeta_i}\wedge dt_{\alpha_\ell}\\
&=
(-1)^{s(\bbeta_i;\aalpha_i)-(|\aalpha_i|-|\bbeta_i|)}\\
&\phantom{=}
\cdot
\omega_{\bbeta_i}\wedge
\omega_{\aalpha_i\setminus\bbeta_i}\wedge
d(c(D_i)-c(D_i\setminus\aalpha_i)-
t_{\bbeta_i}-
t_{\aalpha_i\setminus(\bbeta_i\cup\{\alpha_{\ell}\})})\\
&=
(-1)^{s(\bbeta_i;\aalpha_i)}\cdot(-1)^{|\bbeta_i|-1}\cdot
dt_{\bbeta_i}\wedge
\omega_{\bbeta_i}\wedge\omega_{\aalpha_i\setminus\bbeta_i}
\end{split}
\end{equation*}

Thus, in either case, we find
\begin{equation*}
\begin{split}
\omega_{\H_\eepsilon}
&=
(-1)^{|\bbeta_i|-1}\cdot
(-1)^{s(\bbeta_i;\aalpha_i)}\\
&\phantom{=}
\cdot 
\omega_{\aalpha_1}\odots{\wedge}\omega_{\aalpha_{i-1}}\wedge
dt_{\bbeta_i}\wedge
\omega_{\bbeta_i}\wedge\omega_{\aalpha_i\setminus\bbeta_i}\wedge
\omega_{\aalpha_{i+1}}\odots{\wedge}\omega_{\aalpha_m}\\
&=
(-1)^{(|\aalpha_1|-1)+\cdots+(|\aalpha_{i-1}|-1)}\cdot
(-1)^{|\bbeta_i|-1}\cdot
(-1)^{s(\bbeta_i;\aalpha_i)}\cdot
dt_{\bbeta_i}\wedge
\omega_{\G_{\ol{\eepsilon}}}
\end{split}
\end{equation*}
as required \halmos

\begin{theorem}\label{th:model}
The map
$\H_{\eepsilon}\longrightarrow (P_{\H}^c,\omega_{\H_\eepsilon})$
associating to each oriented nested set the corresponding face of the
polytope $P_D^c$ with orientation given by the volume form $\omega
_{\H_\eepsilon}$ is an isomorphism between $C_*(\ND)$ and the
cellular chain complex of $P_D^c$.
\end{theorem}
\proof Recall (see \eg \cite[\S IX.4]{Ms}) that the cellular chain complex
of a $CW$--complex $X$ is defined by $C\cell_n(X)=H_n(X^n,X^{n-1})$,
where $X^n$ is the $n$--skeleton of $X$, with differential $\partial\cell_
n$ given by the composition
$$H_n(X^n,X^{n-1})\longrightarrow H_{n-1}(X^{n-1})\longrightarrow
H_{n-1}(X^{n-1},X^{n-2})$$
Each $C\cell_n(X)$ is a free abelian group of rank equal to the number
of $n$--cells in $X$. Identifying a given factor with $\IZ$ when $n\geq 1$
amounts to choosing an orientation of the corresponding cell. For a regular
$CW$--complex such as $P_D^c$, that is one where all attaching maps
are homeomorphisms, the boundary $\partial\cell$ has a very simple
description \cite[\S IX.6]{Ms}. Given an oriented $n$--cell $b^n
_\lambda$,
$$\partial\cell b^n_\lambda=
\sum_{\mu}[b^n_\lambda:b^{n-1}_\mu]\cdot b^{n-1}_\mu$$
where the sum ranges over the $(n-1)$--cells of $X$, each taken with a
chosen orientation if $n\geq 2$. The incidence number $[b^n_\lambda:
b^{n-1}_\mu]$ is zero if $b_\mu^{n-1}$ is not contained in the boundary
of $b_\lambda^n$ and $\pm 1$ otherwise. In the latter case, the sign
depends on whether the induced orientation on the boundary of
$b_\lambda^n$ agrees with that on $b_\mu^{n-1}$ if $n\geq 2$
and is otherwise given, for $n=1$ by
$$\partial\cell b^1_\lambda=b^0_\mu-b^0_\nu$$
where, under an orientation preserving identification $b_\lambda^1
\cong[0,1]$, the attaching map send $1$ to $b^0_\mu$ and $0$ to
$b^0_\nu$. Thus, by lemmas \ref{le:oriented edge} and \ref{le:induced},
the map $\H_\eepsilon\longrightarrow (P_{\H_\eepsilon}^c,\omega_{\H_
\eepsilon})$ identifies the differential $\partial$ on $C_*(\ND)$ to the
opposite of the cellular boundary \halmos\\

Theorem \ref{th:model} and the contractibility of $P_D^c$ imply in particular
the following
\begin{corollary}
The complex $C_*(\ND)$ is acyclic.
\end{corollary}

\subsection{The Dynkin complex and the cellular cochain complex of $\AD$}
\label{ss:Dynkin cellular}

Let $A$ be a $D$--algebra. We relate in this subsection the Dynkin complex
of $A$ to the cellular cochain complex of the associahedron $\AD$, when
both are taken with coefficients in a $D$--bimodule $M$ over $A$. We show
in particular that in degrees greater or equal to two, the Dynkin differential
is a geometric boundary operator, albeit in combinatorial guise. We shall need
some terminology.

\begin{definition}\hfill
\begin{enumerate}
\item A nested set $\H\in\ND$ of positive dimension is called irreducible if
it has a unique unsaturated element and reducible otherwise.
\item Two oriented nested sets $\H_\eepsilon$ and $\H'_{\eepsilon'}
$ of positive dimension are equivalent if they have the same unsaturated
elements $D_1,\ldots,D_m$, if $\aalpha^{D_i}_{\H}=\aalpha^{D_i}_
{\H'}$ for any $i=1\ldots m$ and if the orientations $\eepsilon,\eepsilon'$
agree in the obvious sense.
\end{enumerate}
\end{definition}

Note that a nested set of dimension 1 is clearly irreducible. If $\H\in\ND$ is
of dimension 2, \S\ref{ss:2 faces} shows that $\H$ is irreducible when the
corresponding face $\AD^\H$ of $\AD$ is a pentagon or a hexagon and
reducible when $\AD^\H$ is a square. More generally, by proposition
\ref{pr:ND^H}, $\H$ is reducible precisely when $\AD^\H$ is the product
of $p\geq 2$ smaller associahedra.\\

Let $CD^*(A;M)$ be the Dynkin complex of $A$ with coefficients in $M$
and $$\wt{C}^*(\ND;M)=0\longrightarrow M\xrightarrow{\epsilon^*}\Hom
_{\IZ}(C_*(\ND);M)$$ be the augmented cellular cochain complex of $\AD$
with coefficients in $M$. We regard $M$ as sitting in degree $-1$ in $\wt
{C}^*(\ND;M)$. For any $k=0,\ldots,|D|$, define a map $g^k:CD^k(A;M)
\longrightarrow\wt{C}^{k-1}(\ND;M)$ by
\begin{align*}
g^0 m       			&= m_D\\
g^1 m(\H) 			&= \sum_{B\in\H}m_{(B;\alpha^B_\H)}\\
g^k m(\H_{\eepsilon}) &=
\left\{\begin{array}{cl}
m_{(B;\aalpha^B_\H)}&\text{if $\H$ is irreducible with unsaturated set $B$}\\
0				  &\text{if $\H$ is reducible}
\end{array}\right.
\end{align*}
for any $k\geq 2$.

\begin{theorem}\label{th:Dynkin cellular}
$g$ is a chain map from $CD^*(A;M)$ to $\wt{C}^{*-1}(\ND;M)$. For $k\geq 2
$, $g^k$ is an embedding whose image consists of those cochains $c\in\wt{C}
^{k-1}(\ND;M)$ such that
\begin{enumerate}
\item $c(\H_\eepsilon)=0$ for any reducible nested set $\H$.
\item $c(\H_\eepsilon)=c(\H'_{\eepsilon'})$ for any equivalent
oriented nested sets $\H_\eepsilon,\H_\eepsilon'$.
\item $c(\H_\eepsilon)\in M_B^{B\setminus\aalpha^B_{\H}}$
for any irreducible $\H$ with unsaturated element $B$.
\end{enumerate}
\end{theorem}
\proof Denote the differential on $\wt{C}^k(\ND;M)$ by $d^k$, with $d^{-1}=
\epsilon^*$. For $m=\{m_{B}\}_{B\subseteq D}\in CD^0(A;M)$ and $\H$ a
maximal nested set,
$$g^1\dD^0 m(\H)=
\sum_{B\in\H}(m_B-m_{B\setminus\alpha^B_\H})=
m_D=
(d^{-1}g^0 m)(\H)$$
Next, if $m=\{m_{(B;\alpha)}\}_{\alpha\in B\subseteq D}$ is a one--chain, and
$\H_\eepsilon$ an oriented nested set of dimension one with unsaturated element
$B$ and ordered $\aalpha^B_{\H}=\{\alpha_1,\alpha_2\}$, then
$$d^0g^1m(\H_\eepsilon)=g^1m(\H\cup B_1)-g^1m(\H\cup B_2)$$
where $B_1=\con{B}{\alpha_2}{\alpha_1}$ and $B_2=\con{B}{\alpha_1}
{\alpha_2}$. This is equal to
$$m_{(B;\alpha_2)}-m_{(B_2;\alpha_2)}
-m_{(B;\alpha_1)}+m_{(B_1;\alpha_1)}=
g^2\dD^1m(\H_\eepsilon)$$
since $\alpha^{B_i}_{\H\cup B_i}=\alpha_i$ and the only element $B'$ of $\H$
for which $\alpha^{B'}_{\H\cup B_1}\neq\alpha^{B'}_{\H\cup B_2}$ is $B$ with
$\alpha^B_{\H\cup B_1}=\alpha_{2}$ and $\alpha^B_{\H\cup B_2}=\alpha_1$.\\

Let now $m\in CD^{k}(A;M)$ with $k\geq 2$, and let $\H_\eepsilon$ be an
oriented nested set or dimension $k$. Let $p$ be the number of unsaturated
elements of $\H$ and assume first that $p\geq 2$. Then, $\H$ is reducible
so that $$g^{k+1}\dD^k\msp m(\H_\eepsilon)=0$$ If $p\geq 3$, any
boundary facet $\G=\H\cup D_{\bbeta_i}$ of $\H$ is also reducible and
$$d^{k-1}g^k\msp m(\H_\eepsilon)=0$$ as required. If $\H$ only has
two unsaturated elements $D_1$ and $D_2$, a boundary facet $\G=
\H\cup D_{\bbeta_i}$ of $\H$, with $\emptyset\neq\bbeta_i\subsetneq
\aalpha_i=\aalpha^{D_i}_\H$ and $i=1,2$, is reducible unless $|\aalpha
_i|=2$. In that case, the unique unsaturated element of $\G$ is $D_{3-i}$.
It follows that, when $p=2$,
\begin{equation*}
\begin{split}
d^{k-1}g^k\msp m(\H_\eepsilon)
&=
\delta_{|\aalpha_1|=2}\cdot(m(D_2;\aalpha_2)-m(D_2;\aalpha_2))\\
&+
(-1)^{|\aalpha_1|-1}\cdot
\delta_{|\aalpha_2|=2}\cdot(m(D_1;\aalpha_1)-m(D_1;\aalpha_1))\\
&=
0
\end{split}
\end{equation*}
There remains to consider the case when $\H$ is irreducible with unsaturated
element $B$ and $\aalpha=\aalpha^B_\H$ of cardinality $k+1$. If $\emptyset
\neq\bbeta\subsetneq\aalpha$ is of cardinality $2\leq|\bbeta|\leq|\aalpha|-2$,
then $H\cup D_{\bbeta}$ is reducible with unsaturated elements $D_{\bbeta}$
and $B$. It follows that the only non--trivial contributions in $d^{k-1}g^k\msp
m(\H_\eepsilon)$ arise when $\bbeta=\{\alpha_j\}$ or $\bbeta=\aalpha\setminus
\{\alpha_j\}$ for some $j=1,\ldots,k+1$, where $\aalpha=\{\alpha_1,\ldots,\alpha
_{k+1}\}$. The corresponding sign contributions are
\begin{align*}
(-1)^{|\bbeta|-1}\cdot(-1)^{s(\bbeta;\alpha)}&=(-1)^{j-1}\\
(-1)^{|\bbeta|-1}\cdot(-1)^{s(\bbeta;\alpha)}&=(-1)^{k-1}\cdot(-1)^{k-j+1}
\end{align*}
respectively, so that
\begin{equation*}
\begin{split}
d^{k-1}g^k\msp m(\H_\eepsilon)
&=
\sum_{j=1}^{k+1}(-1)^{j-1}\left(
g^k\msp m\left(
(\H\cup\con{B}{(\aalpha\setminus\alpha_j)}{\alpha_j})_{\ol{\eepsilon}}\right)-
g^k\msp m\left(
(\H\cup\con{B}{\alpha_j}{\aalpha\setminus\alpha_j})_{\ol{\eepsilon}}\right)\right)\\
&=
\sum_{j=1}^{k+1}(-1)^{j-1}(
m_{(B;\aalpha\setminus\alpha_j)}-
m_{(\con{B}{\alpha_j}{\aalpha\setminus\alpha_j};\aalpha\setminus\alpha_j)})\\
&=
\dD^k\msp m_{(B;\aalpha)}\\
&=
g^{k+1}\dD^k\msp m(\H_\eepsilon)
\end{split}
\end{equation*}
Thus, $g$ is a chain map as claimed. $g^k$ is an embedding for $k\geq 2$
because for any connected subdiagram $B\subseteq D$ and nonempty subset
$\aalpha\subseteq B$, there exists an irreducible nested set $\H$ with unique
unsaturated element $B$ such that $\aalpha^B_{\H}=\aalpha$. Such an $\H$
may be constructed as follows : let $B_1,\ldots,B_m$ be the connected
components of $B\setminus\aalpha$ and choose a maximal nested set
$\H_i$ on each $B_i$. Let moreover 
$$D_1=B\subset D_2\odots{\subset}D_{p-1}\subset D_p=D$$
be a sequence of encased connected subdiagrams of $D$ such that $|
D_i\setminus D_{i+1}|=1$ for $i=1\ldots p-1$. Then,
$$\H=\H_1\odots{\sqcup}\H_m\sqcup\{D_1,\ldots,D_p\}$$
is a suitable nested set. The image of $g^k$, for $k\geq 2$ is clearly 
characterised by conditions (i)--(iii) \halmos\\

\remark Note that the fact that $g$ is a chain map and that $g^k$ is an
embedding in degrees $\geq 2$ gives another proof of the fact that the
Dynkin differential squares to zero.

\subsection{Deformations of \qc algebras}
\label{ss:Dynkin controls qc}

Let $A$ be a $D$--algebra. Regard $D$ as labelled by attaching an infinite
multiplicity to each edge. By a {\it trivial} \qca structure on $A$ we shall mean
one whose underlying $D$--algebra is $A$ and for which all associators
$\Phi^A_{\G\F}$ are equal to 1. We do not assume that the local monodromies
$S^A_i$ are trivial however. When considering deformations of a trivial \qca
structure on $A$, the elements $S_i^A$ will be assumed to remain undeformed.

\begin{theorem}\label{th:Dynkin controls qc}
The Dynkin complex $CD^*(A;A)$ controls the formal, one--parameter
deformations of trivial \qca structures on $A$. Specifically,
\begin{enumerate}
\item a \qca structure on $A\fml/\hbar^{n+1}A\fml$ which is trivial mod
$\hbar$ canonically determines a Dynkin 3--cocyle $\xi$ and lifts to a
\qca structure on $A\fml/\hbar^{n+2}A\fml$ if, and only if $\xi$ is a
coboundary.
\item Two \qca structures on $A\fml/\hbar^{n+1}A\fml$ which are trivial
mod $\hbar$ and equal mod $\hbar^n$ differ by a Dynkin 2--cocycle
$\varphi$.
They are related by a twist of the form $\{1+\hbar^n a_{(B;\alpha)}\}_
{\alpha\in B\subseteq D}$ if, and only if, $\varphi=\dD a$ and
$Ad(S_i^A)a_{(\alpha_i;\alpha_i)}=a_{(\alpha_i;\alpha_i)}$ for any
$\alpha_i\in D$.
\end{enumerate}
\end{theorem}
\proof The proof of (i) is given in \S \ref{sss:begin obstruction}--\S \ref
{sss:end obstruction}. We prove (ii) first. Let $\Phi^i_{\G\F}$, $i=1,2$, be the associators
of the \qca structures and, for any elementary pair $(\G,\F)$ of \mnss on
$D$, define $\varphi_{\G\F}\in A$ by
$$\Phi^2_{\G\F}=
\Phi^1_{\G\F}+\hbar^n\cdot\varphi_{\G\F}
\mod\hbar^{n+1}$$
Since $\Phi^1_{\G\F}$ and $\Phi^2_{\G\F}$ satisfy the coherence and orientation
axioms mod $\hbar^{n+1}$, $(\G,\F)\rightarrow\varphi_{\G\F}$ is a cellular
one--cocycle on the associahedron $\AD$ with values in $A$. Proposition
\ref{pr:zero on squares}, and the support and forgetfulness properties of
the associators $\Phi^1,\Phi^2$ imply that $\varphi$ satisfies the constraints
(i)--(iii) of theorem \ref{th:Dynkin cellular} respectively so that $\varphi$ is
a Dynkin 2--cocycle. The rest of (ii) is a simple exercise. 

\subsubsection{}\label{sss:begin obstruction}

Let
$$\Phi_{\G\F}=1+\hbar\msp\varphi_{\G\F}^1\odots{+}\hbar^n\varphi_{\G\F}^n$$
be the associators of the \qca structure on $A$.
By assumption, $\Phi_{\G\F}^{-1}=\Phi_{\F\G}$ mod $\hbar^{n+1}$ for any
elementary pair of \mnss $(\G,\F)$ on $D$. We begin by modifying each
$\Phi_{\G\F}$ so that this identity holds mod $\hbar^{n+2}$. Define $\eta
_{\G\F}\in A$ by
$$\Phi_{\G\F}\Phi_{\F\G}=1+\hbar^{n+1}\eta_{\G\F}\mod\hbar^{n+2}$$
Clearly, $\eta_{\G\F}\in A_{B}^{B\setminus\{\alpha,\beta\}}$, where $B=
\supp(\F,\G)$ is
the unique unsaturated element of $\H=\G\cap\F$ and $\{\alpha,\beta\}=
\aalpha^B_\H$, and $\eta_{\G\F}=\eta_{\G'\F'}$ whenever $(\G,\F)$ and
$(\G',\F')$ are equivalent elementary pairs. Moreover, modulo $\hbar^
{n+2}$,
$$\Phi_{\G\F}\Phi_{\F\G}=
\Phi_{\G\F}(\Phi_{\F\G}\Phi_{\G\F}-\hbar^{n+1}\eta_{\F\G})\Phi_{\F\G}=
(\Phi_{\G\F}\Phi_{\F\G})^2-\hbar^{n+1}\eta_{\F\G}$$
whence $\eta_{\F\G}=\eta_{\G\F}$. It follows from this that the associators
$$\wt{\Phi}_{\G\F}=\Phi_{\G\F}-\half{1}\hbar^{n+1}\eta_{\G\F}$$
satisfy $\wt{\Phi}_{\G\F}\wt{\Phi}_{\F\G}=1$ mod $\hbar^{n+2}$, as well
as all the required relations to endow $A\fml/\hbar^{n+2}A\fml$ with a
\qca structure, except possibly for the coherence one.

\subsubsection{}\label{sss:cellular obstruction}

We define next the obstruction $\xi$ as a cellular 2--cochain on $\AD$ with
values in $A$. Let $\H_\eepsilon$ be an oriented nested set of dimension 2,
fix a \mns $\F_0$ on the boundary of $\H$ and let $\F_0,\F_1,\ldots,\F_{k-1},
\F_k=\F_0$ be the vertices of $\H$ listed in their order of appearence on $
\partial\H$ when the latter is endowed with the orientation $\ol{\eepsilon}$.
Define $\xi(\H_\eepsilon;\F_0)\in A$ by
$$\wt{\Phi}_{\F_k\F_{k-1}}\cdots\wt{\Phi}_{\F_1\F_0}=
1+h^{n+1}\xi(\H_\eepsilon;\F_0)\mod\hbar^{n+2}$$
By lemma \ref{le:choices}, $\xi(\H_\eepsilon;\F_0)$ does not depend upon
the choice of $\F_0$ and will be hereafter denoted by $\xi(\H_\eepsilon)$.
Moreover $\xi$ satisfes $\xi(\H_{-\eepsilon})=-\xi(\H_\eepsilon)$, where $-
\veps$ is the opposite orientation of $\H$ since $\wt{\Phi}_{\G\F}^
{-1}=\wt{\Phi}_{\F\G}\mod\hbar^{n+2}$. 

\subsubsection{}
We show next that $\xi$ is a Dynkin 3--cochain. In view of theorem \ref
{th:Dynkin cellular}, it is sufficient to prove the following
\begin{lemma}\hfill
\begin{enumerate}
\item $\xi(\H_\eepsilon)=0$ if $\H$ is a reducible nested set.
\item $\xi(\H_\eepsilon)=\xi(\H'_{\eepsilon'})$ if $\H_\eepsilon,\H'_{\eepsilon'}$
are equivalent.
\item $\xi(\H_\eepsilon)\in A_B^{B\setminus\{\alpha,\beta\}}$ if $\H$ is irreducible
with unsaturated element $B$ and $\aalpha^B_{\H}=\{\alpha,\beta\}$.
\end{enumerate}
\end{lemma}
\proof (i) is a consequence of proposition \ref{pr:zero on squares}. (ii) and (iii)
follow from the analysis of the 2--faces of $\AD$ given in \S\ref{ss:2 faces}
and the forgetfulness and support axioms respectively \halmos

\subsubsection{}\label{sss:3 cocycle}

We claim now that $\xi$ is a Dynkin 3--cocycle. By theorem \ref{th:Dynkin cellular},
it suffices to prove the following
\begin{proposition}
$\xi$ is a cellular 2--cocycle.
\end{proposition}
\proof We shall in fact prove that $\xi$ is a cellular 2--coboundary. Let
$\gamma=\F_0,\F_1,\ldots,\F_k=\F_0$ be an elementary sequence
of \mnss on $D$. Since the associators $\wt{\Phi}_{\G\F}$ satisfy the
coherence axiom, we have
$$\Phi_{\F_k\F_{k-1}}\cdots\Phi_{\F_1\F_0}=
1+\hbar^{n+1}\zeta(\gamma)\mod\hbar^{n+2}$$
for some $\zeta(\gamma)\in A$. Fix a reference \mns $\F_0$ on $D$
and, for any \mns $\F$, use the connectedness of $\AD$ to choose
an edge--path $p_\F$ from $\F_0$ to $\F$. For any oriented 1--edge
$e=(\G,\F)$ in $\AD$, set
$$\eta(e)=\zeta(p_\G\vee e\vee p_\F)\in A$$
where $\vee$ is concatenation. One readily checks that $\eta(\ol{e})=
-\eta(e)$, where $\ol{e}=(\F,\G)$, so that $\eta$ defines a 1--cochain
on $\AD$ with values in $A$ and that $d\eta=\xi$, where $d$ is the
cellular differential \halmos

\subsubsection{}\label{sss:end obstruction}
To complete the proof of (i), we must show the following

\begin{proposition}
The \qca structure on $A\fml/\hbar^{n+1}A\fml$ given by the associators
$\{\Phi_{\G\F}\}$ lifts to one on $A\fml/\hbar^{n+2}A\fml$ if, and only if
$\xi$ is a Dynkin coboundary.
\end{proposition}
\proof
If $\xi=\dD\Theta$ for some $\Theta=\{\Theta_{\G\F}\}$, one readily checks that
the associators $\wt{\Phi}_{\G\F}-\hbar^{n+1}\Theta_{\G\F}$ give the required
lift. Conversely, if
$$\ol{\Phi}_{\G\F}=
\Phi_{\G\F}+\hbar^{n+1}\ol{\Theta}_{\G\F}=
\wt{\Phi}_{\G\F}+\hbar^{n+1}\Theta_{\G\F}$$
endow $A\fml/\hbar^{n+2}A\fml$ with a \qca structure, one readily
checks using theorem \ref{th:Dynkin cellular} that $\Theta$ is a
Dynkin one--cochain and that $\dD\Theta=\xi$ \halmos


%% file: qcqtqb.tex

 \section{Quasi--Coxeter \qt \qbas} \label{se:qcqtqb}

The aim of this section is to define the category of \qc \qt quasibialgebras.
We proceed in stages starting with the notion of $D$--bialgebras, gradually
weakening the bialgebra structure to get $D$--\qt quasibialgebras and
then grafting on a \qc algebra structure.

\subsection{$D$--\bas}

\begin{definition}
A $D$--\ba $(A,\{A_B\},\Delta,\veps)$ is a $D$--algebra $(A,\{A_B\})$
endowed with a bialgebra structure, with coproduct $\Delta$ and counit
$\veps$, such that each $A_B$ is a sub\ba of $A$, that is satisfies
$\Delta(A_B)\subseteq A_B\otimes A_B$.
\end{definition}

\begin{definition}
A morphism of $D$--bialgebras $A,A'$ is a morphism $\{\Psi_\F\}$ of
the underlying $D$--algebras such that each $\Psi_\F$ is a bialgebra
morphism $(A,\Delta,\veps)\rightarrow(A',\Delta',\veps')$.
\end{definition}

If $A$ is a bialgebra, we denote by $\Delta^{(n)}:A\rightarrow
A^{\otimes n}$, $n\geq 0$, the iterated coproduct defined by
$\Delta^{(0)}=\veps$, $\Delta^{(1)}=\id$, and
$$\Delta^{(n+1)}=\Delta\otimes\id^{\otimes n-1}\circ\Delta^{(n)}$$
if $n\geq 1$. Each tensor power $A^{\otimes n}$ of $A$ is an $
A$--bimodule, where $a\in A$ acts by left and right multiplication
by $\Delta^{(n)}(a)$ respectively. If $A$ is a $D$--bialgebra, this
endows each $A^{\otimes n}$ with the structure of a $D$--bimodule
over $A$ by setting $(A^{\otimes n})_B=A_B^{\otimes n}$.
In the notation of \S \ref{ss:Dynkin complex} we then have,
for any $B_1,B_2\subseteq D$, with $B_1$ connected
$$(A^{\otimes n})_{B_1}^{B_2}=
\{\omega\in A_{B_1}^{\otimes n}|
\medspace[\omega,\Delta^{(n)}(a)]=0
\thickspace\text{for all $a\in A_{B_2^i}$}\}$$
where $B_2^i$ are the connected components of $B_2$.

\subsection{$D$--\qbas}
Retain the definitions and notation of \S\ref{sss:qba recap}.

\begin{definition}
A $D$--\qba $(A,\{A_B\},\Delta,\veps,\{\Phi_B\},\{F_{(B;\alpha)}\})$
is a $D$--bialgebra $A$ endowed with the following additional data:
\begin{itemize}
\item {\bf Associators :} for each connected subdiagram $B\subseteq
D$, an invertible element $$\Phi_B\in(A_B^{\otimes 3})^B$$
\item {\bf Structural twists :} for each connected subdiagram $B\subseteq
D$ and vertex $\alpha\in B$, a twist
$$F_{(B;\alpha)}\in(A_B^{\otimes 2})^{B\setminus\alpha}$$
\end{itemize}
satisfying the following axioms:
\begin{itemize}
\item for any connected $B\subseteq D$, $(A_B,\Delta,\veps,\Phi_B)$
is a quasibialgebra.
\item For any connected $B\subseteq D$ and $\alpha\in B$,
\begin{equation}\label{eq:Phi F}
(\Phi_B)_{F_{(B;\alpha)}}=\Phi_{B\setminus\alpha}
\end{equation}
where 
$\Phi_{B\setminus\alpha}=\prod_{B'}\Phi_{B'}$, 
with the product ranging over the connected components of $B\setminus
\alpha$ if $B\neq\alpha$, and $\Phi_\emptyset=1^{\otimes 3}$ otherwise.
\end{itemize}
\end{definition}

Thus, a $D$--\qba has a coassociative coproduct $\Delta$\footnote
{The reader may object that a $D$--\qba is therefore not truly 'quasi'.}.
Moreover, for any \mns $\F$ on $D$, one can coherently twist the family
of \qbas $(A_B,\Delta,\veps,\Phi_B)$, with $B\in\F$, in the following way.
For any connected $B\subseteq D$ and \mns $\F_B$ on $B$, set
\begin{equation}\label{eq:F F}
F_{\F_B}=
\stackrel{\longrightarrow}{\prod_{C\in\F_B}}
F_{(C;\alpha^C_{\F_B})}
\in A_B^{\otimes 2}
\end{equation}
where the product is taken with $F_{(C_1;\alpha^{C_1}_{\F_B})}$ written
to the left of $F_{(C_2;\alpha^{C_2}_{\F_B})}$ whenever $C_1\subset
C_2$. This does not specify the order of the factors uniquely, but two
orders satisfying this requirement are readily seen to yield the same
product. The factorised form of the twist $F_{\F_B}$ implies the
following

\begin{lemma}
Let $\F$ be a \mns on $D$ and $B\in\F$. Then, for any $a\in A_B$,
$$F_\F\cdot\Delta(a)\cdot F_\F^{-1}=
F_{\F_B}\cdot\Delta(a)\cdot F_{\F_B}^{-1}$$
where $\F_B=\{C\in\F|\medspace C\subseteq B\}$ is the \mns on $B$
induced by $\F$.
\end{lemma}

Thus, if $\F$ is a \mns on $D$ and $B\in\F$, the twisted coproduct
\begin{equation}
\Delta_\F(a)=F_\F\cdot\Delta(a)\cdot F_\F^{-1}
\end{equation}
corresponding to $\F$ restricts to $\Delta_{\F_B}$ on $A_B$ so that
$(A_B,\Delta_{\F_B},\veps,(\Phi_B)_{F_{\F_B}})$ is a sub\qba of
$(A,\Delta_{\F},\veps,(\Phi_{D})_{F_\F})$. Turning now to the associators
$\Phi_B$, an inductive application of \eqref{eq:Phi F} readily yields
the following\footnote{this consequence of \eqref{eq:Phi F} arose
during a conversation with R. Nest.}
\begin{lemma}\label{le:twist to 1}
For any connected $B\subseteq D$ and \mns $\F_B$ on $B$,
$$(\Phi_B)_{F_{\F_B}}=1^{\otimes 3}$$
\end{lemma}

In particular, the twisted coproducts $\Delta_{F_{\F_B}},\Delta_{F_\F}$
are in fact coassociative and $(A_B,\Delta_{\F_B},\veps)$ is a sub\ba
of $(A,\Delta_\F,\veps)$.\\

\remark\rklabel{rk:reformulation}Lemma \ref{le:twist to 1} implies
that the associators of a $D$--\qba are not independent variables
since, for any connected $B$ and \mns $\F_B$ on $B$,
\begin{equation}\label{eq:nonab 3 coboundary}
\Phi_B=
\id\otimes\Delta(F_{\F_B}^{-1})\cdot 1\otimes F_{\F_B}^{-1}
\cdot F_{\F_B}\otimes 1\cdot\Delta\otimes\id(F_{\F_B})
\end{equation}
The axioms involving $\Phi_B$ are in fact equivalent to the
requirement that the \rhs of \eqref{eq:nonab 3 coboundary} be
invariant under $A_B$ and independent of the choice of $\F_B$.
It is, however, more convenient to work with the associators
$\Phi_B$.\\

\remark Relation \eqref{eq:Phi F} may be rephrased as follows.
For any subdiagram $B\subseteq D$, let $A_B$ be the algebra
generated by the $A_{B_i}$, where $B_i$ runs over the connected
components of $B$ and set $\Phi_B=\prod_i\Phi_{B_i}$.
Consider the (Drinfeld) tensor category $\Rep_{\Phi_B}(A_B)$
of $A_B$--modules where the associativity constraints are
given by the action of the associator $\Phi_B$. Then, for any
$\alpha\in B\subseteq D$, the twist $F_{(B;\alpha)}$ gives rise to
a tensor structure on the restriction functor $\Rep_{\Phi_B}(A_B)
\rightarrow\Rep_{\Phi_{B\setminus\alpha}}(A_{B\setminus\alpha})$.

\subsection{Morphism of $D$--\qbas}

\begin{definition}\label{de:morphism of Dqba}
A morphism of $D$--\qbas $A,A'$ is a morphism $\Psi=\{\Psi_\F\}$
of the underlying $D$--algebras such that, for any \mns $\F$
on $D$, $\Psi_\F$ is a bialgebra morphism $(A,\Delta_\F,\veps)
\longrightarrow(A',\Delta'_\F,\veps')$.
\end{definition}

\remark Note that a morphism $\Psi:A\to A'$ of $D$--\qbas is not
a morphism of the underlying $D$--bialgebras in general. In other
words, if $\F$ is a maximal nested set, the morphism $\Psi_\F:A
\to A'$ need not satisfy $\Psi_\F^{\otimes 2}\circ\Delta=\Delta'\circ
\Psi_\F$.

\subsection{Twisting of $D$--\qbas}

Let $A$ be a $D$--quasibialgebra.

\begin{definition}\hfill
\begin{enumerate}
\item
A twist of $A$ is a family $F=\{F_B\}_{B\subseteq D}$ labelled
by the connected subdiagrams of $D$, where $F_B\in (A_B^
{\otimes 2})^B$ is an invertible element such that
$$\veps\otimes\id(F_B)=1=\id\otimes\veps(F_B)$$
\item
The twist of $A$ by $F$ is the $D$--\qba
$$A^F=(A,\{A_B\},\Delta,\veps,\{\Phi_B^F\},\{F_{(B;\alpha)}^F\})$$
where, for any connected $B\subseteq D$ and vertex
$\alpha\in B$
\begin{align}
\Phi_B^F
&=
(\Phi_B)_{F_B}
\label{eq:twist Phi}\\
F_{(B;\alpha)}^F
&=
F_{B\setminus\alpha}\cdot F_{(B;\alpha)}\cdot F_B^{-1}
\label{eq:twist F}
\end{align}
with $F_{B\setminus\alpha}=\prod_i F_{B_i}$, where the product
ranges over the connected components of $B\setminus\alpha$ if
$B\neq\alpha$, and $F_\emptyset=1^{\otimes 2}$ otherwise.
\end{enumerate}
\end{definition}

\remark In the notation \eqref{eq:F F}, the twist of $F_{\F_B}$ by
$F$ is given by
$$F^F_{\F_B}=F_{\F_B}\cdot F_B^{-1}$$
Since $F_B$ is invariant under $A_B$, twisting by $F$ does not
change the coproduct $\Delta_{F_{\F_B}}$ on $A_B$, so that $A$
and $A^F$ are isomorphic as $D$--\qbas via the identity map.\footnote
{Notice the difference between morphisms of $D$--\qbas and of
quasibialgebras. As pointed out in \S\ref{sss:qba recap}, a morphism
of \qbas $A\rightarrow A'$ is given by a twist $F$ of $A'$ followed
by a na\"{\i}ve morphism $A\rightarrow A'_F$. Due to the restrictive
nature of the twists we use, the twist of a $D$--\qba $A'$ is still 
'naively' isomorphic to $A'$ in the sense of definition \ref
{de:morphism of Dqba}.}

\subsection{$D$--\qt \qbas}

Retain the definitions of \S\ref{sss:qtqba recap}.

\begin{definition}
A $D$--\qtqba
$$(A,\{A_B\},\Delta,\veps,\{\Phi_B\},\{F_{(B;\alpha)}\},\{R_B\})$$ 
is a $D$--\qba $A$ endowed with an invertible element $R_B\in
A_B^{\otimes 2}$ for each connected subdiagram $B\subseteq
D$ such that $(A_B,\Delta,\veps,\Phi_B,R_B)$ is a \qt quasibialgebra.
\end{definition}

Let $A$ be a $D$--\qt quasibialgebra, $\F$ a \mns on $D$ and
$B\in\F$. By lemma \ref{le:twist to 1}, the twist by $F_{\F_B}$ of
$(A_B,\Delta,\veps,\Phi_B,R_B)$ yields a \qt \ba $(A_B,\Delta_
{\F_B},\veps,(R_B)_{F_{\F_B}})$.

\begin{definition}
A morphism of $D$--\qt \qbas $A,A'$ is a morphism $\{\Psi_\F\}$
of the underlying $D$--\qbas such that, for any \mns $\F$ on $D$
and $B\in\F$, $\Psi_\F$ satisfies
$$\Psi_\F((R_B)_{F_{\F_B}})=(R'_B)_{F'_{\F_B}}$$
and therefore restricts to a morphism of \qt \bas
$$(A_B,\Delta_{\F_B},\veps,(R_B)_{F_{\F_B}})
\longrightarrow
(A'_B,\Delta'_{\F_B},\veps',(R'_B)_{F'_{\F_B}})$$
\end{definition}

\begin{definition}
A twist $F=\{F_B\}_{B\subseteq D}$ of a $D$--\qtqba $A$ is a
twist of the underlying $D$--quasibialgebra. The twisting of $A$
by $F$ is the $D$--\qtqba
$$A^F=(A,\{A_B\},\Delta,\veps,\{\Phi^F_B\},\{F^F_{(B;\alpha)}\},
\{R^F_B\})$$ where $\Phi^F_B$, $F^F_{(B;\alpha)}$ are
given by \eqref{eq:twist Phi}--\eqref{eq:twist F} and
\begin{equation}\label{eq:F twist of R}
R_B^F=
(R_B)_{F_B}=
F_B^{21}\cdot R_B\cdot F_B^{-1}
\end{equation}
\end{definition}

\remark Since for any \mns $\F$ on $D$ and $B\in\F$,
$$(R^F_B)_{F^F_{\F_B}}=
(F_B^{21}\cdot R_B\cdot F_B^{-1})_{F_{\F_B}F_B^{-1}}=
(R_B)_{F_{\F_B}}$$
$A$ and $A^F$ are isomorphic as $D$--\qt \qbas via the
identity map.


\subsection{Quasi--Coxeter \qbas}

Assume henceforth that the diagram $D$ is labelled.

\begin{definition}
A \qc \qba  of type $D$ is a set
$$(A,\{A_B\},\{S_i\},\{\Phi_{(B;\alpha,\beta)}\},
\Delta,\veps,\{F_{(B;\alpha)}\},\{\Phi_B\})$$
where 
\begin{itemize}
\item $(A,\{A_B\},\{S_{i}\},\{\Phi_{(B;\alpha,\beta)}\})$ is a
\qc algebra of type $D$,
\item $(A,\{A_B\},\Delta,\veps,\{F_{(B;\alpha)}\},\{\Phi_B\})$
is a $D$--\qba
\end{itemize}
and, for any pair $(\G,\F)$ of \mnss on $D$, the following holds
\begin{equation}\label{eq:gauge F}
F_\G\cdot\Delta(\Phi_{\G\F})=\Phi_{\G\F}^{\otimes 2}\cdot F_\F
\end{equation}
\end{definition}

Since $\veps\otimes\id(F_\F)=1=\veps\otimes\id(F_\G)$ and
$\veps\otimes\id\circ\Delta=\id$, \eqref{eq:gauge F} implies that
\begin{equation}\label{eq:eps Phi}
\veps(\Phi_{\G\F})=1
\end{equation}
It follows that $\Ad(\Phi_{\G\F})$ is an isomorphism of bialgebras
$(A,\Delta_\F,\veps)\longrightarrow(A,\Delta_\G,\veps)$. Moreover,
by proposition \ref{pr:stability} this isomorphism restricts to an
isomorphism $(A_B,\Delta_{\F_B},\veps)\longrightarrow
(A_B,\Delta_{\G_B},\veps)$ for any element $B\in\F\cap\G$.\\

\remark It was pointed out in remark \rkref{rk:reformulation} that,
for a $D$--\qba the axiom \eqref{eq:Phi F} is equivalent to the
invariance of the \rhs of \eqref{eq:nonab 3 coboundary} and its
independence on the choice of the \mns $\F_B$. Since in a \qc
\qba the twists $F_\F$ are related by gauge transformations, the
invariance of \eqref{eq:nonab 3 coboundary} implies its independence
on the choice of $\F_B$. Thus, for \qc \qba axiom \eqref{eq:Phi F}
is equivalent to the invariance of the \rhs of \eqref{eq:nonab 3 coboundary}.\\

Let us spell out \eqref{eq:gauge F} in diagrammatic notation. By the
connectedness of the associahedron $\AD$, \eqref{eq:gauge F} holds
for any pair $(\G,\F)$ if, and only if it holds for any elementary pair of
\mnss on $D$. Let $(\G,\F)$ be one such pair, $B=\supp(\G,\F)$ the
unique unsaturated element of $\F\cap\G$ and set $\alpha_1=\alpha
^B_\F$ and $\alpha_2=\alpha^B_\G$, so that, in the notation of
\S \ref{ss:qca 3}, $\Phi_{\G\F}=\Phi_{(B;\alpha_2,\alpha_1)}$.

\begin{lemma}\label{le:equivalent gauge}
The relation \eqref{eq:gauge F} is equivalent to
$$F_{(\con{B}{\alpha_2}{\alpha_1};\alpha_1)}\cdot F_{(B,\alpha_2)} 
\cdot\Delta(\Phi_{(B;\alpha_2,\alpha_1)})=
\Phi_{(B;\alpha_2,\alpha_1)}^{\otimes 2}\cdot
F_{(\con{B}{\alpha_1}{\alpha_2};\alpha_2)}\cdot F_{(B,\alpha_1)}$$
\end{lemma}
\proof By definition,
$$F_\G=\stackrel{\longrightarrow}{\prod_{C\in\G}}F_{(C;\alpha^C_\G)}
\qquad\text{and}\qquad
F_\F=\stackrel{\longrightarrow}{\prod_{C\in\F}}F_{(C;\alpha^C_\F)}$$
Let $C\in\F\cap\G$, with $C\neq B$ so that $\alpha^C_\F=\alpha^C
_\G$. If $C\perp B$, then
\begin{equation}\label{eq:comm Phi F}
F_{(C;\alpha^C_{\G})}\cdot\Delta(\Phi_{(B;\alpha_2,\alpha_1)})=
\Delta(\Phi_{(B;\alpha_2,\alpha_1)})\cdot F_{(C;\alpha^C_{\F})}
\end{equation}
since $F_{(C;\alpha^C_{\G})}\in A_C^{\otimes 2}$, $\Delta(\Phi_
{(B;\alpha_2,\alpha_1)})\in A_B^{\otimes 2}$ and $[A_C,A_B]=0$.
If $C\supsetneq B$, then $B\subseteq C\setminus\alpha^C_\G$ and
\eqref{eq:comm Phi F} holds since $F_{(C;\alpha^C_{\G})}$ commutes
with $\Delta(A_{C\setminus\alpha^C_\G})$. Finally, if $C\subsetneq B$,
then $C\subseteq B\setminus\{\alpha_1,\alpha_2\}$ and
$$\Phi_{(B;\alpha_2,\alpha_1)}^{\otimes 2}\cdot F_{(C;\alpha^C_{\G})}=
F_{(C;\alpha^C_{\G})}\cdot \Phi_{(B;\alpha_2,\alpha_1)}^{\otimes 2}$$
since $\Phi_{(B;\alpha_2,\alpha_1)}$ centralises $A_{C\setminus\{\alpha_1,
\alpha_2\}}$. This implies the stated equivalence since, by proposition
\ref{pr:combi<=>practical}, $\G\setminus\F=\con{B}{\alpha_2}{\alpha_1}$
and $\F\setminus\G=\con{B}{\alpha_2}{\alpha_1}$ \halmos

\subsection{Morphism of \qc \qbas}

\begin{definition}
A morphism $\Psi:A\to A'$ of \qc \qbas of type $D$ is a morphism
of the underlying \qcas and $D$--quasibialgebras.
\end{definition}

Thus, $\Psi$ is a collection of algebra morphisms $\Psi_\F:A\to A'$
labelled by \mnss on $D$ such that
\begin{gather*}
\Psi_\F^{\otimes 2}\circ\Delta_\F=\Delta'_\F\circ\Psi_\F\\
\veps'\circ\Psi_\F=\veps\\
\Psi_\F\circ\Ad(\Phi_{\F\G}^A)=\Ad(\Phi_{\F\G}^{A'})\circ\Psi_\G\\
\Psi_\F(S_i^A)=S_i^{A'}
\end{gather*}
for any $\F,\G$ and vertex $\alpha_i\in D$ such that $\{\alpha_i\}\in\F$.

\subsection{Twisting of quasi--Coxeter quasibialgebras}

Let $A$ be a \qc quasibialgebra of type $D$.

\begin{definition}\hfill
\begin{enumerate}
\item A twist of $A$ is a pair $(a,F)$ where
\begin{itemize}
\item $a=\{a_{(B;\alpha)}\}$ is a twist of the underlying \qc algebra such
that $\veps(a_{(B;\alpha)})=1$ for any $\alpha\in B\subseteq D$.
\item $F=\{F_B\}_B$ is a twist of the underlying $D$--quasibialgebra.
\end{itemize}
\item The twisting of $A$ by $(a,F)$ is the \qc \qba
$$(A,\{A_B\},\{S_i^a\},\{\Phi_{(B;\alpha,\beta)}^a\},
\Delta,\veps,
\{F^{(a,F)}_{(B;\alpha)}\},\{\Phi_B^F\})$$
where $S_i^a$, $\Phi^a_{(B;\alpha,\beta)}$ are given by \eqref
{eq:a twist of S}--\eqref{eq:a twist of Phi}, $\Phi_B^F$ is given
by \eqref{eq:twist Phi} and 
\begin{equation*}
\begin{split}
F^{(a,F)}_{(B;\alpha)}
&=
a_{(B;\alpha)}^{\otimes 2}\cdot F_{(B;\alpha)}^F\cdot\Delta(a_{(B;\alpha)})^{-1}\\
&=
F_{B\setminus\alpha}\cdot a_{(B;\alpha)}^{\otimes 2}\cdot
F_{(B;\alpha)}\cdot
\Delta(a_{(B;\alpha)})^{-1}\cdot F_B^{-1}
\end{split}
\end{equation*}
\end{enumerate}
\end{definition}

Note that the twisting of $A$ by $(a_1,F_1)$ followed by a twisting by
$(a_2,F_2)$ is equal to the twisting by $(a_2\cdot a_1,F_2\cdot F_1)$.\\

\remark
In the notation \eqref{eq:F F}, the twist of $F_{\F_B}$ by $(a,F)$ is given 
by
$$F_{\F_B}^{(a,F)}
= a_{\F_B}^{\otimes 2}\cdot F_{\F_B}\cdot\Delta(a_{\F_B})^{-1}\cdot F_B^{-1}
= a_{\F_B}^{\otimes 2}\cdot F_{\F_B}\cdot F_B^{-1}\cdot\Delta(a_{\F_B})^{-1}$$
where $a_{\F_B}=\stackrel{\rightarrow}{\prod}_{C\in\F_B}a_{(C;\alpha^C_\F)}$.

\begin{proposition}
Let $(a,F)$ be a twist of $A$. Then, the assignment $\F\rightarrow
\Ad(a_\F)$ defines an isomorphism of the \qc \qbas $A$ and $A^{
(a,F)}$.
\end{proposition}

\subsection{Quasi--Coxeter \qtqbas}

\begin{definition}
A \qcqtqba of type $D$ is a \qc \qba $A$ of type $D$ endowed
with an invertible element $R_B\in A_B^{\otimes 2}$ for each
connected subdiagram $B\subseteq D$ such that:
\begin{enumerate}
\item $(A_B,\Delta,\veps,\Phi_B,R_B)$ is a \qt quasibialgebra.
\item For any $\alpha_i\in D$, the following holds
\begin{equation}\label{eq:coprod id}
\Delta_{F_{(\alpha_i;\alpha_i)}}(S_i)
=
(R_{\alpha_i})_{F_{(\alpha_i;\alpha_i)}}^{21}\cdot
S_i\otimes S_i
\end{equation}
\end{enumerate}
\end{definition}

\begin{definition}\hfill
\begin{enumerate}
\item
A morphism $\Psi:A\to A'$ of \qc \qt \qbas is a morphism of the
underlying \qc \qbas such that, for any \mns $\F$ on $D$ and
$B\in\F$
$$\Psi_\F((R_B)_{F_{\F_B}})=(R'_B)_{F'_{\F_B}}$$
\item A twist $(a,F)$ of a \qc \qt \qba $A$ is one of the underlying
\qc quasibialgebra. The twisting of $A$ by $(a,F)$ the \qc \qt \qba
$$(A,\{A_B\},\{S_i^a\},\{\Phi_{(B;\alpha,\beta)}^a\},
\Delta,\veps,
\{F^{(a,F)}_{(B;\alpha)}\},\{\Phi_B^F\},\{R^F_B\})$$
where $R_B^F=(R_B)_{F_B}$ is given by \eqref{eq:F twist of R}.
\end{enumerate}
\end{definition}

\remark Since $\veps\otimes\id(R)=1=\id\otimes\veps(R)$ in any \qt \qba
$(A,\Delta,\veps,\Phi,R)$ \cite[\S 3]{Dr3}, applying $\veps\otimes\id$ to
\eqref {eq:coprod id} yields
$$\veps(S_i)=1$$
By \eqref{eq:eps Phi}, this implies that, for any \mns $\F$ on $D$, the
action $\veps\circ\pi_\F$ of the braid group $\BD$ on the trivial $A
$--module is trivial.


%% file: dynkinbi.tex

\section{The \DH bicomplex of a $D$--bialgebra}
\label{se:dynkinbi}

Let $D$ be a connected diagram and $A$ a $D$--bialgebra. By
combining the Dynkin complex of A with the cobar complexes of
its subalgebras $A_B$, we define in this section a bicomplex
which controls the deformations of \qc \qba structures on $A$.

\subsection{}

Let $A$ be a bialgebra and $C^*(A)$ the cobar complex of $A$,
regarded as a coalgebra, defined in \S\ref {ss:Dynkin Hochschild}.
If $C\subseteq A$ a sub--bialgebra, let
$$C^n(A)^C=
\{a\in C^n(A)|\medspace [a,\Delta^{(n)}(c)]=0
\medspace\medspace
\text{for any $c\in C$}\}$$
where $\Delta^{(n)}:C\rightarrow C^{\otimes n}$ is the $n$th iterated
coproduct, be the submodule of $C$--invariants. It is easy to check
that $C^*(A)^C$ is a subcomplex of $C^*(A)$.\\

Assume now that $A$ is a $D$--bialgebra. For $p\in\IN$ and $0\leq
q\leq|D|$, let
$$CD^{q}(A;A^{\otimes p})\subset
\bigoplus_{\substack{\aalpha\subseteq B\subseteq D,\\ |\aalpha|=q}}
(A_B^{\otimes p})^{B\setminus\aalpha}$$
be the group of Dynkin $q$--cochains with values in the $D$--bimodule
$A^{\otimes p}$ over $A$. The Dynkin differential $\dDD$ defines a
vertical differential $CD^{q}(A;A^{\otimes p})\rightarrow CD^{q+1}(A;
A^{\otimes p})$ while the \Ho differential $\dh$ defines a horizontal
differential $CD^{q}(A;A^{\otimes p})\rightarrow CD^{q}(A;A^{\otimes
(p+1)})$. A straightforward computation yields the following.

\begin{theorem}
One has $\dDD\circ\dh=\dh\circ\dDD$. The corresponding cohomology
of the bicomplex $CD^{q}(A;A^{\otimes p})$ is called the \DH cohomology
of the $D$--\ba $A$.
\end{theorem}

\subsection{}

Regard $D$ as labelled by attaching an infinite multiplicity to each edge.
By a trivial \qc \qba structure on $A$ we shall mean one whose underlying
$D$--\ba structure is that of $A$ and for which
$$\Phi_{(B;\beta,\alpha)}=1,\qquad
F_{(B;\alpha)}=1^{\otimes 2}\qquad\text{and}\qquad
\Phi_B=1^{\otimes 3}$$
We do not assume that the local monodromies $S_i$ are trivial however.
When considering deformations of a trivial \qc \qba structure, the local
monodromies will be assumed to remain undeformed.

\begin{theorem}
The \DH bicomplex of $A$ controls the formal, one--parameter
deformations of trivial \qc \qba structures on $A$. Specifically,
\begin{enumerate}
\item A \qc \qba structure on $A\fml/\hbar^{n+1}A\fml$ which is trivial
mod $\hbar$ canonically determines a \DH 4--cocycle $\xi$ and lifts
to a \qc \qba structure on $A\fml/\hbar^{n+2}A\fml$ if, and only if $\xi$
is a coboundary.
\item Two \qc \qba structures on $A\fml/\hbar^{n+1}A\fml$ which are
trivial mod $\hbar$ and agree mod $\hbar^n$ differ by a \DH 3--cocycle
$\eta$ and can be obtained from each other by a twist of the form
$(1+\hbar^n a,1+\hbar^n F)$ if, and only if, $\eta=d(a,F)$ and
$\Ad(S_i^A)a_{(\alpha_i;\alpha_i)}=a_{(\alpha_i;\alpha_i)}$.
\end{enumerate}
\end{theorem}
\proof (i) Let $(\Phi_{(B;\alpha,\beta)},F_{(B;\alpha)},\Phi_B)$ be the
associators and structural twists of the \qc \qba structure on $A\fml/
\hbar^{n+1}A\fml$. We construct in \S \ref{ss:xi}--\ref{ss:theta} a
cochain $(\xi,\eta,\chi,\theta)\in\bigoplus_{i+j=4}CD^{j}(A;A^{\otimes i})$
and check that it is a \DH cocycle.

\subsection{}\label{ss:xi}

Proceeding as in \S \ref{sss:begin obstruction}, we may
assume that $\Phi_{\F\G}=\Phi_{\G\F}^{-1}$ mod $\hbar^{n+2}$ for
any elementary pair $(\G,\F)$ of \mnss on $D$. Let $\H_\veps$ be
an oriented nested set of dimension 2 on $D$, $\F_0$ a \mns on the
boundary of $\H$ and $\F_0,\F_1,\ldots,\F_{k-1},\F_k=\F_0$ the
vertices of $\H$ listed in their order of appearence on $\partial
\H$ when the latter is endowed with the orientation $\ol{\veps}$.
Define $\xi(\H_\veps)\in A$ by
$$\Phi_{\F_k\F_{k-1}}\cdots\Phi_{\F_1\F_0}=
1+h^{n+1}\xi(\H_\veps)\mod\hbar^{n+2}$$
It was proved in \S \ref{sss:cellular obstruction}--\ref{sss:3 cocycle}
that $\xi$ is independent of the choice of $\F_0$, satisfies $\xi
(\H_{-\veps})=-\xi(\H_\veps)$, where $-\veps$ is the opposite
orientation, and $\xi$ is a Dynkin 3--cocycle with values in $A$.

\subsection{} \label{ss:eta}

For any elementary pair $(\G,\F)$ of \mnss on $D$,
define $\eta_{\G\F}\in A^{\otimes 2}$ by
\begin{equation}\label{eq:eta}
F_\G\cdot\Delta(\Phi_{\G\F})-\Phi_{\G\F}^{\otimes 2}\cdot F_\F
=
\hbar^{n+1}\eta_{\G\F}\mod\hbar^{n+2}
\end{equation}

\begin{lemma}
The following holds
\begin{enumerate}
\item $\eta_{\F\G}=-\eta_{\G\F}$ so that $\eta$ is a cellular 1--cochain
on the associahedron $\AD$ with values in $A^{\otimes 2}$.
\item If $d$ is the cellular differential on $\AD$, then $d\eta=-\dH\xi$.
\item $\eta$ is a Dynkin 2--cochain.
\item $\dDD\eta=-\dH\xi$.
\end{enumerate}
\end{lemma}
\proof
(i) Multiplying \eqref{eq:eta} on the left by $\Phi_{\F\G}^{\otimes 2}
=1^{\otimes 2}$ mod $\hbar$ and on the right by $\Delta(\Phi_{\F,
\G})=1^{\otimes 2}$ mod $\hbar$ and using the fact that $\Phi_{\G\F}
\cdot\Phi_{\F\G}=1$ mod $\hbar^{n+2}$, we get, working mod $\hbar
^{n+2}$
$$\hbar^{n+1}\eta_{\G\F}=
\Phi_{\F\G}^{\otimes 2}\cdot F_\G-
F_\F\cdot\Delta(\Phi_{\F\G})=-\hbar^{n+1}\eta_{\F\G}$$
(ii) Since $\Phi_{\F\G}=1$ mod $\hbar$ and $\Phi_{\G\F}\cdot\Phi
_{\F\G}=1$ mod $\hbar^{n+2}$, \eqref{eq:eta} may be rewritten as
$$F_\G=\hbar^{n+1}\eta_{\G\F}+
\Phi_{\G\F}^{\otimes 2}\cdot F_\F\cdot\Delta(\Phi_{\F\G})\mod\hbar^{n+2}$$
Let $\H_\veps$, $\F_0,\ldots,\F_{k-1}$ be as in \S \ref{ss:xi}.
Then, mod $\hbar^{n+2}$,
\begin{equation*}
\begin{split}
F_{\F_0}
&=
\hbar^{n+1}\eta_{\F_0,\F_{k-1}}+
\Phi_{\F_0,\F_{k-1}}^{\otimes 2}\cdot F_{\F_{k-1}}\cdot
\Delta(\Phi_{\F_{k-1},\F_0})\\
&=
\hbar^{n+1}(\eta_{\F_0,\F_{k-1}}+\cdots+\eta_{\F_1,\F_0})\\
&+
\left(\Phi_{\F_0,\F_{k-1}}\cdots\Phi_{\F_1,\F_0}\right)^{\otimes 2}\cdot
F_{\F_0}\cdot
\Delta\left(\Phi_{\F_0,\F_1}\cdots\Phi_{\F_{k-1},\F_0}\right)\\
&=
\hbar^{n+1}d\eta(\H_\veps)
+
\left(1+h^{n+1}\xi(\H_\veps)\right)^{\otimes 2}\cdot
F_{\F_0}\cdot
\Delta\left(1-h^{n+1}\xi(\H_\veps)\right)\\
&=
\hbar^{n+1}\left(d\eta(\H_\veps)
+
1\otimes\xi(\H_\veps)-\Delta(\xi(\H_\veps))+\xi(\H_\veps)\otimes 1\right)
+
F_{\F_0}
\end{split}
\end{equation*}
(iii) Let $(\G,\F)$ be an elementary pair of \mnss on $D$ and set $B=\supp
(\F,\G)$, $\alpha_1=\alpha^B_\F$ and $\alpha_2=\alpha^B_\G$ so that
$\Phi_{\G\F}=\Phi_{(B;\alpha_2,\alpha_1)}$. Reasoning as in the proof of
lemma \ref{le:equivalent gauge}, shows that, modulo $\hbar^{n+2}$,
\begin{equation}\label{eq:new eta}
\begin{split}
\hbar^{n+1}\eta_{\G\F}
&=
F_{(\con{B}{\alpha_2}{\alpha_1};\alpha_1)}\cdot F_{(B;\alpha_2)}\cdot
\Delta(\Phi_{(B;\alpha_2,\alpha_1)})\\
&-
\Phi_{(B;\alpha_2,\alpha_1)}^{\otimes 2}\cdot
F_{(\con{B}{\alpha_1}{\alpha_2};\alpha_2)}\cdot F_{(B;\alpha_1)}
\end{split}
\end{equation}
from which it readily follows that $\eta_{\G\F}\in (A_B^{\otimes 2})^
{B\setminus\{\alpha_1,\alpha_2\}}$ and that $\eta_{\G\F}$ only depends
on the equivalence class of the elementary pair $(\G,\F)$. By theorem
\ref{th:Dynkin cellular}, $\xi$ is therefore a Dynkin 2--cochain.\\
(iv) is a direct consequence of (ii) and theorem \ref{th:Dynkin cellular}
\halmos

\subsection{}\label{ss:chi}

For any $\alpha\in B\subseteq D$, define $\chi_{(B;\alpha)}\in A^{\otimes 3}$
by
\begin{equation}\label{eq:chi}
\begin{split}
\hbar^{n+1}\chi_{(B;\alpha)}
&=
1\otimes F_{(B;\alpha)}\cdot\id\otimes\Delta(F_{(B;\alpha)})\cdot\Phi_B\\
&-
\Phi_{B\setminus\alpha}\cdot F_{(B;\alpha)}\otimes 1\cdot\Delta\otimes\id(F_{(B;\alpha)})
\mod\hbar^{n+2}
\end{split}
\end{equation}

\begin{lemma}
The following holds
\begin{enumerate}
\item $\chi$ is a Dynkin 1--cochain with values in $A^{\otimes 3}$.
\item $\dDD\chi=\dH\eta$
\end{enumerate}
\end{lemma}
\proof
(i) We must show that $\chi_{(B;\alpha)}\in (A_B^{\otimes 2})
^{B\setminus\alpha}$. This readily follows from the support properties
of $\Phi_B,\Phi_{B\setminus\alpha}$ and $F_{(B;\alpha)}$.\\
(ii) Since $F_{(B;\alpha)}=1$ mod $\hbar$, \eqref{eq:chi} may be
rewritten as
\begin{equation}\label{eq:new chi}
(\Phi_B)_{F_{(B;\alpha)}}=
\Phi_{B\setminus\alpha}+\hbar^{n+1}\chi_{(B;\alpha)}
\mod\hbar^{n+2}
\end{equation}
Let $\alpha_1\neq\alpha_2\in B$ and set $D_1=\con{B}{\alpha_2}{\alpha_1}$,
$D_2=\con{B}{\alpha_1}{\alpha_2}$. Then, mod $\hbar^{n+2}$
\begin{equation*}
\begin{split}
(\Phi_B)_{F_{(D_2;\alpha_2)}\cdot F_{(B;\alpha_1)}}
&=
(\Phi_{B\setminus\alpha_1}+\hbar^{n+1}\chi_{(B;\alpha)})_{F_{(D_2;\alpha_2)}}\\
&=
\prod_{D"}\Phi_{D"}\cdot(\Phi_{D_2})_{F_{(D_2;\alpha_2)}}+
\hbar^{n+1}\chi_{(B;\alpha_1)}\\
&=
\prod_{D"}\Phi_{D"}\cdot
(\Phi_{D_2\setminus\alpha_2}+\hbar^{n+1}\chi_{(D_2;\alpha_2)})+
\hbar^{n+1}\chi_{(B;\alpha_1)}\\
&=
\Phi_{B\setminus\{\alpha_1,\alpha_2\}}+
\hbar^{n+1}\left(\chi_{(B;\alpha_1)}+\chi_{(D_2;\alpha_2)})\right)
\end{split}
\end{equation*}
where the product in the third equality ranges over the connected components
$D"$ of $B\setminus\alpha_1$ not containing $\alpha_2$. Permuting $\alpha_1$
and $\alpha_2$, we get mod $\hbar^{n+2}$
\begin{equation}\label{eq:other path}
(\Phi_B)_{F_{(D_1;\alpha_1)}\cdot F_{(B;\alpha_2)}}
=
\Phi_{B\setminus\{\alpha_1,\alpha_2\}}+
\hbar^{n+1}\left(\chi_{(B;\alpha_2)}+\chi_{(D_1;\alpha_1)})\right)
\end{equation}
By \eqref{eq:new eta} however,
\begin{equation*}
\begin{split}
F_{(D_1;\alpha_1)}\cdot F_{(B;\alpha_2)}
&=
(\Phi_{(B;\alpha_2,\alpha_1)})^{\otimes 2}\cdot
F_{(D_2;\alpha_2)}\cdot F_{(B;\alpha_1)}\cdot
\Delta(\Phi_{(B;\alpha_2,\alpha_1)}^{-1})\\
&+
\hbar^{n+1}\eta_{(B;\alpha_2,\alpha_1)}\mod\hbar^{n+2}
\end{split}
\end{equation*}
Since for any $\Phi\in(A_B^{\otimes 3})^B$, $F\in A_B^{\otimes 2}$
and $a\in A_B$,
$$(\Phi)_{a^{\otimes 2}\cdot F\cdot\Delta(a)^{-1}}=
a^{\otimes 3}\cdot(\Phi)_{F}\cdot (a^{\otimes 3})^{-1}$$
the \lhs of \eqref{eq:other path} is also equal mod $\hbar^{n+2}$ to
\begin{equation*}
\begin{split}
(\Phi_B)_{F_{(D_1;\alpha_1)}\cdot F_{(B;\alpha_2)}}
&=
\Ad(\Phi_{(B;\alpha_2,\alpha_1)}^{\otimes 3})
(\Phi_{B\setminus\{\alpha_1,\alpha_2\}}+
\hbar^{n+1}(\chi_{(B;\alpha_1)}+\chi_{(D_2;\alpha_2)}))\\
&+
\hbar^{n+1}\dH\eta_{(B;\alpha_2,\alpha_1)}\\
&=
\Phi_{B\setminus\{\alpha_1,\alpha_2\}}+
\hbar^{n+1}(\chi_{(B;\alpha_1)}+\chi_{(D_2;\alpha_2)}+
\dH\eta_{(B;\alpha_2,\alpha_1)})
\end{split}
\end{equation*}
where the last equality follows from the fact that $\Phi_{(B;\alpha_2,\alpha_1)}$
centralises $A_{B\setminus\{\alpha_1,\alpha_2\}}$ and the fact that $\Phi_
{B\setminus\{\alpha_1,\alpha_2\}}\in A_{B\setminus\{\alpha_1,\alpha_2\}}^
{\otimes 3}$. Comparing the two expressions for $(\Phi_B)_{F_{(D_1;\alpha_1)}
\cdot F_{(B;\alpha_2)}}$ yields
$$\dDD\chi_{(B;\alpha_2,\alpha_1)}=
\chi_{(B;\alpha_1)}-\chi_{(D_1,\alpha_1)}-\chi_{(B;\alpha_2)}+\chi_{(D_2,\alpha_2)}=
\dH\eta_{(B;\alpha_2,\alpha_1)}$$
as claimed \halmos

\subsection{}\label{ss:theta}

For any algebra homomorphism $\wt{\Delta}:A\rightarrow A^{\otimes 2}$
and element $\Phi\in A^{\otimes 3}$, set
$$\Pent_{\wt{\Delta}}(\Phi)=
\id^{\otimes 2}\otimes\wt{\Delta}(\Phi)\cdot
\wt{\Delta}\otimes\id^{\otimes 2}(\Phi)-
1\otimes\Phi\cdot\id\otimes\wt{\Delta}\otimes\id(\Phi)\cdot 1\otimes\Phi$$

Let now $B\subseteq D$ be a connected subdiagram and define $\theta_B
\in A^{\otimes 3}$ by
\begin{equation}\label{eq:bi theta}
\hbar^{n+1}\Theta_B=\Pent_{\Delta}(\Phi_B)\mod\hbar^{n+2}
\end{equation}

\begin{lemma}\hfill
\begin{enumerate}
\item $\theta$ is a Dynkin $0$--cocycle with values in $A^{\otimes 4}$.
\item $\dDD\theta=\dH\eta$.
\item $\dH\theta=0$.
\end{enumerate}
\end{lemma}
\proof
(i) we must prove that $\theta_B\in(A_B^{\otimes 4})^B$ which
readily follows from the fact that $\Phi_B\in(A_B^{\otimes 3})^B$.\\
(ii) Let $\alpha\in B\subseteq D$. One readily checks that
\begin{equation*}
\begin{split}
\Pent_{\Delta_{F_{(B;\alpha)}}}((\Phi_B)_{F_{(B;\alpha)}})
&=
1^{\otimes 2}\otimes F_{(B;\alpha)}\cdot 
1\otimes\id\otimes\Delta(F_{(B;\alpha)})\cdot
\id\otimes\Delta^{(3)}(F_{(B;\alpha)})\\
&\phantom{=}\cdot
\Pent_{\Delta}(\Phi_B)\\
&\phantom{=}\cdot
\Delta^{(3)}\otimes\id(F_{(B;\alpha)}^{-1})\cdot
\Delta\otimes\id(F_{(B;\alpha)}^{-1})\otimes 1\cdot
F_{(B;\alpha)}^{-1}\otimes 1^{\otimes 2}
\end{split}
\end{equation*}
By \eqref{eq:bi theta}, the \rhs of the above equation is equal to $\hbar^{n+1}
\theta_B$ mod $\hbar^{n+2}$. On the other hand, by \eqref{eq:new chi},
the \lhs is equal to
\begin{equation*}
\begin{split}
\Pent_{\Delta_{F_{(B;\alpha)}}}
(\Phi_{B\setminus\alpha}+\hbar^{n+1}\chi_{(B;\alpha)})
&=
\Pent_{\Delta_{F_{(B;\alpha)}}}(\Phi_{B\setminus\alpha})
+\hbar^{n+1}\dH\chi_{(B;\alpha)}\\
&=
\Pent_{\Delta}(\Phi_{B\setminus\alpha})
+\hbar^{n+1}\dH\chi_{(B;\alpha)}\\
&=
\hbar^{n+1}(\theta_{B\setminus\alpha}+\dH\chi_{(B;\alpha)})
\end{split}
\end{equation*}
where the second equality follows from the fact that $\Delta_{F_{(B;\alpha)}}$
restricts to $\Delta$ on $A_{B\setminus\alpha}$. Equating these two expressions,
we therefore get
$$\dDD\theta_{(B;\alpha)}=
\theta_B-\theta_{B\setminus\alpha}=\dH\eta_{(B;\alpha)}$$
(iii) It is known that an obstruction defined by \eqref{eq:bi theta} satisfies
$\dH\theta=0$ \cite[pp. 1448--9]{Dr3}. In the case at hand, a simpler proof
can be given owing to the fact that by remark \rkref{rk:reformulation},
$\Phi_B$ is a non--abelian \Ho coboundary. Let $B\subseteq D$ and
$\alpha_1,\ldots,\alpha_k$ an enumeration of the vertices of $B$. Then
\begin{equation*}
\begin{split}
\theta_B
&=
\sum_{i=0}^{k-1}
(\theta_{B\setminus\{\alpha_1,\ldots,\alpha_i\}}-
\theta_{B\setminus\{\alpha_1,\ldots,\alpha_{i+1}\}})\\
&=
\sum_{i=0}^{k-1}
\dDD\theta_{(B\setminus\{\alpha_1,\ldots,\alpha_i\};\alpha_{i+1})}\\
&=
\sum_{i=0}^{k-1}
\dH\eta_{(B\setminus\{\alpha_1,\ldots,\alpha_i\};\alpha_{i+1})}
\end{split}
\end{equation*}
so that $\theta_B$ is a \Ho coboundary \halmos

\subsection{}

Let now
$$\phi_{(D;\alpha,\beta)}\in A_B^{B\setminus\{\alpha,\beta\}}
\qquad
f_{(B;\alpha)}\in (A_B^{\otimes 2})^{B\setminus\alpha}
\qquad\text{and}\qquad
\psi_B\in (A_B^{\otimes 3})^B$$
with $\veps\otimes\id(f_{(B;\alpha)})=\id\otimes\veps(f_{(B;\alpha)})=0$.
The cocycle $(\xi',\eta',\chi',\theta')$ corresponding to
$$\Phi_{(B;\alpha,\beta)}+\hbar^{n+1}\phi_{(B;\alpha,\beta)}
\qquad
F_{(B;\alpha)}+\hbar^{n+1}f_{(B;\alpha)}
\qquad\text{and}\qquad
\Phi_B+\hbar^{n+1}\psi_B$$
is given by
\begin{align*}
\xi'&=\xi+\dDD\phi\\
\eta'&=\eta+\dDD f-\dH\phi\\
\chi'&=\chi+\dDD\psi+\dH f\\
\theta'&=\theta+\dH\psi
\end{align*}
so that the given \qc \qba structure lifts mod $\hbar^{n+2}$ if, and
only if, $(\xi,\eta,\chi,\theta)$ is a \DH coboundary. This concludes
the proof of (i).

\subsection{}

Let now $(\{S_i\},\{\Phi^j_{(B;\alpha,\beta)}\},\{F^j_{(B;\alpha)}\},
\{\Phi_B\})$, with $j=1,2$, be two \qc \qba structures on $A\fml/
\hbar^{n+1}A\fml$ which are trivial mod $\hbar$ and agree mod
$\hbar^n$. Define $\phi_{(B;\alpha,\beta)},f_{(B;\alpha)}$ and
$\psi_B$ by the following equalities mod $\hbar^{n+1}$
\begin{align*}
\Phi_{(B;\alpha,\beta)}^2
&=
\Phi_{(B;\alpha,\beta)}^1+\hbar^n\phi_{(B;\alpha,\beta)}\\
F_{(B;\alpha)}^2
&=
F_{(B;\alpha)}^1+\hbar^n f_{(B;\alpha)}\\
\Phi_B^2
&=
\Phi_B^1+\hbar^n\psi_B
\end{align*}
Then, linearlising the defining identities of a \qc \qba readily
yields that $(\phi,f,\psi)$ is a \DH 3--cocycle. It is easy to check
that $(\phi,f,\psi)$ is a coboundary if, and only if the two structures
differ by a twist equal to 1 mod $\hbar^n$ \halmos


%% file: transfer.tex

\section{$\Uhg$ as a \qcqtqba}\label{se:transfer}

Let $\g$ be a complex, simple Lie algebra and let $\Dg$ be its
Dynkin diagram. We point out in \ref{ss:DJ}--\ref{ss:DL} that the
quantum group $\Uhg$, when endowed with the \qW operators
and the universal $R$--matrices corresponding to
all subdiagrams of $\Dg$, has the structure of a \qcqtqba of type
$\Dg$ with trivial associators and relative twists. We then transfer
this structure to $\Ug\fml$. This
requires the cohomological construction of non--trivial associators
and structural twists and is similar in spirit to the fact that $\Uhg$
is twist equivalent to a \qtqba of the form $(\Ug\fml,\Delta_0,\exp
(\hbar\Omega),\Phi)$ where $\Delta_0$ is the cocommutative
coproduct on $\Ug$, $\Omega\in\g\otimes\g$ the Casimir operator
of $\g$ and $\Phi$ some associator. The proof is somewhat
lengthier however and occupies the rest of this section.

\subsection{}\label{ss:DJ}

Retain the notation of \S\ref{ss:qWeyl} and regard the quantum group
$\Uhg$ as a topological Hopf algebra over the ring of formal power
series $\ICh$ by endowing it with the coproduct given by
\begin{align*}
\Delta(E_{i})&=E_{i}\otimes 1+q_i^{H_{i}}\otimes E_{i}\\
\Delta(F_{i})&=F_{i}\otimes q_i^{-H_{i}}+1\otimes F_{i}\\
\Delta(H_{i})&=H_{i}\otimes 1+1\otimes H_{i}
\end{align*}

For any subdiagram $D\subseteq\Dg$, the operators $E_{i},
F_{i},H_{i}$, with $i$ such that $\alpha_i\in D$, topologically
generate a subalgebra $\UhgD\subseteq\Uhg$ canonically
isomorphic to the quantum group corresponding to $\gD$
and the restriction of the bilinear form $(\cdot,\cdot)$ to it.
Let
$$R\Dh\in 1^{\otimes 2}+\hbar\UhgD^{\otimes 2}$$
be the universal $R$--matrix of $\UhgD$ \cite{Dr1,Dr2}.
For $D=\alpha_i$, we denote $\UhgD$ by $\Uhsl{2}^{i}$
and $R\Dh$ by $R\ih$.\\

For any $\alpha_i\in\Dg$, let $S\ih\in\wh{\Uhsl{2}^i}$ be
the quantum Weyl group element defined by \eqref{eq:barSih}.
The following result is due to Lusztig, \KR and Soibelman
\cite{Lu,KR,So}.

\begin{proposition}\label{pr:qW ids}\hfill
\begin{enumerate}
\item The following holds in $\wh{\Uhsl{2}^{i}}$
$$(S\ih)^{2}=\exp(\sqrt{-1}\pi H_i)\cdot q^{C_i}$$
where $\exp(\sqrt{-1}\pi H_i)$ and $C_i$ are the
sign and Casimir operators of $\Uhsl{2}^{i}$, that
is the central elements of $\wh{\Uhsl{2}^{i}}$ acting
on the indecomposable representation $\V_{m}$ of
dimension $m+1$ as multiplication by $(-1)^{m}$
and $\half{(\alpha_i,\alpha_i)}\cdot\half{m(m+2)}$
respectively.
\item The following holds in $\wh{\Uhsl{2}^i}^
{\otimes 2}$
$$\Delta(S\ih)=(R\ih)^{21}\cdot S\ih\otimes S\ih$$
\end{enumerate}
\end{proposition}
\proof By \cite[prop. 5.2.2.(b)]{Lu}, $(\ol{S}\ih)^2$ acts on
the subspace of $\V_{m}$ of weight $\ol{\jmath}=-m+2j$,
$j=0\ldots m$, as multiplication by
$$(-1)^{m}q_i^{2j(m-j)+m}=
(-1)^{m}q_i^{\half{1}(m-\ol{\jmath})(m+\ol{\jmath})+m}=
(-1)^{m}q_i^{\half{m(m+2)}+\half{\ol{\jmath}}^{2}}$$
so that $(\ol{S}\ih)^2=\exp(\sqrt{-1}\pi H_i)q^{C_i}q_{i}^
{-H_i^^2/2}$. Since $\Ad(S\ih)(H_i)=-H_i$, (i) holds. (ii)
readily follows from the fact that, by \cite[prop. 5.3.4]{Lu},
$\ol{S}\ih$ satisfies
$$\Delta(\ol{S}\ih)=
(\ol{R}\ih)^{21}\cdot\ol{S}\ih\otimes\ol{S}\ih$$
where $\ol{R}\ih=q_{i}^{-\half{H_{i}\otimes H_{i}}}\cdot
R\ih$ \halmos

\subsection{}\label{ss:DL}

Label the Dynkin diagram $\Dg$ by attaching to each pair
$\alpha_i\neq\alpha_j$ the order $m_{ij}$ of the product
$s_is_j\in W$ of the corresponding simple reflections. The
following is an immediate corollary of proposition \ref{pr:qW ids}
and \cite[\S 13]{Dr1}.

\begin{proposition}
For any $\alpha_i\neq\alpha_j\in D\subseteq\Dg$, set
$$\Phi_{(D;\alpha_i,\alpha_j)}=1,
\qquad
F_{(D;\alpha_i)}=1^{\otimes 2}
\qquad\text{and}\qquad
\Phi_{D}=1^{\otimes 3}$$
Then,
$$\left(
\Uhg,\{\Uhg_D\},\{S\ih\},\{\Phi_{(D;\alpha_i,\alpha_j)}\},
\Delta,\{R\Dh\},\{F_{(D;\alpha_i)}\},\{\Phi_{D}\}
\right)$$
is a \qcqtqba of type $\Dg$. The corresponding braid group
representations are the quantum Weyl group representations
of $\Bg$ on \fd $\Uhg$--modules.
\end{proposition}

\subsection{}\label{ss:mod h}

For any $\alpha_i\in\Dg$, choose root vectors $e_i\in\g_{\alpha_i},
f_i\in\g_{-\alpha_i}$ such that $[e_i,f_i]=h_i$. Then, the assignement
$$E_i\rightarrow e_i,\quad
F_i\rightarrow f_i,\quad
H_i\rightarrow h_i$$
extends uniquely to an isomorphism of Hopf algebras
$$\Psi_{0}:\Uhg/\hbar\Uhg\longrightarrow\Ug$$
We shall say that a $\ICh$--linear map $\Psi:\Uhg\rightarrow
\Ug\fml$ is equal to the identity mod $\hbar$ if its
reduction mod $\hbar$ is equal to $\Psi_0$. Since
\fd $\g$--modules do not possess non--trivial deformations,
the canonical map $\wh{\Ug}\fml\rightarrow\wh{\Ug\fml}$ is
an isomorphism. Any algebra isomorphism $\Psi:\Uhg\rightarrow
\Ug\fml$ therefore extends to an isomorphism
$$\wh{\Uhg}\longrightarrow\wh{\Ug\fml}=\wh{\Ug}\fml$$
which we denote by the same symbol.

\subsection{}\label{ss:transfer}

Extend the bilinear form $(\cdot,\cdot)$ on $\h$ to a non--degenerate,
symmetric, bilinear, ad--invariant form on $\g$. For any subdiagram
$D\subseteq\Dg$, let $\g_{D}\subseteq\lD\subseteq\g$ be the
corresponding simple and Levi subalgebras. Denote by
$$\Omega_{D}=x_a\otimes x^a,
\quad 
C_{D}=x_a\cdot x^a
\quad\text{and}\quad
\rD=
\sum_{\substack{\alpha\succ 0 :\\ \supp(\alpha)\subseteq D}}
\half{(\alpha,\alpha)}\cdot
e_{\alpha}\wedge f_{\alpha}$$
where $\{x_a\}_a,\{x^a\}_a$ are dual basis of $\gD$
with respect to $(\cdot,\cdot)$, the corresponding invariant
tensor, Casimir operator and standard solution of the
modified classical Yang--Baxter equation (MCYBE) for
$\gD$ respectively. Abbreviate $\sl{2}^{\alpha_i}$,
$\Omega_{\alpha_i}$ and $C_{\alpha_{i}}$ to $\sl{2}^{i}$,
$\Omega_{i}$ and $C_{i}$ respectively and let $\wt{s}_i$
be the triple exponentials \eqref{eq:triple exp}.\\

The aim of this section is to prove the following

\begin{theorem}\label{th:transfer}
$\Uhg$ is equivalent to a \qcqtqba of type $\Dg$ of the form
$$\Bigl(
\Ug\fml,\{\Ug_D\fml\},\{S\ic\},\{\Phi_{(D;\alpha_i,\alpha_j)}\},
\Delta_{0},\{\Phi_{D}\},\{R\DKZ\},\{F_{(D;\alpha_i)}\}
\Bigr)$$
where $\Delta_{0}$ is the cocommutative coproduct
on $\Ug$,
\begin{align*}
S\ic	&= \wt{s}_{i}\cdot\exp(\hbar/2\cdot C_i)\\
\Phi_D	&= 1^{\otimes 3}\mod\hbar^2\\
R\DKZ	&= \exp(\hbar\cdot\Omega_{D})\\
\Alt_{2}F_{(D;\alpha_i)}&=
\hbar\cdot(\rD-r_{\g_{D\setminus\{\alpha_i\}}})
\medspace\mod\hbar^{2}
\end{align*}
and $\Phi_{(D;\alpha_i,\alpha_j)}$, $F_{(D;\alpha_i)}$ are
of weight 0.
\end{theorem}
\proof We begin by recursively constructing in
\S \ref{ss:rank 1}--\S \ref{ss:checking markov} two families
$\{\Psi_{(D;\F)}\}$ and $\{\Phi_{(D;\F,\G)}\}$ labelled by
connected subdiagrams $D\subseteq\Dg$ and (elementary
pairs of) \mnss on $D$ satisfying the following properties:
\begin{enumerate}
\item[(i)$_{D}$] For any \mns $\F$ on $D$,
$$\Psi_{(D;\F)}:\UhgD\longrightarrow U\gD\fml$$
is an algebra isomorphism equal to the identity
mod $\hbar$ and restricting to the identity on
$\hD$. Moreover, for any $B\in\F$,
$\Psi_{(D;\F)}$ restricts to $\Psi_{(B;\F_B)}$ on
$U_{\hbar}\g_B$. Lastly, for any
$\alpha_i\in\Dg$,
$$\Psi_{(\alpha_i;\alpha_i)}(S\ih)=S\ic$$
\item[(ii)$_{D}$] For any elementary pair $(\G,
\F)$ of \mnss on $D$, the associator
$$\Phi_{(D;\G,\F)}\in 1+\hbar\Ug_{D_1}\fml^{\ll_{D_2}}$$
where $D_1=\supp(\G,\F)\supset D_2=\zsupp(\G,\F)$,
satisfies
\begin{equation}\label{eq:Phi Psi}
\Psi_{(D;\G)}=\Ad(\Phi_{(D;\G,\F)})\circ\Psi_{(D;\F)}
\end{equation}
and $\Phi_{(D;\F,\G)}=\Phi_{(D;\G,\F)}^{-1}$.
Moreover, if $\alpha^{D}_{\G}=\alpha^{D}_{\F}
=\alpha_i$, then
\begin{equation}\label{eq:ind Phi}
\Phi_{(D;\G,\F)}=
\Phi_{(D\setminus\alpha_i;\G\setminus D,\F\setminus D)}
\end{equation}
\item[(iii)$_{D}$] For any pair of \mnss $\F,\G$ on $D$
and elementary sequences
$$\F=\H_{1},\H_{2},\ldots,\H_{l}=\G
\quad\text{and}\quad
\F=\K_{1},\K_{2},\ldots,\K_{m}=\G$$
one has
$$
\Phi_{(D;\H_{1},\H_{2})}\cdots\Phi_{(D;\H_{l-1},\H_{l})}=
\Phi_{(D;\K_{1},\K_{2})}\cdots\Phi_{(D;\K_{m-1},\K_{m})}
$$
\item[(iv)$_{D}$] For any equivalent elementary pairs of
\mnss $(\G,\F)$ and $(\G',\F')$ on $D$, one has
$$\Phi_{(D;\G,\F)}=\Phi_{(D;\G',\F')}$$
\end{enumerate}

Here and in the sequel, we follow the convention that the
isomorphisms and associators corresponding to non--connected
diagrams are the product of those corresponding to their
connected components. Specifically, let $\alpha_i\in D$
and let $D_1,\ldots,D_k$ be the connected components of
$D\setminus\alpha_i$, so that
$$\Uhg_{D\setminus\alpha_i}
\cong
\Uhg_{D_1}\odots{\otimes}\Uhg_{D_k}
\quad\text{and}\quad
\Ug_{D\setminus\alpha_i}\fml
\cong
\Ug_{D_1}\fml\odots{\otimes}\Ug_{D_k}\fml$$
If $\F$ is a \mns on $D$ with $\alpha^{D}_{\F}=\alpha_i$,
so that $$\F=\{D\}\sqcup\F_1\sqcup\cdots\F_k$$ where $\F_i$
is a \mns on $D_i$, we set
$$\Psi_{(D\setminus\alpha_i;\F\setminus D)}=
  \Psi_{(D_1;\F_1)}\otimes\cdots\otimes\Psi_{(D_k;\F_k)}$$
If $\G$ is another \mns on $D$ with $\alpha^{D}_{\G}=\root
{i}$, so that $$\G=\{D\}\sqcup\G_1\odots{\sqcup}\G_k$$ with
$\G_i$ a \mns on $D_i$, we set
$$\Phi_{(D\setminus\alpha_i;\G\setminus D;\F\setminus D)}=
\Phi_{(D_1;\G_1,\F_1)}\otimes\cdots\otimes\Psi_{(D_k;\G_k,\F_k)}$$

Once constructed, the associators $\Phi_{(\Dg;\G,\F)}$ endow
$\Ug\fml$ with the structure of a \qc algebra $\Q$ which is
equivalent, via the isomorphisms $\Psi_{(\Dg;\F)}$, to the
\qc structure on $\Uhg$ determined by the quantum Weyl group
operators $S\ih$. A suitable collection of associators $\Phi
_{D}$ and structural twists $F_{(D;\alpha_i)}$ promoting $\Q$
to a \qcqtqba structure on $\Ug\fml$ equivalent to that on $
\Uhg$ will then be constructed in \S \ref{ss:graft}--\S \ref
{ss:1-jet of F}.

\subsection{}\label{ss:rank 1}

We first construct, for any $\alpha_i\in\Dg$, an algebra
isomorphism
$$\Psi_{(\alpha_i;\alpha_i)}:
\Uhsl{2}^{i}\longrightarrow\Usl{2}^{i}\fml$$
equal to the identity mod $\hbar$ and mapping $H_{i}$ to
$h_{i}$ and $S\ih$ to $S\ic$.

\begin{lemma}\label{le:square}
Let $A$ be a complete, topological algebra over $\IC\fml$
and $a,b\in A$ two invertible elements such that 
$$a=b\mod\hbar
\qquad\text{and}\qquad
a^{2}=b^{2}$$
Then
$$b=g a g^{-1}
\qquad\text{where}\qquad
g=(ba^{-1})^{1/2}\in 1+\hbar A$$
\end{lemma}
\proof Let $\delta=ba^{-1}\in 1+\hbar A$ so that $b=\delta a$.
Then $b^{2}=a^{2}$ implies that $\delta a\delta=a$ and therefore
that $F(\delta)a=aF(\delta^{-1})$ for any formal power series
$F$. In particular $\delta^{1/2}a=a\delta^{-1/2}$ so that
$$b=\delta^{1/2}\delta^{1/2}a=\delta^{1/2}a\delta^{-1/2}$$
as claimed \halmos\\

Let $$\Psi_i:\Uhsl{2}^{i}\longrightarrow\Usl{2}^{i}\fml$$
be an algebra isomorphism equal to the identity mod $\hbar$
and mapping $H_{i}$ to $h_{i}$ \cite[Prop. 4.3]{Dr2} and
set $S_i=\Psi_i(S\ih)$. Then, $S_i=\wt{s}_i$ mod $\hbar$
and, by proposition \ref{pr:qW ids},
$$S_i^{2}=
\Psi_i((S\ih)^2)=
\Psi_i(\exp(\sqrt{-1}\pi H_i)\cdot q^{C_i})=
\exp(\sqrt{-1}\pi h_i)\cdot q^{C_i}=
S\ic^2$$
Thus, by lemma \ref{le:square},
$$\Psi_{(\alpha_i;\alpha_i)}=
\Ad(S\ic\cdot S_i^{-1})^{1/2}\circ\Psi_i$$
maps $S\ih$ to $S\ic$ and $H_{i}$ to $h_{i}$ since
$$\Ad(S_i)h_i=-h_i=\Ad(S\ic)h_i$$


\subsection{}\label{ss:Psi}

Assume now that, for some $1\leq m\leq|\Dg|-1$, the isomorphisms
$\Psi_{(D;\F)}$ and associators $\Phi_{(D;\F,\G)}$ have been
constructed for all $D$ with $|D|\leq m$ in such a way that
properties (i)$_D$--(iv)$_D$ hold. We now construct isomorphisms
$\Psi_{(D;\F)}$ for all $D$ with $|D|=m+1$ which satisfy
(i)$_D$. We shall need the following.

\begin{proposition}\label{pr:ext}
Let $D\subseteq\Dg$ be a subdiagram. Then, for any algebra
isomorphism
$$\Psi_{D}:\UhgD\longrightarrow\UgD\fml$$
equal to the identity mod $\hbar$, there exists an algebra
isomorphism
$$\Psi:\Uhg\longrightarrow\Ug\fml$$
equal to the identity mod $\hbar$ and restricting to $\Psi
_{D}$ on $\UhgD$. If $\Psi_{D}$ restricts to the identity
on $\hD$, then $\Psi$ may be chosen such that $\left.\Psi
\right|_{\h}=\id$.
\end{proposition}
\proof Let $\wt{\Psi}:\Uhg\longrightarrow\Ug\fml$ be an
algebra isomorphism equal to the identity mod $\hbar$.
Set
$$\tau=\wt{\Psi}\circ\Psi_{D}^{-1}:
\UgD\longrightarrow\Ug\fml$$
so that $\tau=\id+\hbar\tau_{1}\mod\hbar^2$ for some linear
map $\tau_1:\UgD\rightarrow\Ug$. Since $\tau$
is an algebra homomorphism, we readily find that, for
any $x,y\in\gD$,
$$\tau_{1}([x,y])=
[x,\tau_{1}(y)]+[\tau_{1}(x),y]=
\ad(x)\tau_{1}(y)-\ad(y)\tau_{1}(x)$$
so that the restriction of $\tau_{1}$ to $\gD$ is a
1--cocycle with values in $\Ug$ endowed with the adjoint
action of $\gD$. Since $H^{1}(\gD,\Ug)=0$, there exists
$a_{1}\in\Ug$ such that
$$\tau_{1}(x)=\ad(x)a_{1}=-[a_{1},x]$$
for any $x\in\gD$. It follows that
$$\Ad(1+\hbar a_{1})\circ\wt{\Psi}\circ\Psi_{D}^{-1}=
\id+\hbar^{2}\tau_{2}\mod\hbar^{3}$$
for some linear map $\tau_{2}:\UgD\rightarrow\Ug$.
Continuing in this way, we find a sequence of elements
$a_{n}\in\Ug$, $n\geq 2$ such that
$$\Ad(1+\hbar^{n}a_{n})\circ\cdots\circ\Ad(1+\hbar a_{1})
\circ\wt{\Psi}\circ\Psi_{D}^{-1}=
\id\mod\hbar^{n+1}$$
so that, setting
$$a=\lim_{n\rightarrow\infty}
(1+\hbar^{n}a_{n})\cdots(1+\hbar a_{1})
\in 1+\hbar\Ug\fml$$
and
$$\Psi=\Ad(a)\circ\wt{\Psi}:\Uhg\longrightarrow\Ug\fml$$
we find that $\Psi$ is an algebra isomorphism equal
to the identity mod $\hbar$ and extending $\Psi_{D}$.
If $\left.\Psi_{D}\right|_{\hD}=\id$, and $\wt{\Psi}$
is chosen such that $\wt{\Psi}_{|\h}=\id$ \cite
[prop. 4.3]{Dr2}, the obstructions $\tau_{i}$
constructed above are readily seen to be equivariant
for the adjoint actions of $\h$ on $\UgD$ and $\Ug$
so that the $a_{n}$, $n\geq 2$ may be chosen of
weight 0 thus implying that $\left.\Psi\right|_{\h}=
\id$ \halmos\\

Let now $D\subseteq\Dg$ be a connected subdiagram with
$|D|=m+1$. For any $\alpha_i\in D$, choose a reference
\mns $\F_i$ on $D$ such that $\alpha^{D}_{\F_i}=\root
{i}$ and, using proposition \ref{pr:ext}, an algebra
isomorphism
$$\Psi_{(D;\F_i)}:\UhgD\longrightarrow\UgD\fml$$
such that
$$\Psi_{(D;\F_i)}=\id\mod\hbar,
\qquad
\left.\Psi_{(D;\F_i)}
\right|_{\hD}=\id$$
and
$$\left.\Psi_{(D;\F_i)}\right|
_{\Uhg_{D\setminus\alpha_i}}=
\Psi_{(D\setminus\alpha_i;\F_i\setminus D)}$$
For any \mns $\F$ on $D$ with $\alpha^{D}_{\F}=\alpha_i$,
set
\begin{equation}\label{eq:Psi def}
\Psi_{(D;\F)}=
\Ad(\Phi_{(D\setminus\alpha_i;\F\setminus D,\F_i\setminus D)})
\circ
\Psi_{(D;\F_i)}
\end{equation}
We claim that $\Psi_{(D;\F)}$ satisfies (i)$_D$. Since
$\Phi_{(D\setminus\alpha_i;\F\setminus D,\F_i\setminus D)}$
is of weight 0, this amount to showing that the restriction
of $\Psi_{(D;\F)}$ to $\Uhg_{D\setminus\alpha_i}$ is
equal to $\Psi_{(D\setminus\alpha_i;\F\setminus D)}$.
By construction, this restriction is equal to
$$\Ad(\Phi_{(D\setminus\alpha_i;\F\setminus D,\F_i\setminus D)})
\circ
\Psi_{(D\setminus\alpha_i;\F_i\setminus D)}$$
which, by (ii)$_{D\setminus\alpha_i}$, is equal to $\Psi
_{(D\setminus\alpha_i;\F\setminus D)}$. 

\subsection{}\label{ss:Phi}

We next construct associators $\Phi_{(D;\G,\F)}$ and
prove that they satisfy property (ii)$_D$. For any $
\alpha_i,\alpha_j\in D$, the automorphism $\tau_{ji}
=\Psi_{(D;\F_j)}\circ\Psi_{(D;\F_i)}^{-1}$ of $\Ug\fml$
is equal to the identity mod $\hbar$ and fixes $\hD$.
Since $H^{1}(\gD,\UgD)=0$, $\tau_{ji}$ is inner and
there exists an element
$$\Phi_{(D;\F_j,\F_i)}\in 1+\hbar\UgD\fml^{\hD}$$
such that
\begin{equation}\label{eq:Phi_ji def}
\Psi_{(D;\F_j)}=
\Ad(\Phi_{(D;\F_j,\F_i)})\circ
\Psi_{(D;\F_i)}
\end{equation}
We choose the associators $\Phi_{(D;\F_j,\F_i)}$ in
such a way that
$$\Phi_{(D;\F_j,\F_i)}=
\Phi_{(D;\F_i,\F_j)}^{-1}
\quad\text{and}\quad
\Phi_{(D;\F_i,\F_i)}=1$$
For any pair $\F,\G$ of \mnss on $D$ with $\alpha^{D}
_{\F}=\alpha_i$ and $\alpha^{D}_{\G}=\alpha_j$, set
\begin{equation}\label{eq:Phi def}
\Phi_{(D;\G;\F)}
=
\Phi_{(D\setminus\alpha_j;\G\setminus D,\F_j\setminus D)}
\cdot
\Phi_{(D;\F_j,\F_i)}
\cdot
\Phi_{(D\setminus\alpha_i;\F_i\setminus D,\F\setminus D)}
\end{equation}
so that, by (iii)$_{D\setminus\alpha_i}$,
\begin{equation}\label{eq:old Phi}
\Phi_{(D;\G;\F)}=
\Phi_{(D\setminus\alpha_i,\G\setminus D,\F\setminus D)}
\end{equation}
whenever $\alpha^{D}_{\G}=\alpha_i=\alpha^{D}_{\F}$.
We claim that these associators satisfy (ii)$_{D}$.
For any $\F,\G$ with $\alpha^{D}_{\F}=\alpha_i$, $\alpha
^{D}_{\G}=\alpha_j$, we have,
\begin{equation}\label{eq:Phi GF}
\begin{split}
\Psi_{(D;\G)}
&=
\Ad(\Phi_{(D;\G;\F_j)})\circ\Psi_{(D;\F_j)}\\
&=
\Ad(\Phi_{(D;\G;\F_j)}\cdot\Phi_{(D;\F_j,\F_i)})
\circ\Psi_{(D;\F_i)}\\
&=
\Ad(\Phi_{(D;\G;\F_j)}\cdot\Phi_{(D;\F_j,\F_i)}\cdot\Phi_{(D;\F_i,\F)})
\circ\Psi_{(D;\F)}\\
&=
\Ad(\Phi_{(D;\G;\F)})\circ\Psi_{(D;\F)}
\end{split}
\end{equation}
where the first and third equalities follow from \eqref
{eq:Psi def}, the second from \eqref{eq:Phi_ji def} and 
the last one from \eqref{eq:Phi def}. Set now
$$D_1=\supp(\F,\G)\qquad\text{and}\qquad D_2=\zsupp(\F,\G)$$
We claim that $\Phi_{(D;\G,\F)}$ lies in $\Ug_{D_1}\fml$
and is invariant under $\g_{D_{2}}$. It suffices to show
that $\Phi_{(D;\G,\F)}\in\Ug_{D_1}\fml$ since, by \eqref
{eq:Phi GF} and the fact that, by (i)$_D$, $\Psi_{(D,\G)}$
and $\Psi_{(D,\F)}$ have the same restriction on $\Ug_{D_2}$,
$\Phi_{(D;\G,\F)}$ centralises $\Ug_{D_2}$.
Now if $\alpha^{D}_{\G}=\alpha^{D}_{\F}$, then 
\begin{equation}
\supp(\G,\F)=\supp(\G\setminus D,\F\setminus D)
\end{equation}
and the claim follows from the inductive assumption
and \eqref{eq:old Phi}. If, on the other hand $\alpha
^{D}_{\G}\neq\alpha^{D}_{\F}$, then $\supp(\G,\F)=D$
and there is nothing to prove.

\subsection{}\label{ss:enforcing pentagon}

We now modify the associators $\Phi_{(D;\G,\F)}$ so that
they also satisfy (iii)$_{D}$. Introduce to this end some
terminology. Call a \fns $\H$ on $D$ {\it old} (resp. {\it
new}) if
$$|D\setminus\bigcup_{B\in\H\setminus D}B|=1
\quad\text{(resp. $\geq 2$)}$$
For example, if $(\G,\F)$ is an elementary pair of \mnss on $D$,
then $\H=\G\cap\F$ is old precisely
when $\alpha^{D}_{\G}=\alpha^{D}_{\F}$ and therefore when the
associator $\Phi_{(D;\G,\F)}$ is inductively determined by \eqref
{eq:ind Phi}. If, on the other hand, $\H$ is new then $\supp
(\G,\F)=D$ and $\Phi_{(D;\G,\F)}$ is determined by \eqref
{eq:Phi Psi} only up to multiplication by an element of
$$\zeta_{(D;\G,\F)}\in 1+\hbar\cdot Z(\UgD)\fml$$ Our
goal is to modify these new associators by suitable elements
$\zeta_{(D;\G,\F)}$ while keeping the old ones fixed, in
such a way that the generalised pentagon identities corresponding
to the two--faces of the associahedron $\AD$ hold. Note first
the following straightforward

\begin{lemma}
The collection of faces $\AD\old\subseteq\AD$ corresponding to
old \mnss is a subcomplex of $\AD$. For any $\alpha_i\in D$,
let $D^i_1,\ldots,D^i_{k_i}$ be the connected components of $D
\setminus\alpha_i$ and set
$$\A_{D\setminus\alpha_i}=\A_{D_1}\odots{\times}\A_{D_{k_i}}$$
Then, the map $\A_{D\setminus\alpha_i}\longrightarrow\AD\old$
given by 
$$(\H_1,\ldots,\H_{k_i})\longrightarrow\{D\}\sqcup 
\H_1\odots{\sqcup}\H_{k_i}$$
yields an isomorphism
$$\AD\old\cong\bigsqcup_{\alpha_i\in D}\A_{D\setminus\alpha_i}$$
\end{lemma}

Let now $\Sigma$ be an oriented two--face of $\AD$ with vertices
$$\F_{1},\ldots,\F_{k},\F_{k+1}=\F_{1}$$ listed in their order of
appearence along the boundary of $\Sigma$. Thus, for each $i=1
\ldots k$, $(\F_{i+1},\F_i)$ is an elementary pair of \mnss on $D$
and we may set
$$\zeta(\Sigma)=
\Phi_{(D;\F_{k+1},\F_{k})}\cdots\Phi_{(D;\F_{2},\F_{1})}$$

\begin{proposition}\hfill
\begin{enumerate}
\item The element $\zeta(\Sigma)$ lies in $\ZD=1+\hbar Z(\UgD)
\fml$, only depends upon the orientation of $\Sigma$ and
satisfies $\zeta(-\Sigma)=\zeta(\Sigma)^{-1}$.
\item The assignement $\Sigma\rightarrow\zeta(\Sigma)$ defines
a 2-cocycle on $\AD$ relative to the subcomplex $\AD\old$ with
coefficients in the abelian group $\ZD$.
\end{enumerate}
\end{proposition}
\proof (i) By (ii)$_{D}$,
\begin{equation*}
\begin{split}
\Ad(\zeta(\Sigma))
&=
\Ad(\Phi_{(D;\F_{k+1},\F_{k})}\cdots\Phi_{(D;\F_{2},\F_{1})})\\
&=
\Psi_{(D;\F_{k+1})}\circ\Psi_{(D,\F_{k})}^{-1}
\circ\cdots\circ
\Psi_{(D;\F_{2})}\circ\Psi_{(D,\F_{1})}^{-1}\\
&=
\id
\end{split}
\end{equation*}
since $\F_{k+1}=\F_{1}$, so that $\zeta(\Sigma)\in\ZD$ as claimed.
It follows that $\zeta(\Sigma)$ only depends upon the orientation
of $\Sigma$ since $\Phi_{(D;\F_2,\F_1)}$ commutes with $\zeta(\Sigma)$
and we therefore have
$$\Phi_{(D;\F_{k+1},\F_{k})}\cdots\Phi_{(D;\F_{2},\F_{1})}
=
\Phi_{(D;\F_{2},\F_{1})}\cdot\Phi_{(D;\F_{k+1},\F_{k})}
\cdots\Phi_{(D;\F_{3},\F_{2})}$$
(ii)\footnote{I owe this proof to G. Skandalis} We claim that
$\zeta(\Sigma)$ is of the form $\wt{\zeta}(\partial\Sigma)$
where $\wt{\zeta}$ is a homomorphism mapping one--chains
in $\AD$ to $\ZD$, so that $d\zeta=\wt{\zeta}\circ\partial\circ
\partial=0$. Note first that we may attach an element $z(p)
\in\ZD$ to any closed edge--path in $\AD$ \ie a sequence
$p=(\F_2,\F_1),(\F_3,\F_2),\ldots,(\F_{n+1},\F_n)$ of elementary
pairs of \mnss on $D$ such that $\F_{n+1}=\F_1$ by setting
$$z(p)=\Phi_{(D;\F_1,\F_n)}\cdots\Phi_{(D;\F_2,\F_1)}$$
Fix now a \mns $\F_0$ on $D$ and, for each \mns $\F$
on $D$, an edge--path $p_{\F}$ from $\F_0$ to $\F$. For
any oriented 1--face $e=(\G,\F)$ of $\AD$, set
$$\wt{\zeta}(e)=z(p_{\G}^{-1}\vee e\vee p_{\F})$$
where $p_{\G}^{-1}$ is the edge--path from $\G$ to
$\F_0$ obtained by reversing the orientation of $p
_{\G}$ and $\vee$ is the concatenation. It is clear
that $\zeta(\Sigma)=\wt{\zeta}(\partial\Sigma)$ so that
$\zeta$ is a two--cocycle on $\AD$ which, by the inductive
assumption is equal to 1 on the 2--faces of $\AD\old$ \halmos\\

Since $\AD$ and $\AD\old$ are contractible, $H^{2}(\AD,
\AD\old;\ZD)=1$. Thus, there exists a 1--cochain $\xi$
on $\AD$ such that
\begin{equation}
d\xi=\zeta
\qquad\text{and}\qquad
\xi(\G,\F)=1
\end{equation}
whenever $(\G,\F)$ is an elementary pair of \mnss on
$D$ such that $\supp(\G,\F)\subsetneq D$. Replacing
each $\Phi_{(D;\G,\F)}$ by $\Phi_{(D;\G,\F)}\cdot\xi
(\G,\F)^{-1}$ yields a collection of associators
satisfying (ii)$_{D}$ and (iii)$_{D}$.

\subsection{}\label{ss:checking markov}

We now show that the associators $\Phi_{(D;\G,\F)}$ satisfy
property (iv)$_D$. Let $(\G,\F)$ and $(\G',\F')$ be two
equivalent elementary pairs of \mnss on $D$. 
If $\supp(\G,\F)\subsetneq D$, then
$$\alpha^{D}_{\G}=\alpha^{D}_{\F}=\alpha_i=
  \alpha^{D}_{\F'}=\alpha^{D}_{\G'}$$
for some $\alpha_i\in D$ and, by (ii)$_{D\setminus\alpha_i}$
and (iv)$_{D\setminus\alpha_i}$,
$$\Phi_{(D;\G,\F)}=
\Phi_{(D\setminus\alpha_i;\G\setminus D,\F\setminus D)}=
\Phi_{(D\setminus\alpha_i;\G'\setminus D,\F'\setminus D)}=
\Phi_{(D;\G',\F')}$$
Assume now that $\supp(\G,\F)=D=\supp(\G',\F')$ and set
$$\alpha_i=\alpha^{D}_{\F}=\alpha^{D}_{\F'}
\qquad\text{and}\qquad
\alpha_j=\alpha^{D}_{\G}=\alpha^{D}_{\G'}$$

\begin{lemma}\label{le:squares}
There exist two sequences
$$\F=\F_1,\ldots,\F_m=\F'
\qquad\text{and}\qquad
\G=\G_1,\ldots,\G_m=\G'$$
of \mnss on $D$ such that, for any $i=1\ldots m-1$, the
following holds
\begin{enumerate}
\item $(\F_i,\F_{i+1})$ and $(\G_{i},\G_{i+1})$ are equivalent
elementary pairs of \mnss on $D$ such that
$$\supp(\F_i,\F_{i+1})=
\supp(\G_i,\G_{i+1})\subseteq
D\setminus\{\alpha_i,\alpha_j\}$$
\item $(\G_i,\F_i)$ and $(\G_{i+1},\F_{i+1})$ are equivalent
elementary pairs of \mnss on $D$.
\end{enumerate}
\end{lemma}
\proof Let $D_1,\ldots,D_p$ be the connected components of
$D\setminus\{\alpha_i,\alpha_j\}$, so that
\begin{xalignat*}{2}
\F &=\H_1\sqcup\cdots\sqcup\H_p\sqcup\{\con{D}{\alpha_i}{\alpha_j}\}&
\F'&=\H'_1\sqcup\cdots\sqcup\H'_p\sqcup\{\con{D}{\alpha_i}{\alpha_j}\}\\
\G &=\H_1\sqcup\cdots\sqcup\H_p\sqcup\{\con{D}{\alpha_j}{\alpha_i}\}&
\G'&=\H'_1\sqcup\cdots\sqcup\H'_p\sqcup\{\con{D}{\alpha_j}{\alpha_i}\}
\end{xalignat*}
where $\H_k,\H'_k$ are \mnss on $D_k$. By connectedness of $\A_{D_k}$, 
there exists an elementary sequence $\H_k=\H_k^{1},\ldots,\H_k^{
m_k}=\H'_k$ of \mnss on $D_k$. Setting $m=m_1+\cdots+m_k$ and
\begin{align*}
\F_i&=
\H'_1\sqcup\cdots\H'_{k-1}\sqcup
\H_{k}^{i-m_1-\cdots-m_k-1}\sqcup
\H_{k+1}\sqcup\cdots\H_{p}\sqcup
\{\con{D}{\alpha_i}{\alpha_j}\}\\
\G_i&=
\H'_1\sqcup\cdots\H'_{k-1}\sqcup
\H_{k}^{i-m_1-\cdots-m_k-1}\sqcup
\H_{k+1}\sqcup\cdots\H_{p}\sqcup
\{\con{D}{\alpha_j}{\alpha_i}\}
\end{align*}
for any $m_1+\cdots+m_{k-1}+1\leq i\leq m_1+\cdots+m_{k}$ yields
the required sequences \halmos\\

By (iii)$_D$,
\begin{equation*}
\begin{split}
\Phi_{(D;\G,\F)}\cdot
\Phi_{(D;\F_1,\F_2)}\cdots\Phi_{(D;\F_{m-1},\F_{m})}
&=
\Phi_{(D;\G,\G_2)}\cdots\Phi_{(\G_{m-1},\G_{m})}\cdot
\Phi_{(D;\G',\F')}\\
&=
\Phi_{(D;\G',\F')}\cdot
\Phi_{(D;\G_1,\G_2)}\cdots\Phi_{(G_{m-1},\G_{m})}\\
&=
\Phi_{(D;\G',\F')}\cdot
\Phi_{(D;\F_1,\F_2)}\cdots\Phi_{(F_{m-1},\F_{m})}
\end{split}
\end{equation*}
where the second equality follows from the fact that $\Phi_
{(D;\G,\F)}$ commutes with
$$\Ug_{\zsupp(\G,\F)}\fml\ni\Phi_{(D;G_i,G_{i+1})}$$
and the last one from (iv)$_{D\setminus\{\alpha_i,\alpha_j\}
}$ and the fact that $(\G_i,\G_{i+1})$ and $(\F_i,\F_{i+1})$
are equivalent pairs.

\subsection{}\label{ss:graft}

We next graft on to the previously constructed isomorphisms
$\Psi_{(D;\F)}$ and associators $\Phi_{(D;\G,\F)}$ a collection
$\{F_{(D;\alpha_i)}\}$ of relative twists such that, for
any connected subdiagram $D\subseteq\Dg$, the following
properties hold

\begin{enumerate}
\item[(v)$_{D}$] For any $\alpha_i\in D$,
$$F_{(D;\alpha_i)}\in 1^{\otimes 2}+
\hbar(\UgD^{\otimes 2}\fml)^{\ll_{D\setminus\alpha_i}}$$
\item[(vi)$_{D}$] For any \mns $\F$ on $D$,
\begin{equation}\label{eq:twisted Delta}
\Psi_{(D;\F)}\tens{2}\circ\Delta\circ\Psi_{(D;\F)}^{-1}=
\Ad(F_{(D;\F)})\circ\Delta_{0}
\end{equation}
where, as customary
$$F_{(D;\F)}=\stackrel{\longrightarrow}{\prod_{B\in\F}}
F_{(B;\alpha^{B}_{\F})}$$
\item[(vii)$_{D}$] For any pair $(\G,\F)$ of \mnss on $D$
$$F_{(D;\G)}=
\Phi_{(D;\G,\F)}^{\otimes 2}\cdot
F_{(D;\F)}\cdot
\Delta_{0}(\Phi_{(D;\F,\G)})$$
\end{enumerate}

\subsection{}

Assume that, for some $0\leq m\leq |\Dg|-1$, the relative
twists $F_{(D;\alpha_i)}$ have been constructed for all $
D$ with $|D|\leq m$ in such a way that properties (v)$_D
$--(vii)$_D$ hold. Let $D\subseteq\Dg$ be a connected
subdiagram such that $|D|=m+1$. For any \mns $\F$ on $D$,
denote by
$$\Delta_{(D;\F)}:\UgD\fml\longrightarrow\UgD\tens{2}\fml$$
the algebra homomorphism defined by the left--hand side of
\eqref{eq:twisted Delta}. Note that if $\F,\G$ are \mnss
on $D$, property (ii)$_{D}$ of \S \ref{ss:transfer} implies
that
\begin{equation}\label{eq:change of Delta}
\Delta_{(D;\G)}=
\Ad(\Phi_{(D;\G,\F)}^{\otimes 2})\circ
\Delta_{(D;\F)}\circ
\Ad(\Phi_{(D;\F,\G)})
\end{equation}

Fix $\alpha_i\in D$ and a \mns $\F_i$ on $D$ such that
$\alpha^{D}_{\F_i}=\alpha_i$. Since $\Delta_{(D;\F_i)}=
\Delta_{0}$ mod $\hbar$ and $H^{1}(\gD,\UgD\tens{2})=0$,
where $\UgD$ is regarded as a $\gD$--module under the
adjoint action, there exists a twist $$F_i\in 1^{\otimes
2}+\hbar\UgD\tens{2}\fml$$ such that
\begin{equation}\label{eq:Delta Fi}
\Delta_{(D;\F_i)}=\Ad(F_i)\circ\Delta_{0}
\end{equation}
This implies in particular that $F_i$ is invariant under
$\hD$ since $\Delta_{(D;\F_i)}$ and $\Delta_{0}$ coincide
on $\hD$.\\

For any $\alpha_i\neq\alpha_j\in D$, choose a
\mns $\F_j$ on $D$ such that $\alpha^{D}_{\F_j}=\alpha_j
$ and set
\begin{equation}\label{eq:F_j}
F_j=
\Phi_{(D;\F_j,\F_i)}\tens{2}\cdot
F_i\cdot
\Delta_{0}(\Phi_{(D;\F_i,\F_j)})\in
1^{\otimes 2}+\hbar(\UgD\tens{2}\fml)^{\hD}
\end{equation}
It follows from \eqref{eq:Delta Fi} and \eqref{eq:change of Delta}
that, for any $\alpha_j\in D$,
\begin{equation}\label{eq:Delta Fj}
\Delta_{(D;\F_j)}=\Ad(F_j)\circ\Delta_{0}
\end{equation}

For any such $\alpha_j$, set
$$F_{(D;\alpha_j)}=
F_{(D\setminus\alpha_j;\F_j\setminus D)}^{-1}\cdot
F_j\in 1^{\otimes 2}+\hbar(\UgD\tens{2}\fml)^{\hD}$$
We claim that these relative twists satisfy (v)$_{D}
$--(vii)$_{D}$

\begin{proposition}\hfill
\begin{enumerate}
\item $F_{(D;\alpha_j)}$ is invariant under $\ll_
{D\setminus \alpha_j}$.
\item For any pair $(\G,\F)$ of \mnss on $D$
$$F_{(D;\G)}=
\Phi_{(D;\G,\F)}\tens{2}\cdot
F_{(D;\F)}\cdot
\Delta_{0}(\Phi_{(D;\F,\G)})$$
\item For any \mns $\F$ on $D$,
$$\Delta_{(D;\F)}=\Ad(F_{(D;\F)})\circ\Delta_{0}$$
\end{enumerate}
\end{proposition}
\proof 
(i) By \eqref{eq:Delta Fj},
$$\Ad(F_{(D;\alpha_j)})\circ\Delta_{0}=
\Ad(F_{(D\setminus\alpha_j;\F_j\setminus D)}^{-1})\circ
\Delta_{(D;\F_j)}$$
By (i)$_D$ and (vi)$_{D\setminus\alpha_j}$, the right--hand
side restricts to $\Delta_{0}$ on $\Ug_{D\setminus\alpha_j}$.
This implies the invariance of $F_{(D;\alpha_j)}$ under
$\g_{D\setminus\alpha_j}$.
(ii) The stated identity certainly holds if $\F=\F_i$ and
$\G=\F_j$ for some $\alpha_j\in D$ since in that case $F_{(
D;\F)}=F_i$ and $F_{(D;\G)}=F_j$ is given by \eqref{eq:F_j}.
By transitivity of the associators, it therefore suffices
to check it when $\alpha^{D}_{\F}=\alpha_j=\alpha^{D}_{\G}$
for some $\alpha_j\in D$. In that case,
\begin{equation*}
\begin{split}
F_{(D;\G)}
&=
F_{(D\setminus\alpha_j;\G\setminus D)}\cdot F_{(D;\alpha_j)}\\
&=
\Phi_{(D\setminus\alpha_j;\G\setminus D,\F\setminus D)}\tens{2}\cdot
F_{(D\setminus\alpha_j;\F\setminus D)}\cdot
\Delta_{0}(\Phi_{(D\setminus\alpha_j;\F\setminus D,\G\setminus D)}) \cdot F_{(D;\alpha_j)}\\
&=
\Phi_{(D\setminus\alpha_j;\G\setminus D,\F\setminus D)}\tens{2}\cdot
F_{(D\setminus\alpha_j;\F\setminus D)}\cdot
F_{(D;\alpha_j)}\cdot
\Delta_{0}(\Phi_{(D\setminus\alpha_j;\F\setminus D,\G\setminus D)})\\
&=
\Phi_{(D;\G,\F)}\tens{2}\cdot
F_{(D;\F)}\cdot
\Delta_{0}(\Phi_{(D;\F,\G)})
\end{split}
\end{equation*}
where the second equality follows by (vii)$_{D\setminus\root
{j}}$, the third one by the invariance of $F_{(D;\alpha_j)}$
under $\g_{D\setminus\alpha_j}$ and the last one from property
(ii)$_{D}$ of \S \ref{ss:transfer}.
(iii) Let $\alpha_j=\alpha^{D}_{\F}$. By \eqref{eq:Delta Fj},
the stated identity holds if $\F=\F_j$ since in that case
$$F_{(D;\F)}=F_{(D\setminus\alpha_j;\F_j\setminus
D)}\cdot F_{(D;\alpha_j)}=F_j$$ In the general case, we have,
by \eqref{eq:change of Delta} and (ii)
\begin{equation*}
\begin{split}
\Delta_{(D;\F)}
&=
\Ad(\Phi_{(D;\F,\F_j)}\tens{2})\circ
\Delta_{(D;\F_j)}\circ
\Ad(\Phi_{(D;\F_j,\F)})\\
&=
\Ad(\Phi_{(D;\F,\F_j)}\tens{2})\circ
\Ad(F_{(D;\F_j)})\circ
\Delta_{0}\circ
\Ad(\Phi_{(D;\F_j,\F)})\\
&=
\Ad(F_{(D;\F)})\circ
\Ad(\Delta_{0}(\Phi_{(D;\F,\F_j)}))\circ
\Delta_{0}\circ
\Ad(\Phi_{(D;\F_j,\F)})\\
&=
\Ad(F_{(D;\F)})\circ
\Delta_{0}
\end{split}
\end{equation*}
\halmos

\subsection{}\label{ss:Phi_D R_D}

We now construct associators $\Phi_{D}$ and $R$--matrices
$R_D$ such that for any connected subdiagram $D\subseteq
\Dg$, the following holds
\begin{enumerate}
\item[(viii)$_{D}$] $\Phi_{D}\in 1\tens{3}+\hbar(\UgD
\tens{3}\fml)^{\gD}$ satisfies the pentagon equations
with respect to $\Delta_{0}$ and, for any \mns $\F$ on
$D$, $$(\Phi_{D})_{F_{(D;\F)}}=1\tens{3}$$
\item[(ix)$_{D}$] $R_{D}\in 1\tens{2}+\hbar(\UgD\tens{2}
\fml)^{\gD}$ satisfies the hexagon equations with respect
to $\Delta_{0}$ and $\Phi_{D}$ and, for any \mns $\F$
on $D$,
$$(R_D)_{F_{(D;\F)}}=\Psi_{(D;\F)}^{\otimes 2}(R\Dh)$$
\end{enumerate}
For any \mns $\F$ on $D$, set
\begin{equation*}
\begin{split}
\Phi_{(D;\F)}
&=
1\otimes F_{(D;\F)}^{-1}\cdot
\id\otimes\Delta_{(D;\F)}(F_{(D;\F)}^{-1})\cdot
\Delta_{(D;\F)}\otimes\id(F_{(D;\F)})\cdot
F_{(D;\F)}\otimes 1\\
&\in
1\tens{3}+\hbar\UgD^{\otimes 3}\fml\\
\intertext{and}
R_{(D;\F)}
&=
(F_{(D;\F)}^{-1})^{21}\cdot
\Psi_{(D;\F)}^{\otimes 2}(R\Dh)\cdot
F_{(D;\F)}\\
&\in
1\tens{2}+\hbar\UgD\tens{2}\fml
\end{split}
\end{equation*}
so that
\begin{equation*}
\begin{split}
(\UgD\fml,\Delta_{0},\Phi_{(D;\F)},R_{(D;\F)})
&=
(\UgD\fml,\Delta_{(D;\F)},1^{\otimes 3},\Psi_{(D;\F)}^{\otimes 2}(R\Dh))
_{F_{(D;\F)}^{-1}}\\
&\cong
(\UhgD,\Delta,1^{\otimes 3},R\Dh)_
{{{\Psi_{(D;\F)}^{-1}}\tens{2}}(F_{(D;\F)}^{-1})}
\end{split}
\end{equation*}
is a quasitriangular quasibialgebra. In particular, $\Phi_{(D;\F)}$ and $R_
{(D;\F)}$ satisfy the pentagon and hexagon equations with
respect to $\Delta_{0}$ and are invariant under $\gD$ since
$\Delta_{0}$ is coassociative and cocommutative. We claim
that $\Phi_{(D;\F)}$ and $R_{(D;F)}$ are independent of the
choice of $\F$, so that (viii)$_{D}$ and (ix)$_{D}$ hold
with $\Phi_{D}=\Phi_{(D;\F)}$ and $R_D=R_{(D;\F)}$ respectively.
From (vii)$_{D}$, one readily finds that
\begin{align*}
\Phi_{(D;\G)}
&=
\Delta^{(3)}_{0}(\Phi_{(D;\G,\F)})\cdot
\Phi_{(D;\F)}\cdot
\Delta^{(3)}_{0}(\Phi_{(D;\F,\G)})\\
R_{(D;\G)}
&=
\Delta_{0}(\Phi_{(D;\G,\F)})\cdot
R_{(D;\F)}\cdot
\Delta_{0}(\Phi_{(D;\F,\G)})
\end{align*}
where
$$\Delta^{(3)}_{0}=
\Delta_{0}\otimes\id\circ\Delta_{0}=
\id\otimes\Delta_{0}\circ\Delta_{0}:\UgD
\longrightarrow\UgD^{\otimes 3}$$
Thus $\Phi_{(D;\G)}=\Phi_{(D;\F)}$ and $R_{(D;\G)}=R_{(D;\F)}$
since $\Phi_{(D;\F)}$ and $R_{(D;\F)}$ are invariant under $\gD$.

\subsection{The coproduct identity}

Note that, for any $\alpha_i\in\Dg$, we have
\begin{equation*}
\begin{split}
\Ad(F_{(\alpha_i;\alpha_i)})
\Delta_{0}(S\ic)
&=
\Delta_{(\alpha_i;\alpha_i)}(\Psi_{(\alpha_i;\alpha_i)}(S\ih))\\
&=
\Psi_{(\alpha_i;\alpha_i)}\tens{2}
\left((R\ih)^{21}\cdot S\ih\otimes S\ih\right)\\
&=
(R_{\alpha_i})^{21}_{F_{(\alpha_i;\alpha_i)}}\cdot
S\ic\otimes S\ic
\end{split}
\end{equation*}

Thus, the relative twists, associators and $R$--matrices
constructed in \S \ref{ss:graft}--\S \ref{ss:Phi_D R_D}
endow $\Ug\fml$ with the structure of a \qcqtqba $\Q$
which extends the \qca  structure constructed in \S \ref
{ss:rank 1}--\S \ref{ss:checking markov} and is isomorphic,
via the isomorphisms $\Psi_{(D;\F)}$, to the \qcqtqba
structure $\Q_{\hbar}$ on $\Uhg$. In the next two
subsections, we apply suitable $F$--twists to $\Q$
which, while clearly preserving its equivalence to
$\Q_{\hbar}$, bring the $R$--matrices and associators
to the form required by the statement of theorem \ref
{th:transfer}.

\subsection{Symmetrising the $R$-matrices $R_D$}
\label{ss:sym R}

By proposition 3.16. of \cite{Dr3}, there exists, for
each $D\subseteq\Dg$, an invariant twist $$F_D\in 1
\tens{2}+\hbar(\UgD\tens{2}\fml)^{\gD}$$ such that $(R_
{D})_{F_{D}}=R\DKZ$. Performing an $F$--twist of $\Q$
by the collection $\{F_{D}\}_{D\subset\Dg}$, we obtain
an equivalent structure for which $R_{D}=R\DKZ$.

\subsection{Normalising the associators $\Phi_D$}
\label{ss:1-jet of Phi_D}

By lemma \ref{le:1 mod h2}, there exists, for each
$D\subseteq\Dg$, a symmetric invariant twist
$$F_D\in 1\tens{2}+\hbar(\UgD\tens{2}\fml)^{\gD}$$
such that $(\Phi_D)_{F_D}=1\tens{3}\mod\hbar^2$.
Twisting by $\{F_D\}_{D\subseteq\Dg}$ we may therefore
assume that $\Phi_D$ is equal to $1\tens{3}$ mod
$\hbar^2$. This twist does not alter $$R_D=R\DKZ=
\Delta_{0}(e^{\hbar/2 C_D})\cdot
e^{-\hbar/2 C_D}\otimes e^{-\hbar/2 C_D}$$
since $F_D$ is invariant and symmetric. 

\subsection{Computing the 1-jet of $F_{(D;\alpha_i)}$}
\label{ss:1-jet of F}

To complete the proof of theorem \ref{th:transfer},
we need to check that the relative twists satisfy
$$\Alt_{2}(F_{(D;\alpha_i)})
=
\hbar\cdot(\rD-r_{\g_{D\setminus\{\alpha_i\}}})
\mod\hbar^{2}$$

We shall need the following well--known 

\begin{lemma}\label{le:1-jet of Delta}
Let $\Psi:\Uhg\rightarrow\Ug\fml$ be an algebra
isomorphism equal to the identity mod $\hbar$.
Then, the following holds mod $\hbar^2$,
$$\Psi^{\otimes 2}\circ
(\Delta-\Delta^{21})\circ\Psi^{-1}=
2\hbar\cdot\ad(\rg)\circ\Delta_{0}$$
\end{lemma}
\proof It is sufficient to show that both sides
agree on the generators $e_i,f_i,h_i$ of $\g$.
Set $$\Delta_{\Psi}=\Psi\tens{2}\circ\Delta\circ
\Psi^{-1}$$ Then, modulo $\hbar^2$,
\begin{equation*}
\begin{split}
\Delta_{\Psi}(e_i)
&=
\Psi\tens{2}
(E_i\otimes 1+q_i^{H_i}\otimes E_i+\hbar\Delta\varepsilon_i)\\
&=
e_i\otimes 1+\hbar\varepsilon_i\otimes 1+
\hbar\alpha_i\otimes e_i+1\otimes e_i+\hbar 1\otimes\varepsilon'+
\hbar\Psi^{\otimes 2}\circ\Delta(\varepsilon)
\end{split}
\end{equation*}
where $\Psi^{-1}(e_i)=E_i+\hbar\varepsilon_i$, $\Psi(E_i)=e_i
+\hbar\varepsilon_i'$ and $q_i^{H_i}=1+\hbar(\alpha_i,\alpha_i)/2
H_i$ mod $\hbar^2$. Antisymmetring, and using the fact that
$\Delta(x)-\Delta^{21}(x)\in\hbar\Uhg$ for any $x\in\Uhg$,
we find that
$$\Delta_{\Psi}(e_i)-\Delta_{\Psi}^{21}(e_i)=
\hbar\cdot(\alpha_i\otimes e_i-e_i\otimes\alpha_i)=
2\hbar\cdot\alpha_i\wedge e_i\mod\hbar^2$$
A similar calculation yields
$$\Delta_{\Psi}(f_i)-\Delta_{\Psi}^{21}(f_i)=
2\hbar\cdot\alpha_i\wedge f_i\mod\hbar^2$$
and
$$\Delta_{\Psi}(h_i)-\Delta_{\Psi}^{21}(h_i)=0$$
Let now $\n_{\alpha_i}^{\pm}\subset\n^{\pm}$ be the span
of the root vectors $e_{\alpha}$ (resp. $f_{\alpha}$) with
$\alpha\neq\alpha_i$. $\n_{\alpha_i}^{\pm}$ are invariant
under the adjoint action of $\sl{2}^i$ and the inner
product $(\cdot,\cdot)$ yields an $\sl{2}^i$--equivariant
identification $(\n_{\alpha_i}^{+})^{*}\cong\n_{\alpha_i}^{-}$.
Since
$$\rg=
  \sum_{\alpha\succ 0}\frac{(\alpha,\alpha)}{2}
  \cdot e_\alpha\wedge f_\alpha=
  r_{D\setminus\alpha_i}+r_{\alpha_i}$$
where $r_{D\setminus\alpha_i}=\rg-r_{\alpha_i}$ is the
image in $\bigwedge^{2}(\n^+_{\alpha_i}\oplus\n_{\alpha_i}
^{-})$ of
$$\id_{\n^+_{\alpha_i}}\in
  \End(\n^+_{\alpha_i})\cong
  \n^+_{\alpha_i}\otimes\n^-_{\alpha_i}\subset
  (\n^+_{\alpha_i}\oplus\n_{\alpha_i}^{-})$$
we find
$$[\rg,\Delta_{0}(e_i)]=
  [r_{\alpha_i},\Delta_{0}(e_i)]=
  -\sqnorm{\alpha_i}\cdot\ad(e_i) e_i\wedge f_i=
  \alpha_i\wedge e_i$$
and similarly
$$[\rg,\Delta_{0}(f_i)]=\alpha_i\wedge f_i$$
Since $\rg$ is of weight 0, $[\rg,\Delta_0(h_i)]=$ and the
claim is proved \halmos\\

For any connected subdiagram $D\subseteq\Dg$ and \mns $\F$
on $D$, write
$$F_{(D;\F)}=1^{\otimes 2}+\hbar\cdot f_{(D;\F)}\mod\hbar^2$$
where $f_{(D;\F)}\in\UgD^{\otimes 2}$. Taking the coefficient
of $\hbar$ in $$(\Phi_{D})_{F_{(D;\F)}}=1^{\otimes 3}$$ and
using the fact that $\Phi_D=1^{\otimes 3}$ mod $\hbar^{2}$,
we find
$$\dH f_{(D;\F)}=
1\otimes f_{(D;\F)}-
\Delta_{0}\otimes\id(f_{(D;\F)})+
\id\otimes\Delta_{0}(f_{(D;\F)})-
f_{(D;\F)}\otimes 1=0$$
where $\dH:\UgD^{\otimes 2}\rightarrow\UgD^{\otimes 3}$ is
the \Ho differential. It follows that $\Alt_2(f_{(D;\F)})$
lies in $\bigwedge^{2}\gD$. On the other hand, using
(vi)$_D$, we find that, mod $\hbar^{2}$,
$$\Psi_{(D;\F)}\tens{2}\circ(\Delta-\Delta^{21})\circ\Psi_{(D;\F)}=
2\hbar\cdot\ad(\Alt_2(f_{(D;\F)}))\circ\Delta_{0}$$
By lemma \ref{le:1-jet of Delta}, this implies that
$$\Alt_2(f_{(D;\F)})-\rD\in(\bigwedge^{2}\gD)^{\gD}=0$$
as required \halmos


%% file: rigidity.tex

\section{Rigidity of $\Ug$}\label{se:rigidity}

\subsection{}

Retain the notation of section \ref{se:transfer}, particularly
\S\ref{ss:transfer}, but label the Dynkin diagram $\Dg$ by
attaching an infinite multiplicity to each edge.\footnote{A
\qc structure on $\Ug$ \wrt the usual labelling of $\Dg$
is clearly also a \qc structure \wrt the infinite labelling we
are using.
Surprisingly, the proof of theorem \ref{th:rigidity}
does not use the braid relations \eqref{eq:Brieskorn}. This
is why the result is stated in this slightly greater generality.}
The aim of this section is to prove the following

\begin{theorem}\label{th:rigidity}
Up to twisting, there exists a unique \qcqtqba structure of
type $\Dg$ on $\Ug\fml$ of the form
$$\Bigl(
\Ug\fml,\{\Ug_D\fml\},\{S\ic\},\{\Phi_{(D;\alpha_i,\alpha_j)}\},
\Delta,
\{R_D\},\{\Phi_D\},\{F_{(D;\alpha_i)}\}
\Bigr)$$
where $\Delta$ is the cocommutative coproduct on $\Ug$,
\begin{align}
S\ic&=\wt{s}_{i}\cdot\exp(\hbar/2\cdot C_{i}),
\label{eq:Si}\\
R_{D}&=\exp(\hbar\cdot\Omega_{D}),
\label{eq:RD}\\
\Alt_{2}F_{(D;\alpha_i)}&=
\hbar\cdot(\rD-r_{\g_{D\setminus\{\alpha_i\}}})
\medspace\mod\hbar^{2}
\label{eq:FDalpha}
\end{align}
and $\Phi_{(D;\alpha_i,\alpha_j)}$, $F_{(D;\alpha_i)}$ are
of weight 0.
\end{theorem}

\proof Let 
$$\Q^{a}=
\Bigl(
\Ug\fml,\{\Ug_D\fml\},\{S\ic\},\{\Phi_{(D;\alpha_i,\alpha_j)}^{a}\},
\Delta,
\{R_D\},\{\Phi_D^a\},\{F_{(D;\alpha_i)}^a\}\Bigr)$$
$a=1,2$ be two \qcqtqba structures of the above form. We
proceed in four steps

\subsection{Normalising the 1--jets of $\Phi_{D}^{a}$}
\label{ss:1-jet}

We claim first that, up to a suitable twist, we may assume
that $\Phi_{D}^{a}=1^{\otimes 3}$ mod $\hbar^{2}$ for any
$D\subseteq\Dg$ and $a=1,2$.

\begin{lemma}\label{le:1 mod h2}
Let $\Phi\in 1^{\otimes 3}+\hbar(\Ug^{\otimes 3}\fml)^
{\g}$ be a solution of the pentagon and hexagon equations
with respect to $R=e^{\hbar\Omega}$. Then, there exists
a symmetric, invariant twist
$$F\in 1^{\otimes 2}+\hbar(\Ug^{\otimes 2}\fml)^{\g}$$
such that $(\Phi)_{F}=1^{\otimes 3}\mod\hbar^{2}$.
\end{lemma}
\proof Write
$$\Phi=1^{\otimes 3}+\hbar\varphi\mod\hbar^{2}$$
where $\varphi\in(\Ug^{\otimes 3})^{\g}$. The pentagon
equation for $\Phi$ implies that
$$\dH\varphi=
1\otimes\varphi-
\Delta\otimes\id^{\otimes 2}(\varphi)+
\id\otimes\Delta\otimes\id(\varphi)-
\id^{\otimes 2}\otimes\Delta(\varphi)+
\varphi\otimes 1=0$$
where $\dH$ is the \Ho differential. By \cite[Prop. 3.5]
{Dr2}, $(\Phi)^{-1}=(\Phi)^{321}$. Substituting this
into the second of the hexagon relations
\begin{align}
\Delta\otimes\id(R)&=
\Phi^{312}\cdot R^{13}\cdot(\Phi^{132})^{-1}\cdot R^{23}
\cdot\Phi^{123}
\label{eq:hex 1}
\\
\id\otimes\Delta(R)&=
(\Phi^{231})^{-1}\cdot R^{13}\cdot\Phi^{213}\cdot R^{12}
\cdot(\Phi^{123})^{-1}
\label{eq:hex 2}
\end{align}
substracting them, and taking the coefficient of $\hbar$
yields that $\Alt_{3}\varphi=0$. Thus, $\varphi=\dH f$
where $f\in\Ug^{\otimes 2}$ may be chosen invariant under
$\g$. Since
$$\dH f^{21}=-(\dH f)^{321}=-\varphi^{321}=\varphi=\dH f$$
we may further assume, up to replacing $f$ by $(f+f^{21})
/2$, that $f$ is symmetric . Setting $F=1-\hbar f$ yields
the required twist \halmos\\

For any connected $D\subseteq\Dg$, let $F_{D}^{a}\in
1^{\otimes 2}+\hbar(U\gD^{\otimes 2}\fml)^{\gD}$ be
a symmetric invariant twist such that $(\Phi_{D}^{a})
_{F_{D}^{a}}=1^{\otimes 3}\mod\hbar^{2}$. Twisting
$\Q^{a}$ by $F^{a}=\{F_{D}^{a}\}_{D\subseteq\Dg}$
yields the claimed result. Note that
$$(R_{D})_{F_{D}^{a}}=
{F_{D}^{a}}^{21}R_{D}{F_{D}^{a}}^{-1}=
R_{D}$$
since $R_{D}=\Delta(\exp(\hbar/2 C_{D}))\cdot\exp
(-\hbar/2 C_{D})^{\otimes 2}$ and $F_{D}^{a}$ is 
symmetric and invariant under $\gD$, and that
$$\Alt_{2}
(F_{(D;\alpha_i)}^{a})_{F^{a}}=
F_{D\setminus\{\alpha_i\}}^{a}\cdot
\Alt_{2}(F_{(D;\alpha_i)}^{a})\cdot 
{F_{D}^{a}}^{-1}=
\hbar\cdot(\rD-r_{\g_{D\setminus\{\alpha_i\}}})$$
since $F^{a}$ is symmetric. Thus, twisting $\Q^{a}$
by $F^{a}$ preserves the conditions \eqref{eq:Si}--\eqref
{eq:FDalpha}.

\subsection{Matching the associators $\Phi_{D}^{a}$}
\label{ss:associators}

We claim next that, up a twist, we may assume that $\Phi
_{D}^{2}=\Phi_{D}^{1}$ for all $D\subseteq\Dg$. Indeed,
for any such $D$ there exists, by Drinfeld's uniqueness
theorem \cite[prop. 3.12]{Dr3}, a symmetric, invariant
twist
$$F_{D}\in 1^{\otimes 2}+\hbar(U\gD^{\otimes 2}\fml)^{\gD}$$
such that $(\Phi_{D}^{2})_{F_{D}}=\Phi_{D}^{1}$. Twisting
$\Q^{2}$ by $F=\{F_{D}\}_{D\subseteq\Dg}$ yields the claimed
equality of associators and, as in the previous step,
preserves the conditions \eqref{eq:Si}--\eqref{eq:FDalpha}.

\subsection{Matching the twists $F_{(D;\alpha_i)}^{a}$}
\label{ss:twists}

We claim now that, up to a further twist which does
not alter the associators $\Phi_{D}^{1}=\Phi_{D}^{2}$,
we may assume that $F_{(D;\alpha_i)}^{2}=F_{(D;\alpha_i)}
^{1}$ for any $\alpha_i\in D\subseteq\Dg$. We need two
preliminary results.

\begin{lemma}\label{le:non deg}
If $\Phi=1+\hbar^{2}\varphi+\cdots\in 1^{\otimes 3}+\hbar
^{2}(\Ug^{\otimes 3}\fml)^{\g}$ satisfies the pentagon and
hexagon equations with respect to $R=e^{\hbar\Omega}$,
then
$$\Alt_{3}\varphi=\frac{1}{6}[\Omega_{12},\Omega_{23}]$$
\end{lemma}
\proof Since $R$ is symmetric, proposition 3.5 of \cite
{Dr2} implies that $\Phi^{-1}=\Phi^{321}$. Substituting
this into the second hexagon equation \eqref{eq:hex 2},
substracting it from \eqref{eq:hex 1}  and taking the
coefficient of $\hbar^{2}$ shows that
\begin{equation*}
\begin{split}
6\Alt_{3}(\varphi)
&=
\half{1}\left(
\Delta\otimes\id(\Omega^{2})-\id\otimes\Delta(\Omega^{2})+
\Omega^{2}_{12}-\Omega^{2}_{23}\right)
+\Omega_{13}(\Omega_{12}-\Omega_{23})\\
&=
\half{1}\left(
\Omega_{13}\Omega_{23}+\Omega_{23}\Omega_{13}-
\Omega_{12}\Omega_{13}-\Omega_{13}\Omega_{12}\right)
+\Omega_{13}(\Omega_{12}-\Omega_{23})\\
&=
\half{1}\left(
\Omega_{13}(\Omega_{12}+\Omega_{23}+\Omega_{13})-
\Omega_{13}(\Omega_{12}+\Omega_{13})+
\Omega_{23}\Omega_{13}-\Omega_{13}\Omega_{12}\right.\\
&\phantom{=}
-\left.(\Omega_{12}+\Omega_{23}+\Omega_{13})\Omega_{13}
+(\Omega_{23}+\Omega_{13})\Omega_{13}\right)
+\Omega_{13}(\Omega_{12}-\Omega_{23})\\
&=
[\Omega_{23},\Omega_{13}]\\
&=
-[\Omega_{23},\Omega_{12}]
\end{split}
\end{equation*}
where the fourth equality uses the fact that $\Omega_{12}+
\Omega_{23}+\Omega_{13}$ commutes with $\Omega_{ij}$ 
\halmos

\begin{lemma}\label{le:rank 1}
Let $F\in 1+\hbar{\Usl{2}^{i}}^{\otimes 2}\fml$ be a twist
of weight zero and set
$$S\ic=\wt{s}_{i}\cdot\exp(\hbar/2\cdot C_{i}),
\qquad
R_{i}=\exp(\hbar\cdot\Omega^{i})$$
Then, the equation
\begin{equation}\label{eq:coproduct}
\Delta_{F}(S\ic)=(R_{i})_{F}^{21}\cdot S\ic\otimes S\ic
\end{equation}
is equivalent to $F^{\Theta}=F^{21}$ where $\Theta\in\Aut
(\sl{2}^{i})$ is any involution such that $\Theta(h_{\root
{i}})=-h_{\alpha_i}$.
\end{lemma}
\proof Since $\Delta(C_{i})=C_{i}\otimes
1+1\otimes C_{i}+2\Omega^{i}$,
\begin{equation*}
\begin{split}
\Delta_{F}(S\ic)
&=
F\cdot
\exp(\hbar\Omega^{i})
\cdot S\ic\otimes S\ic
\cdot F^{-1}\\
&=
(R_{i})_{F}^{21}\cdot
F^{21}\cdot
S\ic\otimes S\ic\cdot
F^{-1}
\end{split}
\end{equation*}
which is equal to the right--hand side of \eqref{eq:coproduct}
if, and only if
$$\Ad(\wt{s}_{i}^{\otimes 2})F=F^{21}$$
The claim follows since $\Theta=\Ad(\wt{s}_{i})\Ad(c\cdot
h_{\alpha_i})$, for some $c\in\IC$ so that $\Theta$ and
$\Ad(\wt{s}_{i}^{\otimes 2})$ coincide on zero weight
elements of $\Ug^{\otimes 2}$ \halmos\\

Fix now $\alpha_i\in D\subseteq\Dg$ and, for $D'=D,D
\setminus\{\alpha_i\}$, denote $\Phi_{D'}^{2}=\Phi_{D'}
^{1}$ by $\Phi_{D'}$. Write
$$\Phi_{D'}=1+\hbar^{2}\varphi_{D'}\mod\hbar^{3}$$
By lemma \ref{le:non deg}, 
$$\Alt_{3}\varphi_{D}=[\Omega^{D}_{12},\Omega^{D}_{23}]
\qquad\text{and}\qquad
\Alt_{3}\varphi_{D\setminus\{\alpha_i\}}=
[\Omega^{D\setminus\{\alpha_i\}}_{12},
 \Omega^{D\setminus\{\alpha_i\}}_{23}]$$
so that $\Phi_{D}$ and $\Phi_{D\setminus\alpha_i}$ are
non--degenerate in the sense of definition 5.1 of \cite{TL3}
and $\ol{\pi}^{3}(\Alt_{3}\varphi_{D})=\Alt_{3}\varphi
_{D\setminus\alpha_i}$ where
$$\ol{\pi}^{3}:
(U\gD^{\otimes 3})^{\gD}
\longrightarrow
(U\g_{D\setminus\{\alpha_i\}}^{\otimes 3})^
{\g{D\setminus\{\alpha_i\}}}$$
is the generalised \HC homomorphism defined in \S 2 of \cite{TL3}.
Since $(\Phi_{D})_{F_{(D;\alpha_i)}^{a}}=\Phi_{D\setminus\{\alpha_i\}}$
for $a=1,2$, there exists, by \cite[thm. 6.1(iv)]{TL3} a gauge transformation
$$a_{(D;\alpha_i)}\in1+\hbar U\gD\fml^{\ll_{D\setminus\alpha_i}}$$
such that 
$$F_{(D;\alpha_i)}^{1}=
a_{(D;\alpha_i)}^{\otimes 2}\cdot
F_{(D;\alpha_i)}^{2}\cdot
\Delta(a_{(D;\alpha_i)}^{-1})$$
Moreover, by lemma \ref{le:rank 1} and \cite[thm. 6.1(iii)]{TL3}, we may
assume that
\begin{equation}\label{eq:a fixed}
\Ad(\wt{s}_{i})a_{(\alpha_i;\alpha_i)}=
a_{(\alpha_i;\alpha_i)}
\end{equation}
for any $i=1\ldots n$. Twisting $\Q^{2}$ by $a=\{a_{(D;
\alpha_i)}\}_{\alpha_i\in D\subseteq\Dg}$ yields the
required equality of twists while preserving \eqref
{eq:Si}--\eqref{eq:FDalpha} since
$$(S\ic)_{a}=
a_{(\alpha_i;\alpha_i)}\cdot
\wt{s}_{i}\cdot\exp(\hbar/2 C_{i})
a_{(\alpha_i;\alpha_i)}^{-1}=
S\ic$$
by \eqref{eq:a fixed} and, mod $\hbar^{2}$, 
$$\Alt_{2}((F^{2}_{(D;\alpha_i)})_{a})=
\Alt_{2}(F^{2}_{(D;\alpha_i)})+
\hbar\Alt_{2}\dH({a_{(\alpha_i;\alpha_i)}}_{1})=
\Alt_{2}(F^{2}_{(D;\alpha_i)})$$
where $a_{(\alpha_i;\alpha_i)}=1+\hbar{a_{(\alpha_i;\alpha_i)}}_{1}
\mod\hbar^{2}$. 

\subsection
{Matching the associators $\Phi^{a}_{(D;\alpha_i,\alpha_j)}$}
\label{ss:coxeter associators}

We may henceforth assume that
$$\Phi_{D}^{2}=\Phi_{D}^{1}
\qquad\text{and}\qquad
F_{(D;\alpha_i)}^{2}=F_{(D;\alpha_i)}^{1}$$
for any $\alpha_i\in D\subseteq\Dg$, and that $\Phi_{(D;\alpha_i,
\alpha_j)}^{1}=\Phi_{(D;\alpha_i,\alpha_j)}^{2}$ mod $\hbar^{n}$
for some $n\geq 1$ and all $\alpha_i\neq\alpha_j\in D\subseteq
\Dg$. Thus, 
$$\Phi_{(D;\alpha_i,\alpha_j)}^{2}=
\Phi_{(D;\alpha_i,\alpha_j)}^{1}+
\hbar^{n}\varphi_{(D;\alpha_i,\alpha_j)}
\mod \hbar^{n+1}$$
for some $\varphi_{(D;\alpha_i,\alpha_j)}\in U\gD^{\g_{D\setminus
\{\alpha_i,\alpha_j\}}}$. Let $\alpha_i\neq\alpha_j\in D\subseteq
\Dg$ and $\F,\G$ two \fmns such that
$$\F\setminus\G=\conn{D}{\alpha_i}{\alpha_j}
\quad\text{and}\quad
\G\setminus\F=\conn{D}{\alpha_j}{\alpha_i}$$
Substracting the equations
\begin{align*}
F_{\F}\cdot\Delta(\Phi_{(D;\alpha_i,\alpha_j)}^{1})
&=
{\Phi_{(D;\alpha_i,\alpha_j)}^{2}}^{\otimes 2}\cdot F_{\G}\\
F_{\F}\cdot\Delta(\Phi_{(D;\alpha_i,\alpha_j)}^{2})
&=
{\Phi_{(D;\alpha_i,\alpha_j)}^{2}}^{\otimes 2}\cdot F_{\G}
\end{align*}
where, as usual
$$F_{\F}=\stackrel{\longrightarrow}{\prod}_{D'\in\F}
F_{(D',\alpha^{D'}_{\F})}$$
and $F_{(D',\alpha')}=F_{(D',\alpha')}^{1}=F_{(D',\alpha')}^{2}$,
and equating the coefficients of $\hbar^{n+1}$, we find
$$\Delta(\varphi_{(D;\alpha_i,\alpha_j)})-
\varphi_{(D;\alpha_i,\alpha_j)}\otimes 1-
1\otimes \varphi_{(D;\alpha_i,\alpha_j)}=0$$
so that $\varphi_{(D;\alpha_i,\alpha_j)}$ is a primitive
element of $\Ug_{D}$ and therefore lies in $\g_{D}$.
Since $\varphi_{(D;\alpha_i,\alpha_j)}$ is also of weight
0, we find that
$$\varphi_{(D;\alpha_i,\alpha_j)}\in\h_{D}$$
where $\hD\subset\gD$ is the span of the simple roots
$\alpha_{i}\in D$. Since $\Phi^{a}_{(D;\alpha_i,\alpha_j)}$
satisfy the generalised pentagon identities corresponding
to the 2--faces of the \DCP associahedron $\Ag$, we also
find
$$d_{D}\{\varphi_{(D;\alpha_i,\alpha_j)}\}=0$$

\begin{proposition}\label{pr:H20}
Let $\varphi=\{\varphi_{(D;\alpha_i,\alpha_j)}\}$ be a
2--cocycle in the Dynkin complex of $\Ug$ such that
$$\varphi_{(D;\alpha_i,\alpha_j)}\in\h_{D}$$
for any $\alpha_i\neq\alpha_j\in D\subseteq\Dg$. Then,
there exists a Dynkin 1--cochain $a=\{a_{(D;\alpha_i)}\}$
such that
$$a_{(D;\alpha_i)}\in\h_{D}
\quad\text{and}\quad
\dD a=\varphi$$
The element $a$ may be chosen such that $a_{(\alpha_i;
\alpha_i)}=0$ for all $\alpha_i$ and is then unique with
this additional property.
\end{proposition}
\proof It will be convenient to fix an order $\alpha_1,\ldots,\alpha_n$
of the simple roots and identify the group of Dynkin cochains $CD^p
(A;M)$ with elements $m=\{m_{(B;\aalpha)}\}$ where $B$ ranges
over the connected subdiagrams of $\Dg$ and $\aalpha$ over the
subsets $\{\alpha_{i_1},\ldots,\alpha_{i_p}\}\subseteq B$ such that
$i_1<\ldots<i_p$. We wish to solve the equation $\varphi=\dD a$.
In components, this reads
\begin{equation}\label{eq:coh 1}
\varphi_{(D;\alpha_i,\alpha_j)}=
 a_{(D;\alpha_j)}
-a_{(\conn{D}{\alpha_i}{\alpha_j};\alpha_j)}
-a_{(D;\alpha_i)}
+a_{(\conn{D}{\alpha_j}{\alpha_i};\alpha_i)}
\end{equation}
for any connected subdiagram $D\subseteq\Dg$ and $i<j$
such that $\alpha_i,\alpha_j\in D$. The assumptions
$\varphi_{(D;\alpha_i,\alpha_j)},a_{(D;\alpha_i)}\in\h_{D}$
and the fact that $\varphi,a$ lie in the Dynkin complex of
$\g$ imply that
$$\varphi_{(D;\alpha_i,\alpha_j)}
\in\IC\cow{i}\oplus\IC\cow{j}
\quad\text{and}\quad
a_{(D;\alpha_i)}
\in\IC\cow{i}$$
respectively where $\cow{k}$ is the fundamental coweight 
dual to $\alpha_k$. Projecting \eqref{eq:coh 1} on $\cow
{i}$ and $\cow{j}$ we therefore find that it is equivalent
to 
\begin{align*}
a_{(D;\alpha_i)}&=
a_{(\conn{D}{\alpha_j}{\alpha_i};\alpha_i)}
-\varphi_{(D;\alpha_i,\alpha_j)}^{i}\\
a_{(D;\alpha_j)}&=
a_{(\conn{D}{\alpha_i}{\alpha_j};\alpha_j)}
+\varphi_{(D;\alpha_i,\alpha_j)}^{j}
\end{align*}
where
$$\varphi_{(D;\alpha_i,\alpha_j)}=
\varphi_{(D;\alpha_i,\alpha_j)}^{i}\cow{i}+
\varphi_{(D;\alpha_i,\alpha_j)}^{j}\cow{j}$$
and we are identifying $a_{(D;\alpha_i)},a_{(D,\alpha_j)}$
with their components along $\cow{i},\cow{j}$ respectively.
Thus, $\varphi=\dD a$ iff 
\begin{equation}\label{eq:coh 2}
a_{(D;\alpha_i)}=
a_{(\conn{D}{\alpha_j}{\alpha_i};\alpha_i)}
-(-1)^{(i:j)}\varphi_{(D;\alpha_i,\alpha_j)}^{i}
\end{equation}
for all $D$ and $1\leq i\neq j\leq n$ with $\alpha_i,\alpha_j
\in D$, where $(i:j)=0$ if $i<j$ and $1$ otherwise. Induction
on the cardinality of $D$ readily shows that the above equations
possess at most one solution once the values of $a_{(\alpha_i;
\alpha_i)}$ are fixed. To prove that one solution exists, assume
that $a_{(D;\alpha_{i})}$ have been constructed for all $D$
with at most $m$ vertices in such a way that the equations
\eqref{eq:coh 2} hold for all such $D$. We claim that \eqref
{eq:coh 2} may be used to define $a_{(D;\alpha_i)}$ for all
$D$ with $|D|=m+1$ in a consistent way, \ie independently of
$j\neq i$ such that $\alpha_j\in D$. This amounts to showing
that for all such $D$, and distinct vertices $\alpha_i,\root
{j},\root{k}\in D$, one has
\begin{equation}\label{eq:consistent}
a_{(\conn{D}{\alpha_j}{\alpha_i};\alpha_i)}
-(-1)^{(i:j)}\varphi_{(D;\alpha_i,\alpha_j)}^{i}
=
a_{(\conn{D}{\root{k}}{\alpha_i};\alpha_i)}
-(-1)^{(i:k)}\varphi_{(D;\alpha_i,\root{k})}^{i}
\end{equation}

To see this, consider the $(D;\alpha_i,\alpha_j,\root{k})$
component of $\dD\varphi$ \ie the sum
\begin{equation*}
\begin{split}
 &
(-1)^{(i:j)+(i:k)}\left(
\varphi_{(D;\alpha_j,\root{k})}-
\varphi_{(\conn{D}{\alpha_i}{\alpha_j,\root{k}};\alpha_j,\root{k})}
\right)\\
+&
(-1)^{(j:i)+(j:k)}\left(
\varphi_{(D;\alpha_i,\root{k})}-
\varphi_{(\conn{D}{\alpha_j}{\alpha_i,\root{k}};\alpha_i,\root{k})}
\right)\\
+&
(-1)^{(k:i)+(k:j)}\left(
\varphi_{(D;\alpha_i,\alpha_j)}-
\varphi_{(\conn{D}{\root{k}}{\alpha_i,\alpha_j};\alpha_i,\alpha_j)}
\right)
\end{split}
\end{equation*}
Since $\dD\varphi=0$ we get, by projecting on $\cow{i}$,
$$\varphi_{(D;\alpha_i,\alpha_j)}^{i}
=
\varphi_{(\conn{D}{\root{k}}{\alpha_i,\alpha_j};\alpha_i,\alpha_j)}^{i}
+(-1)^{(i:j)+(i:k)}
\left(
\varphi_{(D;\alpha_i,\root{k})}^{i}-
\varphi_{(\conn{D}{\alpha_j}{\alpha_i,\root{k}};\alpha_i,\root{k})}^{i}
\right)$$
so that \eqref{eq:consistent} holds iff
\begin{equation}\label{eq:consistent 2}
a_{(\conn{D}{\alpha_j}{\alpha_i};\alpha_i)}
+(-1)^{(i:k)}
\varphi^{i}_{(\conn{D}{\alpha_j}{\alpha_i,\root{k}};\alpha_i,\root{k})}
=
a_{(\conn{D}{\root{k}}{\alpha_i};\alpha_i)}
+(-1)^{(i:j)}
\varphi^{i}_{(\conn{D}{\root{k}}{\alpha_i,\alpha_j};\alpha_i,\alpha_j)}
\end{equation}
We consider four separate cases.

\subsubsection
{$\conn{D}{\alpha_j}{\alpha_i,\root{k}}=\emptyset$
and $\conn{D}{\root{k}}{\alpha_i,\alpha_j}=\emptyset$}

This case cannot arise since the first condition implies
that any path in $D$ from $\alpha_i$ to $\root{k}$ must
pass through $\alpha_j$ before it reaches $\root{k}$ while
the second one implies that the portion of this path
linking $\alpha_i$ to $\alpha_j$ must first pass through
$\root{k}$.

\subsubsection
{$\conn{D}{\alpha_j}{\alpha_i,\root{k}}=\emptyset$ and
$\conn{D}{\root{k}}{\alpha_i,\alpha_j}\neq\emptyset$}

In this case,
$$\varphi_{(\conn{D}{\alpha_j}{\alpha_i,\root{k}};\alpha_i,\root{k})}=0
\quad\text{and}\quad
\conn{D}{\root{k}}{\alpha_i}=\conn{D}{\root{k}}{\alpha_i,\alpha_j}$$
and \eqref{eq:consistent 2} reads
\begin{equation*}
\begin{split}
a_{(\conn{D}{\alpha_j}{\alpha_i};\alpha_i)}
&=
a_{(\conn{D}{\root{k}}{\alpha_i,\alpha_j};\alpha_i)}
+(-1)^{(i:j)}
\varphi^{i}_{(\conn{D}{\root{k}}{\alpha_i,\alpha_j};\alpha_i,\alpha_j)}\\
&=
a_{(\conn{\conn{D}{\root{k}}{\alpha_i,\alpha_j}}{\alpha_j}{\alpha_i};
\alpha_i)}
\end{split}
\end{equation*}
where the last equality follows from \eqref{eq:coh 2} applied
to the diagram $\conn{D}{\root{k}}{\alpha_i,\alpha_j}$. This
equation however holds since
$$\conn{D}{\alpha_j}{\alpha_i}=
\conn{D}{\alpha_j,\root{k}}{\alpha_i}=
\conn{\conn{D}{\root{k}}{\alpha_i,\alpha_j}}{\alpha_j}{\alpha_i}$$
where the first equality holds because $\conn{D}{\alpha_j}
{\alpha_i,\root{k}}=\emptyset$.

\subsubsection
{$\conn{D}{\alpha_j}{\alpha_i,\root{k}}\neq\emptyset$
and $\conn{D}{\root{k}}{\alpha_i,\alpha_j}=\emptyset$}

This case reduces to the previous one under the interchange
$\alpha_j\leftrightarrow\root{k}$.

\subsubsection
{$\conn{D}{\alpha_j}{\alpha_i,\root{k}}\neq\emptyset$
and $\conn{D}{\root{k}}{\alpha_i,\alpha_j}\neq\emptyset$}

In this case 
$$\conn{D}{\alpha_j}{\alpha_i}=
\conn{D}{\alpha_j}{\alpha_i,\root{k}}
\quad\text{and}\quad
\conn{D}{\root{k}}{\alpha_i}=
\conn{D}{\root{k}}{\alpha_i,\alpha_j}$$
so that \eqref{eq:consistent} reads
\begin{equation*}
\begin{split}
a_{(\conn{\conn{D}{\alpha_j}{\alpha_i,\root{k}}}{\root{k}}{\alpha_i};\alpha_i)}
&=
a_{(\conn{D}{\alpha_j}{\alpha_i,\root{k}};\alpha_i)}
+(-1)^{(i:k)}
\varphi^{i}_{(\conn{D}{\alpha_j}{\alpha_i,\root{k}};\alpha_i,\root{k})}\\
&=
a_{(\conn{D}{\root{k}}{\alpha_i,\alpha_j};\alpha_i)}
+(-1)^{(i:j)}
\varphi^{i}_{(\conn{D}{\root{k}}{\alpha_i,\alpha_j};\alpha_i,\alpha_j)}\\
&=
a_{(\conn{\conn{D}{\root{k}}{\alpha_i,\alpha_j}}{\alpha_j}{\alpha_i};\alpha_i)}
\end{split}
\end{equation*}
where the first and last equalities follow from \eqref
{eq:consistent} for the diagrams $\conn{D}{\alpha_j}
{\alpha_i,\root{k}}$ and $\conn{D}{\root{k}}{\alpha_i,
\alpha_j}$ respectively. This equation however holds
because
$$\conn{\conn{D}{\alpha_j}{\alpha_i,\root{k}}}{\root{k}}{\alpha_i}=
\conn{D}{\alpha_j,\root{k}}{\alpha_i}=
\conn{\conn{D}{\root{k}}{\alpha_i,\alpha_j}}{\alpha_j}{\alpha_i}$$
This concludes the proof of proposition \ref{pr:H20} \halmos\\

We may now conclude the proof of theorem \ref{th:rigidity}.
Twist $\Q^{2}$ by
\begin{equation}
a=\{1-\hbar^{n}\cdot a_{(D;\alpha_i)}\}_{\alpha_i\in D\subseteq\Dg}
\end{equation}
where the $a_{(D;\alpha_i)}$ are given by proposition \ref{pr:H20}
and $a_{(\alpha_i;\alpha_i)}=0$ for any $i$.
Then $\Phi_{(D;\alpha_i,\alpha_j)}^{2}=\Phi_{(D;\alpha_i,\alpha_j)}
^{1}\mod\hbar^{n+1}$ for any $\alpha_i\neq\alpha_j\in D\subseteq\Dg$
and the conditions \eqref{eq:Si}--\eqref{eq:FDalpha} are preserved
since, owing to the fact that $a_{(\alpha_i;\alpha_i)}=0$,
we have $(S\ic)_{a}=S\ic$ for any $\alpha_i\in\Dg$ \halmos\\ 


%% file: acknowledgements.tex
{\bf Acknowledgements.} Above all, I would like to thank R. Nest and M.
Kashiwara for intense and stimulating conversations at the early stages
of this project which proved invaluable at a time when I was struggling
with the definition of \qcas and with my own doubts about their very
existence. I am extremely grateful for their generous time and keen
interest in my work. Conversations and correspondence with a number
of mathematicians greatly helped me to better understand certain aspects
of Drinfeld's work on \qt quasibialgebras. Among these, I would particularly
like to thank P. Etingof, J. Roberts and P. Teichner. It is also a pleasure to
thank E. Frenkel, R. Rouquier and G. Skandalis for their encouragement,
friendly help and for very many valuable conversations.\\

The present project was begun in 2000--2001 while the author was a
postdoctoral fellow at the Mathematical Sciences Research Institute in
Berkeley and carried out at various times in the years 2001--2005 during
visits to the Research Institute of Mathematical Sciences of the University
of Kyoto, MSRI, the mathematics departments of the Universities of Rome
I, Copenhagen, Berkeley, Cambridge and the Erwin Schr\"odinger Institute
in Vienna. I am very grateful to these institutions for their financial support
and wonderful working conditions and am indebted to many mathematicians
for their kind invitations, in particular E. Frenkel, V. Jones, M. Kashiwara,
T. Kohno, R. Nest, K. Saito and C. Teleman. I would also like to thank the
Universit\'e Pierre et Marie Curie for granting me a leave of absence in
2000--2001 and a sabbatical leave in the spring and fall semesters of 2004.
